\documentclass{article}
\usepackage[left = 2.5cm, right = 2.5cm, bottom = 3cm, top = 2.5cm]{geometry}

\usepackage{amsmath,amssymb,amsthm,bm,bbm,dsfont}
\usepackage{thmtools}
\usepackage{thm-restate}

\DeclareMathOperator{\sign}{sign}
\DeclareMathOperator{\expect}{E}
\DeclareMathOperator{\var}{var}

\DeclareMathOperator{\cov}{cov}
\DeclareMathOperator{\trace}{trace}
\DeclareMathOperator{\diag}{diag}

\usepackage{booktabs}
\usepackage{multirow}
\usepackage{rotating}
\usepackage{graphicx}

\usepackage{subcaption}
\usepackage{tabularx}
\usepackage{dcolumn}
\usepackage{tikz}
\usetikzlibrary{positioning,arrows.meta,calc,fit,backgrounds}
\usepackage{float}

\usepackage{pdfpages}

\usepackage{natbib}
\usepackage{bibunits}
\usepackage{hyperref}
\usepackage{autonum}
\defaultbibliographystyle{jss2}
\defaultbibliography{mDYPL}
\usepackage{siunitx}

\usepackage{enumitem}

\usepackage[framemethod=TikZ]{mdframed}

\newtheorem{theorem}{Theorem}[section]

\newtheorem{proposition}[theorem]{Proposition}
\newtheorem{corollary}[theorem]{Corollary}
\newtheorem{lemma}{Lemma}
\theoremstyle{definition}

\surroundwithmdframed[
topline=false,
rightline=false,
bottomline=false,
leftmargin=\parindent,
skipabove=\medskipamount,
skipbelow=\medskipamount,
linecolor=grey!30,
innerlinewidth=1pt
]{proof, definition, remark, theorem, lemma, proposition, corollary}

\newcommand{\restatewithlabelrestore}[1]{%
		\let\hdlrestateLabel\label
		#1%
		\let\label\hdlrestateLabel
}
\makeatletter
\renewcommand{\theparagraph}{\thesubsection.\arabic{paragraph}}
\newcommand{\proofstrategyparagraph}[2]{%
		\refstepcounter{paragraph}%
		\def\@currentlabelname{#1}%
		\label{#2}%
		\paragraph*{\theparagraph\ #1.}%
}
\makeatother

\makeatletter
\newcommand{\vnorm}[1]{%
	\ensuremath{%
		\if@display
		\left\| #1 \right\|
		\else
		\| #1 \|
		\fi
	}%
}
\makeatother

\makeatletter
\newcommand{\abs}[1]{%
	\ensuremath{%
		\if@display
		\left| #1 \right|
		\else
		\lvert #1 \rvert
		\fi
	}%
}
\makeatother

\newcommand{\mnorm}[1]{{\left\vert\kern-0.25ex\left\vert\kern-0.25ex\left\vert #1
		\right\vert\kern-0.25ex\right\vert\kern-0.25ex\right\vert}}
\newcommand{\mnorms}[1]{{\vert\kern-0.25ex\vert\kern-0.25ex\vert #1
		\vert\kern-0.25ex\vert\kern-0.25ex\vert}}
\newcommand*{\bb}{\boldsymbol}
\newcommand{\Op}[1]{\ensuremath{{\mathcal{O}_{\mathrm{p}}(#1)}}}
\newcommand{\op}[1]{\ensuremath{{o_{\mathrm{p}}(#1)}}}

\newcommand{\tnod}{\ensuremath{{\theta_{0}}}}

\newcommand{\logit}[1]{\ensuremath{\textrm{logit}(#1)}}

\newcommand{\prox}[2]{\ensuremath{\textrm{prox}_{#1}\left(#2\right)}}
\newcommand{\betady}{\ensuremath{\hat{\bbeta}^{\textrm{\tiny DY}}}}
\newcommand{\etady}{\ensuremath{\hat{\bb{\eta}}^{\textrm{\tiny DY}}}}

\newcommand{\thetady}{\ensuremath{\hat{\theta}^{\textrm{\tiny DY}}}}

\def\bW {\bb{W}}
\def\bw {\bb{w}}

\def\bE {\bb{E}}
\def\bH {\bb{H}}
\def\bX{\bb{X}}
\def\bx{\bb{x}}
\def\bmu{\bb{\mu}}
\def\b0{\bb{0}}
\def\bu{\bb{u}}
\def\bv{\bb{v}}
\def\bh{\bb{h}}
\def\bg{\bb{g}}
\def\by{\bb{y}}
\def\bY{\bb{y}}
\def\Ytil{\widetilde{Y}}
\def\bYtil{\widetilde{\bb{y}}}
\def\bEtaYtil{\bb{\eta}_{\Ytil}}
\def\bba{\bb{a}}
\def\bA{\bb{A}}
\def\bB{\bb{B}}

\def\br{\bb{r}}

\def\bSigma{\bb{\Sigma}}
\def\bTheta{\bb{\Theta}}
\def\bbeta{\bb{\beta}}
\def\bxi{\bb{\xi}}

\def\bz{\bb{z}}

\def\bH{\bb{H}}
\def\blambda{\bb{\lambda}}

\newcommand{\vspan}{\operatorname{span}}

\newcommand{\rank}{\operatorname{rank}}
\newcommand{\range}{\operatorname{range}}

\usepackage{xcolor}

\usepackage{authblk}
\usepackage{orcidlink}

\author[1]{Philipp Sterzinger~\orcidlink{0009-0007-7348-5810}
}

\title{Proportional-limit asymptotics for Diaconis--Ylvisaker-penalised logistic regression with fitted intercept}

\makeatletter
\let\mainmaketitle\maketitle
\let\main@maketitle\@maketitle
\let\main@title\@title
\let\main@author\@author
\let\main@date\@date
\let\main@thanks\@thanks
\let\mainthanks\thanks
\let\mainand\and
\makeatother

\begin{document}
\begin{bibunit}

	\maketitle

\begin{abstract}
This paper develops estimator-level asymptotic theory for maximum
Diaconis--Ylvisaker prior penalised likelihood for logistic regression with a 
jointly fitted
intercept and nonzero prior slope direction in the proportional-limit regime.  
For $\mathrm{N}(\bb{0}_p, p^{-1}\bb I_p)$ Gaussian covariates and
$p/n\to\kappa\in(0,1)$, a conditional convex Gaussian min--max analysis yields almost-sure
convergence of the fitted intercept and a pseudo-Lipschitz empirical law for
the slope estimator. 
This estimator-level law gives limits for out-of-sample scores, classification error, optimal
thresholding and oracle calibration. 
It also yields oracle-adjusted
fixed-block $Z$-statistics under isotropic Gaussian covariates, and 
yields the main ingredient in establishing the limiting distribution of the 
penalised likelihood-ratio test statistic and identifies the rescaling to 
recover the nominal chi-square distribution.
We extend these results to Gaussian designs with deterministic mean and
positive-definite covariance via affine centering
and whitening and discuss extensions to subgaussian covariates. 
Finally, we propose a consistent response-moment estimator of the oracle parameters entering the state equations that govern the slope limiting law and are required for feasible inference.
\end{abstract}


\section{Introduction}
\label{sec:introduction}

\subsection{Logistic regression}
\label{subsec:logistic-regression}

Consider observations $\{Y_i, \bx_i\}_{i = 1}^n$, independent across $i = 1,\ldots,n$, from the
logistic regression model
\begin{equation}
\label{eq:logistic-model-intro}
    \Pr(Y_i=1\mid\bx_i)
    =
    \rho'\left(\theta_0+\bx_i^\top\bbeta_0\right), 
    \quad
    \rho(x)=\log(1+\exp\{x\})\,,
\end{equation}
where $\theta_0\in\Re$ and $\bbeta_0\in\Re^p$ are the unknown intercept and
slope parameters, respectively. When it exists
(cf. \citealt{albert+anderson:1984}), the maximum likelihood (ML) estimator is
\begin{equation}
\label{eq:log-reg-ml}
    (
        \hat\theta^{\textrm{\tiny ML}},
        \hat{\bbeta}^{\textrm{\tiny ML}}
    )
    =
    \underset{\substack{\theta\in\Re, \, \bbeta\in\Re^p}}{\arg\max} \, \ell(\theta, \bbeta; \bY, \bX), \quad  \ell(\theta, \bbeta; \bY, \bX) =
    \sum_{i=1}^n
    \left\{
        Y_i\left(\theta+\bx_i^\top\bbeta\right)
        -
        \rho\left(\theta+\bx_i^\top\bbeta\right)
    \right\}\,,
\end{equation}
with $\by  = (Y_1, \ldots, Y_n)^\top$ and where $\bX$ is the $n \times p$ matrix with $i$th row $\bx_i^\top$. 
Nonexistence may be undesirable to practitioners, and can be 
problematic if undetected, giving rise to spurious inference and prediction.

\subsection{Maximum Diaconis--Ylvisaker penalised likelihood}

To address nonexistence of ML estimates, this paper studies maximum
Diaconis--Ylvisaker penalised likelihood (MDYPL) (see also \citealt{sterzinger+kosmidis:2026}, \citealt{rigon+aliverti:2023}), 
which maximises $
\ell(\theta,\bbeta;\by,\bX)+\log p(\theta,\bbeta;\bX)
$ using the logarithm of the Diaconis--Ylvisaker prior \citep{diaconis+ylvisaker:1979},
$$
	\log p(\theta, \bbeta; \bX) = \frac{1 - \alpha}{\alpha} \sum_{i = 1}^{n} \left\{ \rho'\left(\theta_P + \bx_i^\top \bbeta_P\right)\left( \theta + \bx_i^\top \bbeta \right) - \rho\left(\theta + \bx_i^\top \bbeta \right) \right\} + C \,, 
$$
where $\alpha\in(0,1)$ and $(\theta_P,\bbeta_P)$ are prior hyperparameters
(see \citealt{diaconis+ylvisaker:1979} for the underlying Bayesian
construction and \citealt{rigon+aliverti:2023} for related bias-reduction
developments) and $C$ is a normalising constant that does not depend on $\theta, \bbeta$.  
From the penalised-likelihood viewpoint adopted here, the
prior parameters $\theta_P, \bbeta_P$ specify the 
target towards which the likelihood fit is regularised, while $\alpha$
controls the strength of regularisation. 
\citet{sterzinger+kosmidis:2026} study the $\bbeta_P = \bb{0}_p$ case, which 
regularises the slope only toward the origin. 
Regularisation towards zero
is natural since ML estimates are typically inflated away from zero (cf. \citealt[Section~8]{cordeiro+mccullagh:1991} and \citealt{sterzinger+kosmidis:2026}, \citealt{zhao+etal:2022}, \citealt{sur+candes:2019}), whilst
allowing
$\bbeta_P\neq\bb0_p$ instead introduces a second slope direction that may
encode prior information. When this direction is aligned with the true
slope, the regularisation target can contain part of the true signal rather
than merely shrinking it towards zero.
The MDYPL estimator is attractive since it always exists
and is unique when the design matrix has full rank
(\citealt[Theorem~1]{rigon+aliverti:2023}), and can be fitted using standard
ML routines with the pseudo-responses $\Ytil_i
    =
    \alpha Y_i
    +(1-\alpha)
    \rho'\left(\theta_P+\bx_i^\top\bbeta_P\right)$. 

\subsection{$p/n \to \kappa \in (0,1)$}

\begin{figure}[ht]
	\centering
	\begin{subfigure}[t]{.32\textwidth}
			\includegraphics[width=\textwidth]{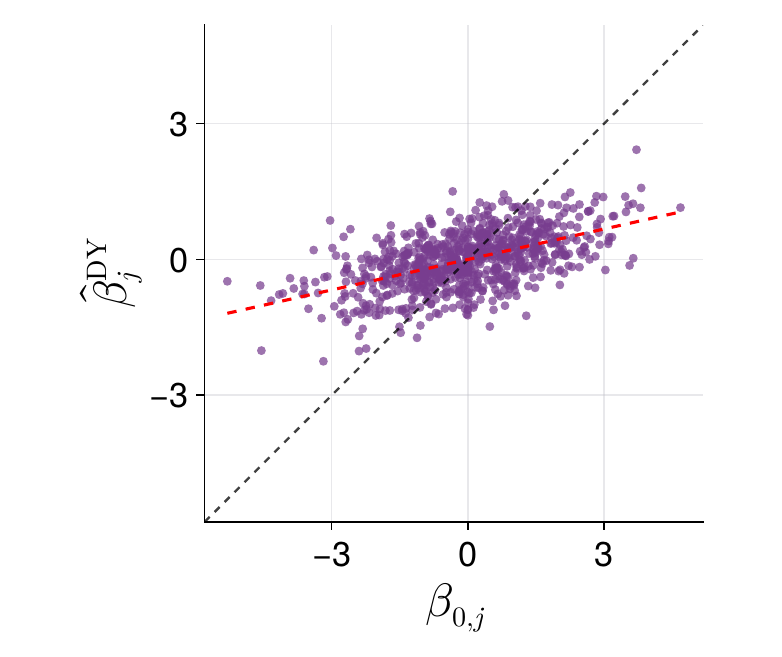}
			\caption{Estimated coordinates vs signal}
			\label{fig:fixed-p-mdypl-decomp}
	\end{subfigure}
	\begin{subfigure}[t]{.32\textwidth}
			\includegraphics[width=\textwidth]{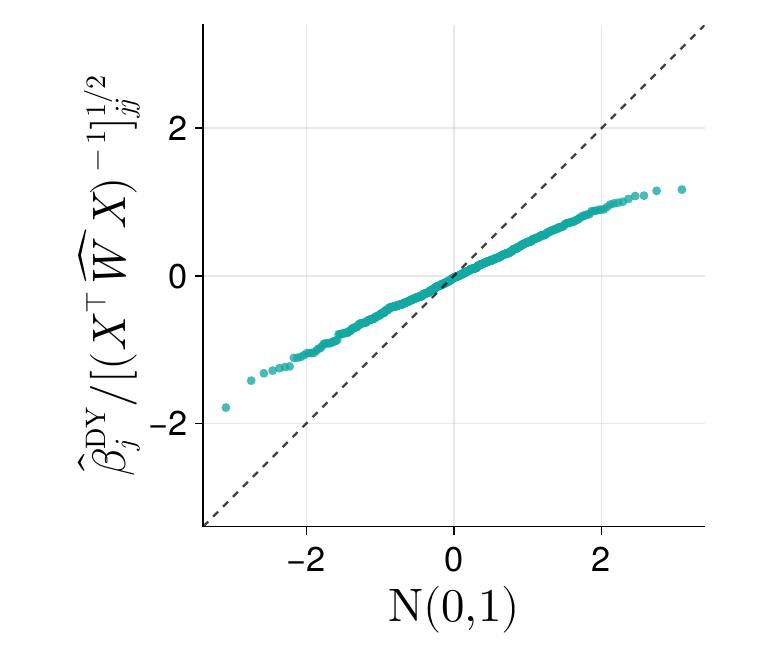}
			\caption{Normal Q--Q plot of classical $Z$-statistics}
			\label{fig:fixed-p-mdypl-qq}
	\end{subfigure}
	\begin{subfigure}[t]{.32\textwidth}
			\includegraphics[width=\textwidth]{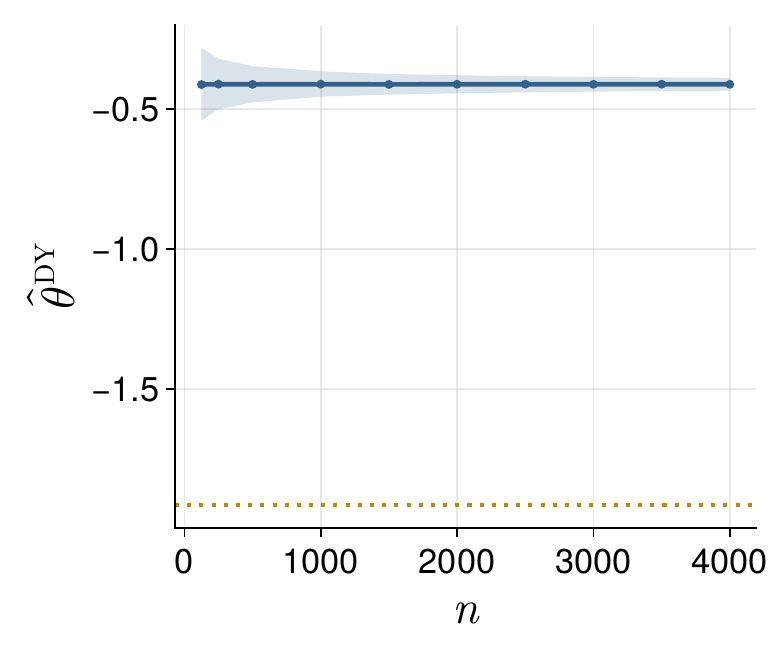}
			\caption{Fitted intercept vs target}
			\label{fig:fixed-p-mdypl-thetady}
	\end{subfigure}
\caption{Illustration of fixed-$p$ asymptotic theory breakdown in the proportional asymptotic regime of 
\eqref{eq:proportional-regime-intro}.
Panel~(a) compares the fitted slope coordinates $\betady_j$ with the
signal
$\bbeta_{0,j}$.
The black dashed line is the identity line and the red dashed line is the
least-squares fit.
Panel~(b) is a normal Q--Q plot of the classical $Z$-statistics 
$Z_j = \betady_j/[(\bX^\top  \widehat{\bb{W}} \bX)^{-1}]^{1/2}_{jj}$, where $\widehat{\bb{W}}$ 
is the diagonal matrix with $i$th entry $\rho''(\bx_i^\top \betady)$. 
Panel~(c) shows the Monte Carlo mean of $\thetady$ and its 
$10$th--$90$th percentile bands. The gold dotted line marks the target $\theta_0$.}
\label{fig:fixed-p} 
\end{figure}

\noindent For fixed $p$, and as $n\to\infty$, standard likelihood theory gives
consistency, asymptotic normality and Cram\'er--Rao efficiency
(see, for example, \citealt{GOURIEROUX198183}),
and under the penalty scaling of
\citet[Section~3.3]{rigon+aliverti:2023} MDYPL 
attains the same optimal limiting distribution of the ML estimator.
In modern datasets, however, 
where the number of features is often comparable to the number of observations, 
classical asymptotic theory with $p/n\to0$ can provide a poor approximation to estimator behaviour.
We thus study the regime
\begin{equation}
\label{eq:proportional-regime-intro}
    p / n \to\kappa\in(0,1),
    \quad
    \bx_i\sim\mathrm N\left(\bb 0_p,p^{-1}\bb I_p\right), 
    \quad
    \theta_0\overset{\mathrm{a.s.}}{\longrightarrow}\theta_0^*,
    \quad
    \var(\bx_i^\top\bbeta_0)\overset{\mathrm{a.s.}}{\longrightarrow}\gamma^2\,,
\end{equation}
which keeps the linear predictor nondegenerate while allowing the number of
parameters to grow at the same order as the sample size.
In this regime, logistic regression behaves fundamentally differently: The
ML estimator may cease to exist \citep{candes+sur:2020} and, when it exists, its sampling
distribution, inferential calibration and predictive behaviour are not
described by classical fixed-$p$ asymptotics theory \citep{sur+candes:2019,sur+et+al:2019,zhao+etal:2022}. An illustration of this breakdown is given in Figure~\ref{fig:fixed-p}. 


\subsection{Related work}
\label{subsec:related}
Logistic regression under the proportional asymptotic regime of
\eqref{eq:proportional-regime-intro} was spearheaded by the twin papers of 
\citet{candes+sur:2020} and \citet{sur+candes:2019}, and the literature has
since developed along both ML and regularised likelihood
lines. 
\citet{candes+sur:2020} establish a sharp phase transition  curve in the $(\kappa,\gamma)$-plane for the existence
of the ML estimator: It exists with probability tending to one
below said curve, i.e. 
$\kappa<h_{\mathrm{\tiny MLE}}(\theta_0^*,\gamma)$, 
and fails to exist with probability tending to one when the inequality is
reversed.
Within the existence region, \citet{sur+candes:2019} characterise the ML
estimator under $\mathrm{N}(\bb{0}_p, n^{-1}\bb I_p)$ covariates and without a fitted
intercept. Their work shows non-vanishing aggregate bias, aggregate variance
inflation and nonclassical likelihood-ratio behaviour, and yields
corresponding corrections for estimation and inference.
\citet{zhao+etal:2022} extend the ML theory to Gaussian
designs with arbitrary covariance and derive corrected coordinatewise
Gaussian test statistics and confidence intervals with nominal asymptotic
calibration. Their formal
results allow only an asymptotically negligible intercept, while the joint
limit of a non-negligible fitted intercept and the slopes is stated as a
conjecture.
For regularised logistic regression with separable convex penalties,
\citet{salehi+et+al:2019} derive an exact asymptotic characterisation of the
estimator through a finite-dimensional system of scalar equations. Their
formulation allows a general convex regulariser, but does not include a
jointly fitted intercept or the second design-aligned signal generated here
by a nonzero MDYPL prior slope, and does not provide results for hypothesis testing. 
\citet{sterzinger+kosmidis:2026} obtain a proportional-asymptotic characterisation of MDYPL logistic regression with
a zero prior slope and without fitted intercept. They use this
characterisation to study adjusted estimation and inference, and conjecture
extensions involving a non-negligible fitted intercept and designs beyond
the Gaussian setting. So, both the ML theory of
\citet{zhao+etal:2022} and the MDYPL theory of
\citet{sterzinger+kosmidis:2026} leave the joint fitted-intercept and slope
limit outside their formal results. 
The closest two-signal statistical precursor is
\citet{li+huang:2024}. In their informative-auxiliary-data setting, the
estimator has a deterministic component in the span of the target and
auxiliary signals. Their formal model, however, imposes a restricted
target--auxiliary correlation structure. The present setting instead permits
an arbitrary positive-semidefinite limiting Gram matrix for
$(\bbeta_0,\bbeta_P)$, including singular limits corresponding to collinear
or vanishing signal directions.

At the methodological level, the CGMT of
\citet{thrampoulidis+etal:2018} provides the main comparison principle in the asymptotic analysis of MDYPL. 
For logistic regression, the key signal--orthogonal decomposition introduced by \citet{salehi+et+al:2019} provides, 
after conditioning on the signal-aligned projection and the randomness generating the responses, 
the independence structure required for the CGMT: the responses are fixed while the signal-orthogonal 
block remains an independent Gaussian matrix.
In the present problem, the conditioning space is
$\vspan\{\bbeta_0,\bbeta_P\}$ and the fitted intercept remains an explicit
optimisation variable. As written, however, the analysis of
\citet{salehi+et+al:2019} does not complete the compact containment,
conditional CGMT comparison and restricted-value localisation needed to
transfer the proposed scalar solution to the original estimator.
The construction of \citet{li+huang:2024} extends this geometry to two
signals and introduces restricted auxiliary optimisation problems intended
to localise the estimator. As written, however, the enlargement of the
projected-ball constraint and the removal of compact scalar domains are not
shown to preserve the optimisation value, and the strict restricted AO gaps
needed for optimiser transfer are not established. The present proof
supplies these missing estimator-level localisation steps.

\subsection{Contribution}

\begin{figure}[ht]
	\centering
	\begin{subfigure}[t]{.32\textwidth}
			\includegraphics[width=\textwidth]{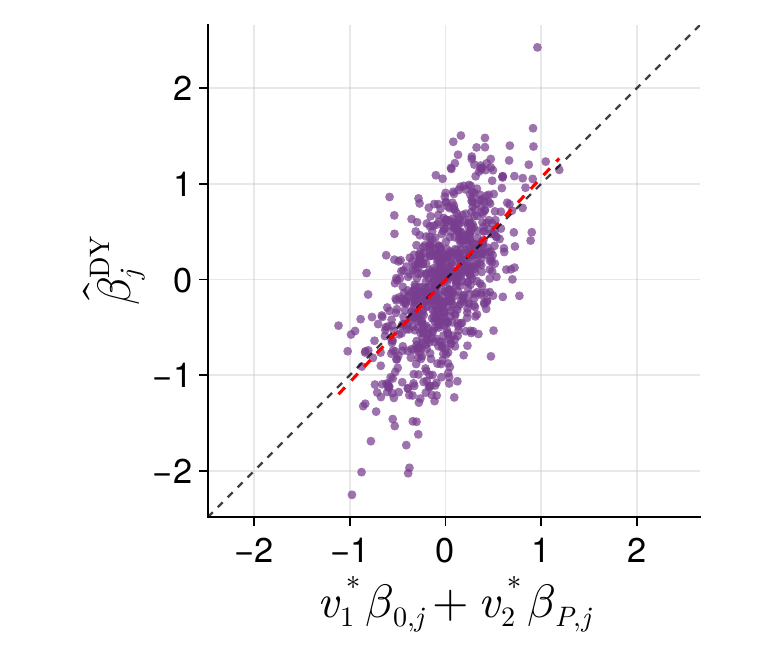}
			\caption{Signal-aligned slope centre}
			\label{fig:mdypl-decomp} 
	\end{subfigure}
	\begin{subfigure}[t]{.32\textwidth}
			\includegraphics[width=\textwidth]{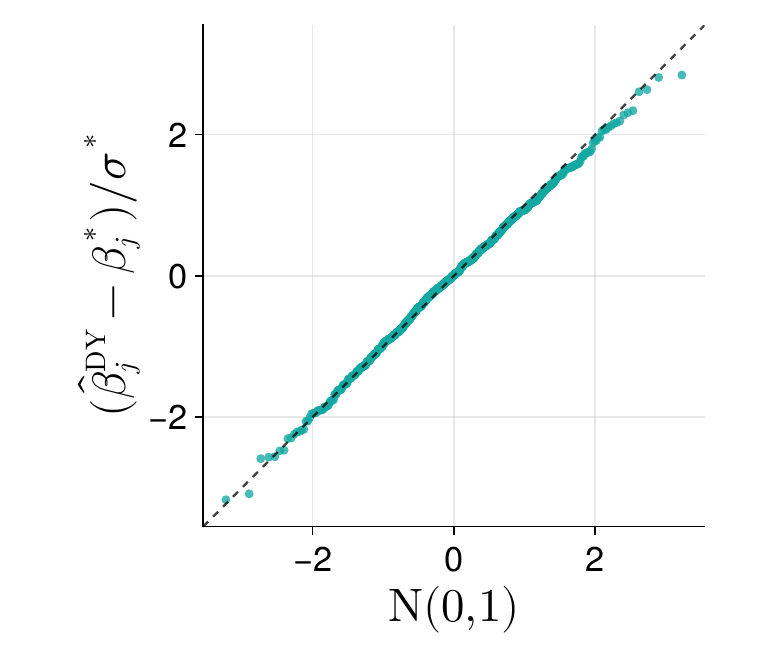}
			\caption{Normal Q--Q plot of standardised slope residuals}
			\label{fig:mdypl-qq} 
	\end{subfigure}
	\begin{subfigure}[t]{.32\textwidth}
			\includegraphics[width=\textwidth]{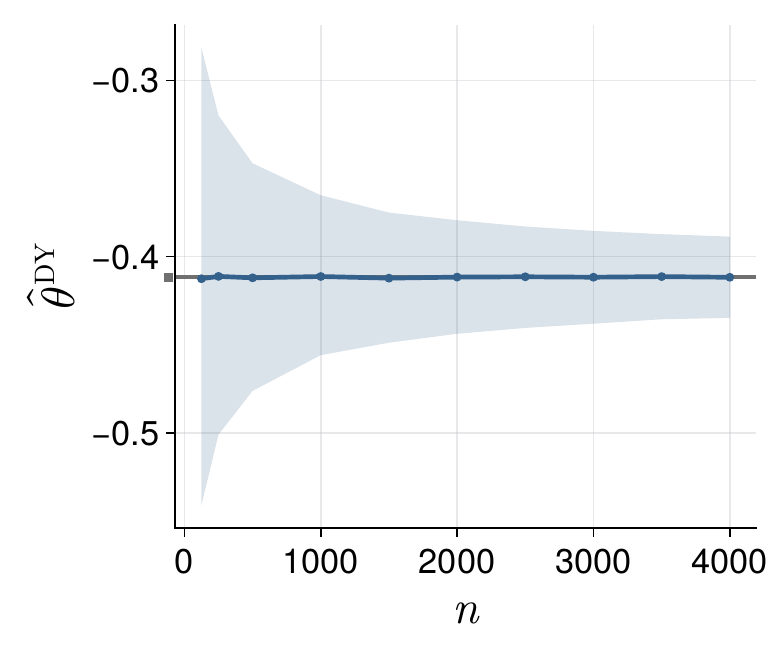}
			\caption{Fitted intercept convergence}
			\label{fig:mdypl-thetady} 
	\end{subfigure}
\caption{Finite-sample illustration of the estimator-level limit in
Theorem~\ref{thm:mdypl-convergence}.
Panel~(a) compares the fitted slope coordinates $\betady_j$ with their
asymptotic centres
$\bbeta_j^*=v_1^*\bbeta_{0,j}+v_2^*\bbeta_{P,j}$.
The grey dashed line is the identity line and the red dashed line is the
least-squares fit.
Panel~(b) is a normal Q--Q plot of the standardised residual coordinates
$(\betady_j-\bbeta_j^*)/\sigma^*$.
Panel~(c) shows the Monte Carlo mean of $\thetady$ and its empirical
$10$th--$90$th percentile band as $n$ varies.
The solid grey line indicates the proportional-limit value $\theta^*$.}
\label{fig:overview} 
\end{figure}

\noindent The main statistical contribution is an estimator-level proportional-limit
law for MDYPL allowing a nonzero prior slope and a fitted
intercept. In particular, we show that there exist deterministic scalars
$\theta^*$, $v_1^*$, $v_2^*$ and $\sigma^*$ and a signal-aligned
vector
$
    \bbeta^*
    =
    v_1^*\bbeta_0+v_2^*\bbeta_P
    \in
    \vspan\{\bbeta_0,\bbeta_P\}
$
such that, informally,
\begin{equation}
\label{eq:informal-result}
    \thetady
    \overset{\mathrm{a.s.}}{\longrightarrow}
    \theta^*,
    \quad
    \betady
    \approx
    \bbeta^*+\sigma^* \sqrt{p} \bb d \,,
\end{equation}
where $\bb d$ is the noise direction that is uniformly distributed on the unit sphere in the signal-orthogonal subspace. 
Our main result, Theorem~\ref{thm:mdypl-convergence}, formalises this approximationas the almost-sure Wasserstein-$2$ convergence of the empirical law of
$(\betady_j-\bbeta_j^*,\bbeta_{0,j},\bbeta_{P,j})$ to the law of
$(\sigma^*G,\bar\beta_0,\bar\beta_P)$, with
$G\sim\mathrm N(0,1)$ independent of $(\bar\beta_0,\bar\beta_P)$.
Figure~\ref{fig:overview} provides an illustration of the three
components of \eqref{eq:informal-result} in Setting~`$\mathrm{M}$' of Section~\ref{subsec:simul}: deterministic centring of the slope
coordinates, Gaussian coordinate fluctuations, and convergence of the fitted
intercept to $\theta^*$ rather than the data-generating value $\theta_0$.

The vector $\bbeta^*$ is defined through the Moore--Penrose
pseudoinverse of the limiting signal Gram matrix, so the result permits any
positive-semidefinite Gram-limit for $(\bbeta_0,\bbeta_P)$, including
collinear or vanishing signal directions. Relative to
\citet{salehi+et+al:2019}, the theorem adds both a jointly fitted intercept
and the second design-aligned signal generated by the nonzero prior
predictor. Relative to \citet{li+huang:2024}, it additionally allows arbitrary positive-semidefinite Gram-limits for the target and source signals. To our knowledge, this is the
first rigorous estimator-level proportional limit for MDYPL with both a
nonzero prior slope and a non-negligible fitted intercept.
The estimator law yields the joint out-of-sample limit of the true and fitted
scores: After removing the prior-aligned slope component and recentering at
the true intercept, the oracle-adjusted score is the true logit observed
through independent additive Gaussian noise. This canonical experiment gives
the optimal score-based classifier and the calibrated and Bregman-optimal
class probability, while its squared correlation with the true logit orders
optimally thresholded classification performance. At the coordinate level, the estimator
decomposition gives oracle-adjusted fixed-block $Z$-statistics under $\mathrm{N}(\bb{0}_p, p^{-1}\bb I_p)$ covariates and, after affine centering and whitening, under Gaussian
covariates with arbitrary deterministic mean and positive-definite
covariance.
For feasible implementation, we introduce a response-moment procedure that
consistently estimates $(\theta_0^*,\gamma^2,\varphi)$ before the remaining
state parameters are computed from plug-in oracle equations, without the
computational and phase-boundary limitations of ProbeFrontier \citep{sur+candes:2019} or the
identification problems of SLOE \citep{yadlowsky+etal:2021} for MDYPL (see \citealt[Figure~8]{sterzinger+kosmidis:2026}).
We also identify the candidate penalised
likelihood-ratio rescaling and missing profile arguments following \citet{sur:2019}, and outline a
universality route for independent-entry subgaussian designs based on
optimal-value comparisons following \citet{montanari+saeed:2022}. 

Our main methodological contribution is a complete conditional-CGMT proof
of the estimator law with fitted intercept and general prior--signal
geometry. We formalise the conditioning and containment arguments that
isolate the independent signal-orthogonal Gaussian block on a compact CGMT
domain. We then give a complete PO--AO reduction, including boundary
arguments not established in comparable analyses and the geometric arguments
required for the general signal--prior geometry. Finally, uniqueness of the
limiting saddle yields strict restricted-value gaps that localise the
estimator, with new restricted PO--AO comparisons handling the fitted
intercept through the profiled objective.



\section{Setup} \label{sec:setup}

Consider covariates $\bx_i \in \Re^p$ that are i.i.d. $\mathrm{N}(\bb{0}_p, p^{-1}\bb{I}_p)$ for $i = 1,\ldots,n$, where $p/n \to \kappa \in (0,1)$ as $n\to \infty$.
Further, let $\bbeta_0, \bbeta_{P} \in \Re^p$, $\theta_P, \theta_0 \in \Re$ be potentially random sequences such that:
\begin{enumerate}[label = (\roman{*})]
	\item \label{ass:betas_conv} For
	$$
		\bb \Gamma_p =  p^{-1} 
        \begin{bmatrix}
            \vnorm{\bbeta_0}_2^2 & \langle \bbeta_0, \bbeta_P \rangle \\ 
            \langle \bbeta_0, \bbeta_P \rangle & \vnorm{\bbeta_P}_2^2 
        \end{bmatrix}, 
        \quad \bb \Gamma = \begin{bmatrix}
			\gamma^2 & \varphi \\
			\varphi & \delta^2
		\end{bmatrix}\,,
	$$
	and $\bb \Gamma$ deterministic and positive semidefinite, there exist universal constants $c_{\Gamma}, C_{\Gamma} >0$ and an $n_0$ such that for all $n > n_0$, and all $t > 0$,
	$$
		\Pr\left(\mnorms{\bb \Gamma_p - \bb \Gamma}_2 > t \right) \leq C_{\Gamma} \exp\{-c_{\Gamma}n \min\{t^2,t\} \} \,.
	$$
    \item \label{ass:betas_dist} The joint empirical distribution of $(\bbeta_0,\bbeta_P)$ converges almost surely to a deterministic distribution function $\pi_{\{\bar{\beta}_0, \bar{\beta}_P\}}(b_0, b_P)$, 
    in Wasserstein $2$-distance, i.e. almost surely
	$$
		\pi_p = \frac{1}{p}\sum_{j = 1}^{p} \delta_{(\bbeta_{0,j}, \bbeta_{P,j})} \overset{W_2}{\longrightarrow} \pi_{\{\bar{\beta}_0, \bar{\beta}_P\}}\,.
	$$
	 \item \label{ass:theta-conv} $\tnod, \theta_P$ have deterministic scalar limits $\theta_0^*, \theta_P^*$ such that for all $n > n_0$ and $t > 0$,
	$$
		\Pr \left(\vnorm{\begin{bmatrix}
			\theta_0 \\ \theta_P \end{bmatrix} - \begin{bmatrix}
				\theta_0^* \\ \theta_P^*
			\end{bmatrix}}_{\infty} > t
			\right) \leq C_{\Theta} \exp \{-c_{\Theta}nt^2\} \,,
	$$
	for some universal constants $c_{\Theta}, C_{\Theta}$ and a $n_0 \in \mathbb{N}$.
    \item $(\theta_0, \theta_P, \bbeta_0, \bbeta_P)$ are jointly independent of the covariates $\bx_i$ ($i = 1, \ldots, n$).
\end{enumerate}
For $\rho(x)$ as in \eqref{eq:logistic-model-intro}, $\rho'(x) = \mathrm d \rho(x) / \mathrm dx$ and $\alpha \in (0,1)$, define the binary-
and pseudo-responses, respectively, by
\begin{equation}
    \label{eq:mdypl}
	Y_i
	=
	\mathds{1}\{ \varepsilon_i < \rho'(\tnod + \bx_i^\top \bbeta_{0})\}, \quad \Ytil_i = \alpha Y_i + (1-\alpha) \rho'(\theta_P + \bx_i^\top\bbeta_{P})\,,
\end{equation}
where $\varepsilon_i \sim \mathrm{Unif}([0,1])$ are i.i.d. uniform random variables on the unit interval, independent of all other randomness.
Define the MDYPL optimisation problem
\begin{equation}
    \label{eq:og}
	(\thetady, \betady) = \underset{\substack{\bbeta \in \Re^p, \, \theta \in \Re}}{\arg\max} \, \ell(\theta, \bbeta ; \bYtil, \bX)
	\,,
\end{equation}
where $\bYtil = (\widetilde{Y}_1, \ldots, \widetilde{Y}_n)^\top$ with $\Ytil_i$ defined as in \eqref{eq:mdypl}.
If $\alpha = 0$ in \eqref{eq:mdypl}, then $\betady = \bbeta_P$, $\thetady = \theta_P$. We exclude this trivial case from our analysis.
If, on the other hand, $\alpha = 1$, then $\thetady, \betady$ are the ML estimates studied in \citet{sur+candes:2019} and \citet{zhao+etal:2022} in this setting.
While \citet{zhao+etal:2022} conjecture behaviour exactly to what is proven here for MDYPL and the results coincide at $\alpha = 1$, we do not formally prove the case $\alpha = 1$ here.
By \citet[][Theorem~1]{rigon+aliverti:2023}, the estimates $(\thetady, \betady)$ of \eqref{eq:og} exist and are unique whenever the matrix $\widetilde{\bX} = [\bb 1,  \bX]$ has full column rank.
For well-definedness, whenever $\widetilde{\bX}$ does not have full column rank, set $\thetady = 0$, $\betady = \bb 0_{p }$.

\section{Proportional asymptotics of the MDYPL estimator}
\label{sec:results}

\subsection{Scalar characterisation of the limit}
\label{sec:limiting-ao-focs}

The quantities $\theta^*$, $\sigma^*$ and $\bbeta^*$ appearing in the
informal representation \eqref{eq:informal-result} are determined by the
unique saddle point of the following deterministic finite-dimensional
min--max problem:
\begin{equation}
\label{eq:limit_ao_natural}
\begin{aligned}
    \underset{\substack{
        \sigma\geq0\\
        \theta\in\Re\\
        t>0\\
        \bu\in\range(\bb\Gamma)
    }}{\min}
    \underset{r\in(0,1)}{\sup} \, 
    \left\{
        \expect
        \left[
            m_\rho
            \left(
                \theta+\bz^\top\bu+\sigma G
                +\frac{t}{r}\Ytil,
                \frac{t}{r}
            \right)
        \right]
        -
        \sigma r\sqrt\kappa
        +
        \frac{rt}{2}
        -
        \theta\expect[\Ytil]
        -
        \expect[\Ytil\bz^\top\bu]
        -
        \frac{t}{2r}\expect[\Ytil^2]
    \right\}\,.
\end{aligned}
\end{equation}
Here $\bz\sim\mathrm N(\bb0_2,\bb I_2)$, $\bb Q
    =[Q_1, Q_2]^\top 
    =
    \bb\Gamma^{1/2}\bz
    \sim
    \mathrm N(\bb0_2,\bb\Gamma)$, $G\sim\mathrm N(0,1)$, $
    \varepsilon\sim\mathrm{Unif}([0,1])$ 
are mutually independent, and $\Ytil
    =
    \alpha
    \mathds{1}
    \left\{
        \varepsilon
        <
        \rho'(\theta_0^*+Q_1)
    \right\}
    +
    (1-\alpha)
    \rho'(\theta_P^*+Q_2)$
is the limiting MDYPL pseudo-response. Further,
$
    m_\rho(x,\lambda)
    =
    \min_{\eta\in\Re} \, 
    \{
        \rho(\eta)
        +
        (\eta-x)^2 / (2 \lambda)
    \}
$
is the Moreau envelope of $\rho$. 
Proposition~\ref{prop:limiting_AO_properties} shows that
\eqref{eq:limit_ao_natural} attains a finite value at an unique
interior saddle point that is equivalent to a solution of the first-order conditions in \eqref{eq:FOCs}. Specifically, on
setting
$
    \lambda^*
    =
    {t^*} / {r^*}$, 
    $
    \xi^*
    =
    \theta^*
    +
    \bz^\top\bu^*
    +
    \sigma^*G
$, 
$(\sigma^*,\lambda^*,\bu^*,\theta^*)$ is the unique solution in
$(0,\infty)^2\times\range(\bb\Gamma)\times\Re$ of the population
equations
\begin{equation}
\label{eq:FOCs}
\begin{aligned}
    \expect
    \left[
        \rho'
        \left(
            \prox{\lambda^*\rho}
            {\xi^*+\lambda^*\Ytil}
        \right)
        -
        \Ytil
    \right]
    &=
    0,
    \\
    \bb P_{\bb \Gamma}
    \expect
    \left[
        \bz
        \left\{
            \rho'
            \left(
                \prox{\lambda^*\rho}
                {\xi^*+\lambda^*\Ytil}
            \right)
            -
            \Ytil
        \right\}
    \right]
    &=
    \bb0,
    \\
    \lambda^*
    \expect
    \left[
        \frac{
            \rho''
            \left(
                \prox{\lambda^*\rho}
                {\xi^*+\lambda^*\Ytil}
            \right)
        }{
            1
            +
            \lambda^*
            \rho''
            \left(
                \prox{\lambda^*\rho}
                {\xi^*+\lambda^*\Ytil}
            \right)
        }
    \right]
    &=
    \kappa,
    \\
     \expect
    \left[
        \left\{
            \rho'
            \left(
                \prox{\lambda^*\rho}
                {\xi^*+\lambda^*\Ytil}
            \right)
            -
            \Ytil
        \right\}^2
    \right]&= \frac{(\sigma^*)^2\kappa}{(\lambda^*)^2}\,.
\end{aligned}
\end{equation}
Here $
    \prox{\lambda\rho}{x}
    =
    {\arg\min}_{\eta\in\Re} \, 
    \{
        \rho(\eta)
        +
        (\eta-x)^2 / {(2\lambda)}
    \}
$ denotes the proximal operator of $\lambda\rho$, and
$\bb P_{\bb \Gamma}$ is the $2 \times 2$ matrix 
that projects onto the range of $\bb\Gamma$. 
Intuitively, the first two equations in \eqref{eq:FOCs} balance the population score in
the intercept and identifiable signal directions. The third matches the
effective curvature to the aspect ratio $\kappa$, while the fourth
determines the signal-orthogonal variance. The projection in the second
equation permits singular $\bb\Gamma$, including collinear or vanishing
signal directions.

\subsection{Proportional limit of the MDYPL estimator}
\label{sec:prop-lim}


The scalar problem of \eqref{eq:limit_ao_natural} determines the parameters in the limiting
decomposition. To express its signal-aligned component, define
\begin{equation}
\label{eq:limiting-alignment}
    \bv^*
    =
    \bb\Gamma^{+1/2}\bu^*,
    \quad 
    \bb B
    =
    [\bbeta_0,\bbeta_P],
    \quad
    \bbeta^*
    =
    \bb B\bv^*\,.
\end{equation}
Here $\bb\Gamma^{+1/2}$ denotes the Moore--Penrose pseudoinverse of
$\bb\Gamma^{1/2}$ (cf. \citealt[Chapter~2]{magnus+neudecker:2019}), so $\bv^*$ gives the canonical alignment coefficients
also when $\bb\Gamma$ is singular. 
\begin{theorem}
    \label{thm:mdypl-convergence}
    Assume the conditions of Section~\ref{sec:setup}, and let
    $\theta^*$, $\sigma^*$ and $\bbeta^*$ be as in
    \eqref{eq:limit_ao_natural} and \eqref{eq:limiting-alignment}. Then the
    MDYPL estimator in \eqref{eq:og} satisfies the following.

    \begin{enumerate}[label=(\roman*)]
        \item
        $$
            \thetady
            \overset{\mathrm{a.s.}}{\longrightarrow}
            \theta^*\,.
        $$

        \item For every pseudo-Lipschitz function
        $\psi:\Re^3\to\Re$ of order two
        ,
        \begin{equation}
        \label{eq:test_fun_conc}
            \frac1p
            \sum_{j=1}^p
            \psi
            \left(
                \betady_j-\bbeta_j^*,
                \bbeta_{0,j},
                \bbeta_{P,j}
            \right)
            \overset{\mathrm{a.s.}}{\longrightarrow}
            \expect
            \left[
                \psi(\sigma^*G,\bar\beta_0,\bar\beta_P)
            \right]\,,
        \end{equation}
        where
        $
            (\bar\beta_0,\bar\beta_P)
            \sim
            \pi_{\{\bar\beta_0,\bar\beta_P\}}
        $
        and $G\sim\mathrm N(0,1)$ is independent of
        $(\bar\beta_0,\bar\beta_P)$.
    \end{enumerate}
\end{theorem}
\noindent A selection of informative aggregate limits, that follow directly from
Theorem~\ref{thm:mdypl-convergence}, are given in Table~\ref{tab:estimation-consequences}.
For the examples in Table~\ref{tab:estimation-consequences}, write
$
    \bb b
    =
    [b_0,b_P]^\top$, $\bb e_1
    =
    [1,0]^\top
$.
The first row is understood componentwise, corresponding to the two
scalar test functions $ab_0$ and $ab_P$.
Rows one and two show that the estimation error around $\bbeta^*$ is
asymptotically orthogonal to both signal directions and has Euclidean norm
$\sigma^*\sqrt p$. 
The remaining rows give the limiting squared norm of
the fitted slope, its mean squared error relative to the true slope, and the
mean squared error after oracle removal and rescaling of the prior-aligned
component when $v_1^* \neq 0$. 
Together, these statements give the aggregate form of the informal
decomposition in \eqref{eq:informal-result}.
\begin{table}[ht]
\centering
\caption{Selected aggregate consequences of
Theorem~\ref{thm:mdypl-convergence}, obtained from the displayed choices of
$\psi$.}
\label{tab:estimation-consequences}
\small
\setlength{\tabcolsep}{5pt}
\begin{tabularx}{\textwidth}{
    @{}
    >{\raggedright\arraybackslash}p{0.27\textwidth}
    >{\raggedright\arraybackslash}X
    >{\raggedright\arraybackslash}p{0.29\textwidth}
    @{}
}
\toprule
$\psi(a,b_0,b_P)$
&
Statistic
&
Almost-sure limit
\\
\midrule
$a\bb b$
&
$\bb B^\top(\betady-\bbeta^*)/p$
&
$\bb0_2$
\\
$a^2$
&
$\vnorm{\betady-\bbeta^*}_2^2/p$
&
$(\sigma^*)^2$
\\
$\left(a+\bv^{*\top}\bb b\right)^2$
&
$\vnorm{\betady}_2^2/p$
&
$(\sigma^*)^2+\bv^{*\top}\bb\Gamma\bv^*$
\\
$\left\{a+(\bv^*-\bb e_1)^\top\bb b\right\}^2$
&
$\vnorm{\betady-\bbeta_0}_2^2/p$
&
$(\sigma^*)^2+
(\bv^*-\bb e_1)^\top
\bb\Gamma
(\bv^*-\bb e_1)$
\\
$a^2/(v_1^*)^2$
&
$\vnorm{(\betady-v_2^*\bbeta_P)/v_1^*-\bbeta_0}_2^2/p$
&
$(\sigma^*)^2/(v_1^*)^2$
\\
\bottomrule
\end{tabularx}
\end{table}
A useful benchmark is provided by the completely informative prior
$(\theta_P,\bbeta_P)=(\theta_0,\bbeta_0)$. In this case, $\alpha=0$ is
optimal and MDYPL recovers the true parameter exactly. In the canonical
rank-one representation of $\vspan\{\bbeta_0,\bbeta_P\}$, the corresponding
boundary solution is
$
    \theta^*=\theta_0^*$, $
    \sigma^*=0$, $
    v_1^*=v_2^*=\frac12$.
These values are also obtained formally from the population equations as
$\alpha\to 0$ from above. By contrast, for the null prior slope
$\bbeta_P=\bb0_p$, Proposition~\ref{prop:limiting_AO_properties} gives
$\sigma_\alpha^*>0$ for every fixed $\alpha\in(0,1)$, and therefore every
such fit has strictly positive limiting slope mean squared error. 
Thus, a completely informative prior strictly improves on optimally tuned
shrinkage towards $\bb0_p$ at any fixed $\alpha$. For an informative but noisy prior, one would
expect this advantage to persist while the prior noise remains sufficiently
small, with a transition depending on the aspect ratio $\kappa$ and the
signal--prior geometry $\bb\Gamma$.

\subsection{Proof strategy}
\label{sec:proof-strategy-main}

Theorem~\ref{thm:mdypl-convergence} is established through a conditional-CGMT
analysis and its proof may be of independent interest: it retains a jointly
fitted intercept throughout, permits arbitrary alignment between the true
and prior slopes, including singular Gram-limits, and provides a complete
route from the scalar AO to the original estimator. We therefore briefly
outline its main steps.
Central to our derivations is the CGMT of
\citet[Theorem~VI.1]{thrampoulidis+etal:2018} that associates the PO and AO problems
\begin{align}
    \Phi(\bb G)
    &=
    \underset{\bw\in\mathcal S_{\bw}}{\min}
    \underset{\bu\in\mathcal S_{\bu}}{\max}
    \left\{
        \bu^\top\bb G\bw+\psi_0(\bw,\bu)
    \right\} \,,
    \tag{PO}\label{eq:cgmt-pair-main}
    \\
    \phi(\bg,\bh)
    &=
    \underset{\bw\in\mathcal S_{\bw}}{\min}
    \underset{\bu\in\mathcal S_{\bu}}{\max}
    \left\{
        \vnorm{\bw}_2\bg^\top\bu
        +
        \vnorm{\bu}_2\bh^\top\bw
        +
        \psi_0(\bw,\bu)
    \right\} \,.
    \tag{AO}
\end{align}
where $\bb G$, $\bg$ and $\bh$ have independent standard Gaussian entries.
CGMT gives the PO--AO value link:
\begin{equation}
\label{eq:cgmt-value-links-main}
\begin{aligned}
    \Pr(\Phi(\bb G)<c)
    &\leq
    2\Pr(\phi(\bg,\bh)\leq c),
    \quad
    \Pr(\Phi(\bb G)>c)
    \leq
    2\Pr(\phi(\bg,\bh)\geq c) \,. 
\end{aligned}
\end{equation}
The lower-tail comparison holds when the constraint sets are compact and
$\psi_0$ is continuous and the upper-tail comparison holds if, in addition,
the constraint sets are convex and $\psi_0$ is convex--concave.
We translate the MDYPL optimisation into a suitable primary optimisation and
use these PO--AO value links to derive the limiting behaviour of its
optimisers.

\proofstrategyparagraph{Signal decomposition and conditioning}
{par:proof-strategy-signal-decomposition}
To apply CGMT, the Gaussian matrix in the bilinear term must be
independent of the remaining objective and the profile $\psi$ non-random.
Writing $\bbeta=\bb B\bv+\bb E\bw$, where $\bb E$ has orthonormal columns spanning
$\range(\bb B)^\perp$ and $\bX = \bH / \sqrt{p}$ gives $\bH\bbeta=\bH_1\bv+\bH_2\bw$, 
with $\bH_1 = \bH \bb B$, $\bH_2 = \bH \bb E$ and the
variables $Y_i, \widetilde{Y}_i$ depend only on 
$\bH_1$, but not on $\bH_2$. Hence, conditional on
$\mathcal C_n=\sigma(\theta_0,\theta_P,\bb B,\bH_1,\bb\varepsilon)$, the only randomness
stems from $\bb H_2$, which is
standard Gaussian, and on the 
containment event $\{ \vnorm{\betady}_2 / \sqrt p \leq C_\beta\}$, MDYPL is represented exactly by the
primary optimisation
\begin{equation}
\label{eq:mdypl-po-main}
\begin{aligned}
    \Phi(\bH_2)
    =
    \underset{\bw\in W_p}{\min}
    \underset{\blambda\in\mathcal B_1^n}{\max}
        -\frac{\blambda^\top\bH_2\bw}{n\sqrt p}
        +\psi(\blambda)
	, \,
    \psi(\blambda)
    =
    \frac1n
    \underset{(\bv,\theta,\bb\eta)\in\mathcal K}{\min} \,
        \bb1^\top\bb\rho(\bb\eta)-\bYtil^\top\bb\eta
        +\blambda^\top
		\left(
            \bb\eta-\theta\bb1-\frac1{\sqrt p}\bH_1\bv
        \right)
	\,,
\end{aligned}
\end{equation}
where
$\mathcal K=V_p\times[-C_\theta,C_\theta]\times\mathcal B_{C_\eta}^n$, 
$\mathcal B_C^d=\{\bx\in\Re^d:\vnorm{\bx}_2/\sqrt d\le C\}$, 
$V_p, W_p$ are compact sets defined in Section~\ref{sec:po}.
Conditional on
$\mathcal C_n$, the AO associated with \eqref{eq:mdypl-po-main} is
\begin{equation}
\label{eq:mdypl-ao-main}
    \phi(\bg,\bh)
    =
    \underset{\bw\in W_p}{\min}
    \underset{\blambda\in\mathcal B_1^n}{\max} \,
        -
        \frac{1}{n\sqrt p} \{\vnorm{\bw}_2\bg^\top\blambda
            +
            \vnorm{\blambda}_2\bw^\top\bh \} 
        +
        \psi(\blambda) \,,
\end{equation}
where $\bg\in\Re^n$ and $\bh\in\Re^{p-s}$ are independent standard
Gaussian vectors.

\proofstrategyparagraph{AO scalarisation \& limiting AO}
{par:proof-strategy-ao-scalarisation}

The next step is to reduce the high-dimensional AO min--max problem in
\eqref{eq:mdypl-ao-main} to an optimisation over finite-dimensional
coordinates and establish its convergence to the deterministic limiting AO
in \eqref{eq:limit_ao_natural}. The scalarisation follows the general
blueprint of \citet{salehi+et+al:2019}, but formalises a number of
nontrivial steps missing from their argument. The radial boundary $r=0$ is
removed by a joint small-$r$ environment: a projection argument gives a
uniform gap at $\bEtaYtil$, while coordinatewise clipping and an excess-loss
bound control the moving minimiser and yield strictly positive profile
increments on a deterministic interval. Uniform convergence of the
scalarised AO to its limiting counterpart follows largely from uniform
Lipschitz control of the objective and pointwise concentration.
The main difficulty in the convergence argument is that at finite $n$, the signal
coordinate is constrained to $\range(\bb\Gamma_p)$, whereas the limiting AO
is constrained to $\range(\bb\Gamma)$. Although
$\bb\Gamma_p\to\bb\Gamma$, these ranges need not converge when
$\bb\Gamma$ is singular, so uniform convergence of the objectives alone does
not compare the corresponding optimisation values.
We instead compare both problems through the full fixed ball.
Lemma~\ref{lemma:nullspace_no_help_barFn} shows that, for the conditional mean
finite objective, projecting $\bu$ onto $\range(\bb\Gamma_p)$ cannot
increase the objective, so enlarging the finite signal domain to the full
ball does not change its conditional optimisation value.
Lemma~\ref{lemma:nullspace_no_help_F} gives the analogous projection result for
the limiting objective, with strict inequality whenever $\bu$ has a
nonzero component orthogonal to $\range(\bb\Gamma)$.
The value-Lipschitz property then transfers the finite conditional
comparison to the random finite AO up to its uniform fluctuation around the
conditional mean.
Proposition~\ref{prop:uniform_AO_control} combines these
deterministic domain comparisons with uniform concentration on the
ambient domain and yields
$
    \phi(\bg,\bh)
    \overset{\mathrm p}{\longrightarrow}
    \bar\phi
$, 
where $\bar\phi$ is the value of the limiting AO in
\eqref{eq:limit_ao_natural}. 
Proposition~\ref{prop:limiting_AO_properties} then gives attainment, the
unique interior saddle
$(\sigma^*,r^*,\bu^*,\theta^*,t^*)$, the natural and ambient domain
relaxations, and the population equations \eqref{eq:FOCs}.

\proofstrategyparagraph{Value comparisons and optimiser localisation}
{par:proof-strategy-localisation}

The next step is to transfer the unique limiting-AO saddle to every
PO optimiser. For any such optimiser
$(\hat\bw,\hat\bv,\hat\theta,\hat{\bb\eta})$, define
$
    \hat\sigma
    =
    \vnorm{\hat\bw}_2 / \sqrt p$, $
    \hat\bu
    =
    \bb\Gamma_p^{1/2}\hat\bv\,,
$
and, for $\epsilon>0$, consider the localising set
$
    \mathcal S^\epsilon
    =
	\{
        (\bw,\bv,\theta):
        \abs{
            \vnorm{\bw}_2 / {\sqrt p}-\sigma^*
        }
        <
        \epsilon, \, 
        \|
            \bb\Gamma_p^{1/2}\bv-\bu^*
        \|_2
        <
        \epsilon, \, 
        \abs{\theta-\theta^*}
        <
        \epsilon
	\}
$. 
We rule out an optimiser outside $\mathcal S^\epsilon$ by comparing the
unrestricted PO value with values obtained under radial, signal, and
intercept restrictions.
The radial deviation is imposed directly through
$
    \mathcal S_{\bw}^{\epsilon,c}
    =
	\{
        \bw\in W_p:
        \abs{
			\vnorm{\bw}_2 / {\sqrt p}-\sigma^*
        }
        \geq
        \epsilon
	\}
$. 
The radial complement is handled by restricting $\bw$ directly. The
signal--intercept complement is not convex, so
Section~\ref{sec:cgmt-cost-comparisons} fixes a finite $1/2$-net and
covers it by the compact convex pieces $C_{p,j}^\epsilon$,
$j\in\mathcal J$. A signal piece cuts
$V_p\times[-C_\theta,C_\theta]$ by the half-space
$
    \bb a_j^\top
    (
        \bb\Gamma_p^{1/2}\bv-\bu^*
    )
    \geq
    \epsilon / 2
$, 
whereas the two intercept pieces leave the signal coordinate unrestricted
and restrict $\theta$ to $\Theta_+^\epsilon$ or
$\Theta_-^\epsilon$. Restricting the profile to any nonempty
$C_{p,j}^\epsilon$ preserves compactness and convexity and gives the
restricted profile $\psi_j^\epsilon$ and PO value $\Psi_j^\epsilon$.
Let $\Psi_{\bw}^{\epsilon,c}$ denote the radially restricted PO value and
let $\Psi_j^\epsilon$ denote the value for a nonempty signal or intercept
piece. Lemma~\ref{lemma:restricted_piece_scalarisation} represents all
corresponding AOs through the same finite-dimensional objective as the
unrestricted AO, with only the relevant feasible domain restricted.
Consequently, the uniform convergence already established on the fixed
ambient domain, together with the value-Lipschitz property and monotonicity
under domain restriction, controls all restricted values simultaneously.
Each restricted limiting domain excludes the unique limiting saddle.
Together with finiteness of the cover,
Lemma~\ref{lemma:limiting_AO_gap} therefore gives
a common gap $\eta_\epsilon>0$. Uniform AO control transfers these gaps to
the finite restricted problems, and the conditional CGMT comparisons yield,
with exponentially high probability,
\begin{equation}
\label{eq:proof-strategy-po-separation}
    \Psi
    \leq
    \bar\phi+\eta_\epsilon
    <
    \bar\phi+2\eta_\epsilon
    \leq
    \min
	\{
        \Psi_{\bw}^{\epsilon,c},
        \underset{j\in\mathcal J_n^\epsilon}{\min}\,
        \Psi_j^\epsilon
	\}\,.
\end{equation}
Here the upper-tail comparison bounds the unrestricted PO from above, while
the lower-tail comparison bounds each restricted PO from below. An optimiser
outside $\mathcal S^\epsilon$ would be feasible for one of the restricted
problems, forcing its value to be no larger than $\Psi$ and contradicting
\eqref{eq:proof-strategy-po-separation}.
Proposition~\ref{prop:PO_localisation} consequently gives
$
    \hat\sigma
    \overset{\mathrm p}{\longrightarrow}
    \sigma^*$, $
    \hat\bu
    \overset{\mathrm p}{\longrightarrow}
    \bu^*$, $
    \hat\theta
    \overset{\mathrm p}{\longrightarrow}
    \theta^*
$.

\proofstrategyparagraph{Convergence on test functions}
{par:proof-strategy-test-functions}

By Lemma~\ref{lemma:compact_PO_uniqueness}, on $\mathcal E_n$ the compact
PO has a unique primal minimiser, which we denote by
$(\hat\bv,\hat\bw,\hat\theta,\hat{\bb\eta})$and set $\hat\bbeta=\bb B\hat\bv+\bb E\hat\bw$ and $\hat\sigma=\vnorm{\hat\bw}_2 / \sqrt p$. 
Lemma~\ref{lemma:compact_PO_spherical_representation} gives the exact
decomposition $ \hat\bbeta
    =
    \bb B\hat\bv
    +
    \hat\sigma \sqrt p \bb d$, 
where, conditionally on $\mathcal C_n$, the direction $\bb d$ is uniform
on the unit sphere of $\range(\bb B)^\perp$ and independent of the radial
and signal-aligned optimiser coordinates. We replace the noise $\sqrt p \bb d$ with a 
Gaussian surrogate $\bb z_p \sim \mathrm N(\bb 0_p,\bb I_p)$ and show that this replacement 
is negligible on the class of test functions we consider. 
The final proof is then an assembly of four ingredients: (i) containment and full rank, which identify MDYPL with the compact PO, 
(ii) AO concentration and restricted CGMT comparisons, which localise the compact-PO optimiser, 
(iii) Gaussian replacement of the signal-orthogonal slope and its propagation through the pseudo-Lipschitz average, and 
(iv) conditional concentration of this Gaussian empirical average around its expectation, followed by the $W_2$ 
limit of the signal array. 
Summable failure probabilities, and Borel--Cantelli then yield Theorem~\ref{thm:mdypl-convergence}.

\begin{figure}[ht]
    \centering
    \includegraphics[width=\textwidth]
    {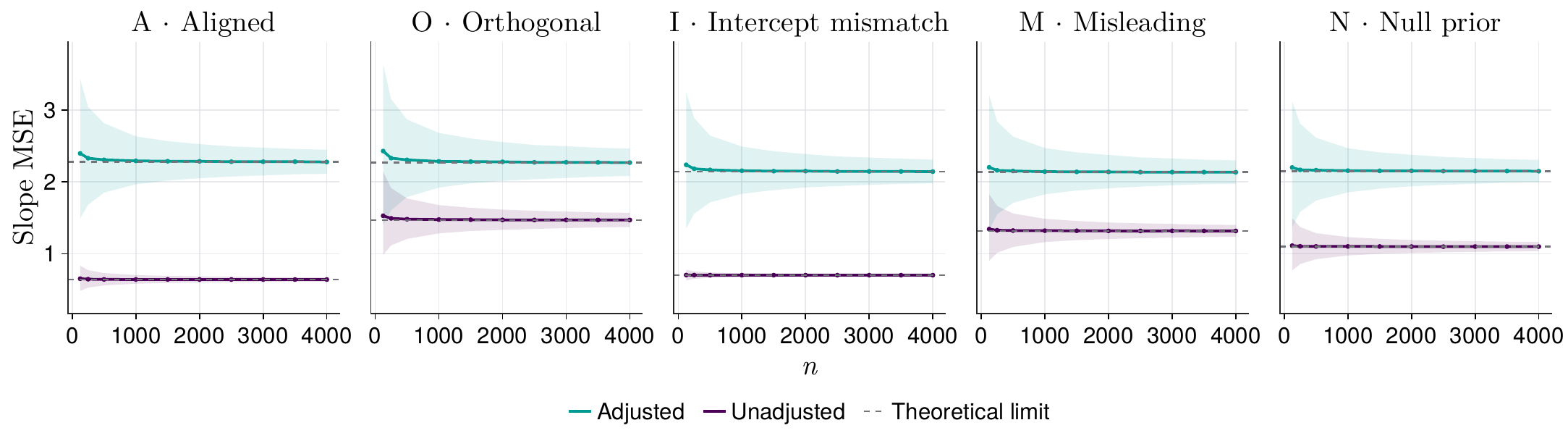}
    \caption{Finite-sample and limiting slope MSEs of the unadjusted MDYPL
    estimator $\betady$ and the oracle-adjusted estimator $\check{\bbeta}$
    across the five simulation settings as $n$ varies with $p/n=0.2$.
    Each estimator uses the setting-specific value of $\alpha$ selected by
    its own limiting slope-MSE criterion.
    Solid lines show Monte Carlo means, shaded bands show the empirical
    $10$--$90$th percentiles, and dashed lines mark the corresponding
    state-equation limits.}
    \label{fig:slope-mse-comparison}
\end{figure}

\subsection{Synthetic data illustration}
\label{subsec:simul}

We conclude this section with a synthetic-data illustration of
Theorem~\ref{thm:mdypl-convergence} that we use throughout the paper. 
Throughout, we fix $\gamma = \delta = 1.5$ and $\kappa = 0.2$, whilst varying
$
    n\in
    \{125,250,500,1000,\ldots,4000\}
$, and 
$p=\kappa n$.
The signal slope $\bbeta_0$ is constructed as
$\bbeta_0=\sqrt p \gamma\bu_1$, where $\bu_1$ is a unit vector of normalised
Gaussian noise, while the true intercept is chosen to give marginal event
probability $0.2$, yielding $\theta_0\approx-1.9163$.
We consider five configurations for the parameters $\theta_P, \bbeta_P$.
For $\bbeta_P \neq \bb 0_p$, we take 
$
    \bbeta_P(r)
    =
    \sqrt p \delta
   \{
        r\bu_1+\sqrt{1-r^2}\bu_2
    \}$,
where $\bu_2$ is a unit vector obtained from independent Gaussian noise and
orthogonalised against $\bu_1$.
Thus, the non-null prior slopes
have equal strength and differ only in their alignment with $\bbeta_0$.
The five configurations are `$\mathrm A$' (aligned and calibrated,
$r=0.8$ and $\theta_P=\theta_0$), `$\mathrm O$' (orthogonal and calibrated,
$r=0$ and $\theta_P=\theta_0$), `$\mathrm I$' (aligned with an incorrectly
calibrated intercept, $r=0.8$, $\theta_P=-\theta_0$), `$\mathrm M$' (misleading, $r=-0.8$, $\theta_P=-\theta_0$), and
`$\mathrm N$' (null, $\bbeta_P=\bb0_p$ and $\theta_P=0$). 
Each setting uses three numerically selected, setting-specific values of
$\alpha$: $\alpha_{s,\mathrm{slope}}^{\mathrm{adj}}$,
$\alpha_{s,\mathrm{slope}}^{\mathrm{DY}}$, and
$\alpha_{s,\mathrm{pred}}^{\mathrm{DY}}$, selected using the adjusted
slope-MSE, unadjusted slope-MSE, and unadjusted classification-error
criteria, respectively.
The complete simulation design is
given in Section~\ref{sec:simul-setup} of the supplementary material.
For our five simulation settings,
Figure~\ref{fig:slope-mse-comparison} compares the finite-sample and limiting
slope risks, with each estimator evaluated at the setting-specific value of $\alpha$
selected using its own limiting slope-MSE criterion.
Alignment lowers the MSE of the unadjusted estimator in settings `$\mathrm{A}$'
and `$\mathrm{I}$', whereas orthogonal or misleading prior slopes increase its
risk relative to the null prior. Oracle adjustment removes this deterministic
prior-aligned component and produces a limiting risk that is similar across
the five prior geometries. In the present settings, however, the associated
variance amplification makes the adjusted MSE larger than the unadjusted MSE.

\section{Prediction}
\label{sec:pred}

\subsection{Misclassification error rate}
\label{subsec:mer}

\begin{figure}[ht]
	\centering
	\includegraphics[width=\textwidth]{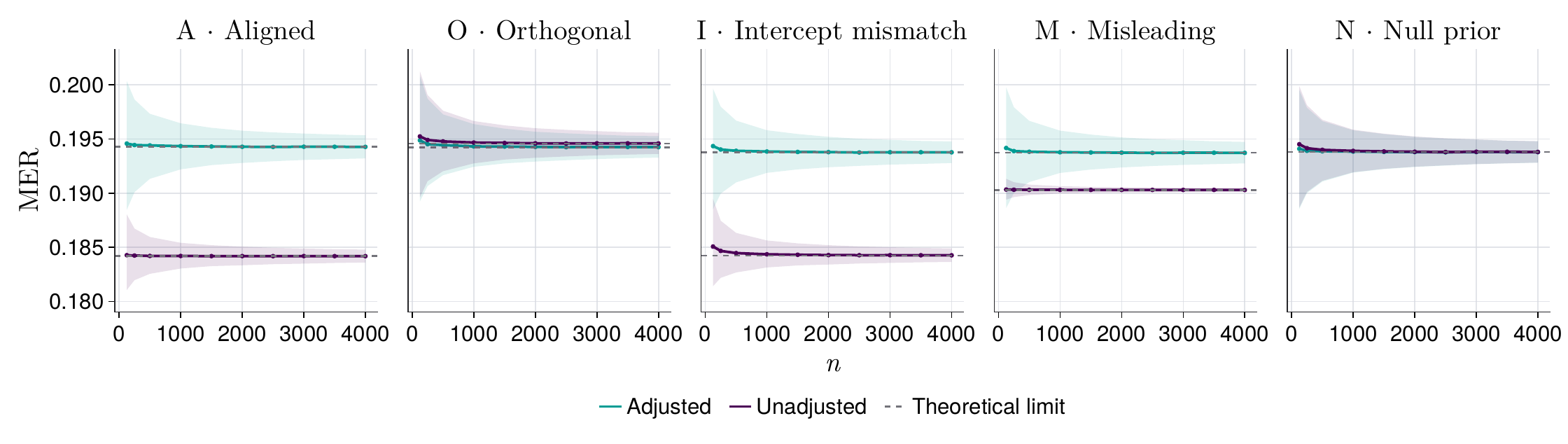}
	\caption{Finite-sample misclassification risks for an independent test observation in the settings of
Section~\ref{subsec:simul}. Classifiers use the unadjusted MDYPL or
oracle-adjusted score with its limiting optimal cutoff and orientation. The adjusted and unadjusted classifiers use the fits at
$\alpha_{s,\mathrm{slope}}^{\mathrm{adj}}$ and
$\alpha_{s,\mathrm{pred}}^{\mathrm{DY}}$, respectively, together
with their corresponding limiting optimal cutoff and orientation. Solid curves show Monte Carlo medians, shaded
bands the empirical $10$--$90$th percentiles, and dashed lines the corresponding
asymptotic risks.}
	\label{fig:mer}
\end{figure}
The convergence result of Theorem~\ref{thm:mdypl-convergence} enables the characterisation of 
the asymptotic out-of-sample prediction performance. 
In this section, we study the out-of-sample misclassification error for the unadjusted and
oracle-adjusted MDYPL logits, with particular attention to optimal
thresholding. 
Since these logits and their limits depend on the tuning parameter $\alpha$, 
that we wish to choose to optimise performance, 
when comparing different values of $\alpha$, we attach an $\alpha$ subscript
to the estimator and to the corresponding limiting quantities, and write
$
    (\sigma_\alpha^*,\bu_\alpha^*,\theta_\alpha^*)$, $
    \bv_\alpha^*
    =
    \bb\Gamma^{+1/2}\bu_\alpha^*
    =
    (v_{1,\alpha}^*,v_{2,\alpha}^*)^\top$, $
    \bbeta_\alpha^*
    =
    \bb B\bv_\alpha^* 
$. 
For a covariate vector $\bx\in\Re^p$, define the true logit, the unadjusted
fitted logit, and the oracle-adjusted logit, respectively, by $\eta_0(\bx)
    =
    \theta_0+\bx^\top\bbeta_0$, $\hat\eta_\alpha^{\mathrm{DY}}(\bx)
    =
    \thetady_\alpha+\bx^\top\betady_\alpha$, $\hat\eta_\alpha^{\mathrm{adj}}(\bx)
    =
    \theta_0+\bx^\top\check{\bbeta}_\alpha$ and $ \check{\bbeta}_\alpha
    =
     (   \betady_\alpha-v_{2,\alpha}^*\bbeta_P) / 
        v_{1,\alpha}^*$. 
The adjusted logit is defined only when $v_{1,\alpha}^*\neq0$. For a
threshold $h\in\Re$, an orientation $a\in\{-1,1\}$, and
$j\in\{\mathrm{DY},\mathrm{adj}\}$, define the oriented threshold
classifier by
\begin{equation}
\label{eq:prediction-threshold-classifier}
    \widehat Y_\alpha^{j,a}(\bx;h)
    =
    \mathds{1}\left\{
        a\left(\hat\eta_\alpha^j(\bx)-h\right)>0
    \right\} \,,
\end{equation}
so that $a=1$ gives an upper-threshold rule and $a=-1$ a lower-threshold
rule.
To evaluate these classifiers on a fresh observation, let
$
    \mathcal F_n
    =
    \sigma(
        \bbeta_0,\bbeta_P,\theta_0,\theta_P,
        \{(\bx_i,\varepsilon_i)\}_{i=1}^n
    )
$
and let $(\bx_{\rm new},\varepsilon_{\rm new})$ be independent of
$\mathcal F_n$, with
$\bx_{\rm new}\sim\mathrm N(\bb0_p,p^{-1}\bb I_p)$ and
$\varepsilon_{\rm new}\sim\mathrm{Unif}([0,1])$ mutually independent. Set
$
    Y_{\rm new}
    =
    \mathds{1}\{
        \varepsilon_{\rm new}
        <
        \rho'(\eta_0(\bx_{\rm new}))
    \}
$.
For the rule in \eqref{eq:prediction-threshold-classifier}, define its
conditional test error by
\begin{equation}
\label{eq:prediction-finite-risk}
    \mathcal E_{n,\alpha}^{j,a}(h)
    =
    \Pr\left(
        Y_{\rm new}
        \neq
        \widehat Y_\alpha^{j,a}(\bx_{\rm new};h)
        \mid
        \mathcal F_n
    \right) \,.
\end{equation}
Its expectation is the corresponding unconditional misclassification error.
To characterise the limiting behaviour of
\eqref{eq:prediction-finite-risk}, recall the definition of $\bb Q$ from
Section~\ref{sec:limiting-ao-focs}, let $\Phi_{\rm N}$ denote the
standard normal distribution function and let $G\sim\mathrm N(0,1)$ and
$\varepsilon\sim\mathrm{Unif}([0,1])$ be mutually independent and independent also
of $\bb Q$. 
Then, we define the limiting logits 
\begin{equation}
\label{eq:prediction-limiting-logits}
\begin{aligned}
    S_0
    =
    \theta_0^*+Q_1,
    \quad
    S_\alpha^{\textrm{\tiny DY}}
    =
    \theta_\alpha^*
    +\bb Q^\top\bv_\alpha^*
    +\sigma_\alpha^*G, \quad 
    S_\alpha^{\textrm{adj}}
    =
    S_0
    +\frac{\sigma_\alpha^*}{v_{1,\alpha}^*}G\,, 
\end{aligned}
\end{equation}
with limiting response $Y=\mathds{1}\{\varepsilon<\rho'(S_0)\}$
For either limiting fitted logit, denote its conditional mean and residual
standard deviation, given $\bb Q$, by $( \overline S_\alpha^{\textrm{\tiny DY}}(\bb Q),
        \tau_\alpha^{\textrm{\tiny DY}}) = (\theta_\alpha^*+\bb Q^\top\bv_\alpha^*,
        \sigma_\alpha^*)$, $(\overline S_\alpha^{\textrm{adj}}(\bb Q),
        \tau_\alpha^{\textrm{adj}}) = (\theta_0^*+Q_1,
        \tau_\alpha)$,
with $\tau_\alpha = \sigma^*_\alpha / \abs{v_{1,\alpha}^*}$. 
Thus
$S_\alpha^j\mid\bb Q\sim
\mathrm N\{\overline S_\alpha^j(\bb Q),(\tau_\alpha^j)^2\}$. 
\begin{theorem}
\label{thm:prediction-classification}
Assume the conditions of Section~\ref{sec:setup}, and fix
$\alpha\in(0,1)$. The statements below hold for $j=\mathrm{DY}$ and, when
$v_{1,\alpha}^*\neq0$, also for $j=\mathrm{adj}$.
\begin{enumerate}[label=(\roman*)]
    \item Conditionally on $\mathcal F_n$, almost surely,
    \begin{equation}
    \label{eq:prediction-joint-limit}
        \left(
            \eta_0(\bx_{\rm new}),
            \hat\eta_\alpha^j(\bx_{\rm new}),
            Y_{\rm new}
        \right)
        \overset{\mathrm d}{\longrightarrow}
        \left(
            S_0,
            S_\alpha^j,
            Y
        \right) \,.
    \end{equation}
    The convergence also holds unconditionally.
    \item Suppose that $\var(S_\alpha^j)>0$. Then, for every fixed
    $h\in\Re$ and $a\in\{-1,1\}$,
    \begin{equation}
    \label{eq:prediction-risk-convergence}
        \mathcal E_{n,\alpha}^{j,a}(h)
        \overset{\mathrm{a.s.}}{\longrightarrow}
        \mathcal E_\alpha^{j,a}(h) \,,
    \end{equation}
    where
    \begin{equation}
    \label{eq:prediction-limiting-risk}
        \mathcal E_\alpha^{j,a}(h)
        =
        \expect\left[
            \rho'(S_0)
            \Phi_{\rm N}\left\{
                \frac{
                    a\{h-\overline S_\alpha^j(\bb Q)\}
                }{
                    \tau_\alpha^j
                }
            \right\}
            +
            \{1-\rho'(S_0)\}
            \Phi_{\rm N}\left\{
                \frac{
                    a\{\overline S_\alpha^j(\bb Q)-h\}
                }{
                    \tau_\alpha^j
                }
            \right\}
        \right] \,,
    \end{equation}
    and which is understood by continuity if $\tau_\alpha = 0$. Moreover,
    \begin{equation}
    \label{eq:prediction-unconditional-risk-convergence}
        \Pr\left(
            Y_{\rm new}
            \neq
            \widehat Y_\alpha^{j,a}(\bx_{\rm new};h)
        \right)
        =
        \expect\left\{
            \mathcal E_{n,\alpha}^{j,a}(h)
        \right\}
        \longrightarrow
        \mathcal E_\alpha^{j,a}(h) \,.
    \end{equation}
\end{enumerate}
\end{theorem}
\noindent
Theorem~\ref{thm:prediction-classification} shows that 
for each fixed $\alpha$, the limiting prediction experiment is the
Gaussian channel in \eqref{eq:prediction-limiting-logits}, that maintains a non-vanishing gaussian noise component. 

\subsection{Optimal thresholding}
\label{subsec:opt-threshold}

The natural question arising from
Theorem~\ref{thm:prediction-classification} is how the threshold $h$ and
orientation $a$ should be chosen optimally. Since this calculation is the same
for the true and fitted logits and depends only on their Gaussian law, we
consider a generic score $S$ satisfying
\begin{equation}
\label{eq:prediction-generic-parameters}
    \begin{bmatrix}
        S_0\\
        S
    \end{bmatrix}
    \sim
    \mathrm N\left(
        \begin{bmatrix}
            \theta_0^*\\
            \mu
        \end{bmatrix},
        \begin{bmatrix}
            \gamma^2 & \chi\\
            \chi & \nu^2
        \end{bmatrix}
    \right), 
    \quad
    \nu>0 \,.
\end{equation}
Let $\varepsilon\sim\mathrm{Unif}([0,1])$ be independent of $(S_0,S)$, set
$q=\rho'(S_0)$, and let $Y=\mathds{1}\{\varepsilon<q\}$. Gaussian
conditioning gives
\begin{equation}
\label{eq:prediction-conditional-logit}
    S_0\mid S=s
    \sim
    \mathrm N\left(
        \theta_0^*+\frac{\chi}{\nu^2}(s-\mu),
        \gamma^2-\frac{\chi^2}{\nu^2}
    \right) \,.
\end{equation}
Consequently, the conditional class probability based on the score is
\begin{equation}
\label{eq:prediction-calibration-map}
    \pi_S(s)
    =
    \Pr(Y=1\mid S=s)
    =
    \expect_Z\left[
        \rho'\left\{
            \theta_0^*
            +\frac{\chi}{\nu^2}(s-\mu)
            +\left(
                \gamma^2-\frac{\chi^2}{\nu^2}
            \right)^{1/2}Z
        \right\}
    \right],
    \quad
    Z\sim\mathrm N(0,1) \,.
\end{equation}
Suppose first that $\chi\neq0$, and define $ h_S
    =
    \mu-({\nu^2}/ {\chi})\theta_0^*$
The conditional mean in \eqref{eq:prediction-conditional-logit} can then be
written as
$
    ({\chi} / {\nu^2}) (s-h_S)
$. 
For every $c\geq0$, the function
$t\mapsto\expect_Z[\rho'(t+cZ)]$ is strictly increasing. Moreover,
$\rho'(-t)=1-\rho'(t)$ and the symmetry of $Z$ imply that this function
equals $1/2$ exactly at $t=0$. Thus
\begin{equation}
\label{eq:prediction-generic-level-set}
    \pi_S(s)>\frac12
    \quad\Longleftrightarrow\quad
    \chi(s-h_S)>0 \,.
\end{equation}
Under zero--one loss, the conditional error from predicting class one at
score value $s$ is $1-\pi_S(s)$, whereas the conditional error from
predicting class zero is $\pi_S(s)$. Hence the optimal score-based rule
predicts class one exactly when $\pi_S(s)>1/2$.

\begin{proposition}
\label{prop:gaussian-score-prediction}
Under the Gaussian-score setup in
\eqref{eq:prediction-generic-parameters},
suppose that $\chi\neq0$. Then
\begin{equation}
\label{eq:prediction-optimal-score-rule}
    \mathds{1}\left\{
        \chi(S-h_S)>0
    \right\}
	\, ,
\end{equation}
is the unique minimiser, up to almost-sure equivalence, of
$
    \Pr(Y\neq g(S))
$
over all $\sigma(S)$-measurable functions $g:\Re\to\{0,1\}$. 
\end{proposition}
\noindent
If $\chi=0$, joint Gaussianity implies that $S$ and $S_0$ are independent,
so the score contains no class information beyond the marginal class
probability, and an optimal score-based rule is a constant majority-class
rule. Accordingly, we call a score $S$ satisfying
\eqref{eq:prediction-generic-parameters} informative if
$\chi=\cov(S_0,S)\neq0$.
For the adjusted and unadjusted limiting logits, the respective optimal
score-based cutoffs are
\begin{equation}
\label{eq:prediction-thresholds}
\begin{aligned}
    h_\alpha^{\textrm{adj}}
    &=
    -
    \frac{\tau_\alpha^2}{\gamma^2}
    \theta_0^*,
    \quad
    h_\alpha^{\textrm{\tiny DY}}
    =
    \theta_\alpha^*
    -
    \frac{
        \vnorm{\bu_\alpha^*}_2^2+(\sigma_\alpha^*)^2
    }{
        \gamma^2v_{1,\alpha}^*
        +\varphi v_{2,\alpha}^*
    }
    \theta_0^* \,.
\end{aligned}
\end{equation}
The adjusted case requires $\gamma^2>0$ and $v_{1,\alpha}^*\neq0$ and is an
upper-threshold rule because
$\cov(S_0,S_\alpha^{\mathrm{adj}})=\gamma^2>0$. The unadjusted case
requires
$\gamma^2v_{1,\alpha}^*+\varphi v_{2,\alpha}^*\neq0$ and uses the upper or
lower inequality according as this covariance is positive or negative. Thus,
for
$
    a_\alpha^{\textrm{adj}}
    =
    1$, $
    a_\alpha^{\textrm{\tiny DY}}
    =
    \sign(
        \gamma^2v_{1,\alpha}^*
        +\varphi v_{2,\alpha}^*
    )
$, 
the optimal rules are
$\mathds{1}\{a_\alpha^j(S_\alpha^j-h_\alpha^j)>0\}$. If
$\gamma^2v_{1,\alpha}^*+\varphi v_{2,\alpha}^*=0$, then
$S_\alpha^{\mathrm{DY}}$ is independent of $S_0$ and is not informative.
An asymptotically deterministic shift of a fitted logit changes its mean and
therefore the location of its optimal cutoff, but not its minimised
classification error. For example, replacing $\theta_0$ in the adjusted
logit by any estimator converging almost surely to $\theta_\alpha^*$ changes
the adjusted cutoff in \eqref{eq:prediction-thresholds} to
$
    \theta_\alpha^*
    -
    (
        1+\tau_\alpha^2 / {\gamma^2}
    )
    \theta_0^* 
$, 
without changing the minimised error or the optimal choice of $\alpha$. 

Having identified the optimal thresholds and orientations, we next ask when
oracle adjustment improves the resulting optimally oriented classifier. For
the generic Gaussian score in \eqref{eq:prediction-generic-parameters},
suppose that $\gamma^2>0$ and $\chi\neq0$, and consider the affine
transformation of $S$
\begin{equation}
\label{eq:prediction-gaussian-channel-main}
    W = \frac{\gamma^2}{\chi}(S-\mu)
    =
    S_0-\theta_0^*+\epsilon_R,
    \quad
    \epsilon_R\perp S_0,
    \quad
    \var(\epsilon_R)
    =
    \gamma^2\left(\frac{1}{R^2}-1\right), 
    \quad
    R^2
    =
    \frac{\chi^2}{\gamma^2\nu^2}
    =
    \operatorname{corr}(S_0,S)^2 \,.
\end{equation}
Thus $S$ is informative when $R^2>0$, and, after choosing the
orientation in Proposition~\ref{prop:gaussian-score-prediction} and
removing irrelevant location and scale, every informative score is
equivalent to observing the same latent logit $S_0$ with an amount of
independent Gaussian noise that decreases in $R^2$. 
For the two limiting fitted logits,
\begin{equation}
\label{eq:prediction-information-logits}
\begin{aligned}
    \left(R_\alpha^{\textrm{\tiny DY}}\right)^2
    =
    \frac{
        (
            \gamma^2v_{1,\alpha}^*
            +\varphi v_{2,\alpha}^*
        )^2
    }{
        \gamma^2
        (
            \vnorm{\bu_\alpha^*}_2^2
            +(\sigma_\alpha^*)^2
        )
    },\quad 
    \left(R_\alpha^{\textrm{adj}}\right)^2
    =
    \frac{\gamma^2}{
        \gamma^2
        +(\sigma_\alpha^*)^2/(v_{1,\alpha}^*)^2
    }
    =
    \frac{\gamma^2}{\gamma^2+\tau_\alpha^2} \,, 
\end{aligned}
\end{equation}
so that 
$S_\alpha^{\mathrm{adj}}$ is informative whenever it is defined, because
$\cov(S_0,S_\alpha^{\mathrm{adj}})=\gamma^2>0$, while
$S_\alpha^{\textrm{\tiny DY}}$ is informative if and only if
$\gamma^2v_{1,\alpha}^*+\varphi v_{2,\alpha}^*\neq0$.
Consequently,
increasing $R^2$, that is, reducing the amount of Gaussian
noise corrupting the score, should translate into a reduction in the
optimally oriented misclassification error, which we show below.
\begin{proposition}
\label{prop:prediction-comparison}
Assume the conditions of Theorem~\ref{thm:prediction-classification} with
$\gamma^2>0$. For $j\in\{\mathrm{DY},\mathrm{adj}\}$ and every
$\alpha\in(0,1)$ for which $S_\alpha^j$ is defined and informative, consider
\begin{equation}
\label{eq:prediction-mer-logits}
    \underset{a \in \{-1,1\}}{\min}\underset{h\in\Re}{\inf} \, \mathcal E_\alpha^{j,a}(h)
    =
    \mathcal E^*((R_\alpha^j)^2),
    \, 
    \mathcal E^*(r)
    =
    \expect\left[
        \rho'(S_0)\,\Phi_{\rm N}\left(\frac{h_r-S_0}{\varsigma_r}\right)
        +\left\{1-\rho'(S_0)\right\}
        \Phi_{\rm N}\left(\frac{S_0-h_r}{\varsigma_r}\right)
    \right] \,,
\end{equation}
for $\varsigma_r^2=\gamma^2(1-r)/r$ and
$h_r=-\theta_0^*(1-r)/r$, $r \in (0,1)$ and where for $r=1$, we set
$\mathcal E^*(1)=\expect[\min\{\rho'(S_0),1-\rho'(S_0)\}]$.

Then $\mathcal E^*$ is strictly decreasing on $(0,1]$, and when the respective scores are defined and informative, 
\begin{enumerate}[label=(\roman*)]
    \item up to ties, and if attained, 
    \begin{equation}
    \label{eq:prediction-optimal-alpha}
        \underset{\alpha \in (0,1)}{\arg\min}\, \mathcal E^*((R_\alpha^j)^2)
        =
        \underset{\alpha \in (0,1)}{\arg\max}\,\left(R_\alpha^j\right)^2 \,, 
    \end{equation}
    which for the adjusted logit, is the $\alpha$ that minimises slope error, i.e. 
    $$
        \underset{\alpha \in (0,1)}{\arg\min}\,
        \mathcal E^*((R_\alpha^{\rm adj})^2)
        =
        \underset{\alpha \in (0,1)}{\arg\min}\,
        \frac{(\sigma_\alpha^*)^2}{(v_{1,\alpha}^*)^2}\,.
    $$
    \item oracle
    adjustment lowers the optimised error if and only if
    $(R_\alpha^{\mathrm{adj}})^2>(R_\alpha^{\textrm{\tiny DY}})^2$, which
    holds whenever $\varphi=0$, $v_{1,\alpha}^*\neq0$ and
    $\delta^2(v_{2,\alpha}^*)^2>0$.
\end{enumerate}
\end{proposition}
\noindent
Proposition~\ref{prop:prediction-comparison} is illustrated in
Figure~\ref{fig:mer} in the simulation setup of Section~\ref{sec:prop-lim} with the adjusted and unadjusted logits evaluated
at their respective numerically selected values of $\alpha$. By Proposition~\ref{prop:prediction-comparison}(i),
the adjusted slope-MSE and classification-error criteria have the same
minimisers.
In the settings we considered, prediction based on the unadjusted logit, which retains the component aligned with $\bbeta_P$, is useful when 
$\bbeta_P$ contains information about $\bbeta_0$, and comparable to the oracle adjusted logit when this component is 
orthogonal to or misaligned with the true signal. Notably, after oracle adjustment the optimally 
thresholded classification error is nearly invariant to substantial changes in the prior geometry and remains 
close to the null-prior benchmark. 

\subsection{Calibration}
\label{subsec:calibration}

So far, we have been concerned with label prediction, for 
which we found the optimal prediction rule to be defined by the 
$1/2$-level set
$
    \{s\in\Re:\pi_S(s)>1/2\} 
$ of the conditional class probability $\pi_S$. 
For probability
prediction, the relevant quantity is $\pi_S(s)$ itself (cf. \citealt{li+sur:2025}). Since
$
    \pi_S(S)
    =
    \Pr(Y=1\mid S)
    =
    \expect[\rho'(S_0)\mid S] \,, 
$
the tower property gives
\begin{equation}
\label{eq:prediction-calibration-identity}
    \Pr(Y=1\mid\pi_S(S))
    =
    \pi_S(S) \,,
\end{equation}
almost surely. 
Thus $\pi_S(S)$ is calibrated in the limiting prediction experiment, in the
sense that, conditional on the reported probability, the probability of
class one equals the reported value.
The same conditional-expectation representation gives an optimality
property. For a differentiable strictly convex function
$\mathfrak f:(0,1)\to\Re$, define the Bregman divergence
$
    D_{\mathfrak f}(a,b)
    =
    \mathfrak f(a)-\mathfrak f(b)-\mathfrak f'(b)(a-b) 
$. 
The Bregman-centroid identity
\citep[Proposition~1]{banerjee+et+al:2005}, applied conditionally on $S$,
implies that every measurable $f:\Re\to(0,1)$ for which the displayed
expectations are finite satisfies
\begin{equation}
\label{eq:prediction-bregman-optimality}
    \expect\left[
        D_{\mathfrak f}\{q,\pi_S(S)\}
    \right]
    \leq
    \expect\left[
        D_{\mathfrak f}\{q,f(S)\}
    \right] \,,
\end{equation}
with equality if and only if $f(S)=\pi_S(S)$ almost surely. Hence
$\pi_S(S)$ is the unique Bregman-optimal approximation to the true class
probability among all probability predictions measurable with respect
to $S$.
For the unadjusted limiting logit $S=S_\alpha^{\mathrm{DY}}$, the parameters
in \eqref{eq:prediction-calibration-map} are
$
    \mu=\theta_\alpha^*$, $
    \nu^2=\vnorm{\bu_\alpha^*}_2^2+(\sigma_\alpha^*)^2$, $
    \chi=\gamma^2v_{1,\alpha}^*+\varphi v_{2,\alpha}^* 
$.
For the adjusted limiting logit $S=S_\alpha^{\mathrm{adj}}$, they are
$
    \mu=\theta_0^*$, $
    \nu^2=\gamma^2+(\sigma_\alpha^*)^2 / (v_{1,\alpha}^*)^2$, $
    \chi=\gamma^2
$.
Substitution of these quantities into
\eqref{eq:prediction-calibration-map} gives the respective oracle
calibration maps.

\section{Hypothesis testing}
\label{sec:testing}

\subsection{Oracle-adjusted $Z$-statistics}
\label{subsec:adjusted-z-stats}

We now turn to testing, for which the proof of
Theorem~\ref{thm:mdypl-convergence} provides the main
ingredients. We first consider a deterministic coordinate block
$I\subseteq\{1,\ldots,p\}$ with $\abs{I}=k$ fixed. 
Let
$\bb J_I\in\Re^{p\times k}$ be the coordinate-selection
matrix that chooses the elements indexed by $I$, and write
$
    \bb B_I
    =
    \bb J_I^\top\bb B 
$ 
for the tested block. 
To analyse the tested block, we decompose $\betady$ into a signal-aligned and signal-orthogonal component, i.e. 
\begin{equation}
\label{eq:isotropic-fixed-block-exact-split}
    \betady
    =
    \bb B\hat\bv
    +
    \sqrt p \hat\sigma \bb d \,, 
\end{equation}
with $\hat\bv  
    =
    (
        \bb B^\top\bb B
    )^{+}
    \bb B^\top\betady$, $\hat\sigma
    =
    \vnorm{\bb P^\perp\betady}_2 / 
        \sqrt p
    $, $\bb d
    =
        \bb P^\perp\betady / 
        \vnorm{\bb P^\perp\betady}_2$, 
where $\bb P^\perp$ is the orthogonal projector onto $\range(\bb B)^\perp$.
For $\betady_I = \bb J_I^\top \betady$, we then have
\begin{equation}
\label{eq:isotropic-fixed-block-decomposition}
    \betady_I -\bb B_I\bv^*
    =
    \bb B_I
    (
        \hat\bv-\bv^*
    )
    +
    \sqrt p \hat\sigma 
    \bb J_I^\top\bb d \,. 
\end{equation}
Generalising the spherical representation of
\citet[Proposition~2.1]{zhao+etal:2022}, we show that, conditionally on
$\mathcal C_n$ from \ref{par:proof-strategy-signal-decomposition}, the direction $\bb d$ is uniform on the unit sphere of 
$\range(\bb B)^\perp$ and independent of $(\hat\bv,\hat\sigma)$, while the
proof of Theorem~\ref{thm:mdypl-convergence} further gives
$
\hat\sigma
\overset{\mathrm p}{\longrightarrow}
\sigma^*
$
and
$
\vnorm{
\bb\Gamma_p^{1/2}
    (
\hat\bv-\bv^*
    )
    }_2
\overset{\mathrm p}{\longrightarrow}
0
$.
The leverage condition \eqref{eq:isotropic-adjusted-z-block-leverage} then makes the
bias $\bb B_I(\hat\bv-\bv^*)$ in \eqref{eq:isotropic-fixed-block-decomposition} vanish, and it ensures that the tested block is
asymptotically orthogonal to $\range(\bb B)$, in which case the projection
of the uniform noise becomes asymptotically Gaussian and provides the
desired pivot.
\begin{theorem}
\label{thm:isotropic-adjusted-z}
Assume the conditions of Section~\ref{sec:setup}. Let
$I\subseteq\{1,\ldots,p\}$ be a deterministic sequence of index sets
with $\abs{I}=k$ fixed, write
$
    \bb B_I
    =
    \bb J_I^\top\bb B
$,
assume the leverage condition
\begin{equation}
\label{eq:isotropic-adjusted-z-block-leverage}
    \mnorms{
        \bb B_I
        \bb\Gamma_p^{+1/2}
    }_2
    =
    \mathcal O_{\mathrm p}(1) \,,
\end{equation}
and let
$\bv^*=[v_1^*,v_2^*]^\top$
be defined in \eqref{eq:limiting-alignment}. Then
\begin{equation}
\label{eq:isotropic-fixed-block-limit}
    \frac1{\sigma^*}
    \left(
        \betady_I-\bb B_I\bv^*
    \right)
    \overset{\mathrm d}{\longrightarrow}
    \mathrm N(\bb0_k,\bb I_k) \,.
\end{equation}
Consequently, for $v_1^* \neq 0$ and a deterministic $b_I^0\in\Re^k$, under
$
    H_0:
    \bbeta_{0,I}=b_I^0 
$, 
we have
\begin{equation}
\label{eq:isotropic-adjusted-z-limit}
    \bb Z_I^{\mathrm{adj}}(b_I^0)
    =
    \frac{\check{\bbeta}_I-b_I^0}{(\sigma^* / v_1^*)}
    \overset{\mathrm d}{\longrightarrow}
    \mathrm N(\bb0_k,\bb I_k) \,.
\end{equation}
\end{theorem}
\noindent
The derivations above using the decomposition in
\eqref{eq:isotropic-fixed-block-decomposition} and the conditional spherical
representation of $\bb d$ establish the stronger result that
Theorem~\ref{thm:isotropic-adjusted-z} applies to any contrast
$\bb a\in\Re^p$ with $\vnorm{\bb a}_2=1$ that may depend on the signal, but
not on the design or responses, provided that
$$
    \mnorms{
        \bb a^\top
        \bb B
        \bb\Gamma_p^{+1/2}
    }_2
    =
    \mathcal O_{\mathrm p}(1)\,.
$$
We further note that the marginal asymptotic standard error of a coordinate $i \in I$ in Theorem~\ref{thm:isotropic-adjusted-z} is
$
    \sigma^*/\abs{v_1^*}
$,
while $(\sigma^* / \abs{v_1^*})^2$ is also the limiting MSE of the
oracle-adjusted slope from Table~\ref{tab:estimation-consequences}. 
Thus, among the admissible values of $\alpha$,
minimising the adjusted slope MSE also minimises the asymptotic width of the
corresponding oracle confidence intervals.
\begin{figure}[ht]
	\centering
	\includegraphics[width=\textwidth]{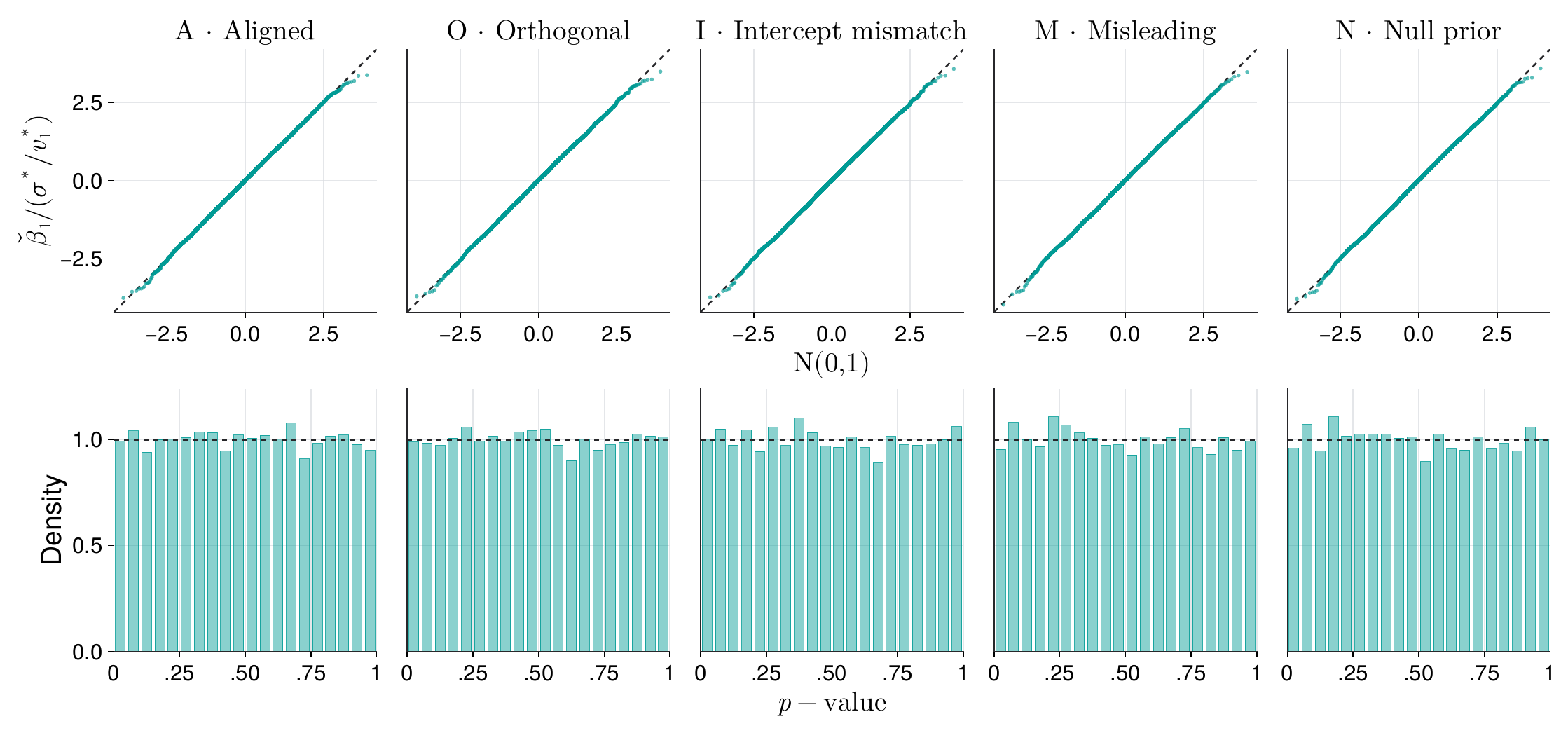}
	\caption{Finite-sample calibration of the oracle-adjusted statistic
$\bb Z_{1}^{\mathrm{adj}}(0)$ for the first null coordinate in the settings of
Section~\ref{subsec:simul}, based on $10{,}000$ replications with $n=4000$
and $p=800$. Top: normal Q--Q plots against $\mathrm N(0,1)$. Bottom:
densities of the corresponding two-sided $p$-values. Black dashed lines mark
the identity and $\mathrm{Unif}([0,1])$ reference density, respectively.}
	\label{fig:z}
\end{figure}
The calibration of the oracle adjusted $Z$ statistics of Theorem~\ref{thm:isotropic-adjusted-z} is illustrated in Figure~\ref{fig:z}, where we investigate the finite-sample calibration
of the oracle-adjusted $Z$-statistic pertaining to the first null coordinate over $10{,}000$ independent runs.
Detailed tables of empirical coverage and quantile calibration are provided in Section~\ref{subsec:appendix-inference-simul} of the supplementary material.

\subsection{Penalised likelihood-ratio test statistic}
\label{subsec:plr}

We next identify the high-dimensional calibration of the MDYPL
penalised likelihood-ratio test statistic. 
For this, 
Theorem~\ref{thm:isotropic-adjusted-z} provides the
fixed-block Gaussian limit that forms the distributional ingredient of
an MDYPL analogue of the high-dimensional likelihood-ratio theory for
logistic maximum likelihood developed by \citet{sur+candes:2019} (see
\citealt[Theorem~6 and Section~3.11]{sur:2019} for the detailed argument).
Extending that argument to MDYPL would largely replicate its
full-versus-reduced analysis, with the likelihood replaced by the MDYPL
objective and the fitted intercept included among the nuisance parameters.
We thus omit the full argument and instead outline the required substitutions and record the implied calibration result.
Under the null hypothesis
$
    H_0:\bbeta_{0,I}=b_I^0
$
the oracle null centre for the MDYPL block is
$
    m_I^0
    =
    v_1^*b_I^0+v_2^*\bbeta_{P,I}
$. Thus, define the penalised likelihood-ratio statistic
\begin{equation}
\label{eq:oracle-PLR-definition}
    \Lambda_I
    =
    \underset{\substack{
        \theta\in\Re\\
        \bbeta\in\Re^p
    }}{\max} \, 
    \ell(\theta,\bbeta;\bYtil,\bX)
	-
    \underset{\substack{
        \theta\in\Re\\
        \bbeta\in\Re^p:  \bbeta_I=m_I^0
    }}{\max} \, 
    \ell(\theta,\bbeta;\bYtil,\bX)
    \,,
\end{equation}
where $\ell(\theta,\bbeta;\bY,\bX)$ is defined in
\eqref{eq:log-reg-ml}, and $\bYtil$ is the $n$-vector with entries
defined in \eqref{eq:mdypl}.
The limiting distribution of $\Lambda_I$ rests on three ingredients. First,
a Taylor expansion around the unrestricted maximiser, followed by profiling
the nuisance slopes and the fitted intercept, should give
\begin{equation}
\label{eq:adjusted-PLR-profile-expansion}
    2\Lambda_I
    =
    \left(
        \betady_I-m_I^0
    \right)^\top
    H_{I}
    \left(
        \betady_I-m_I^0
    \right)
    +
    o_{\mathrm p}(1)
    \,,
\end{equation}
where $H_{I}$ is the Schur complement of the
nuisance block in the negative Hessian. In the coordinatewise setting,
\citet[Theorem~12, equations~(3.147)--(3.152)]{sur:2019} obtain this
expansion from the unrestricted score equations, a Taylor expansion, and the
leave-one-predictor-out approximation in
\citet[Theorem~11, equations~(3.106)--(3.111)]{sur:2019}. The MDYPL
extension requires the corresponding blockwise comparison, with the
pseudo-response and fitted intercept incorporated into the argument.
Second, Theorem~\ref{thm:isotropic-adjusted-z} gives the required centred
MDYPL block limit:
\begin{equation}
\label{eq:adjusted-PLR-block-limit}
    \frac{1}{\sigma^*}
    \left(
        \betady_I-m_I^0
    \right)
    \overset{\mathrm d}{\longrightarrow}
    \mathrm N(\bb0_k,\bb I_k)
    \,,
\end{equation}
which is the fixed-block counterpart of the null-coordinate limit in
\citet[Theorem~5]{sur:2019}.
Third, the random profile curvature in
\eqref{eq:adjusted-PLR-profile-expansion} must satisfy
\begin{equation}
\label{eq:adjusted-PLR-profile-curvature}
    H_{I}
    \overset{\mathrm p}{\longrightarrow}
    \frac{1}{\lambda^*}\bb I_k
    \,.
\end{equation}
In the scalar argument of \citet{sur:2019}, Theorem~12 produces the random
curvature $\kappa/\lambda_{[-j]}$, and Theorem~13, using
Lemmas~24--25, proves that
$\lambda_{[-j]} \to \lambda^*$ in probability.
Here $\lambda^*=t^*/r^*$ is derived from \eqref{eq:FOCs}.
\eqref{eq:adjusted-PLR-profile-curvature} requires the
corresponding matrix-valued calculation for a fixed block, with the fitted
intercept added the profiled nuisance parameters.
Combining \eqref{eq:adjusted-PLR-profile-expansion} and
\eqref{eq:adjusted-PLR-profile-curvature} gives
\begin{equation}
\label{eq:adjusted-PLR-expansion}
    2\Lambda_I
    =
    \frac{1}{\lambda^*}
    \vnorm{
        \betady_I-m_I^0
    }_2^2
    +
    o_{\mathrm p}(1)
    \,.
\end{equation}
Together with \eqref{eq:adjusted-PLR-block-limit}, this implies
\begin{equation}
	\label{eq:plr-limits} 
	    2\Lambda_I
    \overset{\mathrm d}{\longrightarrow}
    \frac{(\sigma^*)^2}{\lambda^*}
    \chi_k^2
,
    \quad\textrm{and hence}\quad
    \frac{\lambda^*}{(\sigma^*)^2}
    2\Lambda_I
    \overset{\mathrm d}{\longrightarrow}
    \chi_k^2
    \,.
\end{equation}
\begin{table}[H]
\centering
\begin{tabular}{@{}lrrrrr@{}}
\toprule
Nominal coverage & \multicolumn{5}{c}{Setting} \\
\cmidrule(l){2-6}
 & `$\mathrm{A}$' & `$\mathrm{O}$' & `$\mathrm{I}$' & `$\mathrm{M}$' & `$\mathrm{N}$' \\
\midrule
90\% & 0.898 & 0.898 & 0.897 & 0.897 & 0.899 \\[-2pt]
 & (0.003) & (0.003) & (0.003) & (0.003) & (0.003) \\
95\% & 0.947 & 0.948 & 0.951 & 0.950 & 0.951 \\[-2pt]
 & (0.002) & (0.002) & (0.002) & (0.002) & (0.002) \\
99\% & 0.989 & 0.989 & 0.989 & 0.989 & 0.989 \\[-2pt]
 & (0.001) & (0.001) & (0.001) & (0.001) & (0.001) \\
\bottomrule
\end{tabular}
\caption{Empirical coverage of nominal
confidence regions for the oracle-rescaled PLR statistic $2\lambda^*\Lambda_I/(\sigma^*)^2$ with null block $I=\{1,\ldots,10\}$, 
along with Monte Carlo standard errors, based on $10{,}000$ replications.}
\label{tab:inference-plr-coverage}
\end{table}
\noindent
Table~\ref{tab:inference-plr-coverage} reports the empirical coverage of
confidence regions obtained by calibrating the candidate oracle-rescaled PLR
statistic in \eqref{eq:plr-limits}, with
$I=\{1,\ldots,10\}$, against the $\chi_{10}^2$ reference distribution,
together with the corresponding Monte Carlo standard errors $\sqrt{\hat c(l) (1 - \hat c(l)) / R}$, where $\hat c(l)$
records the empirical coverage at confidence level $l$.
Across all settings and nominal levels, the empirical coverage is close to
its nominal value, with the observed differences lying within approximately
one Monte Carlo standard errors.
Further results are provided in
Section~\ref{subsec:appendix-inference-simul} of the supplementary material.

\section{Beyond isotropic Gaussian designs} 
\label{sec:beyond-isotropic} 

\subsection{Gaussian designs with arbitrary mean and covariance}
\label{sec:extensions-covariates}

We now seek to relax some of the distributional assumptions about the
covariates. The random design setup in Section~\ref{sec:setup} requires
mean-zero, isotropic Gaussian covariates. The first two restrictions can be
removed by an affine reparameterisation, where whitening standardises the
covariance and the mean shift is handled by the fitted intercept. 
Indeed,
assume that
$$
\bx_i
\overset{\mathrm{i.i.d.}}{\sim}
\mathrm{N}(\bmu_p,\bSigma_p) \,, 
$$
where $(\bmu_p,\bSigma_p)$ is a deterministic triangular array,
$\bSigma_p = \bb L_p \bb L_p^\top$ is positive definite, and $\bb L_p=\bSigma_p^{1/2}$ is its
symmetric positive-definite square root. Consider the affine parameterisation 
\begin{equation}
\label{eq:affine-reparameterisation}
\widetilde\bx_i
=
p^{-1/2}\bb L_p^{-1}(\bx_i-\bmu_p), \quad \widetilde\bx_i\sim\mathrm N(\bb0_p,p^{-1}\bb I_p)
\,,
\end{equation}
Now, the equivariance of the MDYPL estimator under affine transformations of the covariates
(see, also, \citet[][Proposition~2.1]{zhao+etal:2022})
yields that the estimates $(\thetady, \betady)$ in the original coordinates $\bx_i$ and the estimates $(\widetilde{\theta}^{\textrm{\tiny DY}}, \widetilde{\bbeta}^{\textrm{\tiny DY}})$ in the transformed coordinates $\widetilde{\bx}_i$ satisfy
\begin{equation}
    \label{eq:invariance} 
    \widetilde{\theta}^{\textrm{\tiny DY}}
=
\thetady+\bmu_p^\top\betady,
\quad
\widetilde{\bbeta}^{\textrm{\tiny DY}}
=
\sqrt p\bb L_p^\top\betady\,.
\end{equation}
Thus, to apply Theorem~\ref{thm:mdypl-convergence} to
\eqref{eq:affine-reparameterisation}, define
$
\widetilde\theta_0
=
\theta_0+\bmu_p^\top\bbeta_0$, 
$\widetilde\theta_P
=
\theta_P+\bmu_p^\top\bbeta_P$, 
$\widetilde{\bb B}_p
=
\sqrt p\,\bb L_p^\top\bb B$, 
$\widetilde{\bb\Gamma}_p
=
\frac{1}{p}
\widetilde{\bb B}_p^\top
\widetilde{\bb B}_p
=
\bb B^\top\bSigma_p\bb B
$, 
so that $\widetilde{\bb\Gamma}_p$ is the covariance matrix of the centred
true and prior linear predictors
$\widetilde\bx_i^\top\widetilde\bbeta_0$ and
$\widetilde\bx_i^\top\widetilde\bbeta_P$.
If the transformed parameters
$\widetilde\theta_0,\widetilde\theta_P,
\widetilde\bbeta_0,\widetilde\bbeta_P$ satisfy the conditions of
Section~\ref{sec:setup}, Theorem~\ref{thm:mdypl-convergence} applies directly
in the centred and whitened coordinates. 
Let $\widetilde{\bb\Gamma}$ denote the limiting transformed signal Gram matrix,
let
$
    (\widetilde\theta^*,\sigma^*,\bu^*)
$
be the corresponding limiting state parameters, and set
$
    \bv^*
    =
    \widetilde{\bb\Gamma}^{+1/2}\bu^*
    =
    (v_1^*,v_2^*)^\top\
$. 
Corollary~\ref{cor:gaussian-covariance-empirical-law} below gives the original-coordinate intercept limit and,
under separate inverse-whitening conditions, convergence of the slope on pseudo-Lipschitz test functions. Define
\begin{equation}
\label{eq:gaussian-covariance-standardisation}
\begin{aligned}
    \bTheta_p&=\bSigma_p^{-1},\quad
    \bb\Delta_p=\diag(\bTheta_p),\quad
    \bb R_p=\bb\Delta_p^{-1/2}\bTheta_p\bb\Delta_p^{-1/2}, \quad 
    \tau_j^2&=\frac1{(\bTheta_p)_{jj}}
    =\var(x_{i,j}\mid\bx_{i,-j})\,,
\end{aligned}
\end{equation}
and let
$$
    \bb m_p=\bb B^\top\bmu_p,
    \quad
    \bmu_p^\perp
    =\bmu_p-\bSigma_p\bb B\widetilde{\bb\Gamma}_p^+\bb m_p, \quad 
    q_p=\frac1p(\bmu_p^\perp)^\top
    \bSigma_p^{-1}\bmu_p^\perp\,.
$$
\begin{corollary}
\label{cor:gaussian-covariance-empirical-law}
Assume that the transformed parameters
$\widetilde\theta_0,\widetilde\theta_P,
\widetilde\bbeta_0,\widetilde\bbeta_P$ satisfy the
conditions of Section~\ref{sec:setup}, and let
$\bv^*=[v_1^*,v_2^*]^\top$,
$\widetilde\theta^*$ and $\sigma^*$ denote the corresponding limiting state
parameters of \eqref{eq:FOCs}.
\begin{enumerate}[label=(\roman*)]
    \item Suppose that, almost surely,
    \begin{equation}
    \label{eq:original-intercept-shift-conditions}
        \bb m_p\longrightarrow\bb m,
        \quad
        \limsup_{p\to\infty} \, 
        \bb m_p^\top\widetilde{\bb\Gamma}_p^+\bb m_p<\infty,
        \quad
        q_p\log p\longrightarrow0\,,
    \end{equation}
    for a deterministic
    $\bb m=(m_0,m_P)^\top\in\Re^2$. Then
    \begin{equation}
    \label{eq:gaussian-covariance-original-intercept-limit}
        \thetady\overset{\mathrm{a.s.}}{\longrightarrow}
        \widetilde\theta^*-\bb m^\top\bv^*
        =\widetilde\theta^*-v_1^*m_0-v_2^*m_P\,.
    \end{equation}

    \item Suppose that
    $\sup_p \, \mnorms{\bb R_p}_2<\infty$, 
    and, almost surely,
    \begin{equation}
    \label{eq:gaussian-covariance-standardised-signal-law}
        \frac1p\sum_{j=1}^p
        \delta_{(\sqrt p\,\tau_j\bbeta_{0,j},
        \sqrt p \tau_j\bbeta_{P,j})}
        \overset{W_2}{\longrightarrow}\pi_\Sigma\,,
    \end{equation}
    for a deterministic probability law
    $\pi_\Sigma$. Then, for any pseudo-Lipschitz function $\psi$
    of order two,
    \begin{equation}
    \label{eq:gaussian-covariance-test-function-law}
    \begin{aligned}
        \frac1p\sum_{j=1}^p
        \psi\left(
            \sqrt p\tau_j
            \{\betady_j-v_1^*\bbeta_{0,j}-v_2^*\bbeta_{P,j}\},
            \sqrt p\tau_j\bbeta_{0,j},
            \sqrt p\tau_j\bbeta_{P,j}
        \right)
       \overset{\mathrm{a.s.}}{\longrightarrow}
        \expect\left[
            \psi(\sigma^*G,\bar\beta_{0,\Sigma},\bar\beta_{P,\Sigma})
        \right]\,,
    \end{aligned}
    \end{equation}
    where
    $
        (\bar\beta_{0,\Sigma},\bar\beta_{P,\Sigma})\sim\pi_\Sigma$, 
        $G\sim\mathrm N(0,1)
    $,
    independent from each other.
\end{enumerate}
\end{corollary}
\noindent
A uniformly bounded condition number of $\bSigma_p$ is sufficient for
$\sup_p \mnorms{\bb R_p} \leq \infty$, while 
$$
    \bb m_p
    \overset{\mathrm{a.s.}}{\longrightarrow}
    \bb m,
    \quad
    \sup_p \, 
    \bmu_p^\top\bSigma_p^{-1}\bmu_p
    <
    \infty
$$
imply all three conditions in
\eqref{eq:original-intercept-shift-conditions}.
The convergence result Corollary~\ref{cor:gaussian-covariance-empirical-law}
immediately gives convergence of the fitted logits at an independent test point, which replicates the 
prediction results from Section~\ref{sec:pred}
after
replacing
$
    (\theta_0^*,\theta_P^*,\bb\Gamma,\theta_\alpha^*)
$
by
$
    (
        \widetilde\theta_0^*,
        \widetilde\theta_P^*,
        \widetilde{\bb\Gamma},
        \widetilde\theta_\alpha^*
    )
$. 
Finally, for inference, given a fixed block $I$, let
$$
\bb C_{I,p}
=
\left(
\bb J_I^\top\bTheta_p\bb J_I
\right)^{-1}
=
\cov(\bx_{i,I}\mid\bx_{i,-I})\,.
$$
\begin{corollary}

\label{cor:gaussian-covariance-adjusted-z}
Assume that the centred and whitened triangular array satisfies the
conditions of Section~\ref{sec:setup}. Let
$I\subseteq\{1,\ldots,p\}$ be a deterministic sequence of index sets
with $\abs{I}=k$ fixed, and suppose that
\begin{equation}
\label{eq:covariance-adjusted-z-leverage}
    \mnorms{
        \sqrt p
        \bb C_{I,p}^{1/2}
        \bb J_I^\top
        \bb B
        \widetilde{\bb\Gamma}_p^{+1/2}
    }_2
    =
    \mathcal O_{\mathrm p}(1) \,.
\end{equation}
Let
$\bv^*=[v_1^*,v_2^*]^\top$
be the canonical alignment vector of the centred and whitened problem, and
suppose that $v_1^*\neq0$. For deterministic
$b_I^0\in\Re^k$, under
$
    H_0:
    \bbeta_{0,I}=b_I^0
$,
\begin{equation}
\label{eq:covariance-adjusted-z-limit}
    \bb Z_I^{\mathrm{adj}}(b_I^0)
    =
    \sqrt p  \bb C_{I,p}^{1/2} \frac{\check{\bbeta}_I-b_I^0}{(\sigma^* / v_1^*)}
    \overset{\mathrm d}{\longrightarrow}
    \mathrm N(\bb0_k,\bb I_k)\,.
\end{equation}

\end{corollary}
\noindent
For singletons $I=\{j\}$, we have that 
$\bb C_{I,p}=1/(\bTheta_p)_{jj}=\tau_j^2$. As with Theorem~\ref{thm:isotropic-adjusted-z}, the fixed block pivots 
can be extended to unit vector contrasts that satisfy the leverage condition \eqref{eq:covariance-adjusted-z-leverage} with $\bb J_I$ replaced by $\bb a$. 
If $\bSigma_p$ is unknown, the fixed block $\bb C_{I,p}$ can be estimated by
the fixed-block analogue of the scalar residual-regression construction in
\citet[Section~5.1, eq.~(5.1)]{zhao+etal:2022}. Let
$\bb M_{-I}=\bb I_n-\bb P_{[\bb 1,\bX_{-I}]}$ and $d_n=n-p+k-1$, where
$\bb P_{[\bb 1,\bX_{-I}]}$ is the orthogonal projector onto the indicated
column space. Gaussian conditional regression (see also \citealt[Section~5.1, eq.~(5.1)]{zhao+etal:2022}) gives
$$
\widehat{\bb C}_{I,p}
=
\frac{1}{d_n}
\bX_I^\top\bb M_{-I}\bX_I,
\quad
\bb C_{I,p}^{-1/2}
\widehat{\bb C}_{I,p}
\bb C_{I,p}^{-1/2}
\overset{\mathrm p}{\longrightarrow}
\bb I_k\,.
$$
For the PLR test statistic, no
adjustment of the PLR statistic beyond the rescaling of the isotropic
case is required: The affine reparameterisation in \eqref{eq:affine-reparameterisation}
preserves $\Lambda_I$, with centering absorbed by the unrestricted
intercept. After whitening, the contrast
$\bb C_{I,p}^{1/2}\bb J_I^\top\bb L_p^{-\top}$ has orthonormal rows,
so Gaussian rotational invariance reduces the restriction to the
coordinate-block case. 

\subsection{Independent-entry subgaussian designs}
\label{sec:extensions-subgaussian}

There is strong empirical evidence that the Gaussian proportional-limit
predictions remain accurate for independent-entry subgaussian designs (see
\citealt[Section~6]{zhao+etal:2022} and
\citealt[Sections~8 and~S2.5]{sterzinger+kosmidis:2026}). However, the
conditional-CGMT proof of Theorem~\ref{thm:mdypl-convergence} cannot simply
be repeated with a CGMT universality result like that of \citet{han+shen:2023}.
Indeed, under subgaussian designs the conditioned signal-orthogonal block of Section~\ref{par:proof-strategy-signal-decomposition} 
is generally neither independent of the pseudo-responses nor isotropic subgaussian.
A more natural route is the
optimal-value universality theory of \citet{montanari+saeed:2022} applied to the
unconditioned primal problem, and then analyse its matched Gaussian
counterpart using the existing compactification and CGMT argument:
Let $\bW$ have rows $\bxi_i^\top/\sqrt p$, where the entries of $\bxi_i$
are independent, centred, of unit variance and uniformly subgaussian, and
let $\bX$ be the matched Gaussian design. For design
$\bb D\in\{\bW,\bX\}$, let $\bYtil_{\bb D}$ denote the pseudo-response
vector in \eqref{eq:mdypl}, regenerated under $\bb D$ using the same signal
arrays and response-generating uniforms, and, for compact
$\mathcal A_p\subseteq\Re\times\Re^p$, define
\begin{equation}
\label{eq:subgaussian-restricted-value}
    \Psi_{\bb D}(\mathcal A_p)
    =
    -\frac1n
    \max_{(\theta,\bbeta)\in\mathcal A_p}
    \ell(\theta,\bbeta;\bYtil_{\bb D},\bb D)\,.
\end{equation}
For sufficiently large fixed $C_\theta,C_\beta$ and a deterministic
sequence $a_p\to0$ to be chosen sufficiently slowly, set
\begin{equation}
\label{eq:subgaussian-containment-set}
    \mathcal K_p
    =
    \left\{
        (\theta,\bbeta):
        \abs{\theta}\leq C_\theta,\quad
        \frac{\vnorm{\bbeta}_2}{\sqrt p}\leq C_\beta,\quad
        \frac{\vnorm{\bbeta}_\infty}{\sqrt p}\leq a_p
    \right\} \,. 
\end{equation}
Then, given a smoothness argument to handle the discontinuous binary responses, their comparison 
gives, for every sequence of compact sets
$\mathcal A_p\subseteq\mathcal K_p$ and every bounded Lipschitz function
$f:\Re\to\Re$,
\begin{equation}
\label{eq:subgaussian-value-universality}
    \expect [f(\Psi_{\bW}(\mathcal A_p))]
    -
    \expect [f(\Psi_{\bX}(\mathcal A_p))]
    \longrightarrow0\,.
\end{equation}
Once the containment of \eqref{eq:subgaussian-containment-set} is established,
\eqref{eq:subgaussian-value-universality} may be applied separately to
$\mathcal K_p$ and to the radial, signal-alignment and intercept deviation
sets used in Section~\ref{par:proof-strategy-localisation}. 
The strict Gaussian value
gaps in \eqref{eq:proof-strategy-po-separation} would then transfer to the
subgaussian problem and yield
$$
    \Psi_{\bW}(\mathcal K_p)
    \overset{\mathrm p}{\longrightarrow}
    \bar\phi,
    \quad
    (
        \hat\sigma_{\bW},
        \hat\bu_{\bW},
        \hat\theta_{\bW}
    )
    \overset{\mathrm p}{\longrightarrow}
    (
        \sigma^*,
        \bu^*,
        \theta^*
    ) \,.
$$
This gives the subgaussian analogue of the finite-dimensional localisation of Section~\ref{par:proof-strategy-localisation}, while full convergence on test functions would additionally require a universality step linking the subgaussian signal-orthogonal residual to its Gaussian counterpart on pseudo-Lipschitz test functions.
The substantive missing estimator-specific requirement is consequently the
delocalised containment
$$
    \abs{\thetady_{\bW}}=\Op{1},
    \quad
    \vnorm{\betady_{\bW}}_2=\Op{\sqrt p},
    \quad
    \vnorm{\betady_{\bW}}_\infty=\op{\sqrt p} \,.
$$
The first two bounds should follow by the same argument as
Lemma~\ref{lemma:boundedness}, with its Gaussian operator-norm and uniform
$\ell_1$-lower bounds replaced by their subgaussian counterparts, the latter of which may be obtained from
\citet[Theorem~3.1]{litvak+etal:2005} after adaptation to the intercept column in the design.
The sup-norm delocalisation is more delicate. 
On the Gaussian side it follows from Lemma~\ref{lemma:mdypl-gaussian-l2-coupling}
and the Gaussian maximum
bound \citep[Exercise~2.5.10]{vershynin:2018}.  
For the subgaussian side, a promising route is the
leave-one-predictor-out analysis of
\citet[Theorem~11, (3.106)--(3.111) and
(3.121)--(3.122)]{sur:2019}, removing predictor $j$ also from the true and
prior projections defining the pseudo-responses. Uniform score concentration
and profile-curvature control should then yield
$
    \vnorm{\betady_{\bW}}_\infty
    =
    \op{\sqrt p}
$.
\section{Estimating unknowns}
\label{sec:estimating-unknowns}

\begin{figure}[ht]
	\centering 
	\includegraphics[width=.65\textwidth]{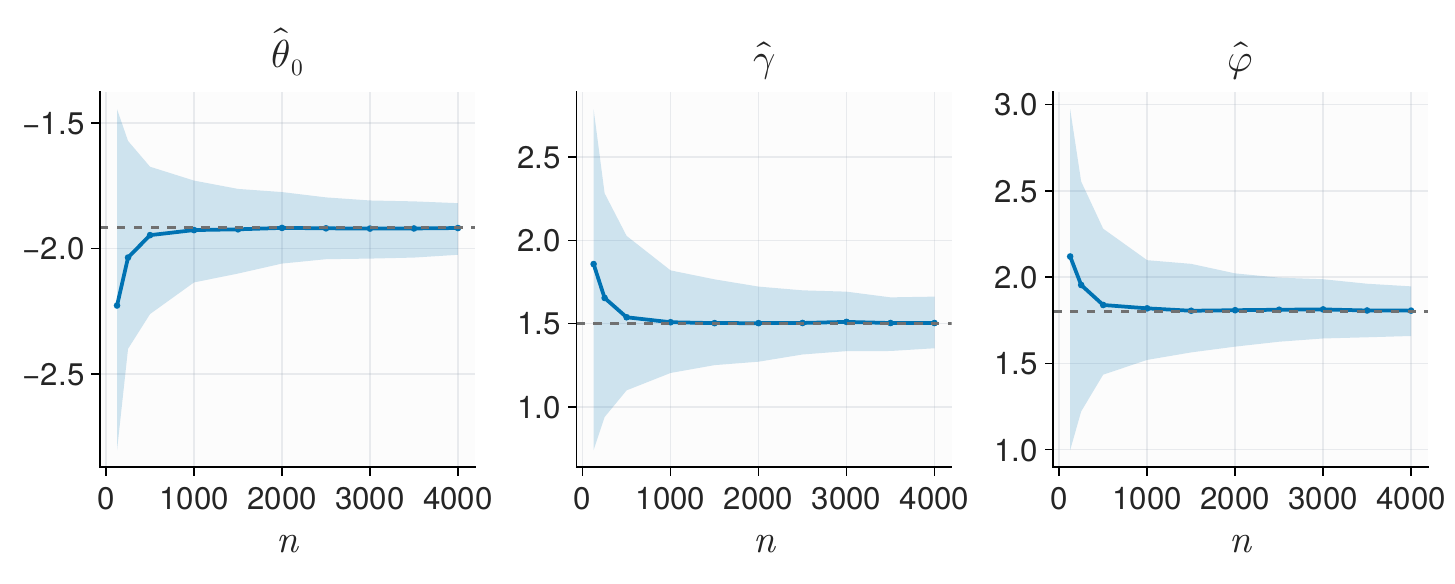}
	\caption{Convergence of response-moment estimates from \eqref{eq:unknowns-o-sample-moments} 
	for oracle parameters $\theta_0, \gamma, \varphi$ in simulation 
	setting `$\mathrm{A}$' of Section~\ref{subsec:simul} with 
	$n \in \{125, 250, 500, 1000, \ldots, 4000 \}$, $p/ n = 0.2$. 
	The solid blue line records the Monte Carlo mean from $1000$ 
	independent runs, while the blue shaded region gives 
	the $10$- and $90$-percentile bounds. 
	The true value is indicated as a grey dashed line.} 
	\label{fig:response-moments}
\end{figure}

\begin{figure}
	\includegraphics[width=\textwidth]{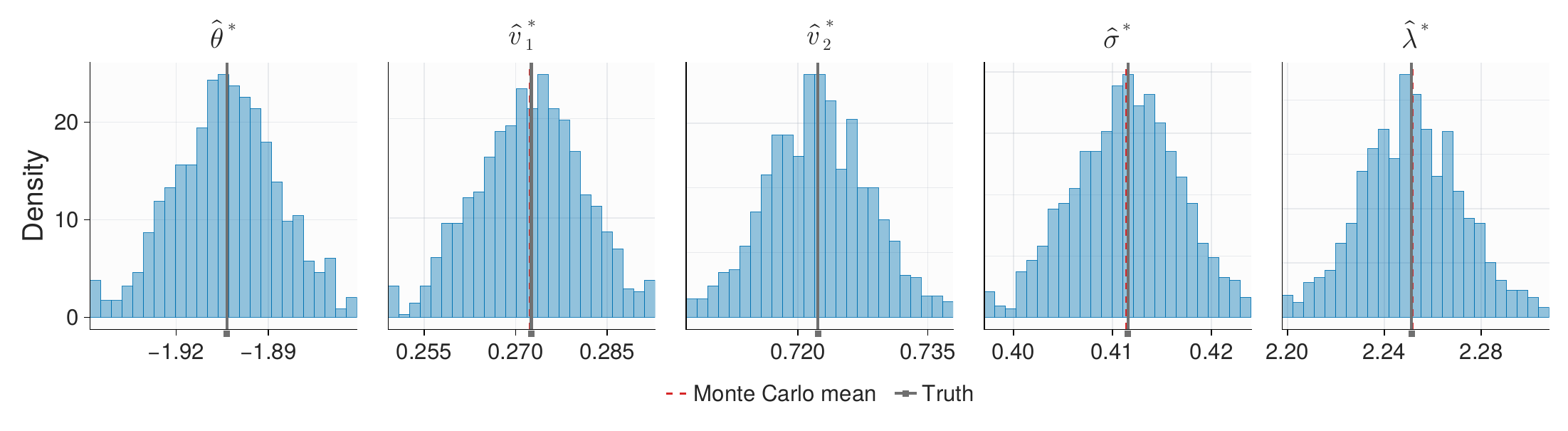}
	\caption{Histograms of estimates of limiting state parameters from 
	\eqref{eq:FOCs} using response-moment estimates from \eqref{eq:unknowns-o-sample-moments} for 
	$\theta_0^*, \gamma, \varphi$ in simulation setting `$\mathrm{A}$' of Section~\ref{subsec:simul} with $n = 4000$, $p=800$ over $1000$ independent 
	runs. 
	The red dashed line indicates the Monte Carlo mean, whilst the grey solid line indicates the truth.} 
	\label{fig:state-pars}
\end{figure}

\noindent The oracle corrections to $\betady$, and the test statistics based on them,
are useful in practice only if the limiting parameters $(\kappa, \theta_0^*, \theta_P^*, \bb \Gamma)$ that enter 
in \eqref{eq:FOCs}, can be
consistently estimated. 
We first note that the parameters $\kappa, \delta, \theta_P^*$ 
can be straightforwardly estimated through the observed quantities
$p/n$, $\delta_p^2=\vnorm{\bbeta_P}_2^2/p$ and $\theta_P$, respectively.
Once estimates of the missing parameters 
$\theta_0^*$, $\bb\Gamma$ are supplied, the limiting state parameters
$(\theta^*,\bu^*,\sigma^*,\lambda^*)$ are determined by the unique
solution of the oracle equations in \eqref{eq:FOCs}, and can be
computed numerically.
\citet{sur+candes:2019} first proposed to estimate the signal strength $\gamma^2$ with their {ProbeFrontier} algorithm, 
that repeatedly subsamples the data, locates the empirical
separability frontier and inverts the phase-transition curve of
\citet{candes+sur:2020}. Beyond computational constraints, this method is tied to the existence boundary of the phase 
transition curve, and does not provide an extension to $(\kappa, \gamma)$ pairs beyond that boundary, where the
ML estimator fails to exist asymptotically but MDYPL remains well defined. 
The SLOE method of
\citet{yadlowsky+etal:2021} instead reparameterises the ML
first-order equations using the corrupted signal strength
$
    \nu^2
    =
    \lim
    \var(\bx_{\rm new}^\top\hat\bbeta)
$
rather than the true signal strength $\gamma^2$, and estimates $\nu^2$ by a
leave-one-out approximation. \citet[Section~7]{sterzinger+kosmidis:2026} adapt this idea
to MDYPL pseudo-responses without a fitted intercept and with
$\bbeta_P=\bb 0_p$.
While well defined beyond the phase transition, \citet[Figure~8]{sterzinger+kosmidis:2026} observe that
for fixed $\kappa$ and $\alpha<1$, the same
corrupted signal strength $\nu^2$ can arise from different values of
the true signal strength $\gamma$. Consequently, the inverse SLOE equations
need not be uniquely identified and can admit multiple admissible roots.
To avoid these identification issues, we estimate $\theta_0$, $\gamma$, and $\varphi$ directly from simple moment 
conditions, before solving the MDYPL state equations \eqref{eq:FOCs} with the resulting plug-in estimates. 
The general idea is not entirely new: \citet[Section~7.2]{zhao+etal:2022} combine the marginal response proportion in \eqref{eq:unknowns-response-channel}
with ProbeFrontier to simultaneously estimate the intercept and signal strength. Here, we use response-moments alone.
For $Z\sim\mathrm N(0,1)$ and $(\theta,g)\in\Re\times[0,\infty)$, define
\begin{equation}
\label{eq:unknowns-response-channel}
\begin{aligned}
    \pi(\theta,g)
    =
    \expect[\rho'(\theta+gZ)], \quad
    \zeta(\theta,g)
    =
    \expect[Z\rho'(\theta+gZ)]
    =
    g c(\theta, g), \quad c(\theta, g) = \expect[\rho''(\theta + g Z)] 
    \,.
\end{aligned}
\end{equation}
To connect \eqref{eq:unknowns-response-channel} to observable moments, let
$\mathcal G_n=\sigma(\bbeta_0,\bbeta_P,\theta_0,\theta_P)$,
$\gamma_p^2=\vnorm{\bbeta_0}_2^2/p$. 
Stein's Lemma \citep[Lemma~1]{stein:1981} yields
\begin{equation}
\label{eq:unknowns-response-population-moments}
\begin{aligned}
    \expect[Y_i\mid\mathcal G_n]
    =
    \pi(\theta_0,\gamma_p), \quad
    p\left\|
        \expect[\bx_iY_i\mid\mathcal G_n]
    \right\|_2^2
    =
    \zeta(\theta_0,\gamma_p)^2 \,, 
\end{aligned}
\end{equation}
For each fixed $g\geq0$, the map
$\theta\mapsto\pi(\theta,g)$ is strictly increasing from zero to one.
Hence, for every $q\in(0,1)$, there is a unique $\theta_q(g)$ such that
$\pi(\theta_q(g),g)=q$. Furthermore, implicit differentiation \citep[Theorem~9.28]{rudin+etal:1976} and
Cauchy--Schwarz \citep[Theorem~1.37(d)]{rudin+etal:1976} give $\partial  \zeta(\theta_q(g),g) / \partial g > 0$ so that 
the profiled map is strictly increasing and
$(\theta,g)\mapsto\{\pi(\theta,g),\zeta(\theta,g)\}$ is injective on
$\Re\times[0,\infty)$. Since
$\zeta(\theta,g)=g c(\theta,g)\geq0$, the conditional moments in \eqref{eq:unknowns-response-population-moments} identify $(\theta_0,\gamma_p)$ and have the 
corresponding sample moments 
\begin{equation}
\label{eq:unknowns-o-sample-moments}
\begin{aligned}
    \widehat\pi
    =
    \bar Y,
    \quad
    \widehat{\zeta}^2
    =
    \frac{p}{n(n-1)}
    \sum_{i\neq j}
    Y_iY_j\bx_i^\top\bx_j, 
    \quad
    \widehat\zeta
    =
    \left\{\max(\widehat{\zeta}^2,0)\right\}^{1/2}\,,
\end{aligned}
\end{equation}
for $\bar Y = \sum_{i = 1}^{n}Y_i / n$. 
Let
$(\widehat\theta_0,\widehat\gamma)$ be the unique solution of 
\begin{equation}
\label{eq:unknowns-response-inversion}
    \pi(\widehat\theta_0,\widehat\gamma)
    =
    \widehat\pi,
    \quad
    \zeta(\widehat\theta_0,\widehat\gamma)
    =
    \widehat\zeta \,,
\end{equation}
when it exists, and set $(\widehat\theta_0,\widehat\gamma)=(0,0)$ otherwise.
To estimate $\varphi$, note that conditional on $\mathcal G_n$, 
\begin{equation}
    \label{eq:cov-y-eta_p}
    \cov(Y_i, \bx_i^\top \bbeta_P \mid \mathcal{G}_n) = \varphi_p c(\theta_0, \gamma_p) \,,
\end{equation}
for $\varphi_p = \langle \bbeta_0, \bbeta_P \rangle / p$. 
Thus, define the sample covariance estimator for \eqref{eq:cov-y-eta_p} as 
$$
    \widehat b = \frac{1}{n - 1} \sum_{i = 1}^{n} (Y_i - \bar Y) \left(\bbeta_P^\top \{\bx_i - \bar \bx \} \right) \,,
$$
where $\bar \bx = \sum_{i = 1}^n \bx_i / n$. Thus, given estimates $\widehat{\gamma}, \widehat{\theta}_0$ for $\gamma_p, \theta_0$, we estimate $\varphi$ as 
\begin{equation}
\label{eq:unknowns-response-geometry-estimates}
    \widetilde{\varphi} = \frac{\widehat b}{c(\widehat{\theta}_0, \widehat{\gamma})}, \quad \widehat{\varphi} = \sign(\widetilde{\varphi}) \min \{\abs{\widetilde{\varphi}}, \widehat{\gamma} \delta_p\} \,,
\end{equation}
where the signed truncation enforces positive semidefiniteness of the estimated
signal Gram matrix.
\begin{proposition}
\label{prop:unknowns-o-identification-consistency}
Assume the conditions of Section~\ref{sec:setup}. Then
\begin{equation}
\label{eq:unknowns-o-consistency}
    (\widehat\theta_0,\widehat\gamma^2,\widehat\varphi)
    \overset{\mathrm p}{\longrightarrow}
    (\theta_0^*,\gamma^2,\varphi).
\end{equation}
\end{proposition}
\noindent
The resulting estimated signal Gram matrix is
\begin{equation}
\label{eq:unknowns-o-gram-estimator}
    \widehat{\bb\Gamma}
    =
    \begin{bmatrix}
        \widehat\gamma^2 & \widehat\varphi\\
        \widehat\varphi & \delta_p^2
    \end{bmatrix}
    \succeq
    \bb0_2,
    \quad
    \widehat{\bb\Gamma}
    \overset{\mathrm p}{\longrightarrow}
    \bb\Gamma \,. 
\end{equation}
Figure~\ref{fig:response-moments} illustrates the convergence of the response-moment estimates
from Proposition~\ref{prop:unknowns-o-identification-consistency} for the simulation setting `$\mathrm{A}$' from Section~\ref{subsec:simul}.
Feasible estimates are obtained by replacing
$(\theta_0^*,\bb\Gamma)$ in \eqref{eq:FOCs} with
$(\widehat\theta_0,\widehat{\bb\Gamma})$. Consistency of the resulting state estimates additionally requires continuity
of this unique oracle solution in
$(\theta_0,\theta_P, \bb\Gamma,\kappa,\alpha)$. We further illustrate the finite sample distribution 
of the limiting state parameters in Figure~\ref{fig:state-pars} for the
simulation setting `$\mathrm{A}$' from Section~\ref{subsec:simul}, which demonstrates good recovery. 
Finally, we note that in contrast to SLOE-type procedures (cf.
Section~\ref{sec:estimating-unknowns-sloe} of the supplementary material and
\citealt[Section~7]{sterzinger+kosmidis:2026}), the first stage response-moments do
not depend on $\alpha$, so that the first stage is computed once and \eqref{eq:FOCs} may then be solved over candidate values of $\alpha$, for example 
to find $\alpha$ that minimises the oracle adjusted slope-MSE $(\sigma^*)^2 / (v_1^*)^2$. 
Further simulation results and a comparison to SLOE are given in
Section~\ref{sec:estimating-unknowns-sloe} of the supplementary material. 

\section{Summary}
\label{sec:summary}

This paper develops a proportional-asymptotic theory for MDYPL logistic
regression with a fitted intercept and a nonzero prior signal via a rigorous conditional CGMT analysis. 
We establish convergence of the fitted intercept and a pseudo-Lipschitz empirical law for
the slope estimator, which yields optimal score-based classification and
calibration rules and oracle-adjusted fixed-block inference, with an affine
extension to noncentred, correlated Gaussian designs. We also propose a
simple response-moment procedure that consistently estimates
$(\theta_0^*,\gamma^2,\varphi)$ before solving the oracle equations for the
remaining state parameters. 
The main remaining restriction is Gaussianity of the design. We view the
value-universality framework of \citet{montanari+saeed:2022} as a promising
route to independent-entry subgaussian designs: the conditional-CGMT
pipeline developed here already provides the unrestricted and restricted
Gaussian value comparisons required for estimator localisation.
The remaining estimator-specific requirement is the delocalised containment needed 
to invoke that framework.

\section{Supplementary materials} 

The repository \url{https://github.com/psterzinger/MDYPL-2} provides the supplementary material document
and scripts to reproduce all analyses and outputs in that document and the main text.

\section{Declarations}

The author is grateful to Ioannis Kosmidis for his helpful comments that greatly improved this work.
For the purpose of open access, the author opts for a Creative Commons Attribution (CC BY)
licence to any Author Accepted Manuscript version arising from this submission.

\putbib
\end{bibunit}

\clearpage
\thispagestyle{empty}

\makeatletter
\edef\main@footnotevalue{\the\c@footnote}
\setcounter{footnote}{0}
\let\@maketitle\main@maketitle
\let\@title\main@title
\let\@author\main@author
\let\@date\main@date
\let\@thanks\main@thanks
\let\thanks\mainthanks
\let\and\mainand
\def\@title{Supplementary material document for:\\
\textit{\main@title}}
\mainmaketitle
\setcounter{footnote}{\main@footnotevalue}
\makeatother
\thispagestyle{empty}

\clearpage

\begin{bibunit}

\appendix

\section{Simulations}
\label{sec:simuls}

\subsection{Simulation setup}
\label{sec:simul-setup}

We use one common data-generating model and five contrasting prior
specifications. The aspect ratio, true signal strength, true event
probability and prior configurations are held fixed, while the slope and
intercept information supplied by the prior slope varies across settings.
Throughout we set 
$$
    n\in
    \{125,250,500,1000,\ldots,4000\},
    \quad
    p=\kappa n,
    \quad
    \kappa=0.2, \quad 
	\gamma = \delta = 3/2 \,.
$$
We construct the signal $\bbeta_0$ as follows: 
Let $I_0=\{1,\ldots,20\}$ denote a fixed block of true-null coordinates. We
generate a standard Gaussian vector $\bb g_1\in\Re^{p-20}$ and define
$\bu_1$ by normalising $\bb g_1$ and prepending $20$ zero coordinates.
Independently, we generate $\bb g_2\sim\mathrm N(\bb0_p,\bb I_p)$ and set
$$
    \widetilde{\bu}_2
    =
    \bb g_2
    -
    \langle\bb g_2,\bu_1\rangle\bu_1,
    \quad
    \bu_2
    =
    \frac{
        \widetilde{\bu}_2
    }{
        \vnorm{\widetilde{\bu}_2}_2
    }\,,
$$
so that $\bu_1$ and $\bu_2$ are orthonormal, $u_{1,j}=0$, $j\in I_0$, 
while $u_{2,j}\neq0$ for every
$j\in I_0$ almost surely. At each value of $p$, the two directions are
generated once, and are then held fixed over all
replications, prior settings and all three tuning streams at that sample size.

We then set
$$
    \bbeta_0
    =
    \sqrt p\,\gamma\bu_1,
    \quad
    \bbeta_P(r)
    =
    \sqrt p\,\delta
    \left\{
        r\bu_1+\sqrt{1-r^2}\bu_2
    \right\}\,,
$$
where $\delta=\gamma=1.5$ for every non-null prior. Hence,
$$
    \frac{\vnorm{\bbeta_0}_2^2}{p}=\gamma^2,
    \quad
    \frac{\vnorm{\bbeta_P(r)}_2^2}{p}=\delta^2,
    \quad
    \frac{\langle\bbeta_0,\bbeta_P(r)\rangle}{p}=r\gamma\delta\,.
$$
All non-null prior slopes therefore have the same strength and differ only in
their alignment with the true slope. Since the set of correlations used in
the simulations is finite and satisfies $\abs{r}<1$, every non-null
$\bbeta_P(r)$ is dense and on the block $I_0$,
$$
    \beta_{0,j}=0,
    \quad
    \beta_{P,j}(r)
    =
    \sqrt p\,\delta\sqrt{1-r^2}\,u_{2,j}
    \neq0,
    \quad
    j\in I_0\,,
$$
almost surely. Thus $I_0$ is a null block for the true slope, but not for the
non-null prior slopes. 
Setting `$\mathrm{N}$' remains the sole exception, with an
identically zero prior slope.

We choose the intercepts through marginal event probabilities. For
$s\geq0$ and $\pi\in(0,1)$, let $a(s,\pi)$ be the unique solution of
$$
    \expect\left[
        \rho'(a(s,\pi)+sZ)
    \right]
    =
    \pi,
    \quad
    Z\sim\mathrm N(0,1)\,.
$$
The true marginal event probability is fixed at $\pi_0=0.2$, giving
$
    \theta_0
    =
    a(1.5,0.2)
    \approx
    -1.9163
$ and $a(1.5,0.8)\approx1.9163$. For the non-null prior settings,
$\theta_P=a(1.5,\pi_P)$. The five settings are given in
Table~\ref{tab:simul-setups}.


\begin{table}[t]
\centering
\caption{Simulation settings.
The common true model has signal strength $\gamma=1.5$ and marginal event
probability $\pi_0=0.2$.
For the non-null prior settings, $\delta=1.5$,
$r=\operatorname{corr}(\bx^\top\bbeta_0,\bx^\top\bbeta_P)$, and
$\pi_P=\expect[\rho'(\theta_P+\bx^\top\bbeta_P)]$.}
\label{tab:simul-setups}
\small
\setlength{\tabcolsep}{3pt}
\begin{tabularx}{\textwidth}{
    @{}
    >{\raggedright\arraybackslash}p{0.24\textwidth}
    c
    c
    c
    c
    >{\raggedright\arraybackslash}X
    @{}
}
\toprule
Setting
&
$\delta$
&
$r$
&
$\pi_P$
&
$\theta_P$
&
Prior information
\\
\midrule
`$\mathrm{A}$': aligned and calibrated
&
$1.5$
&
$0.8$
&
$0.2$
&
$-1.9163$
&
Aligned slope and correctly calibrated intercept
\\
`$\mathrm{O}$': orthogonal and calibrated
&
$1.5$
&
$0$
&
$0.2$
&
$-1.9163$
&
Orthogonal slope and correctly calibrated intercept
\\
`$\mathrm{I}$': aligned, incorrect intercept
&
$1.5$
&
$0.8$
&
$0.8$
&
$1.9163$
&
Aligned slope and oppositely calibrated intercept
\\
`$\mathrm{M}$': misleading
&
$1.5$
&
$-0.8$
&
$0.8$
&
$1.9163$
&
Anti-aligned slope and oppositely calibrated intercept
\\
`$\mathrm{N}$': null
&
$0$
&
---
&
$0.5$
&
$0$
&
Identically zero prior slope
\\
\bottomrule
\end{tabularx}
\end{table}
Settings `$\mathrm{A}$' and `$\mathrm{O}$' isolate the effect of prior-slope
alignment, settings `$\mathrm{A}$' and `$\mathrm{I}$' isolate prior-intercept
misspecification, and settings `$\mathrm{I}$' and `$\mathrm{M}$' compare aligned
and anti-aligned prior slopes under the same incorrect intercept. Setting
`$\mathrm{N}$' provides a neutral baseline and a rank-deficient prior geometry.

For each setting
$s\in\{$`$\mathrm A$',`$\mathrm O$',`$\mathrm I$',`$\mathrm M$',`$\mathrm N$'$\}$,
we use three setting-specific numerically selected tuning parameters for
the limiting adjusted slope MSE, unadjusted slope MSE, and unadjusted
classification error.
For this, let $(\bv_{\alpha,s}^*,\sigma_{\alpha,s}^*)$ denote the relevant solution
components of the population equations in \eqref{eq:FOCs}, let
$\bb e_1=[1,0]^\top$, and define the numerical search interval
$
    I_\alpha=$\{0.01, 0.05,0.10,\ldots,0.95,0.975,0.99, 0.995\}$
$.
We define the numerical-search criteria
$$
\begin{aligned}
    \mathcal Q_{s,\mathrm{adj}}(\alpha)
    &=
    \frac{
        (\sigma_{\alpha,s}^*)^2
    }{
        (v_{1,\alpha,s}^*)^2
    } \\
    \mathcal Q_{s,\mathrm{unadj}}(\alpha)
    &=
        (\sigma_{\alpha,s}^*)^2
        +
        \left(
            \bv_{\alpha,s}^*-\bb e_1
        \right)^\top
        \bb\Gamma_s
        \left(
            \bv_{\alpha,s}^*-\bb e_1
        \right)
    \\
    \mathcal Q_{s,\mathrm{pred}}(\alpha)
    &=
    \frac{
        \left(
            \bb e_1^\top\bb\Gamma_s\bb e_1
        \right)
        \left\{
            (\bv_{\alpha,s}^*)^\top
            \bb\Gamma_s
            \bv_{\alpha,s}^*
            +
            (\sigma_{\alpha,s}^*)^2
        \right\}
    }{
        \left(
            \bb e_1^\top
            \bb\Gamma_s
            \bv_{\alpha,s}^*
        \right)^2
    }\,.
\end{aligned}
$$
We denote by $\alpha_{s,\mathrm{slope}}^{\mathrm{adj}}$,
$\alpha_{s,\mathrm{slope}}^{\mathrm{DY}}$, and
$\alpha_{s,\mathrm{pred}}^{\mathrm{DY}}$ the values returned by the
numerical search below for the adjusted slope-MSE, unadjusted slope-MSE,
and unadjusted classification-error criteria, respectively.
We evaluate each criterion over $I_\alpha$. 

We use $R=10{,}000$ independent Monte Carlo replications at each value of $n$.
Unless a display explicitly shows the sample-size grid, the main-text
illustrations use $n=4000$ and $p=800$.
For replication $b=1,\ldots,R$, we generate one training sample from the common
true model and reuse it across all five prior settings and all three tuning streams.
For out-of-sample prediction, we evaluate the finite-sample conditional test
error in \eqref{eq:prediction-finite-risk}. The adjusted rule uses the fit at
$\alpha_{s,\mathrm{slope}}^{\mathrm{adj}}$, whereas the unadjusted rule uses the
fit at $\alpha_{s,\mathrm{pred}}^{\mathrm{DY}}$. In each case we use the
corresponding population cutoff and orientation in
\eqref{eq:prediction-thresholds}. 
For testing, the common true-null block $I_0$ is used for the fixed-coordinate Gaussian
experiment. In the adjusted-tuning stream, we evaluate the statistic in
\eqref{eq:isotropic-adjusted-z-limit} for $I=\{1\}$ and $b_I^0=0$ in every
replication. 
For the penalised likelihood-ratio experiment, we also use the fit at
$\alpha_{s,\mathrm{slope}}^{\mathrm{adj}}$ and take
$I=\{1,\ldots,10\}$. The oracle restriction centre is
$$
    \bb m_{I,s}^{0,\mathrm{adj}}
    =
    v_{1,s,\mathrm{adj}}^*\bbeta_{0,I}
    +
    v_{2,s,\mathrm{adj}}^*\bbeta_{P,s,I}
    =
    v_{2,s,\mathrm{adj}}^*\bbeta_{P,s,I}\,, 
$$
and record
$$
    \frac{
        \lambda_{s,\mathrm{adj}}^*
    }{
        (\sigma_{s,\mathrm{adj}}^*)^2
    }
    \,2\Lambda_{I,s}^{(b)} \,. 
$$
Limiting parameter estimation for Section~\ref{sec:estimating-unknowns} is performed on 
replications
$b\in\{10,20,\ldots,R\}$. For each fit, we retain the
SLOE variance estimate $\widehat\nu_{\mathrm{SLOE},s}^2$ described in
Section~\ref{sec:estimating-unknowns-sloe} of the supplementary material,
which removes the fitted-intercept contribution from the estimated variance
of the fitted score. The
response-moment estimator, estimates
$\gamma$ using $\widehat{\zeta}^2$ in
\eqref{eq:unknowns-o-sample-moments}, and estimates the remaining response
geometry using the covariance-centred aligned moment $\widehat b$ and the
signed-truncated ratio $\widehat\varphi$ in
\eqref{eq:unknowns-response-geometry-estimates}. 
For each successful estimate
$(\widehat\theta_0,\widehat\gamma^2,\widehat\varphi)$ produced by either
procedure, we form
$$
    \widehat{\bb\Gamma}
    =
    \begin{bmatrix}
        \widehat\gamma^2 & \widehat\varphi\\
        \widehat\varphi & \widehat\delta^2
    \end{bmatrix}\,,
$$
and solve the population first-order equations at the common value
$\alpha_{s,\mathrm{slope}}^{\mathrm{DY}}$ used for the underlying
unadjusted MDYPL fit, recording the resulting estimates of
$(\theta^*,\bu^*,\bv^*,\bb\chi^*,\sigma^*,\lambda^*)$.

\subsubsection{SLOE-based estimation of the limiting parameters}
\label{sec:estimating-unknowns-sloe}

The original SLOE idea \citep{yadlowsky+etal:2021} is to reparameterise the ML analogue of the limiting
system in \eqref{eq:FOCs} in terms of moments of the fitted limiting score.
Applied to the present setting, Theorem~\ref{thm:mdypl-convergence} gives
$$
    \frac{1}{p}
    \vnorm{\betady}_2^2
    \overset{\mathrm{a.s.}}{\longrightarrow}
    \nu^2
    =
    (\sigma^*)^2
    +
    {\bv^*}^{\top}
    \bb\Gamma
    \bv^* \,, 
$$
and
$$
    \frac{1}{p}
    \langle
        \betady,
        \bbeta_P
    \rangle
    \overset{\mathrm{a.s.}}{\longrightarrow}
    \chi_P
    =
    \varphi v_1^*
    +
    \delta^2v_2^* \,. 
$$
Since
$$
    {\bv^*}^{\top}
    \bb\Gamma
    \bv^*
    =
    \gamma^2(v_1^*)^2
    +
    2\varphi v_1^*v_2^*
    +
    \delta^2(v_2^*)^2 \,, 
$$
the two limiting identities can, when $v_1^*\neq0$, be inverted as
$$
    \varphi
    =
    \frac{
        \chi_P-\delta^2v_2^*
    }{
        v_1^*
    } \,, 
$$
and
$$
    \gamma^2
    =
    \frac{
        \nu^2
        -
        (\sigma^*)^2
        -
        2\chi_Pv_2^*
        +
        \delta^2(v_2^*)^2
    }{
        (v_1^*)^2
    } \,. 
$$
Substituting these expressions into \eqref{eq:FOCs} gives a system for
$
    \left(
        \theta_0^*,
        v_1^*,
        v_2^*,
        \sigma^*,
        \lambda^*
    \right)
$
parameterised by
$
    \left(
        \kappa,
        \theta_P^*,
        \theta^*,
        \delta^2,
        \nu^2,
        \chi_P
    \right)
$, 
for which the direct plug-in estimates are
$
    p / n$, 
    $
    \theta_P$, 
    $
    \thetady$, 
    $
    \delta_p^2
    =
    \vnorm{\bbeta_P}_2^2 / p$, 
    $
    \widehat\nu^2
    =
    \vnorm{\betady}_2^2 / p$, 
    $
    \widehat\chi_P
    =
    \langle
        \betady,
        \bbeta_P
    \rangle / p
$.
For numerical stability, the implementation instead replaces
$(\bv^*,\sigma^*)$ by $(\bb\chi,\nu^2)$, where
$\bb\chi=(\cov(S^{\textrm{\tiny DY}}, Q_1),\chi_P)^\top$, treats
$(\widehat\chi_P,\widehat\nu^2)$ as plug-ins, and solves for
$(\theta_0^*,\gamma^2,\varphi,\lambda^*,\chi)$. On the full-rank region
with $v_1^*\neq0$, this is algebraically equivalent to the direct
parameterisation above.

\citet{yadlowsky+etal:2021} further provide a leave-one-out alternative to the direct estimation of
$\widehat\nu^2$. Let
$\widehat\eta_i=\widehat\eta^{\textrm{\tiny DY}}(\bx_i)$, set
$
    \widehat{\bb W}
    =
    \operatorname{diag}
    \{
        \rho''(\widehat\eta_i)
    \}_{i=1}^n
$,
and let $\widehat h_i$ be the $i$th diagonal element of the corresponding
weighted hat matrix for the augmented design. Adapting
\citet{yadlowsky+etal:2021} to the MDYPL pseudo-responses, as in
\citealt[equation~(11)]{sterzinger+kosmidis:2026}, define
$$
    \widehat\nu_{\textrm{\tiny SLOE}}^2
    =
    \frac1n
    \sum_{i=1}^n
    \left(
        \widehat\eta_i^{\textrm{\tiny SLOE}}
        -
        \bar\eta^{\textrm{\tiny SLOE}}
    \right)^2,
    \quad
    \widehat\eta_i^{\textrm{\tiny SLOE}}
    =
    \widehat\eta_i
    -
    \frac{\widehat h_i}{1-\widehat h_i}
    \frac{
        \Ytil_i-\rho'(\widehat\eta_i)
    }{
        \rho''(\widehat\eta_i)
    },
    \quad
    \bar\eta^{\textrm{\tiny SLOE}}
    =
    \frac1n
    \sum_{i=1}^n
    \widehat\eta_i^{\textrm{\tiny SLOE}} \,.
$$
The centring removes the fitted-intercept contribution, while
Theorem~\ref{thm:mdypl-convergence} supplies the required full-sample
intercept and slope-norm limits. If the corresponding leave-one- and
leave-two-out stability results hold for the augmented MDYPL problem, the
argument of \citet[Proposition~2]{yadlowsky+etal:2021} yields $\widehat\nu_{\textrm{\tiny SLOE}}^2 \overset{\mathrm p}{\longrightarrow} \nu^2 \,.$

While the SLOE reparameterisation is algebraically valid, it does not by itself
establish invertibility of the nuisance map. In particular, even with a fitted
intercept, a fixed-intercept slice of the map from the true signal strength
$\gamma$ to the corrupted fitted-score strength $\nu$ can fold. Distinct
values of $\gamma$ may therefore produce the same fitted-score moments, and
the inverse SLOE system may have more than one admissible root. 
Moreover, consistency of a selected root also
requires local identification, uniform convergence of the plug-in equations,
existence of a nearby sample root and a rule selecting that root when several
solutions are present.

\begin{figure}[H]
	\centering
	\includegraphics[width=0.8\textwidth]{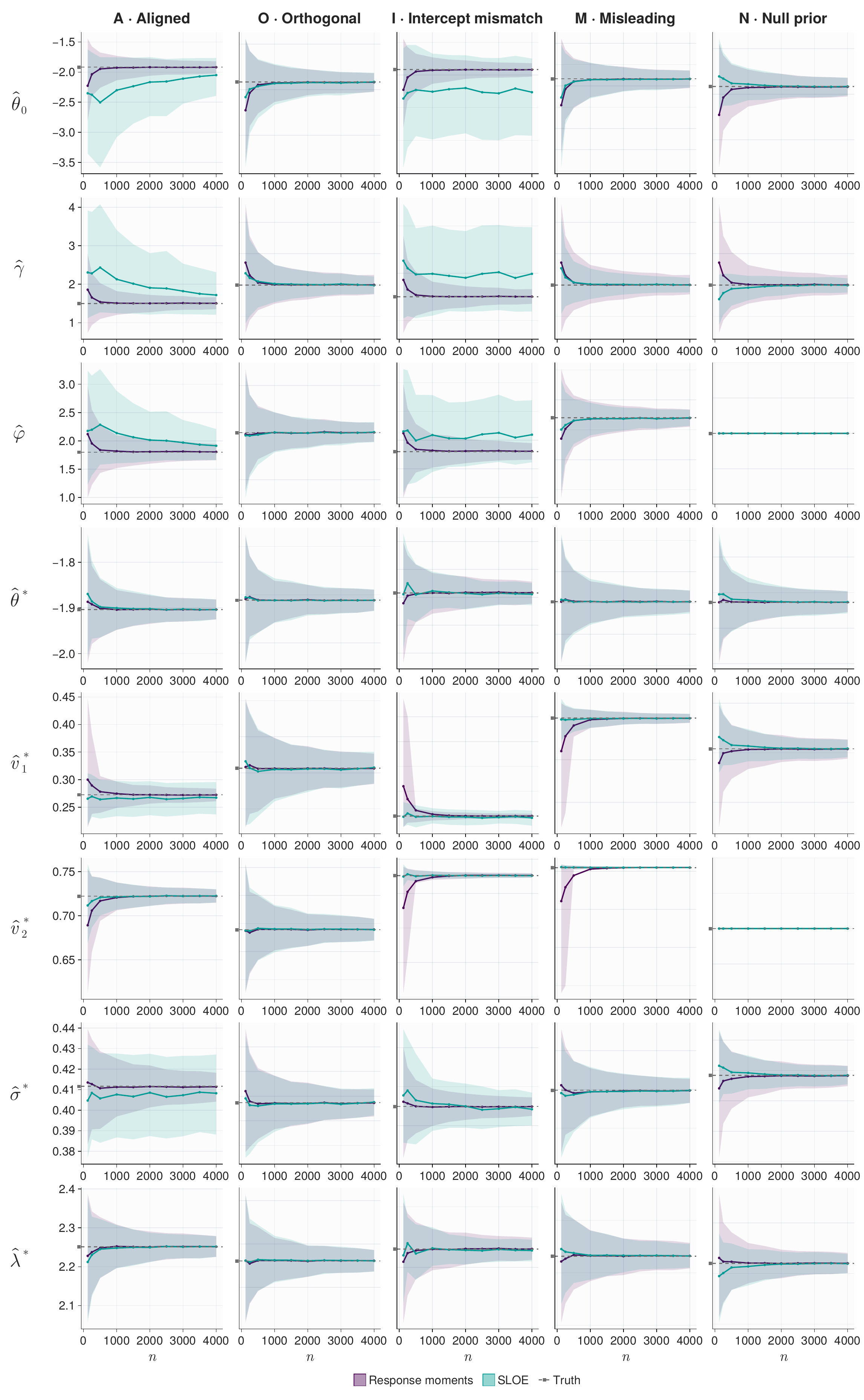}
	\caption{Finite-sample behaviour of the response-moment and fitted-score
	SLOE estimates of the unknown model parameters and induced limiting state
	parameters for varying $n$.
	Solid curves and shaded bands show means and pointwise $10$--$90$th
	percentiles over successful unknown-parameter estimates.
	Grey dashed lines mark the corresponding population targets.
    }
	\label{fig:appendix-estimating-unknowns}
\end{figure}

\begin{figure}[H]
	\centering
	\includegraphics[width=0.8\textwidth]{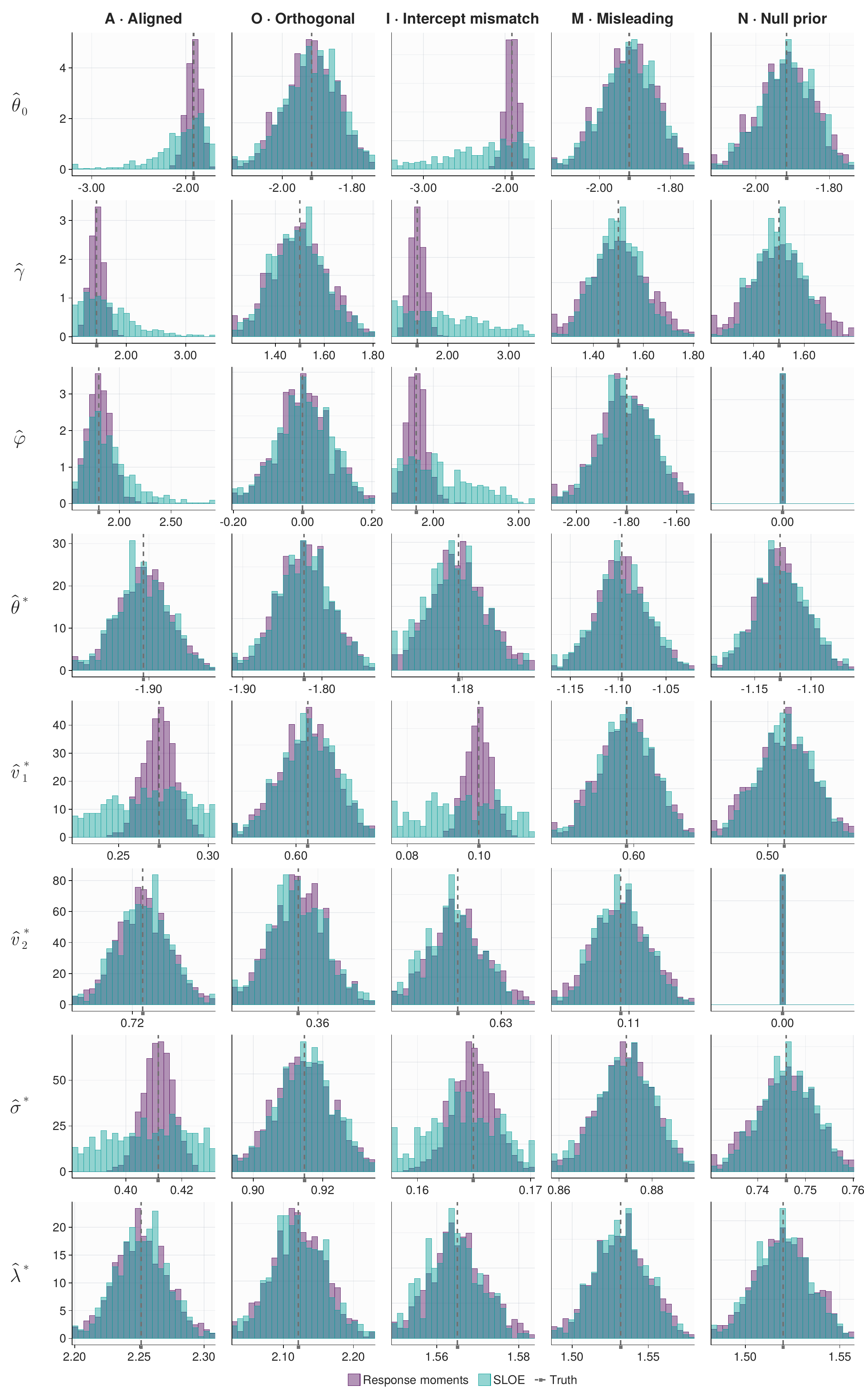}
	\caption{Monte Carlo distributions of the response-moment and fitted-score
	SLOE estimates of the unknown model parameters and induced limiting state
	parameters at $n=4000$ and $p=800$, which are obtained by re-solving the oracle equations
	at $\alpha_{s,\mathrm{slope}}^{\mathrm{DY}}$.
	Grey dashed vertical lines mark the corresponding population targets.}
	\label{fig:appendix-estimating-vars}
\end{figure}

\subsection{Additional simulation output}

\subsubsection{Additional output for Section~\ref{subsec:simul}}
\label{subsec:appendix-limit-simul}

\begin{figure}[H]
	\centering
	\includegraphics[width=\textwidth]{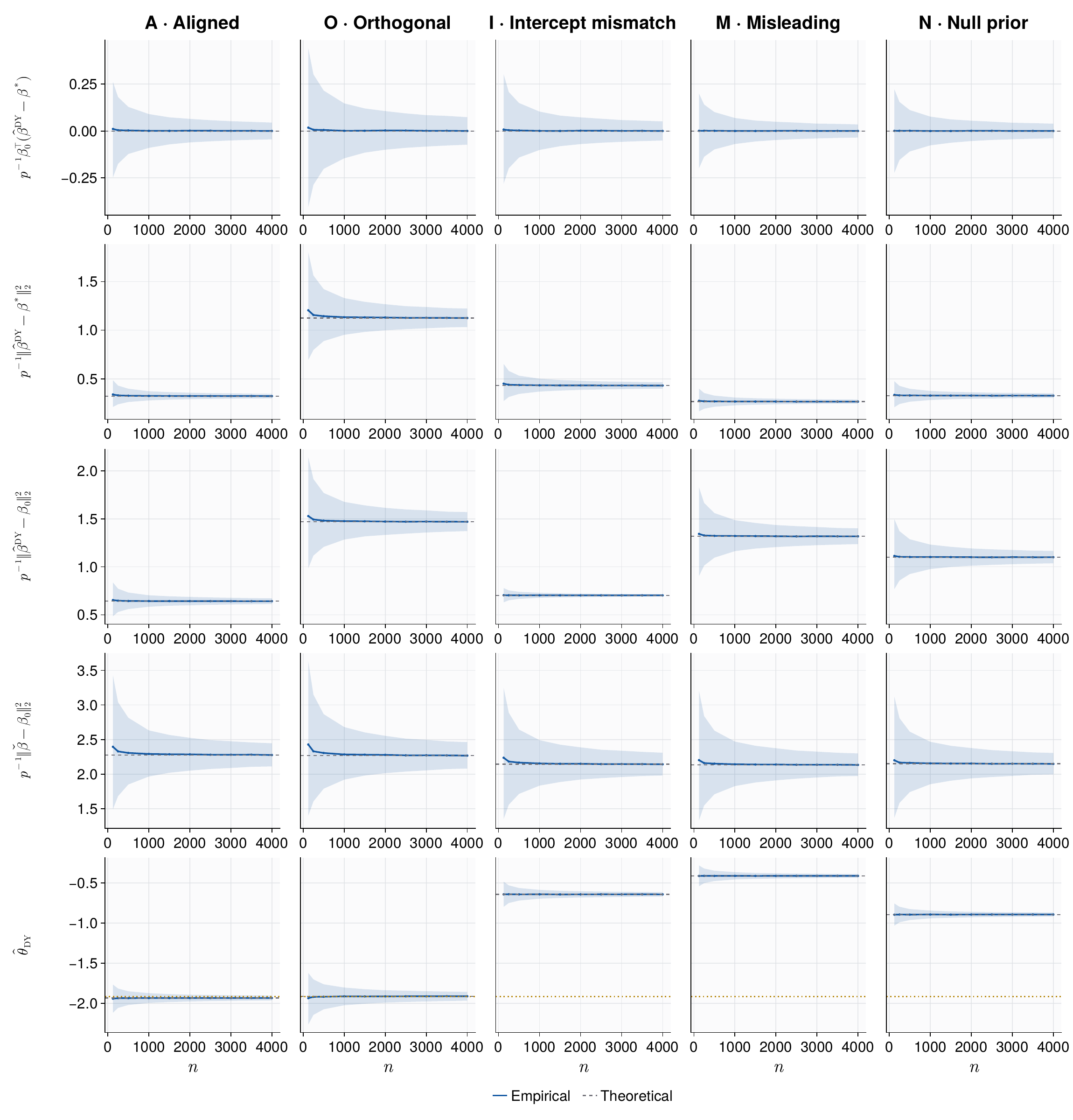}
	\caption{Aggregate finite-sample convergence toward the limits in
	Theorem~\ref{thm:mdypl-convergence} for the consequences summarised in
	Table~\ref{tab:estimation-consequences}.
	The third row uses $\alpha_{s,\mathrm{slope}}^{\mathrm{DY}}$, while the
	remaining rows use $\alpha_{s,\mathrm{slope}}^{\mathrm{adj}}$.
	Solid blue curves show Monte Carlo means, shaded bands span the empirical
	$10$--$90$th percentiles, and grey dashed lines mark the corresponding
	state-equation limits.
	In the fitted-intercept panels, gold dotted lines mark the data-generating
	value $\theta_0$.}
	\label{fig:mdypl-convergence-all}
\end{figure}

\subsubsection{Additional output for Section~\ref{sec:testing}}
\label{subsec:appendix-inference-simul}

\begin{table}[H]
\centering
\begin{tabular}{@{}lrrrrrrrrr@{}}
\toprule
Setting & 1\% & 5\% & 10\% & 25\% & 50\% & 75\% & 90\% & 95\% & 99\% \\
\midrule
$\mathrm{N}(0,1)$ reference& -2.326 & -1.645 & -1.282 & -0.674 & 0.000 & 0.674 & 1.282 & 1.645 & 2.326 \\
$\mathrm{A}$  $\cdot$ Aligned & -2.397 & -1.617 & -1.256 & -0.648 & 0.001 & 0.670 & 1.260 & 1.609 & 2.320 \\
$\mathrm{O}$ $\cdot$ Orthogonal & -2.381 & -1.623 & -1.244 & -0.651 & 0.000 & 0.673 & 1.258 & 1.608 & 2.290 \\
$\mathrm{I}$ $\cdot$ Wrong intercept & -2.405 & -1.631 & -1.253 & -0.644 & 0.001 & 0.681 & 1.240 & 1.621 & 2.309 \\
$\mathrm{M}$ $\cdot$ Misleading & -2.381 & -1.620 & -1.247 & -0.643 & 0.006 & 0.678 & 1.258 & 1.629 & 2.303 \\
$\mathrm{N}$ $\cdot$ Null prior & -2.397 & -1.639 & -1.247 & -0.656 & -0.002 & 0.676 & 1.259 & 1.627 & 2.338 \\
\bottomrule
\end{tabular}
\caption{Empirical quantiles of the oracle-adjusted statistic
$\bb Z_{1}^{\mathrm{adj}}(0)$ for the first true-null coordinate,
compared with the corresponding $\mathrm N(0,1)$ quantiles.
Results are based on $10{,}000$ replications with $n=4000$ and $p=800$.}
\label{tab:inference-adjusted-z-quantiles}
\end{table}

\begin{table}[H]
\centering
\begin{tabular}{@{}lrrrrr@{}}
\toprule
Nominal coverage & \multicolumn{5}{c}{Setting} \\
\cmidrule(l){2-6}
 & `$\mathrm{A}$' & `$\mathrm{O}$' & `$\mathrm{I}$' & `$\mathrm{M}$' & `$\mathrm{N}$' \\
\midrule
90\% & 0.907 & 0.907 & 0.905 & 0.904 & 0.902 \\[-2pt]
 & (0.003) & (0.003) & (0.003) & (0.003) & (0.003) \\
95\% & 0.952 & 0.952 & 0.951 & 0.951 & 0.951 \\[-2pt]
 & (0.002) & (0.002) & (0.002) & (0.002) & (0.002) \\
99\% & 0.990 & 0.990 & 0.989 & 0.990 & 0.990 \\[-2pt]
 & (0.001) & (0.001) & (0.001) & (0.001) & (0.001) \\
\bottomrule
\end{tabular}
\caption{Empirical coverage of nominal $90\%$, $95\%$, and $99\%$
two-sided confidence intervals for the null coordinate
$\beta_{0,1}=0$, based on the oracle-adjusted statistic
$\bb Z_{1}^{\mathrm{adj}}(0)$, along with Monte Carlo standard errors.
Results are based on $10{,}000$ replications with $n=4000$ and $p=800$.}
\label{tab:inference-adjusted-z-coverage}
\end{table}

\begin{table}[H]
\centering
\begin{tabular}{@{}lrrrrrrrrr@{}}
\toprule
Setting & 1\% & 5\% & 10\% & 25\% & 50\% & 75\% & 90\% & 95\% & 99\% \\
\midrule
$\chi_{10}^2$ reference& 2.558 & 3.940 & 4.865 & 6.737 & 9.342 & 12.549 & 15.987 & 18.307 & 23.209 \\
$\mathrm{A}$ $\cdot$ Aligned & 2.644 & 3.966 & 4.870 & 6.788 & 9.317 & 12.607 & 16.034 & 18.507 & 23.353 \\
$\mathrm{O}$ $\cdot$ Orthogonal & 2.557 & 3.912 & 4.812 & 6.812 & 9.359 & 12.592 & 16.076 & 18.400 & 23.509 \\
$\mathrm{I}$ $\cdot$ Wrong intercept & 2.598 & 3.949 & 4.838 & 6.784 & 9.360 & 12.548 & 16.113 & 18.211 & 23.560 \\
$\mathrm{M}$ $\cdot$ Misleading & 2.597 & 3.990 & 4.895 & 6.781 & 9.336 & 12.525 & 16.097 & 18.293 & 23.459 \\
$\mathrm{N}$ $\cdot$ Null prior & 2.657 & 3.956 & 4.908 & 6.761 & 9.368 & 12.553 & 16.015 & 18.267 & 23.628 \\
\bottomrule
\end{tabular}
\caption{Empirical quantiles of the candidate oracle-rescaled penalised
likelihood-ratio statistic
$2\lambda^*\Lambda_I/(\sigma^*)^2$, for
$I=\{1,\ldots,10\}$, compared with the corresponding $\chi_{10}^2$
quantiles.
Results are based on $10{,}000$ replications with $n=4000$ and $p=800$.}
\label{tab:inference-plr-quantiles}
\end{table}

\begin{figure}[H]
	\centering
	\includegraphics[width=\textwidth]{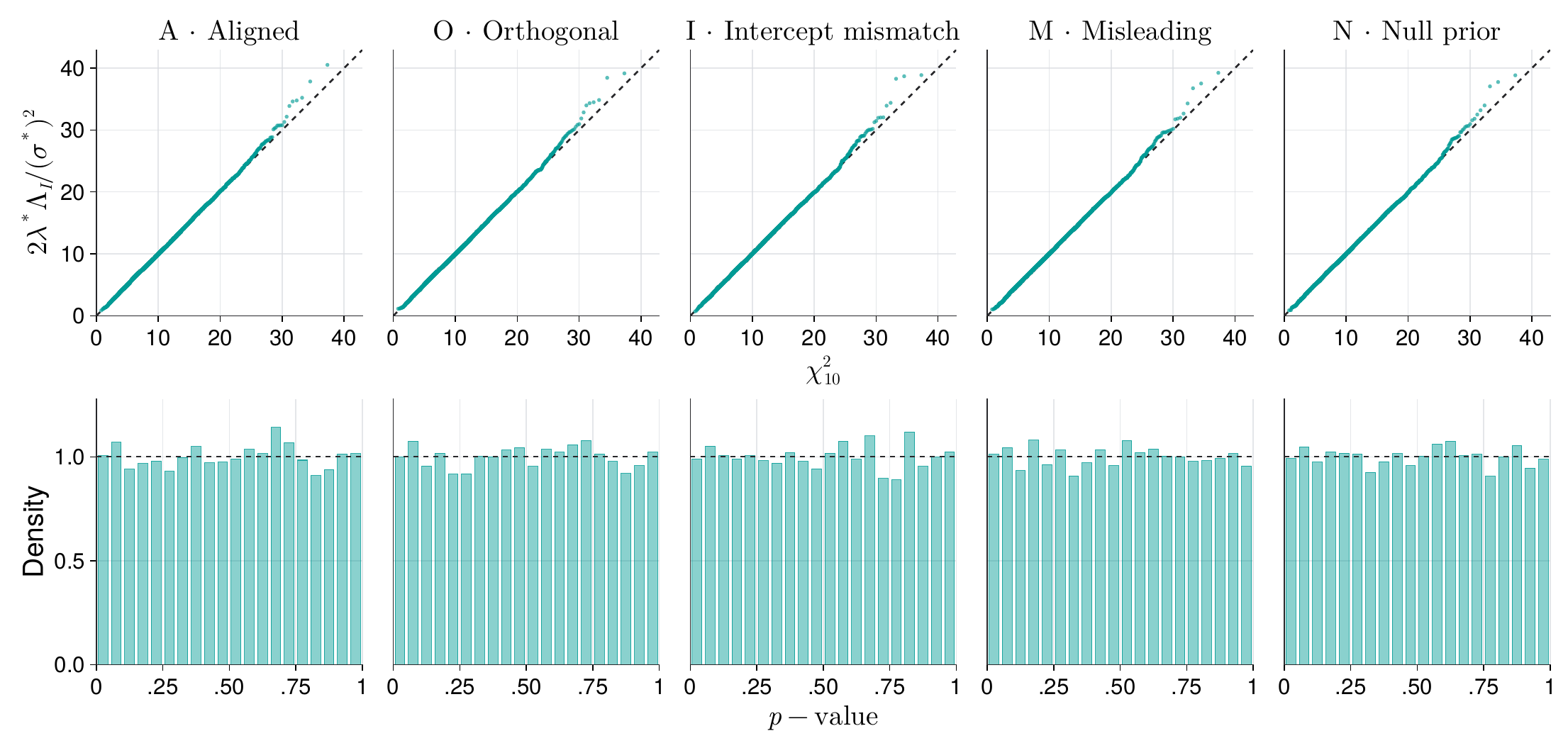}
	\caption{Finite-sample calibration of the oracle-adjusted penalised likelihood-ratio statistic for the true-null block
	$I=\{1,\ldots,10\}$ in the settings of
Section~\ref{subsec:simul}, based on $10{,}000$ replications with $n=4000$
and $p=800$. Top: Q--Q plots of $2\lambda^*\Lambda_I/(\sigma^*)^2$ against
	the $\chi_{10}^2$
	distribution. Bottom:
densities of the corresponding one-sided $p$-values. 
	Black dashed lines mark the identity line and the
	$\mathrm{Unif}([0,1])$ reference density, respectively.}
	\label{fig:plr}
\end{figure}

\newpage

\section{Proof of Theorem~\ref{thm:mdypl-convergence}}
\label{sec:proof-mdypl-conv}
\begingroup
\renewcommand{\thetheorem}{\arabic{theorem}}
\setcounter{theorem}{0}
\begin{theorem}
    \label{thm:mdypl-convergence}
    Assume the conditions of Section~\ref{sec:setup}, and let
    $\theta^*$, $\sigma^*$ and $\bbeta^*$ be as in
    \eqref{eq:limit_ao_natural} and \eqref{eq:limiting-alignment}. Then the
    MDYPL estimator in \eqref{eq:og} satisfies the following.

    \begin{enumerate}[label=(\roman*)]
        \item
        $$
            \thetady
            \overset{\mathrm{a.s.}}{\longrightarrow}
            \theta^*\,.
        $$

        \item For every pseudo-Lipschitz function
        $\psi:\Re^3\to\Re$ of order two
        ,
        \begin{equation}
        \label{eq:test_fun_conc}
            \frac1p
            \sum_{j=1}^p
            \psi
            \left(
                \betady_j-\bbeta_j^*,
                \bbeta_{0,j},
                \bbeta_{P,j}
            \right)
            \overset{\mathrm{a.s.}}{\longrightarrow}
            \expect
            \left[
                \psi(\sigma^*G,\bar\beta_0,\bar\beta_P)
            \right]\,,
        \end{equation}
        where
        $
            (\bar\beta_0,\bar\beta_P)
            \sim
            \pi_{\{\bar\beta_0,\bar\beta_P\}}
        $
        and $G\sim\mathrm N(0,1)$ is independent of
        $(\bar\beta_0,\bar\beta_P)$.
    \end{enumerate}
\end{theorem}
\endgroup

\subsection{Proof roadmap}\label{sec:roadmap}
By Lemma~\ref{lemma:compact_PO_uniqueness}, on $\mathcal E_n$ the compact
PO has a unique primal minimiser
$(\hat\bv,\hat\bw,\hat\theta,\hat{\bb\eta})$. Set
$$
    \hat\bbeta
    =
    \bb B\hat\bv+\bb E\hat\bw,
    \quad
    \hat\sigma
    =
    \frac{\vnorm{\hat\bw}_2}{\sqrt p}\,.
$$
The direction $\bb d$ of \eqref{eq:compact_PO_direction} is defined also
when $\hat\sigma=0$ and satisfies
\begin{equation}
\label{eq:betady-comp}
    \hat\bbeta
    =
    \underbrace{
        \bb B\hat\bv
    }_{\textrm{signal-aligned component}}
    +
    \underbrace{
        \sqrt p\,\hat\sigma\,\bb d
    }_{\textrm{signal-orthogonal component}}\,.
\end{equation}
On the compactification event of
Proposition~\ref{prop:mdypl-po-link},
$(\hat\theta,\hat\bbeta)=(\thetady,\betady)$.
The localisation argument identifies the limits of
$\hat\bu=\bb\Gamma_p^{1/2}\hat\bv$,
$\hat\sigma$ and $\hat\theta$, while
Lemma~\ref{lemma:mdypl-gaussian-l2-coupling} converts
\eqref{eq:betady-comp} into a standard Gaussian coupling for the actual
MDYPL slope.

The route from the estimator \eqref{eq:og} to these limits comes from a CGMT comparison with two branches. The left
branch embeds MDYPL into a PO, whose optimiser is the object we want
to localise. The right branch studies the corresponding AO, where the actual
scalar and limiting calculations are performed. The two branches are then
brought back together through restricted value comparisons, as summarised in
Figure~\ref{fig:roadmap}.

\definecolor{roadmapConstruct}{HTML}{EAF3FF}
\definecolor{roadmapConstructBorder}{HTML}{3E6C9A}
\definecolor{roadmapProb}{HTML}{FFF3D9}
\definecolor{roadmapProbBorder}{HTML}{A36A00}
\definecolor{roadmapLimit}{HTML}{E8F5EA}
\definecolor{roadmapLimitBorder}{HTML}{3D7A4A}
\definecolor{roadmapConclusion}{HTML}{F1EAFE}
\definecolor{roadmapConclusionBorder}{HTML}{65499C}

\begin{figure}[H]
\centering
\begin{tikzpicture}[
    >={Stealth[length=1.8mm,width=1.25mm]},
    base/.style={
        draw=black!85,
        line width=0.45pt,
        rounded corners=2pt,
        align=center,
        inner xsep=3pt,
        inner ysep=2.3pt,
        minimum height=7mm,
        font=\footnotesize
    },
    leftbox/.style={base,text width=36mm},
    midbox/.style={base,text width=35mm},
    rightbox/.style={base,text width=49mm},
    widebox/.style={base,text width=111mm},
    construction/.style={
        fill=roadmapConstruct,
        draw=roadmapConstructBorder
    },
    probabilistic/.style={
        densely dashed,
        fill=roadmapProb,
        draw=roadmapProbBorder
    },
    deterministic/.style={
        rounded corners=0.5pt,
        fill=roadmapLimit,
        draw=roadmapLimitBorder
    },
    conclusion/.style={
        line width=0.85pt,
        fill=roadmapConclusion,
        draw=roadmapConclusionBorder
    },
    branch/.style={
        font=\footnotesize\bfseries,
        align=center
    },
    flow/.style={->,line width=0.48pt},
    compare/.style={<->,densely dashed,line width=0.52pt},
    support/.style={->,densely dotted,line width=0.48pt}
]

\node[leftbox,construction] (mdypl)
{\textbf{Original MDYPL estimator}\\
$(\thetady,\betady)$, \eqref{eq:og}};

\node[leftbox,construction,below=3.5mm of mdypl] (contain)
{\textbf{Containment and compactification}\\
$\mathcal M_n\cap\mathcal E_n$, \eqref{event:M_n}, \eqref{event:E_n}\\
Proposition~\ref{prop:mdypl-po-link}};

\node[leftbox,construction,below=3.5mm of contain] (po)
{\textbf{PO}\\
$\Psi(\bH_2)$, \eqref{eq:PO}\\
PO optimiser to be localised};

\node[midbox,probabilistic,right=8mm of po] (bridge)
{\textbf{Conditional PO--AO CGMT bridge}\\
conditional on $\mathcal C_n$\\
Proposition~\ref{thm:pathwise_cgmt_localisation}};

\node[rightbox,construction,right=8mm of bridge] (ao)
{\textbf{Conditional AO}\\
$\phi(\bg,\bh)$, \eqref{eq:AO}};

\node[rightbox,construction,below=2.8mm of ao] (scalarise)
{\textbf{AO scalarisation}\\
vector AO $\to(\sigma,r,\bu,\theta,t)$\\
Section~\ref{sec:ao-scalar}};

\node[rightbox,probabilistic,below=2.8mm of scalarise] (finiteao)
{\textbf{Scalarisation event and finite scalar AO}\\
$\mathcal R_n$, \eqref{event:R_n};
Proposition~\ref{prop:scalarisation_events}\\
$\phi(\bg,\bh)=\operatorname{val}(F_n;\bb\Omega_n)$,
\eqref{eq:AO_as_val_Fn}};

\node[rightbox,probabilistic,below=2.8mm of finiteao] (control)
{\textbf{Uniform convergence and moving-domain control}\\
$\mathcal U_n(s)$ and $\mathcal N_n(s)$\\
Proposition~\ref{prop:uniform_AO_control}};

\node[rightbox,deterministic,below=2.8mm of control] (limit)
{\textbf{Limiting AO}\\
$\bar\phi=\operatorname{val}(F;\bb\Omega^*)$,
\eqref{eq:bar_phi_val}};

\node[rightbox,deterministic,below=2.8mm of limit] (saddle)
{\textbf{Limiting-AO properties}\\
attainment, unique saddle, domain relaxations, FOCs\\
Proposition~\ref{prop:limiting_AO_properties}};

\node[rightbox,deterministic,below=2.8mm of saddle] (gaps)
{\textbf{Limiting value gaps on bad sets}\\
Lemma~\ref{lemma:limiting_AO_gap}};

\node[rightbox,probabilistic,below=2.8mm of gaps] (aobounds)
{\textbf{Restricted AO value bounds}\\
finite convex bad-set pieces\\
Lemma~\ref{lemma:restricted_AO_value_bounds}};

\node[midbox,probabilistic] (transfer) at (bridge |- aobounds)
{\textbf{Restricted CGMT transfer}\\
AO barriers $\Longrightarrow$ PO barriers\\
Lemma~\ref{lemma:restricted_PO_value_barriers}};

\node[leftbox,conclusion] (local) at (po |- transfer)
{\textbf{PO localisation}\\
$(\hat\bw,\hat\bv,\hat\theta)\in\mathcal S^\epsilon$\\
Proposition~\ref{prop:PO_localisation}};

\coordinate (bottomcenter) at ($(local.west)!0.5!(aobounds.east)$);

\node[widebox,conclusion,anchor=north] (test)
at ($(bottomcenter)+(0,-13mm)$)
{\textbf{Gaussian coupling and convergence on test functions}\\
spherical PO direction, pathwise Gaussian surrogate,\\
pseudo-Lipschitz transfer and return to MDYPL};

\node[widebox,conclusion,below=3mm of test] (theorem)
{\textbf{Theorem~\ref{thm:mdypl-convergence}}\\
$\thetady\to\theta^*$ and the limiting empirical law of
$\betady-\bbeta^*$};

\node[branch,above=3mm of mdypl] (lefttitle)
{MDYPL / PO branch};
\node[branch] at (bridge |- lefttitle)
{CGMT bridge};
\node[branch] at (ao |- lefttitle)
{AO analysis branch};

\draw[flow] (mdypl) -- (contain);
\draw[flow] (contain) -- (po);

\draw[compare] (po) -- (bridge);
\draw[compare] (bridge) -- (ao);

\draw[flow] (ao) -- (scalarise);
\draw[flow] (scalarise) -- (finiteao);
\draw[flow] (finiteao) -- (control);
\draw[flow] (control) -- (limit);
\draw[flow] (limit) -- (saddle);
\draw[flow] (saddle) -- (gaps);
\draw[flow] (gaps) -- (aobounds);

\draw[support] (control.east) -- ++(5mm,0) |- (aobounds.east);

\draw[flow,densely dashed] (bridge) -- (transfer);
\draw[flow] (aobounds) -- (transfer);
\draw[flow] (transfer) -- (local);

\draw[flow] (po) -- (local);

\draw[flow] (local.south) -- (local.south |- test.north);
\draw[support] (contain.west) -- ++(-6mm,0) |- (test.west);
\draw[flow] (test) -- (theorem);

\end{tikzpicture}
\caption{Proof roadmap for Theorem~\ref{thm:mdypl-convergence}.
The left branch embeds MDYPL, given a containment event, into a compact primary
optimisation (PO).
The middle branch applies the conditional CGMT to unrestricted and
restricted PO--AO values.
The right branch scalarises the auxiliary optimisation (AO), controls the
finite and moving domains, identifies the limiting AO and establishes
attainment, its unique saddle, the required domain relaxations, and the
associated first-order conditions, and obtains
bad-set gaps and restricted AO bounds.
These bounds localise the unique PO minimiser.
Its exact spherical representation and Gaussian surrogate feed the
test-function assembly, which combines them with the containment event
into the theorem.
Blue nodes denote constructions and reformulations, amber dashed nodes
denote probabilistic finite-$n$ comparisons and controls, green nodes
denote deterministic limiting-AO arguments, and violet nodes
denote conclusions feeding into
Theorem~\ref{thm:mdypl-convergence}.}
\label{fig:roadmap}
\end{figure}

In words, the hierarchy in Figure~\ref{fig:roadmap} is as follows. First,
Lemma~\ref{lemma:containment} and Proposition~\ref{prop:mdypl-po-link} put the
MDYPL estimator on the event $\mathcal M_n\cap\mathcal E_n$ inside a compact convex--concave
primary optimisation. This is the PO branch of the argument: it produces the
PO optimiser that we ultimately want to show belongs to
$\mathcal S^\epsilon$. Second, after conditioning on
$\mathcal C_n=\sigma(\bb B,\theta_0,\theta_P,\bH_1,\bb\varepsilon)$, the pair
\eqref{eq:PO}--\eqref{eq:AO} is a valid pathwise CGMT instance, with $\bH_2$
random in the PO and $(\bg,\bh)$ random in the AO. Third, the AO branch is
analysed on its own: Proposition~\ref{prop:scalarisation_events} reduces it to
the scalar problem in $\bb\omega=(\sigma,r,\bu,\theta,t)$.
Proposition~\ref{prop:uniform_AO_control} connects the scalar problem to
$\bar\phi=\operatorname{val}(F;\bb\Omega^*)$, and
Proposition~\ref{prop:limiting_AO_properties} establishes attainment, the
unique limiting saddle, the required domain relaxations, and the population
equations \eqref{eq:FOCs}. Finally, the deterministic gaps of
Lemma~\ref{lemma:limiting_AO_gap} and the restricted AO bounds of
Lemma~\ref{lemma:restricted_AO_value_bounds} are inserted back into the CGMT
comparison. This rules out PO minimisers outside
$\mathcal S^\epsilon$ and gives localisation by
Proposition~\ref{prop:PO_localisation}. The spherical representation in
Lemma~\ref{lemma:compact_PO_spherical_representation} and the Gaussian
surrogate in \eqref{eq:mdypl-gaussian-surrogate}, with law and approximation
bound \eqref{eq:gaussian-surrogate-law} and
\eqref{eq:gaussian-surrogate-direction-bound}, supply the Gaussian
replacement of the localised PO slope, and
Proposition~\ref{prop:test_function_convergence} assembles it,
conditionally on $\mathcal G_n$, with the containment and full-rank event
$\mathcal M_n\cap\mathcal E_n$, the conditional triangular-array law, and
the $W_2$ signal replacement into Theorem~\ref{thm:mdypl-convergence}.
Lemma~\ref{lemma:mdypl-gaussian-l2-coupling} records the resulting vector
coupling for the MDYPL slope.

\subsection{CGMT} \label{sec:CGMT}

For completeness, we restate the form of the CGMT used in the conditional
PO--AO comparison. Its application to the present problem is recorded in
Proposition~\ref{thm:pathwise_cgmt_localisation}.

\begin{theorem}[CGMT, \cite{thrampoulidis+etal:2018} Theorem~VI.1, (i)--(ii)]
\label{thm:cgmt}
Consider the primary optimisation (PO) and auxiliary optimisation (AO)
\begin{equation}
\label{eq:CGMT_PO}
\Phi(\bb G)
=
\underset{\bw \in \mathcal S_{\bw}}{\min}
\underset{\bu \in \mathcal S_{\bu}}{\max}
\, \left\{ \bu^\top \bb G \bw + \psi(\bw,\bu) \right\} \,,
\end{equation}
and
\begin{equation}
\label{eq:CGMT_AO}
\phi(\bg,\bh)
=
\underset{\bw \in \mathcal S_{\bw}}{\min}
\underset{\bu \in \mathcal S_{\bu}}{\max}
\, \left\{ \vnorm{\bw}_2 \, \bg^\top \bu + \vnorm{\bu}_2 \, \bh^\top \bw + \psi(\bw,\bu) \right\} \,,
\end{equation}
where $\bb G \in \Re^{m \times n}$, $\bg \in \Re^{m}$, $\bh \in \Re^{n}$, $\mathcal S_{\bw} \subset \Re^{n}$,
$\mathcal S_{\bu} \subset \Re^{m}$, and $\psi : \Re^{n} \times \Re^{m} \to \Re$.
Let $\bw_{\Phi} = \bw_{\Phi}(\bb G)$ and $\bw_{\phi} = \bw_{\phi}(\bg,\bh)$ denote any optimal minimisers of
\eqref{eq:CGMT_PO} and \eqref{eq:CGMT_AO}, respectively.

Assume that $\mathcal S_{\bw}$ and $\mathcal S_{\bu}$ are compact, that $\psi$ is continuous on
$\mathcal S_{\bw} \times \mathcal S_{\bu}$, and that $\bb G,\bg,\bh$ have i.i.d.\ standard normal entries.

\begin{enumerate}[label=(\roman*)]
	\item For every $c \in \Re$,
	$$
	\Pr\left( \Phi(\bb G) < c \right)
	\leq 2 \Pr\left( \phi(\bg,\bh) \leq c \right) \,.
	$$

	\item If, in addition, $\mathcal S_{\bw}$ and $\mathcal S_{\bu}$ are convex and $\psi$ is convex--concave on
	$\mathcal S_{\bw} \times \mathcal S_{\bu}$, then for every $c \in \Re$,
	$$
	\Pr\left( \Phi(\bb G) > c \right)
	\leq 2 \Pr\left( \phi(\bg,\bh) \geq c \right) \,.
	$$
	In particular, for every $\mu \in \Re$ and $t > 0$,
	$$
	\Pr\left( \abs{\Phi(\bb G) - \mu} > t \right)
	\leq 2 \Pr\left( \abs{\phi(\bg,\bh) - \mu} \geq t \right) \,.
	$$
\end{enumerate}
\end{theorem}

\subsection{Conditioning architecture}
\label{sec:good-events}

To employ the CGMT machinery of Theorem~\ref{thm:cgmt}, we work with a collection of good events.
The proof is organised in two stages. First, we derive the PO--AO comparison and the scalar AO reduction on these events.
Second, after the limiting AO has been defined, we prove that these events hold with sufficiently large probability.

The guiding convention is the following. Let
\begin{equation}
	\label{eq:sigma-fields}
		\mathcal{G}_n = \sigma(\bb B, \theta_0, \theta_P),
	\quad
	\mathcal{C}_n = \sigma(\bb B, \theta_0, \theta_P, \bH_1, \bb \varepsilon) \,.
\end{equation}
The field $\mathcal C_n$ is used for the pathwise CGMT comparison, because conditional on $\mathcal C_n$ the function $\psi$ in the PO is fixed and $\bH_2$ remains a fresh Gaussian matrix.
The field $\mathcal G_n$ is used for the probability bounds on the AO events, because conditional on $\mathcal G_n$ the randomness in $(\bH_1,\bb\varepsilon,\bg,\bh)$ is still available.

\subsection{Setting up the PO} \label{sec:po}

The purpose of this section is to replace the original MDYPL problem
\eqref{eq:og}, on a high-probability containment event, by a compact
convex--concave primary optimisation suitable for CGMT.
The output of the section is
Proposition~\ref{prop:mdypl-po-link}.

To bring the optimisation problem of \eqref{eq:og} into a CGMT-friendly min--max form, let us first rewrite it as a constrained optimisation problem
\begin{equation}
	\label{eq:og-constr}
	\underset{\substack{\bbeta \in \Re^p \\ \theta \in \Re \\ \bb \eta \in \Re^n}}{\min} \, \frac{1}{n} \left\{ \bb 1^\top \bb \rho(\bb{\eta}) - \bYtil^\top \bb{\eta} \right\}, \quad \textrm{such that } \bb{\eta} = \theta \bb 1 + \frac{1}{\sqrt{p}} \bH \bbeta, \quad \bH = \sqrt{p} \bX \,.
\end{equation}
The Lagrangian of \eqref{eq:og-constr} is then
\begin{equation}
	\label{eq:og-lagrangian}
	\underset{\substack{\bbeta \in \Re^p \\ \theta \in \Re \\ \bb \eta \in \Re^n}}{\min} \underset{\bb \lambda \in \Re^n}{\max} \, \frac{1}{n} \left\{ \bb 1^\top \bb \rho(\bb{\eta}) - \bYtil^\top \bb{\eta} + \bb \lambda^\top \left(\bb{\eta} - \theta \bb 1 - \frac{1}{\sqrt{p}}\bH \bbeta \right)\right\} \,.
\end{equation}
Next, let us introduce a containment event on the MDYPL estimates that allows us to constrain the min- and max-variables of \eqref{eq:og-lagrangian} to CGMT-conformal, convex sets.
To start, by \citet[][Theorem~1]{rigon+aliverti:2023}, the MDYPL estimates exist and are unique, whenever the matrix $\widetilde{\bH} = \begin{bmatrix}
	\bb 1 & \bH
	\end{bmatrix}$ has full column rank. Thus define the event
\begin{equation}
	\label{event:E_n}
	\mathcal{E}_n = \left\{ \rank\left(\begin{bmatrix}
	\bb 1 & \bH
	\end{bmatrix}\right) = p + 1 \right\} \,.
\end{equation}
Since $p/n\to\kappa\in(0,1)$, one has $n\geq p+1$ for all sufficiently
large $n$. For every such $n$, $\Pr(\mathcal E_n)=1$, since, conditionally
on the preceding columns, each Gaussian column of $\bH$ has a density and
the span of $\bb1$ and those columns is a proper subspace of $\Re^n$.
Further, define $\etady = \thetady \bb 1 + \bH \betady / \sqrt{p}$ and $\hat{\blambda}^{\textrm{\tiny DY}} = \bYtil - \bb \rho'(\etady)$.
Define the compact containment event
\begin{equation}
	\label{event:M_n}
		\mathcal{M}_n
	=
	\left\{
		\frac{\vnorm{\betady}_2}{\sqrt p} \leq C_\beta,\,
		\abs{\thetady} \leq C_\theta,\,
		\frac{\vnorm{\etady}_2}{\sqrt n} \leq C_\eta,\,
		\frac{\vnorm{\hat{\blambda}^{\textrm{\tiny DY}}}_2}{\sqrt n} \leq 1
	\right\} \,,
\end{equation}
where the constants $C_\beta,C_\theta,C_\eta$ are chosen as in
Lemma~\ref{lemma:containment}, and may be enlarged later without invalidating
the containment bound.
We henceforth enlarge $C_\beta$ and $C_\theta$, if
necessary, so that $C_\beta, C_\theta>B_a$
where $B_a$ is the deterministic constant from
Lemma~\ref{lemma:limiting_AO_closed_saddle}.

Now, denoting the $d$-dimensional Euclidean ball centred at zero, with radius $C \sqrt{d}$, as
$\mathcal{B}^d_C = \{\bv \in \Re^d: \vnorm{\bv}_2 / \sqrt{d} \leq C \}$, on $\mathcal E_n \cap \mathcal M_n$, the optimiser of
\eqref{eq:og-lagrangian} coincides with the optimiser of
\begin{equation}
    \label{eq:lagrangian2}
    \underset{\substack{\bbeta \in \mathcal{B}^p_{C_\beta} \\ \abs{\theta} \leq C_\theta \\ \bb{\eta} \in \mathcal{B}^n_{C_\eta}}}{\min} \underset{\blambda \in \mathcal{B}^n_1}{\max} \,  \frac{1}{n} \left\{\bb{1}^\top \bb{\rho}(\bb{\eta})-  \bYtil^\top \bb{\eta}+ \bb{\lambda}^\top \left(\bb{\eta} - \theta \bb{1} - \frac{1}{\sqrt{p}} \bb{H} \bbeta \right) \right\} \,.
\end{equation}

Now note that the bilinear form $\bb{\lambda}^\top \bb{H} \bbeta$ is not independent of the vector $\bYtil$ whose elements \eqref{eq:mdypl},
depend on $\bb X = \bb{H} / \sqrt{p}$. Hence, following the idea of \citet{salehi+et+al:2019}, we decompose $\bbeta$, and consequently $\bH\bbeta$, into a part that depends on $\bYtil$ and one that is independent.
Towards this, for $ \bB = [\bbeta_0, \bbeta_P]$, define the orthogonal projector on the column space of $\bB$:
\begin{equation}
	\label{eq:b-proj}
	 \bb{P} = \bb{B}  \left(\bb{B}^\top \bb{B} \right)^{+} \bb{B}^\top \,,
\end{equation}
where $\bA^{+}$ is the Moore--Penrose pseudoinverse of $\bA$ (cf. \citealt[Chapter~2]{magnus+neudecker:2019}). Let $\bb{P}^\perp = \bb{I}_p - \bb{P}$ 
and for $s=\rank(\bb B)$, fix once and for all, separately on each
rank stratum $\{\rank(\bb B)=s\}$, a Borel-measurable map
$$
	\bb B
	\mapsto
	\bb E(\bb B)\in\Re^{p\times(p-s)}: \quad 
	\bb E\bb E^\top=\bb P^\perp,
	\quad
	\bb E^\top\bb E=\bb I_{p-s}\,.
$$
For example, a Borel-measurable choice is obtained by applying 
Gram--Schmidt to $\bb P^\perp\bb e_1,\ldots,\bb P^\perp\bb e_p$,
selecting the smallest admissible index at each step and retaining the
first $p-s$ nonzero residuals.
Then, for any $\bbeta \in \Re^p$, we can uniquely decompose the linear predictor
\begin{equation}
\label{eq:H_beta}
    \bb{H} \bbeta = \bb{H} \bb{P} \bbeta +\bb{H}  \bb{P}^\perp  \bbeta = \bb{H}\bB (\bB^\top \bB)^{+} \bB^\top \bbeta  + \bb{H} \bb{E}\bb{E}^\top \bbeta = \bb{H}_1 \bv + \bb{H}_2 \bw \,,
\end{equation}
for $\bb{H}_1 = \bb H \bB, \bb H_2 = \bb H \bb E$, and $\bv =  \left(\bb{B}^\top \bb{B} \right)^{+} \bb{B}^\top  \bbeta \in \operatorname{range}(\bb B^\top) \subseteq \Re^2, \bw = \bb E^\top \bbeta \in \Re^{p-s}$.

By construction of $\bb E$, the rows of $\bH_2$ are of the form $\bb h_i^\top \bb E $ and are thus i.i.d. $\mathrm{N}(\bb 0_{p - s}, \bb I_{p - s})$. Similarly, the rows of $\bH_1$ are of the form
$\bb h_i^\top \bB$, which are i.i.d. $\mathrm{N}(\bb{0}_2, \bB^\top \bB)$. Finally, note that conditional on $\bb B$,
$$
	\cov \left(\bb B^\top \bh_i, \bb E^\top \bh_j\right) = \bb B^\top \bb E = \bb{0}_{2 \times (p - s)} \,,
$$
by Lemma~\ref{lemma:orth_decomp} and therefore $\bb H_1$, $\bb H_2$ are independent.
Furthermore, since $\bH\bbeta_0$ and $\bH\bbeta_P$ are the two columns of $\bH_1=\bH \bb B$,
the response vector $\bYtil$ is measurable with respect to $\sigma(\bb B,\theta_0,\theta_P,\bH_1,\bb \varepsilon)=\mathcal C_n \,,$
and conditional on $\mathcal{C}_n$, $\bH_2$ is independent of $\bYtil$.

Since, by Lemma~\ref{lemma:orth_decomp}, the decomposition $\bbeta = \bB \bv + \bb E \bw$ is componentwise unique for $\bv \in \range(\bb B^\top)$, we can rewrite \eqref{eq:lagrangian2} equivalently as
\begin{equation}
\label{eq:lagrangian3}
    \underset{\substack{(\bv, \bw) \in \mathcal K^p_{C_\beta} \\ \abs{\theta} \leq C_\theta \\ \bb{\eta} \in \mathcal{B}^n_{C_\eta}}}{\min} \underset{\blambda \in \mathcal{B}^n_1}{\max} \,  \frac{1}{n} \left\{\bb{1}^\top \bb{\rho}( \bb{\eta})-  \bYtil^\top\bb{\eta}+ \bb{\lambda}^\top \left(\bb{\eta} - \theta \bb{1} - \frac{1}{\sqrt{p}} \bb{H}_1 \bv - \frac{1}{\sqrt{p}}\bb{H}_2 \bw \right) \right\} \,,
\end{equation}
where $\mathcal{K}_{C_\beta}^p = \{\bv \in \range(\bb B^\top) \subseteq \Re^2, \bw \in \Re^{p-s}: \vnorm{\bB \bv + \bb E \bw}_2^2 \leq {C_\beta}^2 p \}$.
To obtain a product domain, choose $C_\sigma\geq C_\beta$,
enlarging it if necessary so that $C_\sigma>B_a$, and define
\begin{equation}
\label{eq:W_p}
    V_p
    =
    \left\{
        \bv\in\range(\bB^\top):
        \vnorm{\bB\bv}_2\leq C_\beta\sqrt p
    \right\},
    \quad
    W_p
    =
    \left\{
        \bw\in\Re^{p-s}:
        \vnorm{\bw}_2\leq C_\sigma\sqrt p
    \right\} \,.
\end{equation}
Since $\bB^\top\bE=\bb0$ and
$\bE^\top\bE=\bb I_{p-s}$, every
$(\bv,\bw)\in\mathcal K^p_{C_\beta}$ satisfies
$$
    \vnorm{\bB\bv}_2^2+\vnorm{\bw}_2^2
    =
    \vnorm{\bB\bv+\bE\bw}_2^2
    \leq
    C_\beta^2p \,,
$$
and hence
$
    \mathcal K^p_{C_\beta}
    \subseteq
    V_p\times W_p
$.
Both $V_p$ and $W_p$ are nonempty, compact, and convex: $W_p$ is a Euclidean ball, while $V_p$ is the closed ball of radius $C_\beta\sqrt p$ for the norm $\bv\mapsto\vnorm{\bB\bv}_2$ on $\range(\bB^\top)$, since $\range(\bB^\top)\cap\ker(\bB)=\{\bb0_2\}$.
The particular form of $V_p$ is required
to ensure compactness and to allow for a uniform treatment of all
constellations of prior and signal alignment.
By Lemma~\ref{lemma:product_domain_relaxation}, replacing the coupled coordinate set $\mathcal K^p_{C_\beta}$ in
\eqref{eq:lagrangian3} by the Cartesian product $V_p\times W_p$ does not change the optimiser in
$(\theta,\bbeta,\bb{\eta})$ on $\mathcal M_n\cap\mathcal E_n$.
Hence, on $\mathcal M_n \cap \mathcal E_n$,
\eqref{eq:lagrangian3} is equivalent to
\begin{equation}
    \label{eq:lagrangian4}
    \underset{\substack{\bv \in V_p \\ \bw \in W_p \\
    \abs{\theta} \leq C_\theta \\ \bb{\eta} \in \mathcal{B}^n_{C_\eta}} }{\min}
    \underset{\blambda \in \mathcal{B}^n_1}{\max} \,  \frac{1}{n}
    \left\{\bb{1}^\top \bb{\rho}( \bb{\eta})-  \bYtil^\top\bb{\eta}+ \bb{\lambda}^\top \left(\bb{\eta} -
    \theta \bb{1} - \frac{1}{\sqrt{p}} \bb{H}_1 \bv - \frac{1}{\sqrt{p}}\bb{H}_2 \bw \right) \right\} \,.
\end{equation}
For each fixed $\bw\in W_p$, the objective in
\eqref{eq:lagrangian4} is continuous on the compact product domain,
convex in $(\bv,\theta,\bb\eta)$, and affine, hence concave, in
$\blambda$, because $\rho$ is convex and all remaining terms are
linear. Hence Sion's min--max theorem
(see, for example, \citealt[][Theorem~3]{simons:1995}) applies
to the compact convex sets
$V_p\times[-C_\theta,C_\theta]\times\mathcal B^n_{C_\eta}$
and $\mathcal B_1^n$.
Thus, we get that
\begin{equation}
\label{eq:lagrangian5}
    \begin{aligned}
        &\underset{\substack{\bv \in V_p\\ \bw \in W_p\\ \abs{\theta} \leq C_\theta \\ \bb{\eta} \in
        \mathcal{B}^n_{C_\eta}}}{\min} \underset{\blambda \in \mathcal{B}^n_1}{\max} \,  \frac{1}{n}
        \left\{\bb{1}^\top \bb{\rho}( \bb{\eta})-  \bYtil^\top\bb{\eta}+ \bb{\lambda}^\top \left(\bb{\eta} -
        \theta \bb{1} - \frac{1}{\sqrt{p}} \bb{H}_1 \bv - \frac{1}{\sqrt{p}}\bb{H}_2 \bw \right) \right\} \\
        &= \underset{\substack{ \bw \in W_p}}{\min} \underset{\blambda \in \mathcal{B}^n_1}{\max} \,
        -\frac{1}{n \sqrt{p}} \blambda^\top \bH_2 \bw  + \underset{\substack{\bv \in V_p\\ \abs{\theta}
        \leq C_\theta \\ \bb{\eta} \in \mathcal{B}^n_{C_\eta}}}{\min} \, \frac{1}{n} \left\{\bb 1^\top \bb \rho(\bb{\eta}) -
        \bYtil^\top \bb{\eta} +\blambda^\top \bb{\eta} - \theta \bb 1^\top \blambda - \frac{1}{\sqrt{p}}\blambda^\top \bH_1 \bv  \right\} \\
        &= \underset{\substack{ \bw \in W_p}}{\min} \underset{\blambda \in \mathcal{B}^n_1}{\max} \,
        -\frac{1}{n \sqrt{p}} \blambda^\top \bH_2 \bw  + \psi(\blambda) \,,
    \end{aligned}
\end{equation}
for
$$
	\psi(\blambda) = \underset{\substack{\bv \in V_p\\ \abs{\theta}
        \leq C_\theta \\ \bb{\eta} \in \mathcal{B}^n_{C_\eta}}}{\min} \, \frac{1}{n} \left\{\bb 1^\top \bb \rho(\bb{\eta}) -
        \bYtil^\top \bb{\eta} +\blambda^\top \bb{\eta} - \theta \bb 1^\top \blambda - \frac{1}{\sqrt{p}}\blambda^\top \bH_1 \bv  \right\}  \,.
	$$

By Lemma~\ref{lemma:psi_conv_conc_concavity}, applied with $Z=Z_p$,
$S_{\blambda}=\mathcal B_1^n$, and $S_{\bw}=W_p$, $\psi$ is concave in
$\blambda$ and continuous.
Hence, given our PO-candidate in \eqref{eq:lagrangian5}, we are finally in a place to formulate an AO in line with the CGMT of Theorem~\ref{thm:cgmt}.

We thus conclude this section with the following result about the equivalence of MDYPL from \eqref{eq:og} and the PO of \eqref{eq:lagrangian5}:
\begin{proposition}[MDYPL compactification and PO link]
\label{prop:mdypl-po-link}
There exist constants
$C_\sigma,C_\beta,C_\theta,C_\eta,C,c>0$ and
$N\in\mathbb N$ such that, for all $n>N$,
$$
    \Pr(\mathcal M_n\cap\mathcal E_n)
    \geq
    1-C\exp\{-cn\}\,,
$$
where $\mathcal M_n$,
$\mathcal E_n$ are defined in \eqref{event:M_n} and
\eqref{event:E_n}, respectively, and
$C_\sigma,C_\beta,C_\theta,C_\eta,C,c$ depend only on
$$
    \kappa,\alpha,\bb\Gamma,\theta_0^*,\theta_P^*\,,
$$
and on the constants in the probability bounds of
Section~\ref{sec:setup}.

Furthermore, on $\mathcal E_n$, let $\hat\bw$ be the unique minimiser of
the profiled PO in the final line of \eqref{eq:lagrangian5}, and let
$(\hat\bv,\hat\theta,\hat{\bb\eta})$ be the inner minimising variables at
any saddle point of the corresponding fixed-$\hat\bw$ problem. By
Lemma~\ref{lemma:compact_PO_uniqueness}, these variables are unique.
Extend
$(\hat\bv,\hat\bw,\hat\theta,\hat{\bb\eta})$
by zero on $\lnot\mathcal E_n$, and define
$$
    \hat\bbeta
    =
    \bb B\hat\bv+\bb E\hat\bw\,.
$$
Then, on $\mathcal M_n\cap\mathcal E_n$,
$$
    \hat\theta=\thetady,
    \quad
    \hat\bbeta=\betady,
    \quad
    \hat{\bb\eta}=\etady\,.
$$
Consequently, for all $n>N$,
$$
    \Pr
    \left(
        \{\hat\theta=\thetady\}
        \cap
        \{\hat\bbeta=\betady\}
        \cap
        \{\hat{\bb\eta}=\etady\}
    \right)
    \geq
    1-C\exp\{-cn\}\,.
$$
\end{proposition}

\begin{proof}
By Lemma~\ref{lemma:containment},
$$
    \Pr(\lnot\mathcal M_n)
    \leq
    C\exp\{-cn\}\,.
$$
Since $p/n\to\kappa<1$, the Gaussian augmented design has full column rank
almost surely for all sufficiently large $n$, and hence
$\Pr(\mathcal E_n)=1$. After increasing $N$ if necessary,
$$
    \Pr(\mathcal M_n\cap\mathcal E_n)
    =
    \Pr(\mathcal M_n)
    \geq
    1-C\exp\{-cn\}\,.
$$

Work now on $\mathcal M_n\cap\mathcal E_n$. By
Lemma~\ref{lemma:compact_PO_uniqueness},
$(\hat\bv,\hat\bw,\hat\theta,\hat{\bb\eta})$
is the unique tuple of minimising variables of the compact saddle problem
and hence is a product-domain optimiser.
Lemma~\ref{lemma:product_domain_relaxation} therefore gives
$$
    \hat\theta=\thetady,
    \quad
    \hat\bbeta=\betady,
    \quad
    \hat{\bb\eta}=\etady
$$
on $\mathcal M_n\cap\mathcal E_n$. The final probability bound follows
from the first part.
\end{proof}

\subsection{Conditional PO--AO pair} \label{sec:po-ao-pair}

To apply the CGMT PO--AO value-link to \eqref{eq:lagrangian5}, we require that the function $\psi(\bb \lambda)$ is non-random. Hence, we condition on the sigma field $\mathcal C_n = \sigma(\bb B, \theta_0, \theta_P, \bH_1, \bb \varepsilon)$ from \eqref{eq:sigma-fields}.
Conditional on $\mathcal C_n$, the function $\psi$, the feasible sets, and the response vector $\bYtil$ are deterministic, while $\bb H_2$ is a random matrix with standard normal entries. Our conditional PO is given by
\begin{equation}
    \label{eq:PO}
    \Psi(\bb H_2) = \underset{\substack{ \bw \in W_p}}{\min} \underset{\blambda \in \mathcal{B}^n_1}{\max} \, -\frac{1}{n \sqrt{p}} \blambda^\top \bH_2 \bw  + \psi(\blambda) \,,
\end{equation}
and the corresponding conditional AO is given by
\begin{equation}
    \label{eq:AO}
    \phi(\bb g, \bb h) = \underset{\substack{ \bw \in W_p}}{\min} \underset{\blambda \in \mathcal{B}^n_1}{\max} \, -\frac{1}{n \sqrt{p}} \left\{\vnorm{\bw}_2 \bb g^\top\blambda + \vnorm{\blambda}_2 \bw^\top \bb h \right\} + \psi(\blambda) \,,
\end{equation}
where $\bg \in \Re^n, \bh \in \Re^{p-s}$, and conditional on $\mathcal C_n$, the only remaining PO randomness is $\bb H_2$, and the only remaining AO randomness is $(\bg,\bh)$.
Thus, conditional on $\mathcal C_n$, the PO/AO pair in
\eqref{eq:PO}--\eqref{eq:AO} is a valid pathwise CGMT instance. The
following proposition records the corresponding conditional value
comparison, in a form that can later be applied to the unrestricted
PO and to the restricted values used in Section
\ref{sec:cgmt-cost-comparisons}.

\begin{proposition}[Conditional PO--AO CGMT bridge]
	\label{thm:pathwise_cgmt_localisation}
	Let $W_0$ and $D_0$ be $\mathcal C_n$-measurable compact-valued
	random sets such that almost surely
	$$
		W_0\subseteq W_p,
		\quad
		D_0\subseteq V_p\times[-C_\theta,C_\theta] \,,
	$$
	both sets are nonempty, and $D_0$ is convex. Set
	$$
		Z_0
		=
		Z_p(D_0)
		=
		D_0
		\times
		\mathcal B^n_{C_\eta} \,.
	$$
	For $\bb\zeta=(\bv,\theta,\bb\eta)\in Z_0$, write
	$$
		L(\bb\zeta;\blambda)
		=
		\frac1n
		\left\{
			\bb1^\top\bb\rho(\bb\eta)
			-\bYtil^\top\bb\eta
			+\blambda^\top\bb\eta
			-\theta\bb1^\top\blambda
			-\frac{1}{\sqrt p}\blambda^\top\bH_1\bv
		\right\} \,,
	$$
	and define
	$$
		\psi_{Z_0}(\blambda)
		=
		\underset{\bb\zeta\in Z_0}{\min}
		L(\bb\zeta;\blambda) \,.
	$$
	Let
	$$
		\Psi_{W_0,Z_0}(\bb H_2)
		=
		\underset{\bw\in W_0}{\min}
		\underset{\blambda\in\mathcal B^n_1}{\max} \,
		\left\{
			-\frac{1}{n\sqrt p}\blambda^\top\bH_2\bw
			+
			\psi_{Z_0}(\blambda)
		\right\} \,, 
	$$
	and
	\begin{equation}
		\label{eq:ao_scalarisation_start}
		\phi_{W_0,Z_0}(\bg,\bh)
		=
		\underset{\bw\in W_0}{\min}
		\underset{\blambda\in\mathcal B^n_1}{\max} \,
		\left\{
			-\frac{1}{n\sqrt p}
			\left(
				\vnorm{\bw}_2\bg^\top\blambda
				+
				\vnorm{\blambda}_2\bw^\top\bh
			\right)
			+
			\psi_{Z_0}(\blambda)
		\right\}\,,
	\end{equation}
	where $\bg\in\Re^n$ and $\bh\in\Re^{p-s}$ have i.i.d. standard normal
	entries and are independent of $\mathcal C_n$. Then, almost surely for every
	$a\in\Re$, 
	$$
		\Pr\left(\Psi_{W_0,Z_0}(\bH_2)<a \mid \mathcal C_n\right)
		\leq
		2\Pr\left(\phi_{W_0,Z_0}(\bg,\bh)\leq a \mid \mathcal C_n\right)\,.
	$$
	If, in addition, $W_0$ is convex, then, almost
	surely for every $a\in\Re$, 
	$$
		\Pr\left(\Psi_{W_0,Z_0}(\bH_2)>a \mid \mathcal C_n\right)
		\leq
		2\Pr\left(\phi_{W_0,Z_0}(\bg,\bh)\geq a \mid \mathcal C_n\right)\,.
	$$
\end{proposition}

\begin{proof}
	Work conditionally on $\mathcal C_n$ from \eqref{eq:sigma-fields}.
	Outside a $\mathcal C_n$-measurable null set, $W_0$, $D_0$, $Z_0$,
	$\bYtil$, and $\bH_1$ are deterministic, while $\bH_2$ is an
	independent Gaussian matrix with i.i.d. standard normal entries.

	By Lemma~\ref{lemma:psi_conv_conc_concavity}, applied with
	$Z=Z_0$, $S_{\blambda}=\mathcal B^n_1$, and $S_{\bw}=W_0$,
	the function $\psi_{Z_0}$ is continuous and concave on
	$\mathcal B^n_1$. Set
	$$
		\widetilde W_0
		=
		\left\{
			\frac{\bw}{n\sqrt p}:\bw\in W_0
		\right\}\,.
	$$
	Then $\widetilde W_0$ is compact, and
	$$
		\Psi_{W_0,Z_0}(\bH_2)
		=
		\underset{\tilde\bw\in\widetilde W_0}{\min}
		\underset{\blambda\in\mathcal B^n_1}{\max} \, 
		\left\{
			\blambda^\top(-\bH_2)\tilde\bw
			+
			\psi_{Z_0}(\blambda)
		\right\}\,.
	$$
	Since $-\bH_2$ has the same law as $\bH_2$, Theorem
	\ref{thm:cgmt} (i) gives
	$$
		\Pr\left(\Psi_{W_0,Z_0}(\bH_2)<a\mid\mathcal C_n\right)
		\leq
		2\Pr\left(\widetilde\phi_{W_0,Z_0}(\bg,\bh)\leq a
		\mid\mathcal C_n\right)\,,
	$$
	where
	$$
		\widetilde\phi_{W_0,Z_0}(\bg,\bh)
		=
		\underset{\tilde\bw\in\widetilde W_0}{\min}
		\underset{\blambda\in\mathcal B^n_1}{\max} \, 
		\left\{
			\vnorm{\tilde\bw}_2\bg^\top\blambda
			+
			\vnorm{\blambda}_2\tilde\bw^\top\bh
			+
			\psi_{Z_0}(\blambda)
		\right\}\,.
	$$
	Changing back to $\bw=n\sqrt p \tilde\bw$ and using the symmetry
	$(\bg,\bh)\stackrel{d}{=}(-\bg,-\bh)$ shows that
	$\widetilde\phi_{W_0,Z_0}(\bg,\bh)$ has the same law as
	$\phi_{W_0,Z_0}(\bg,\bh)$. This proves the lower-tail comparison.

	For the upper tail, note that if $W_0$ is convex, then
	$\widetilde W_0$ is convex. Since $\mathcal B^n_1$ is convex and
	Lemma~\ref{lemma:psi_conv_conc_concavity} gives the required
	convex--concave structure, the upper-tail part of
	Theorem~\ref{thm:cgmt} applies and yields
	$$
		\Pr\left(\Psi_{W_0,Z_0}(\bH_2)>a\mid\mathcal C_n\right)
		\leq
		2\Pr\left(\widetilde\phi_{W_0,Z_0}(\bg,\bh)\geq a
		\mid\mathcal C_n\right) \,.
	$$
	Again using the equality in law between $\widetilde\phi_{W_0,Z_0}$ and
	$\phi_{W_0,Z_0}$ gives the stated upper-tail comparison.
\end{proof}
Note that for $W_0=W_p$ and
$D_0=V_p\times[-C_\theta,C_\theta]$, the definition of $\psi_{Z_0}$ coincides
with the reduced function $\psi$ in \eqref{eq:lagrangian5}. The
corresponding values are therefore precisely those in \eqref{eq:PO}
and \eqref{eq:AO}.

Now that we have the desired PO--AO link, we analyse the limiting behaviour
of the AO. The restricted comparisons needed to localise the PO are
carried out in Section~\ref{sec:cgmt-cost-comparisons}, after the scalar
limiting value and the limiting-AO properties have been established in
Sections~\ref{sec:ao-limit} and~\ref{sec:limiting-AO-properties}.

\subsection{AO scalarisation} \label{sec:ao-scalar}

The next step is to reduce the optimisations over vector valued arguments to an equivalent min--max problem with scalar optimisers.
Let $Z_p
	=
	V_p
	\times
	[-C_\theta,C_\theta]
	\times
	\mathcal B^n_{C_\eta}$,
and, for $\bb{\zeta}=(\bv,\theta,\bb\eta)\in Z_p$, define
\begin{equation}
	\label{eq:ao_scalarisation_psi_def}
	\psi(\blambda)
	=
	\underset{\bb{\zeta}\in Z_p}{\min} \,
	\left\{
		G_n(\bb\eta)
		+
		\frac{1}{n}
		\blambda^\top
		\bb a (\bb \zeta)
	\right\} , \quad
	G_n(\bb\eta)
	=
	\frac{1}{n}
	\left\{
		\bb 1^\top\bb\rho(\bb\eta)
		-
		\bYtil^\top\bb\eta
	\right\},
	\quad
	\bb a (\bb \zeta)
	=
	\bb\eta
	-
	\theta\bb 1
	-
	\frac{1}{\sqrt p}
	\bH_1\bv \,.
\end{equation}
Denote by $\phi(\bg,\bh)$ the AO of \eqref{eq:ao_scalarisation_start} with 
$W_0 = W_p$ and $D_0 = V_p\times[-C_\theta, C_\theta]$.
Moreover, since $s\le 2$ and $p\to\infty$, we may increase $N$ if necessary
and assume throughout this scalarisation that $p-s\ge 1$.

Now, as a first scalarisation step, we decompose the noise-aligned component $\bw$ into a length and a direction. Since
$$
	W_p
	=
	\left\{
		\bw\in\Re^{p-s}:
		\vnorm{\bw}_2
		\leq
		C_\sigma\sqrt p
	\right\} \,,
$$
every nonzero $\bw\in W_p$ can be written as
\begin{equation}
	\label{eq:ao_scalarisation_w_decomp}
	\bw
	=
	\sqrt p\,\sigma\bb{q},
	\quad
	\sigma
	=
	\frac{\vnorm{\bw}_2}{\sqrt p}
	\in[0,C_\sigma],
	\quad
	\bb{q}\in\mathbb S^{p-s-1} \,.
\end{equation}
The case $\bw=\bb 0$ is represented by $\sigma=0$, with arbitrary $\bb{q}\in\mathbb S^{p-s-1}$. Substituting \eqref{eq:ao_scalarisation_w_decomp} into \eqref{eq:ao_scalarisation_start} gives
\begin{equation}
	\label{eq:ao_scalarisation_after_w_decomp}
	\begin{aligned}
		\phi(\bg,\bh)
		=
		\underset{\sigma\in[0,C_\sigma]}{\min}
		\underset{\bb{q}\in\mathbb S^{p-s-1}}{\min}
		\underset{\blambda\in\mathcal B^n_1}{\max} \,
		\left\{
			\psi(\blambda)
			-
			\frac{\sigma}{n}
			\bg^\top\blambda
			-
			\frac{\sigma}{n}
			\vnorm{\blambda}_2
			\bb{q}^\top\bh
		\right\} \,.
	\end{aligned}
\end{equation}
For $p > s$ and every $\bb{q}\in\mathbb S^{p-s-1}$, $\bb{q}^\top\bh\leq\vnorm{\bh}_2$ by Cauchy--Schwarz \citep[Theorem~1.37(d)]{rudin+etal:1976}.
Since $\sigma\geq0$ and $\vnorm{\blambda}_2\geq0$, it follows pointwise in $\blambda$ that
\begin{equation}
	\label{eq:ao_scalarisation_q_bound}
	\psi(\blambda)
	-
	\frac{\sigma}{n}
	\bg^\top\blambda
	-
	\frac{\sigma}{n}
	\vnorm{\blambda}_2
	\bb{q}^\top\bh
	\geq
	\psi(\blambda)
	-
	\frac{\sigma}{n}
	\bg^\top\blambda
	-
	\frac{\sigma}{n}
	\vnorm{\blambda}_2
	\vnorm{\bh}_2 \,.
\end{equation}
Equality in \eqref{eq:ao_scalarisation_q_bound} is attained by taking $\bb{q}=\bh/\vnorm{\bh}_2$ if $\bh\neq\bb 0$, and is immediate for every $\bb{q}$ if $\bh=\bb 0$. Hence,
taking maxima over $\blambda$ on both sides yields
$$
\underset{\blambda \in \mathcal{B}_1^n}{\max} \, \psi(\blambda)
	-
	\frac{\sigma}{n}
	\bg^\top\blambda
	-
	\frac{\sigma}{n}
	\vnorm{\blambda}_2
	\bb{q}^\top\bh
	\geq \underset{\blambda \in \mathcal{B}_1^n}{\max} \,
	\psi(\blambda)
	-
	\frac{\sigma}{n}
	\bg^\top\blambda
	-
	\frac{\sigma}{n}
	\vnorm{\blambda}_2
	\vnorm{\bh}_2 \,,
$$
and thus
\begin{equation}
	\label{eq:q-min}
	\underset{\bb q \in \mathbb S^{p-s-1}}{\min}
\underset{\blambda \in \mathcal{B}_1^n}{\max} \, \psi(\blambda)
	-
	\frac{\sigma}{n}
	\bg^\top\blambda
	-
	\frac{\sigma}{n}
	\vnorm{\blambda}_2
	\bb{q}^\top\bh
	\geq \underset{\bb q \in \mathbb S^{p-s-1}}{\min} 	\underset{\blambda \in \mathcal{B}_1^n}{\max} \,
	\psi(\blambda)
	-
	\frac{\sigma}{n}
	\bg^\top\blambda
	-
	\frac{\sigma}{n}
	\vnorm{\blambda}_2
	\vnorm{\bh}_2 \,.
\end{equation}
Since the expression on the right-hand side of \eqref{eq:q-min} does not depend on $\bb q$ and the left-hand side can attain this value, it follows that
\begin{equation}
	\label{eq:ao_lambda_max}
	\phi(\bg,\bh)
	=
	\underset{\sigma\in[0,C_\sigma]}{\min}
	\underset{\blambda\in\mathcal B^n_1}{\max} \,
	\left\{
		\psi(\blambda)
		-
		\frac{\sigma}{n}
		\bg^\top\blambda
		-
		\frac{\sigma}{n}
		\vnorm{\blambda}_2
		\vnorm{\bh}_2
	\right\} \,.
\end{equation}
Now, we seek to decompose $\blambda$ into a radius $\tilde r \in [0, \sqrt n]$ and a direction $\bb d \in \mathbb S^{n -1}$, the latter of which we will profile out in a similar fashion to \eqref{eq:q-min}. 
For this, fix $\sigma\in[0,C_\sigma]$ and expand $\psi$ using \eqref{eq:ao_scalarisation_psi_def}. Define
\begin{equation}
	\label{eq:ao_scalarisation_A_sigma}
	\bb a_\sigma(\bb{\zeta})
	=
	\bb a (\bb \zeta)
	-
	\sigma\bg
	=
	\bb\eta
	-
	\theta\bb 1
	-
	\frac{1}{\sqrt p}
	\bH_1\bv
	-
	\sigma\bg \,.
\end{equation}
Then the inner maximisation in \eqref{eq:ao_lambda_max} is
\begin{equation}
	\label{eq:ao_scalarisation_lambda_min}
	\begin{aligned}
		\underset{\blambda\in\mathcal B^n_1}{\max} \,
		\left\{
			\psi(\blambda)
			-
			\frac{\sigma}{n}
			\bg^\top\blambda
			-
			\frac{\sigma}{n}
			\vnorm{\blambda}_2
			\vnorm{\bh}_2
		\right\}
		=
		\underset{\blambda\in\mathcal B^n_1}{\max}
		\underset{\bb{\zeta}\in Z_p}{\min} \,
		\left\{
			G_n(\bb\eta)
			+
			\frac{1}{n}
			\blambda^\top
			\bb a_\sigma(\bb{\zeta})
			-
			\frac{\sigma\vnorm{\bh}_2}{n}
			\vnorm{\blambda}_2
		\right\} \,.
	\end{aligned}
\end{equation}
Now, we apply Lemma~\ref{lemma:psi_conv_conc_ball_to_radius} (with $R=\sqrt n, \bb b=\sigma\bg, c=\sigma\vnorm{\bh}_2$) pathwise for each $\sigma$.
In particular, we apply two min--max swaps and profile out the direction $\bb d$ using the identity $\max_{\bu\in\mathbb S^{n-1}} \,
	\bu^\top
	\bb a_\sigma(\bb{\zeta})
	=
	\vnorm{
		\bb a_\sigma(\bb{\zeta})
	}_2$, that is just the Cauchy--Schwarz inequality \citep[Theorem~1.37(d)]{rudin+etal:1976}, to get 
\begin{equation}
	\label{eq:ao_scalarisation_lambda_to_rtilde}
	\begin{aligned}
		\underset{\blambda\in\mathcal B^n_1}{\max}
		\underset{\bb{\zeta}\in Z_p}{\min} \,
		\left\{
			G_n(\bb\eta)
			+
			\frac{1}{n}
			\blambda^\top
			\bb a_\sigma(\bb{\zeta})
			-
			\frac{\sigma\vnorm{\bh}_2}{n}
			\vnorm{\blambda}_2
		\right\} 
		=
		\underset{\tilde r\in[0,\sqrt n]}{\max}
		\underset{\bb{\zeta}\in Z_p}{\min} \,
		\left\{
			G_n(\bb\eta)
			+
			\frac{\tilde r}{n}
			\vnorm{
				\bb a_\sigma(\bb{\zeta})
			}_2
			-
			\frac{\sigma\tilde r}{n}
			\vnorm{\bh}_2
		\right\} \,.
	\end{aligned}
\end{equation}
Since \eqref{eq:ao_scalarisation_lambda_to_rtilde} holds for every fixed
$\sigma\in[0,C_\sigma]$, combining it with \eqref{eq:ao_lambda_max}
and then taking the outer minimum over $\sigma$ gives
\begin{equation}
	\label{eq:ao_scalarisation_rtilde}
	\begin{aligned}
		\phi(\bg,\bh)
		=
		\underset{\sigma\in[0,C_\sigma]}{\min}
		\underset{\tilde r\in[0,\sqrt n]}{\max} \,
		\left\{
			-
			\frac{\sigma\tilde r}{n}
			\vnorm{\bh}_2
			+
			\underset{\bb{\zeta}\in Z_p}{\min} \,
			\left[
				\frac{\tilde r}{n}
				\vnorm{
					\bb a_\sigma(\bb{\zeta})
				}_2
				+
				G_n(\bb\eta)
			\right]
		\right\} \,.
	\end{aligned}
\end{equation}
Finally, reparameterise
\begin{equation}
	\label{eq:ao_scalarisation_r_reparam}
	r
	=
	\frac{\tilde r}{\sqrt n}
	\in[0,1] \,.
\end{equation}
Using \eqref{eq:ao_scalarisation_A_sigma} and \eqref{eq:ao_scalarisation_r_reparam} in \eqref{eq:ao_scalarisation_rtilde} yields the scalar AO
\begin{equation}
	\label{eq:ao_r}
	\begin{aligned}
		\phi(\bg,\bh)
		=
		\underset{\sigma\in[0,C_\sigma]}{\min}
		\underset{r\in[0,1]}{\max} \,
			-
			\sigma r
			\frac{\vnorm{\bh}_2}{\sqrt n}
			+
			\underset{
				\substack{
					\bv\in V_p\\
					\abs{\theta}\leq C_\theta\\
					\bb\eta\in\mathcal B^n_{C_\eta}
				}
			}{\min} \,
				\frac{r}{\sqrt n}
				\vnorm{
					\bb\eta
					-
					\theta\bb 1
					-
					\frac{1}{\sqrt p}
					\bH_1\bv
					-
					\sigma\bg
				}_2
				+
				\frac{1}{n}
				\left\{
					\bb 1^\top\bb\rho(\bb\eta)
					-
					\bYtil^\top\bb\eta
				\right\}
		 \,.
	\end{aligned}
\end{equation}
For the analysis of the radial variable $r$ it is convenient to profile
out $(\bv,\theta)$. Recall the empirical loss $G_n(\bb\eta)$ from
\eqref{eq:ao_scalarisation_psi_def} and set
$$
	\bEtaYtil
	=
	\logit{\bYtil}
$$
coordinatewise, so that $\rho'(\eta_{\Ytil,i})=\Ytil_i$ and
$\bEtaYtil$ is the unconstrained minimiser of $G_n$. Define the compact
convex predictor set
$$
	\mathcal R_0
	=
	\left\{
		\theta\bb1
		+
		\frac1{\sqrt p}\bH_1\bv:
		\abs{\theta}\leq C_\theta,
		\ \bv\in V_p
	\right\}\,,
$$
and, for $\sigma\in[0,C_\sigma]$ and $\bb\eta\in\Re^n$,
$$
	M_n(\sigma,\bb\eta)
	=
	\frac1{\sqrt n}
	\operatorname{dist}_2
	\left(
		\bb\eta-\sigma\bg,
		\mathcal R_0
	\right)
	=
	\underset{\substack{\abs{\theta}\leq C_\theta\\ \bv\in V_p}}{\inf}\,
	\frac1{\sqrt n}
	\vnorm{
		\bb\eta-\theta\bb1-\frac1{\sqrt p}\bH_1\bv-\sigma\bg
	}_2\,.
$$
For $r\in[0,1]$, define the radial profile
\begin{equation}
\label{eq:radial_AO_profile}
	A_n(\sigma,r)
	=
	-\sigma r\frac{\vnorm{\bh}_2}{\sqrt n}
	+
	\underset{\bb\eta\in\mathcal B^n_{C_\eta}}{\min}
	\left\{
		G_n(\bb\eta)
		+
		rM_n(\sigma,\bb\eta)
	\right\}\,.
\end{equation}
Since the minimum over $(\bv,\theta)$ in \eqref{eq:ao_r} is exactly
$rM_n(\sigma,\bb\eta)$, the scalar AO is
\begin{equation}
\label{eq:ao_r_profile}
	\phi(\bg,\bh)
	=
	\underset{\sigma\in[0,C_\sigma]}{\min}
	\underset{r\in[0,1]}{\max}\,
	A_n(\sigma,r)\,.
\end{equation}
For $a\in[0,1]$, write
$$
	\phi^r(a;\bg,\bh)
	=
	\underset{\sigma\in[0,C_\sigma]}{\min}
	\underset{r\in[a,1]}{\max}\,
	A_n(\sigma,r)\,,
$$
and, given a constant $c_r>0$, define the event
$$
	\mathcal E_n^r(c_r)
	=
	\left\{
		\phi^r(c_r;\bg,\bh)
		=
		\phi^r(0;\bg,\bh)
		=
		\phi(\bg,\bh)
	\right\}\,.
$$
A sufficient, stronger property is that, for every
$\sigma\in[0,C_\sigma]$,
$$
	\underset{r\in[0,1]}{\max}\,A_n(\sigma,r)
	=
	\underset{r\in[c_r,1]}{\max}\,A_n(\sigma,r)\,,
$$
which holds as soon as $r\mapsto A_n(\sigma,r)$ is nondecreasing on
$[0,c_r]$. This pointwise property is what
Proposition~\ref{prop:small_r_exclusion} establishes, on the good event
introduced before Proposition~\ref{prop:scalarisation_events} below.
Given $\mathcal E_n^r(c_r)$, the AO can equivalently be written as
\begin{equation}
	\label{eq:ao_r_pos}
    \begin{aligned}
        \phi(\bb g, \bb h) = \underset{\sigma \in [0, C_\sigma]}{\min}\underset{r \in [c_r, 1]}{\max} \, -\sigma r\frac{\vnorm{\bb h}_2}{\sqrt{n}}
        +\underset{\substack{\bv \in V_p\\ \abs{\theta} \leq C_\theta\\
        \bb{\eta} \in \mathcal{B}^n_{C_\eta}}}{\min} \,  \frac{r}{\sqrt{n}} \vnorm{\bb{\eta}- \theta \bb 1 - \frac{1}{\sqrt{p}}\bH_1 \bv - \sigma \bb g}_2
        + \frac{1}{n} \left\{\bb 1^\top \bb \rho(\bb{\eta}) - \bYtil^\top \bb{\eta} \right\} \,.
    \end{aligned}
\end{equation}
Next, we use the following identity
\begin{equation}
	\label{eq:norm_variational_identity}
    \vnorm{\bb a}_2 = \underset{t>0}{\inf}\, \frac{ \vnorm{\bb a}_2^2 }{2t}+\frac{t}{2} \,,
\end{equation}
to re-write
$$
    \begin{aligned}
    \frac{r}{\sqrt{n}} \vnorm{\bb{\eta}- \theta \bb 1 - \frac{1}{\sqrt{p}}\bH_1 \bv -
    \sigma\bb g}_2
    &= r \vnorm{\frac{1}{\sqrt{n}}\left(\bb{\eta}- \theta \bb 1 - \frac{1}{\sqrt{p}}\bH_1 \bv -
    \sigma\bb g\right)}_2 \\
    &= r \underset{t  > 0}{\inf}\,  \frac{1}{2t} \vnorm{\frac{1}{\sqrt{n}}\left(\bb{\eta}- \theta \bb 1 - \frac{1}{\sqrt{p}}\bH_1 \bv -
    \sigma\bb g\right)}_2^2 + \frac{t}{2} \\
    &= \underset{t > 0}{\inf}\, \frac{r}{2t} \frac{1}{n}\vnorm{\bb{\eta}- \theta \bb 1 - \frac{1}{\sqrt{p}}\bH_1 \bv -
    \sigma\bb g}_2^2 + \frac{rt}{2} \,.
    \end{aligned}
$$
Now, define
\begin{equation}
	\label{eq:psi-t-def}
		\begin{aligned}
			\phi^{r,t}(a,b; \bg,\bh)
			&=
			\underset{\sigma \in [0, C_\sigma]}{\min}
			\underset{r \in [a, 1]}{\max} \,
			-\sigma r\frac{\vnorm{\bb h}_2}{\sqrt{n}}
			+
			\underset{\substack{\bv \in V_p\\ \abs{\theta} \leq C_\theta\\
			\bb{\eta} \in \mathcal{B}^n_{C_\eta}}}{\min}
			\underset{t \in (0,b]}{\inf} \,
			\frac{r}{2t}\frac{1}{n}
			\vnorm{\bb{\eta}-\theta\bb 1-\frac{1}{\sqrt p}\bH_1\bv-\sigma\bb g}_2^2
			+
			\frac{rt}{2}
			\\
			&+
			\frac{1}{n}
			\left\{\bb 1^\top\bb\rho(\bb{\eta})-\bYtil^\top\bb{\eta}\right\} \,.
		\end{aligned}
\end{equation}
For $b = \infty$ in \eqref{eq:psi-t-def}, the same definition is understood with $t \in (0, \infty)$.
Given constants $c_r, C_t$, introduce the event
$$
	\mathcal E_n^t(c_r, C_t) = \left\{\phi^{r,t}(c_r, C_t; \bg,\bh) = \phi^{r,t}(c_r, \infty; \bg,\bh) = \phi^r(c_r; \bg, \bh) \right\} \,.
$$
Upon $\mathcal E_n^r(c_r) \cap \mathcal E_n^t(c_r,C_t)$, and since the
$t$- and $(\bv,\theta,\bb\eta)$-optimisations are all minimising
operations over independent variables, we can re-express the AO as
\begin{equation}
	\label{eq:ao_t}
	    \begin{aligned}
	        \phi(\bb g, \bb h)
	        &=
	        \underset{\sigma \in [0, C_\sigma]}{\min}
	        \underset{r \in [c_r, 1]}{\max} \,
	        -\sigma r\frac{\vnorm{\bb h}_2}{\sqrt{n}}
	        +
			\underset{t \in (0, C_t]}{\inf}
			\underset{\substack{\bv \in V_p\\ \abs{\theta} \leq C_\theta\\
	        \bb{\eta} \in \mathcal{B}^n_{C_\eta}}}{\min} \,
	        \frac{r}{2t}
	        \vnorm{\frac{1}{\sqrt n}\left(\bb{\eta}-\theta\bb 1-\frac{1}{\sqrt p}\bH_1\bv-\sigma\bb g\right)}_2^2
	        +
	        \frac{rt}{2}
	        \\
	        &+
	        \frac{1}{n}
	        \left\{\bb 1^\top\bb\rho(\bb{\eta})-\bYtil^\top\bb{\eta}\right\} \,.
	    \end{aligned}
\end{equation}
A simple completion of squares yields that
$$
    \begin{aligned}
        \frac{r}{2t} \vnorm{\bb a - \frac{t}{r} \bYtil}_2^2 &= \frac{r}{2t} \vnorm{\bb a}_2^2 - \bb a^\top \bYtil + \frac{t}{2r} \vnorm{\bYtil}_2^2 \,.
    \end{aligned}
$$
Hence, for $\bb a = \bb{\eta}- \theta \bb 1 - \frac{1}{\sqrt{p}}\bH_1 \bv -  \sigma\bb g$, we
get, upon rearrangement, that
\begin{equation}
\label{eq:completion-of-squares}
    \begin{aligned}
        \frac{r}{2t} \vnorm{\bb{\eta}- \theta \bb 1 - \frac{1}{\sqrt{p}}\bH_1 \bv -  \sigma\bb g}_2^2
    - \bb{\eta}^\top \bYtil &=
    \frac{r}{2t} \vnorm{\bb{\eta}- \theta \bb 1 - \frac{1}{\sqrt{p}}\bH_1 \bv -  \sigma\bb g - \frac{t}{r} \bYtil}_2^2 \\
    &- \bYtil^\top\left(\theta \bb 1 + \frac{1}{\sqrt{p}}\bH_1 \bv + \sigma\bb g  \right) \\
    &- \frac{t}{2r} \vnorm{\bYtil}^2_2 \,.
    \end{aligned}
\end{equation}
Hence, the full AO becomes
$$
    \begin{aligned}
        \phi(\bb g, \bb h) 
        &= \underset{\sigma \in [0, C_\sigma]}{\min} \underset{r \in [c_r, 1]}{\max}
        \underset{\substack{\bv \in V_p\\ \abs{\theta} \leq C_\theta\\
        \bb{\eta} \in \mathcal{B}^n_{C_\eta}}}{\min} \underset{t \in (0, C_t]}{\inf} \,
        -\sigma r \frac{\vnorm{\bb h}_2}{\sqrt{n}} + \frac{rt}{2} 
        + \frac{1}{n} \left\{ \frac{r}{2t}\vnorm{\bb{\eta}- \theta \bb 1 - \frac{1}{\sqrt{p}}\bH_1 \bv -  \sigma\bb g - \frac{t}{r} \bYtil}_2^2 
    - \bYtil^\top\left(\theta \bb 1 + \frac{1}{\sqrt{p}}\bH_1 \bv +  \sigma\bb g \right)
    - \frac{t}{2r} \vnorm{\bYtil}^2_2  + \bb 1^\top \bb \rho (\bb{\eta})\right\}  \,.
    \end{aligned}
$$
Lastly, we wish to scalarise the $n$-vector minimisation over $\bb \eta$. 
For fixed $(\sigma,r,\bv,\theta,t)$, the $\bb\eta$-dependent part is
the constrained minimisation
\begin{equation}
	\label{eq:constr_moreau}
		\begin{aligned}
		\underset{\bb{\eta} \in \mathcal{B}^n_{C_\eta}}{\min} \, \frac{1}{n} \left( \bb 1^\top \bb \rho(\bb{\eta}) + \frac{r}{2t} \vnorm{\bb{\eta} - \theta \bb 1 - \frac{1}{\sqrt{p}} \bH_1 \bv - \sigma \bg - \frac{t}{r} \bYtil}_2^2 \right) \,,
	\end{aligned}
\end{equation}
which we want to express as an unconstrained, separable, Moreau envelope. 
For this, define the constrained Moreau envelope
$$
	M_{\rho}^{\mathrm{const}}(\bb b,\lambda)
	=
	\underset{\bb\eta\in\mathcal B^n_{C_\eta}}{\min}
	\left\{
		\bb1^\top\bb\rho(\bb\eta)
		+
		\frac{1}{2\lambda}
		\vnorm{\bb\eta-\bb b}_2^2
	\right\} \,,
$$
and the unconstrained coordinatewise Moreau envelope
$$
	M_{\rho}(\bb b,\lambda)
	=
	\sum_{i=1}^n m_\rho(b_i,\lambda),
	\quad
	m_\rho(b,\lambda)
	=
	\underset{\eta\in\Re}{\min}
	\left\{
		\rho(\eta)
		+
		\frac{1}{2\lambda}(\eta-b)^2
	\right\} \,,
$$
where the case $\lambda=0$, is understood through the continuous extension $M_\rho(\bb b,0)=\bb1^\top\bb\rho(\bb b)$.
Then define the AO value with the constrained Moreau envelope by
$$
\begin{aligned}
	\phi^{M,\mathrm{const}}(c_r,C_t;\bg,\bh)
	=
	\underset{\sigma\in[0,C_\sigma]}{\min}
	\underset{r\in[c_r,1]}{\max}
	\underset{\substack{\bv\in V_p\\ \abs{\theta}\leq C_\theta}}{\min}
	\underset{t\in(0,C_t]}{\inf} \,
	&\frac1n
	M_\rho^{\mathrm{const}}
	\left(
			\theta\bb1
	+
	\frac1{\sqrt p}\bH_1\bv
	+
	\sigma\bg
	+
	\frac tr\bYtil ,
		\frac tr
	\right)
	-\sigma r\frac{\vnorm{\bh}_2}{\sqrt n} +
	\frac{rt}{2} 
	\\
	&-
	\frac1n
	\left\{
		\bYtil^\top
		\left(
			\theta\bb1+\frac1{\sqrt p}\bH_1\bv+\sigma\bg
		\right)
		+
		\frac{t}{2r}\vnorm{\bYtil}_2^2
	\right\} \,.
\end{aligned}
$$
By the completion of squares in \eqref{eq:completion-of-squares}, $\phi^{M,\mathrm{const}}(c_r,C_t;\bg,\bh) = \phi^{r,t}(c_r,C_t;\bg,\bh) \,.$
Similarly, define the AO value with the unconstrained Moreau envelope by
$$
\begin{aligned}
	\phi^{M}(c_r,C_t;\bg,\bh)
	=
	\underset{\sigma\in[0,C_\sigma]}{\min}
	\underset{r\in[c_r,1]}{\max}
	\underset{\substack{\bv\in V_p\\ \abs{\theta}\leq C_\theta}}{\min}
	\underset{t\in(0,C_t]}{\inf} \,
	&\frac1n
	M_\rho
	\left(
			\theta\bb1
	+
	\frac1{\sqrt p}\bH_1\bv
	+
	\sigma\bg
	+
	\frac tr\bYtil ,
		\frac tr
	\right)
	-\sigma r\frac{\vnorm{\bh}_2}{\sqrt n} +
	\frac{rt}{2} 
	\\
	&
	-
	\frac1n
	\left\{
		\bYtil^\top
		\left(
			\theta\bb1+\frac1{\sqrt p}\bH_1\bv+\sigma\bg
		\right)
		+
		\frac{t}{2r}\vnorm{\bYtil}_2^2
	\right\} \,.
\end{aligned}
$$
Now, define the Moreau feasibility event as
$$
	\mathcal E_n^{\mathrm M}(c_r,C_t,C_\eta)
	=
	\left\{
		\phi^{M}(c_r,C_t;\bg,\bh)
		=
		\phi^{M,\mathrm{const}}(c_r,C_t;\bg,\bh) 
		= \phi^{r,t}(c_r,C_t;\bg,\bh)
	\right\} \,.
$$
For constants $c_r>0$, $C_t<\infty$, and $C_\eta<\infty$, define the
scalarisation-valid event
\begin{equation}
	\label{event:R_n}
	\mathcal R_n(c_r,C_t,C_\eta)
	=
	\mathcal E_n^r(c_r)
	\cap
	\mathcal E_n^t(c_r,C_t)
	\cap
	\mathcal E_n^{\mathrm M}(c_r,C_t,C_\eta) \,, 
\end{equation}
where the quantities $\phi^r$, $\phi^{r,t}$, and
$\phi^{M,\mathrm{const}}$, and hence the events
$\mathcal E_n^r$ and $\mathcal E_n^t$, also depend on the fixed truncation
radius $C_\eta$. 
When the constants have been fixed, we write simply
$\mathcal R_n$ for $\mathcal R_n(c_r,C_t,C_\eta)$.

Next, upon $\mathcal R_n$, we wish extend the infimum over $t\in(0,C_t]$ to the compact interval $[0,C_t]$ by continuously extending
$
M_\rho(\bx+\lambda\bba,\lambda)/n
$
to $\lambda=0$, with boundary value
$
\bb1^\top\bb\rho(\bx) / n 
$ for $\bx \in \Re^n$, $\bb a \in (0,1)^n$. To that end, 
note that since $0<\rho'<1$, the function $\rho$ is $1$-Lipschitz. Note further then that for every
$b\in\Re$ and $\lambda>0$,
$$
    0
    \leq
    \rho(b)-m_\rho(b,\lambda)
    \leq
    \frac{\lambda}{2}\,.
$$
Indeed, the first inequality follows by testing $\eta=b$, while
$\rho(\eta)\geq\rho(b)-\abs{\eta-b}$ gives
$$
    m_\rho(b,\lambda)
    \geq
    \rho(b)
    +
    \underset{u\geq0}{\inf}
    \left\{
        -u+\frac{u^2}{2\lambda}
    \right\}
    =
    \rho(b)-\frac{\lambda}{2}\,.
$$
Consequently, for every $\bb x\in\Re^n$ and
$\bba\in(0,1)^n$,
$$
    \frac1n
    \abs{
        M_\rho(\bb x+\lambda\bba,\lambda)
        -
        \bb1^\top\bb\rho(\bb x)
    }
    \leq
    \frac{3\lambda}{2}\,.
$$
Since $\lambda=t/r\leq t/c_r$, the Moreau term, and hence the full
objective defining $\phi^M(c_r,C_t;\bg,\bh)$, extends continuously to
$t=0$. 
Therefore the infimum over $t\in(0,C_t]$ equals the minimum
over $t\in[0,C_t]$.
With that, upon $\mathcal R_n(c_r,C_t,C_\eta)$, our scalarised AO becomes
\begin{equation}
	\label{eq:AO_scalar}
    \begin{aligned}
        \phi(\bb g, \bb h) = \underset{\sigma \in [0, C_\sigma]}{\min} \underset{r \in [c_r, 1]}{\max}
        \underset{\substack{\bv \in V_p\\ \abs{\theta} \leq C_\theta \\ t \in [0, C_t]}}{\min} \,
         &\frac{1}{n} M_{\rho}\left(\theta \bb 1 + \frac{1}{\sqrt{p}} \bH_1 \bv + \sigma \bg + \frac{t}{r} \bYtil, \frac{t}{r} \right)-\sigma r \frac{\vnorm{\bb h}_2}{\sqrt{n}}  
        + \frac{rt}{2} \\ &- \frac{1}{n} \left(  \bYtil^\top\left(\theta \bb 1 + \frac{1}{\sqrt{p}}\bH_1 \bv +  \sigma\bb g \right)  + \frac{t}{2r} \vnorm{\bYtil}^2_2  \right) \,.
    \end{aligned}
\end{equation}
Before stating the concluding result of this section, it will be useful to introduce a structural event about the signal parameters $\bbeta_0, \bbeta_P, \theta_0, \theta_P$.
For fixed $t_0>0$, define
\begin{equation}
	\label{event:S_n}
	\mathcal{S}_n(t_0)
	=
	\left\{
		\mnorm{\bb \Gamma_p - \bb \Gamma}_2 \leq t_0,\,
		\vnorm{
			\begin{bmatrix}
				\theta_0\\
				\theta_P
			\end{bmatrix}
			-
			\begin{bmatrix}
				\theta_0^*\\
				\theta_P^*
			\end{bmatrix}
		}_{\infty}
		\leq t_0
	\right\} \,.
\end{equation}
By Assumptions~\ref{ass:betas_conv} and~\ref{ass:theta-conv}, for every fixed
$t_0>0$, after adjusting constants $C, c > 0$,
$$
	\Pr\left(\lnot \mathcal S_n(t_0)\right)
	\leq
	C\exp\{-cn\} \,.
$$
All constants in the $\mathcal G_n$-conditional AO bounds below are deterministic on $\mathcal S_n(t_0)$.

We also need the good event on which the radial boundary $r=0$ is
excluded. Recall $G_n$, $\bEtaYtil$, $\mathcal R_0$, $M_n$, and the
radial profile $A_n$ from \eqref{eq:radial_AO_profile}. For
$i=1,\ldots,n$, define the envelope
$$
	L_i
	=
	\abs{\eta_{\Ytil,i}}
	+
	C_\theta
	+
	C_\sigma\abs{g_i}
	+
	\underset{\bv\in V_p}{\sup}
	\abs{
		\frac1{\sqrt p}(\bH_1\bv)_i
	}\,,
$$
and, for $a\in(0,1/2)$ and $K>0$, the interior index set
$$
	\mathcal I_n(a,K)
	=
	\left\{
		i\in\{1,\ldots,n\}:
		\Ytil_i\in[a,1-a],
		\ L_i\leq K
	\right\}\,.
$$
For $d>0$, $a\in(0,1/2)$, and $K,C_L>0$, define the small-$r$ environment
\begin{equation}
\label{event:small_r_environment}
\begin{aligned}
	\mathcal E_n^{\mathrm{sr}}(d,a,K,C_L)
	=
	\Bigg\{&
		\underset{\sigma\in[0,C_\sigma]}{\inf}
		\left[
			M_n(\sigma,\bEtaYtil)
			-
			\sigma\frac{\vnorm{\bh}_2}{\sqrt n}
		\right]
		\geq d,
		\\
	&
		\frac1n\sum_{i=1}^nL_i^2
		\leq C_L,
		\qquad
		\frac1n\sum_{i\notin\mathcal I_n(a,K)}L_i^2
		\leq
		\frac{d^2}{64}
	\Bigg\}\,.
\end{aligned}
\end{equation}
Once $C_\eta\geq\sqrt{C_L}$, as ensured below, the first component is
a uniform lower bound on the initial slope
$M_n(\sigma,\bEtaYtil)-\sigma\vnorm{\bh}_2/\sqrt n$ of
$r\mapsto A_n(\sigma,r)$ at $r=0$; the two envelope components control how
far the minimiser of $G_n+rM_n(\sigma,\cdot)$ can move from $\bEtaYtil$ for
small $r>0$. Lemma~\ref{lemma:small_r_environment} shows that this event
has exponentially high conditional probability, and
Proposition~\ref{prop:small_r_exclusion} shows that, on it,
$A_n(\sigma,\cdot)$ is nondecreasing on a deterministic interval
$[0,r_0]$, uniformly in $\sigma$.

\begin{proposition}[Scalarisation event bounds]
\label{prop:scalarisation_events}
Fix $t_0>0$. After possibly increasing $C_\eta$, there exist deterministic
constants
$
	d_r>0$, 
	$
	a\in(0,1/2)$, 
	$
	K,C_L,C_t,C_M\in(0,\infty)$, 
	$
	c_r\in(0,1/2]
$,
such that
$
	C_M\geq\sqrt{C_L}$, 
	$
	C_\eta
	\geq
	C_M+
	{C_t} / {c_r},
$
and
\begin{equation}
\label{eq:scalarisation_restricted_Ct_bound}
	C_t
	>
	\max\left\{
		1,
		\frac{
			(C_\theta+3C_\beta+2C_\sigma)^2
			+
			1
			+
			4\log2
		}{c_r}
	\right\}\,.
\end{equation}
There are also constants $C,c>0$ and $N\in\mathbb N$ such that, for all
$n>N$,
\begin{equation}
\label{eq:scalarisation_event_sr_prob}
	\mathds{1}_{\mathcal S_n(t_0)}
	\Pr\left(
		\lnot\mathcal E_n^{\mathrm{sr}}(d_r,a,K,C_L)
		\mid
		\mathcal G_n
	\right)
	\leq
	C\exp\{-cn\}\,,
\end{equation}
\begin{equation}
\label{eq:scalarisation_event_r_prob}
	\mathds{1}_{\mathcal S_n(t_0)}
	\Pr\left(
		\lnot\mathcal E_n^r(c_r)
		\mid
		\mathcal G_n
	\right)
	\leq
	C\exp\{-cn\}\,,
\end{equation}
\begin{equation}
\label{eq:scalarisation_event_t_prob}
	\mathds{1}_{\mathcal S_n(t_0)}
	\Pr\left(
		\lnot\mathcal E_n^t(c_r,C_t)
		\mid
		\mathcal G_n
	\right)
	\leq
	C\exp\{-cn\}\,,
\end{equation}
and
\begin{equation}
\label{eq:scalarisation_event_M_prob}
	\mathds{1}_{\mathcal S_n(t_0)}
	\Pr\left(
		\lnot\mathcal E_n^{\mathrm M}(c_r,C_t,C_\eta)
		\mid
		\mathcal G_n
	\right)
	\leq
	C\exp\{-cn\}\,.
\end{equation}
For these constants, abbreviate
$$
	\mathcal R_n
	=
	\mathcal R_n(c_r,C_t,C_\eta)\,,
$$
where $\mathcal R_n(c_r,C_t,C_\eta)$ is defined in
\eqref{event:R_n}. Then
\begin{equation}
\label{eq:Rn_cond_prob}
	\mathds{1}_{\mathcal S_n(t_0)}
	\Pr\left(
		\lnot\mathcal R_n
		\mid
		\mathcal G_n
	\right)
	\leq
	C\exp\{-cn\}\,.
\end{equation}
Moreover,
\begin{equation}
\label{eq:Rn_uncond_prob}
	\Pr\left(
		\lnot\mathcal R_n
	\right)
	\leq
	C\exp\{-cn\}
	+
	\Pr\left(
		\lnot\mathcal S_n(t_0)
	\right)
	\leq
	C\exp\{-cn\}\,.
\end{equation}
The constants depend only on
$$
	\kappa,
	\alpha,
	\bb\Gamma,
	\theta_0^*,
	\theta_P^*,
	t_0,
	C_\beta,
	C_\theta,
	C_\sigma\,,
$$
and on the constants in the probability bounds of
Section~\ref{sec:setup}.
\end{proposition}

\begin{proof}
Lemma~\ref{lemma:small_r_environment} gives deterministic constants
$d_r,a,K,C_L$, constants $C,c>0$, and
$N\in\mathbb N$ such that
\eqref{eq:scalarisation_event_sr_prob} holds with these constants. Set
$$
	c_{a,K}
	=
	\rho''(\logit{1-a}+K),
	\qquad
	r_0
	=
	\min\left\{
		\frac12,
		\frac{c_{a,K} d_r}{16}
	\right\}\,.
$$
Let $\bar r$ and $c_0$ be the deterministic constants from
Lemma~\ref{lemma:limiting_AO_closed_saddle} and
Lemma~\ref{lemma:limiting_AO_lower_radius}, respectively, and set
$$
	c_r
	=
	\min\left\{
		r_0,
		\bar r,
		c_0
	\right\}\,.
$$

Let $K_t$ be the constant of Lemma~\ref{lemma:t_ub} and choose $C_t$ so that
\begin{equation}
\label{eq:scalarisation_Ct_choice}
	C_t
	>
	\max\left\{
		B_t,
		1,
		\frac{2(K_t+\log2)}{c_r}
	\right\}\,,
\end{equation}
where $B_t$ is the constant from
Lemma~\ref{lemma:limiting_AO_closed_saddle}; this is the threshold required by Lemma~\ref{lemma:t_ub} and by the restricted scalarisation below. Lemma~\ref{lemma:t_ub} gives constants
$C,c>0$ and $N\in\mathbb N$ such that
\eqref{eq:scalarisation_event_t_prob} holds.

By Lemma~\ref{lemma:moreau_center_probability}, there exist
$C_M<\infty$, constants $C,c>0$, and
$N\in\mathbb N$ such that, for all $n>N$,
$$
	\Pr\left(
		\lnot\mathcal E_n^{\mathrm{ctr}}(C_M)
		\mid
		\mathcal G_n
	\right)
	\leq
	C\exp\{-cn\}\,;
$$
since $\mathcal E_n^{\mathrm{ctr}}(C_M)$ increases with $C_M$, we may and do take $C_M\geq\sqrt{C_L}$. Increase $C_\eta$, if
necessary, once and for all so that
\begin{equation}
\label{eq:scalarisation_final_Ceta}
	C_\eta
	\geq
	C_M+
	\frac{C_t}{c_r}\,.
\end{equation}
In particular, $C_\eta\geq\sqrt{C_L}$.

On $\mathcal E_n^{\mathrm{sr}}(d_r,a,K,C_L)$,
Proposition~\ref{prop:small_r_exclusion}, applied with
$\mathcal R=\mathcal R_0$, gives the pointwise profile equality in
\eqref{eq:small_r_profile_exclusion} for the final cutoff $c_r$. Taking the
outer minimum in $\sigma$ yields $\mathcal E_n^r(c_r)$. Hence
$$
	\mathcal E_n^{\mathrm{sr}}(d_r,a,K,C_L)
	\subseteq
	\mathcal E_n^r(c_r)\,,
$$
and \eqref{eq:scalarisation_event_r_prob} follows from
\eqref{eq:scalarisation_event_sr_prob}.

The bound in Lemma~\ref{lemma:t_ub} is independent of $C_\eta$: its proof
uses the feasible point $\bb\eta=\bb0$ and a lower bound over all
$\bb\eta\in\Re^n$. Thus
\eqref{eq:scalarisation_event_t_prob} remains valid for the final value of
$C_\eta$. By \eqref{eq:scalarisation_final_Ceta} and
Lemma~\ref{lemma:moreau_feasibility},
$$
	\mathcal E_n^{\mathrm{ctr}}(C_M)
	\subseteq
	\mathcal E_n^{\mathrm M}(c_r,C_t,C_\eta)\,,
$$
which proves \eqref{eq:scalarisation_event_M_prob}.

Taking the largest of the thresholds $N$, the smallest of the exponents $c$ and the largest of the prefactors $C$ obtained above, and increasing $N$ if necessary so that $p-s\geq1$ for every $n>N$, the four preceding estimates hold with common
constants. Since
$$
	\lnot\mathcal R_n
	\subseteq
	\lnot\mathcal E_n^r(c_r)
	\cup
	\lnot\mathcal E_n^t(c_r,C_t)
	\cup
	\lnot\mathcal E_n^{\mathrm M}(c_r,C_t,C_\eta)\,,
$$
a union bound proves \eqref{eq:Rn_cond_prob}.

Finally, $\mathcal S_n(t_0)\in\mathcal G_n$, so the tower property gives
$$
\begin{aligned}
	\Pr(\lnot\mathcal R_n)
	&\leq
	\expect\left[
		\mathds{1}_{\mathcal S_n(t_0)}
		\Pr(
			\lnot\mathcal R_n
			\mid
			\mathcal G_n
		)
	\right]
	+
	\Pr(\lnot\mathcal S_n(t_0))
	\\
	&\leq
	C\exp\{-cn\}
	+
	\Pr(\lnot\mathcal S_n(t_0))\,.
\end{aligned}
$$
Assumptions~\ref{ass:betas_conv} and~\ref{ass:theta-conv} bound the final
probability by $C\exp\{-cn\}$ after adjusting constants. This proves
\eqref{eq:Rn_uncond_prob}.
\end{proof}

\subsection{Limiting AO}
\label{sec:ao-limit}

The scalarisation in Section~\ref{sec:ao-scalar} reduces the AO to a
min--max problem over finite-dimensional variables. Although $\bv$ is
two-dimensional, its feasible set $V_p$ need not be uniformly scaled in
$n$: when a positive eigenvalue of $\bb\Gamma_p$ is small, the corresponding
coefficient in $\bv$ may be large even though the induced signal-aligned
predictor remains bounded. We therefore reparameterise the AO at the intrinsic predictor
coordinate
$$
	\bu
    =
    \bb\Gamma_p^{1/2}\bv\,.
$$
With that, let
$$
    \bb Q_n
    =
    \frac1{\sqrt p}\bH_1\,.
$$
Conditional on $\mathcal G_n$, the rows $\bb q_i^\top$ of $\bb Q_n$ are
independent with
$$
    \bb q_i
    \sim
    \mathrm N(\bb0_2,\bb\Gamma_p)\,,
$$
and $\bb Q_n$ is independent of
$(\bb\varepsilon,\bg,\bh)$ conditional on $\mathcal G_n$.
Hence, for every $i$,
$$
    \var\left(
        [\bb Q_n\bv]_i
        \mid
        \mathcal G_n
    \right)
    =
    \bv^\top\bb\Gamma_p\bv
    =
    \vnorm{\bu}_2^2\,.
$$
Thus $\vnorm{\bu}_2$ is the natural intrinsic scale of the signal-aligned
predictor, rather than the potentially ill-conditioned Euclidean scale of
$\bv$. Moreover,
$$
    \vnorm{\bu}_2^2
    =
    \bv^\top\bb\Gamma_p\bv
    =
    \frac1p\vnorm{\bB\bv}_2^2
    \leq
    C_\beta^2\,,
$$
and
$
    \bu\in\range(\bb\Gamma_p)
$.
Hence
$$
    \bu
    \in
    \{\bu\in\Re^2:\vnorm{\bu}_2\leq C_\beta\}
    \cap
    \range(\bb\Gamma_p)\,.
$$
Conversely, if
$$
    \bu
    \in
    \{\bu\in\Re^2:\vnorm{\bu}_2\leq C_\beta\}
    \cap
    \range(\bb\Gamma_p)\,,
$$
set
$$
    \bv
    =
    \bb\Gamma_p^{+1/2}\bu\,.
$$
Then
$$
    \bv
    \in
    \range(\bb\Gamma_p)
    =
    \range(\bB^\top),
    \quad
    \bb\Gamma_p^{1/2}\bv
    =
    \bu\,,
$$
and
$$
    \frac1p\vnorm{\bB\bv}_2^2
    =
    \bv^\top\bb\Gamma_p\bv
    =
    \vnorm{\bu}_2^2
    \leq
    C_\beta^2\,.
$$
Thus $\bv\in V_p$. Consequently, the map
$
    \bv\mapsto\bb\Gamma_p^{1/2}\bv
$
sends $V_p$ onto
\begin{equation}
\label{eq:finite_signal_domain}
    U_n
    =
    U
    \cap
    \range(\bb\Gamma_p),
    \quad
    U
    =
    \{\bu\in\Re^2:\vnorm{\bu}_2\leq C_\beta\}\,.
\end{equation}
The change of variables therefore replaces the potentially ill-conditioned
set $V_p$ by a moving section $U_n$ of the fixed ambient ball $U$.
Now, let $\bb P_{\bb \Gamma_p}$ denote the orthogonal projector onto
$\range(\bb\Gamma_p)$ and define the direct pseudoinverse coordinates
$$
    \bb Z_n^0
    =
    \bb Q_n\bb\Gamma_p^{+1/2}\,.
$$
Conditional on $\mathcal G_n$, the rows of $\bb Z_n^0$ are centred Gaussian
with covariance
$$
    \bb\Gamma_p^{+1/2}
    \bb\Gamma_p
    \bb\Gamma_p^{+1/2}
    =
    \bb P_{\bb \Gamma_p}\,.
$$
Moreover, for every $\bv\in V_p$, with
$\bu=\bb\Gamma_p^{1/2}\bv$,
$$
\begin{aligned}
    \bb Z_n^0\bu
    &=
    \bb Q_n
    \bb\Gamma_p^{+1/2}
    \bb\Gamma_p^{1/2}
    \bv
    =
    \bb Q_n
    \bb P_{\bb \Gamma_p}
    \bv
    =
    \bb Q_n\bv\,.
\end{aligned}
$$
Thus $\bb Z_n^0$ already represents the finite AO on its actual signal
domain $U_n$. Its row covariance is, however, $\bb P_{\bb \Gamma_p}$, which is discontinuous
in $\bb\Gamma_p$ across rank strata: $\bb\Gamma_p\to\bb\Gamma$ does not
imply $\bb P_{\bb \Gamma_p}\to\bb P_{\bb \Gamma}$. An ambient extension built on
$\bb Z_n^0$ would inherit that discontinuity. Since the finite and limiting
objectives will be compared on $U$, we complete the Gaussian coordinates in
the missing directions.
Conditional on $\mathcal G_n$, let
$\widetilde{\bb Z}\in\Re^{n\times2}$ have i.i.d.\ standard Gaussian
entries, independently of
$(\bb Q_n,\bb\varepsilon,\bg,\bh)$, and define
\begin{equation}
\label{eq:ao_H_completion}
    \bb Z
    =
    \bb Z_n^0
    +
    \widetilde{\bb Z}
    \left(
        \bb I_2-\bb P_{\bb \Gamma_p}
    \right)
    =
    \bb Q_n\bb\Gamma_p^{+1/2}
    +
    \widetilde{\bb Z}
    \left(
        \bb I_2-\bb P_{\bb \Gamma_p}
    \right)\,.
\end{equation}
Conditional on $\mathcal G_n$, the two summands are independent, and the
rows $\bz_i^\top$ of $\bb Z$ are independent centred Gaussian vectors with
covariance
$$
\begin{aligned}
    \cov(
        \bz_i
        \mid
        \mathcal G_n
    )
    =
    \bb\Gamma_p^{+1/2}
    \bb\Gamma_p
    \bb\Gamma_p^{+1/2}
    +
    \left(
        \bb I_2-\bb P_{\bb \Gamma_p}
    \right)^2
   =
    \bb P_{\bb \Gamma_p}
    +
    \bb I_2-\bb P_{\bb \Gamma_p}
    =
    \bb I_2\,.
\end{aligned}
$$
Thus $\bb Z$ has i.i.d.\ standard Gaussian entries conditional on
$\mathcal G_n$. Since the rows of $\bb Q_n$ belong to
$\range(\bb\Gamma_p)$ almost surely,
$$
    \bb Q_n\bb P_{\bb \Gamma_p}
    =
    \bb Q_n\,,
$$
whereas
$$
    \left(
        \bb I_2-\bb P_{\bb \Gamma_p}
    \right)
    \bb\Gamma_p^{1/2}
    =
    \bb0_{2\times2}\,.
$$
It follows that
\begin{equation}
\label{eq:ao_H_repar}
    \frac1{\sqrt p}\bH_1
    =
    \bb Q_n
    =
    \bb Z\bb\Gamma_p^{1/2}\,,
\end{equation}
almost surely. Moreover, $\bb Z$ is independent of
$(\bb\varepsilon,\bg,\bh)$ conditional on $\mathcal G_n$.
Writing
$$
    \bb q_i
    =
    [q_{i,1},q_{i,2}]^\top
    =
    \bb\Gamma_p^{1/2}\bz_i\,,
$$
the pseudo-responses \eqref{eq:mdypl} satisfy
$$
    \Ytil_i
    =
    \alpha\mathds{1}
    \left\{
        \varepsilon_i
        <
        \rho'(\theta_0+q_{i,1})
    \right\}
    +
    (1-\alpha)
    \rho'(\theta_P+q_{i,2})\,.
$$
For $\bu\in U_n\subseteq\range(\bb\Gamma_p)$,
$$
    \left(
        \bb I_2-\bb P_{\bb \Gamma_p}
    \right)
    \bu
    =
    \bb0_2\,,
$$
and hence
$$
    \bb Z\bu
    =
    \bb Z_n^0\bu\,.
$$
Thus the completion leaves
$(\bH_1,\bYtil,\bg,\bh)$ and the finite AO unchanged on its actual
signal domain. The added coordinates define the signal-coordinate part of
the objective on the full ambient ball $U$ in a base law that is standard
on every rank stratum, so that, by \eqref{eq:ao_H_repar}, its dependence on
$\bb\Gamma_p$ enters through $\bb\Gamma_p^{1/2}$.

Combining \eqref{eq:finite_signal_domain} and
\eqref{eq:ao_H_repar}, for every $\bv\in V_p$,
$$
    \bb Q_n\bv
    =
    \bb Z\bb\Gamma_p^{1/2}\bv
    =
    \bb Z\bu,
    \quad
    \bu
    =
    \bb\Gamma_p^{1/2}\bv
    \in
    U_n\,.
$$
Accordingly, \eqref{eq:AO_scalar} becomes
\begin{equation}
\label{eq:AO_scalar_u}
\begin{aligned}
    \phi(\bb g,\bb h)
    =
    \underset{\sigma\in[0,C_\sigma]}{\min}
    \underset{r\in[c_r,1]}{\max}
    \underset{
        \substack{
            \bu\in U_n\\
            \abs{\theta}\leq C_\theta\\
            t\in[0,C_t]
        }
    }{\min}
    \,
    &\frac1n
    M_\rho
    \left(
        \theta\bb1
        +
        \bb Z\bu
        +
        \sigma\bg
        +
        \frac tr\bYtil,
        \frac tr
    \right)
    -
    \sigma r
    \frac{\vnorm{\bh}_2}{\sqrt n}
    +
    \frac{rt}{2}
    \\
    &
    -
    \frac1n
    \left\{
        \bYtil^\top
        \left(
            \theta\bb1
            +
            \bb Z\bu
            +
            \sigma\bg
        \right)
        +
        \frac{t}{2r}
        \vnorm{\bYtil}_2^2
    \right\}\,.
\end{aligned}
\end{equation}
Now, fix once and for all a baseline structural tolerance
$t_\star\in(0,1]$, and let
$
	d_r$, 
    $
	a$, 
    $
	K$, 
    $
	C_L$, 
    $
	c_r$, 
    $
	C_t$, 
    $
	C_M$, 
    $
	C_\eta
$
and $\mathcal R_n$ be the resulting constants and event from
Proposition~\ref{prop:scalarisation_events} with structural tolerance
$t_\star$. Throughout this subsection and its auxiliary proofs, every
subsequent structural tolerance denoted by $t_0$ is chosen in
$(0,t_\star]$. Since for $t_0\leq t_\star$,
$
	\mathcal S_n(t_0)
	\subseteq
	\mathcal S_n(t_\star)
$, 
the bounds of Proposition~\ref{prop:scalarisation_events} remain valid
on each such event $\mathcal S_n(t_0)$.

Recall from \eqref{eq:AO_scalar_u} that, on the scalarisation-valid
event $\mathcal R_n$, the AO can be written in terms of the variables $\bb\omega=(\sigma,r,\bu,\theta,t)$.
Since value convergence will be obtained from uniform convergence of the
objective, which compares optimisation values only over a common fixed
domain, and since the finite and limiting signal domains $U_n$ and $U^*$
need not converge to one another, we embed
both problems into a fixed ambient domain.
We use the following domains:
\begin{equation}
	\label{eq:Omegas}
	\begin{aligned}
		\bb\Omega
		&=
		[0,C_\sigma]\times[c_r,1]\times U
		\times[-C_\theta,C_\theta]\times[0,C_t],\\
		\bb\Omega_n
		&=
		[0,C_\sigma]\times[c_r,1]\times
		\left(U\cap\range(\bb\Gamma_p)\right)
		\times[-C_\theta,C_\theta]\times[0,C_t],\\
		\bb\Omega^*
		&=
		[0,C_\sigma]\times[c_r,1]\times
		\left(U\cap\range(\bb\Gamma)\right)
		\times[-C_\theta,C_\theta]\times[0,C_t]\,.
	\end{aligned}
\end{equation}
The domain $\bb\Omega_n$ is the actual finite-$n$ domain of the scalar AO after the change of variables $\bu=\bb\Gamma_p^{1/2}\bv$.
The ambient domain $\bb\Omega$ is fixed and will be used for uniform convergence.
The domain $\bb\Omega^*$ is the limiting domain.

For $\bb\omega=(\sigma,r,\bu,\theta,t)\in\bb\Omega$, write
$$
	\lambda=\frac{t}{r},
	\quad
	\xi_i=\theta+\bz_i^\top\bu+\sigma g_i\,,
$$
where $\bz_i^\top$ is the $i$th row of the matrix $\bb Z$ defined in
\eqref{eq:ao_H_completion}.
In the finite-$n$ model, set
$$
	\bb q_i=[q_{i,1}, q_{i,2}]^\top = \bb\Gamma_p^{1/2}\bz_i\,,
$$
and
$$
	\Ytil_i
	=
	\alpha\mathds{1}
	\left\{
		\varepsilon_i<\rho'(\theta_0+q_{i,1})
	\right\}
	+
	(1-\alpha)\rho'(\theta_P+q_{i,2})\,.
$$
Define the objective function
\begin{equation}
	\label{eq:F_n_def}
	\begin{aligned}
		F_n(\bb\omega)
		&=
		\frac1n
		\sum_{i=1}^n
		\left\{
			m_\rho\left(\xi_i+\lambda \Ytil_i,\lambda\right)
			-
			\Ytil_i\xi_i
			-
			\frac{\lambda}{2}\Ytil_i^2
		\right\}
		-
		\sigma r\frac{\vnorm{\bh}_2}{\sqrt n}
		+
		\frac{rt}{2}\,.
	\end{aligned}
\end{equation}
The continuous extension at $t=0$ is understood through $m_\rho(b,0)=\rho(b)\,.$
The conditional mean of the objective in \eqref{eq:F_n_def} given $\mathcal{G}_n$ is
\begin{equation}
	\label{eq:F_bar_n_def}
	\bar F_n(\bb\omega)
	=
	\expect\left[
		F_n(\bb\omega)
		\mid
		\mathcal G_n
	\right]\,.
\end{equation}
The limiting objective is defined in the same standard-Gaussian base
coordinate $\bz$ as follows. Let
$$
	\bb Q=[Q_1,Q_2]^\top = \bb\Gamma^{1/2}\bz,
	\quad
	\bz\sim\mathrm N(\bb0_2,\bb I_2),
	\quad
	G\sim\mathrm N(0,1),
	\quad
	\varepsilon\sim\mathrm{Unif}([0,1]) \,,
$$
with all variables mutually independent, and let
$$
	\Ytil
	=
	\alpha\mathds{1}
	\left\{
		\varepsilon<\rho'(\theta_0^*+Q_1)
	\right\}
	+
	(1-\alpha)\rho'(\theta_P^*+Q_2)\,.
$$
For $\xi=\theta+\bz^\top\bu+\sigma G\,,$
define
\begin{equation}
	\label{eq:F_limit_def}
	\begin{aligned}
		F(\bb\omega)
		&=
		\expect\left[
			m_\rho\left(\xi+\lambda \Ytil,\lambda\right)
			-
			\Ytil\xi
			-
			\frac{\lambda}{2}\Ytil^2
		\right]
		-
		\sigma r\sqrt\kappa
		+
		\frac{rt}{2}\,,
	\end{aligned}
\end{equation}
where all expectations in \eqref{eq:F_limit_def} are taken with respect to $(\bz,G,\varepsilon)$.

The next step is to compare the finite-$n$ objective $F_n$ with the
limiting objective $F$. We will establish uniform convergence on the fixed
ambient domain $\bb\Omega$ and separately control the replacement of the
finite signal domain in $\bb\Omega_n$ by the limiting signal domain in
$\bb\Omega^*$.

To state the resulting value comparisons compactly, let
$$
    \bb O
    =
    \Sigma\times R\times \bb K
    \subseteq
    \bb\Omega\,,
$$
where $\Sigma$ is a domain for $\sigma$, $R$ is a domain for $r$, and
$\bb K$ is a joint domain for $(\bu,\theta,t)$. For
$f:\bb O\to\Re$, define
\begin{equation}
\label{eq:val_operator}
    \operatorname{val}(f;\bb O)
    =
    \underset{\sigma\in\Sigma}{\min}\,
    \underset{r\in R}{\max}\,
    \underset{(\bu,\theta,t)\in \bb K}{\min}\,
    f(\sigma,r,\bu,\theta,t)\,.
\end{equation}
All value operators below use this coordinate order.
With this notation, on the scalarisation-valid event $\mathcal R_n$,
\begin{equation}
	\label{eq:AO_as_val_Fn}
	\phi(\bg,\bh)
	=
	\operatorname{val}(F_n;\bb\Omega_n)\,.
\end{equation}
The limiting AO value is
\begin{equation}
	\label{eq:bar_phi_val}
	\bar\phi
	=
	\operatorname{val}(F;\bb\Omega^*)\,.
\end{equation}
Now, consider the uniform convergence event for $F_n(\bb \omega)$ to $F(\bb \omega)$ on $\bb \Omega$. For $s>0$, define
\begin{equation}
	\label{event:U_n}
	\mathcal U_n(s)
	=
	\left\{
		\underset{\bb\omega\in\bb\Omega}{\sup}
		\abs{
			F_n(\bb\omega)-F(\bb\omega)
		}
		\leq s
	\right\}\,.
\end{equation}
The value operator in \eqref{eq:val_operator} is monotone and translation
equivariant. Hence, for any bounded $f,g:\bb O\to\Re$ for which the nested
extrema are attained,
\begin{equation}
\label{eq:value_lipschitz}
    \abs{
        \operatorname{val}(f;\bb O)
        -
        \operatorname{val}(g;\bb O)
    }
    \leq
    \underset{\bb\omega\in\bb O}{\sup}
    \abs{
        f(\bb\omega)-g(\bb\omega)
    }\,.
\end{equation}
Indeed, with
$\Delta=\sup_{\bb\omega\in\bb O}\abs{f(\bb\omega)-g(\bb\omega)}$,
the bounds $g-\Delta\leq f\leq g+\Delta$, together with monotonicity and
translation equivariance, give \eqref{eq:value_lipschitz}.
Thus, on $\mathcal U_n(s)$,
$\abs{\operatorname{val}(F_n;\bb O)-\operatorname{val}(F;\bb O)}
\leq s$
for every scalar domain $\bb O\subseteq\bb\Omega$ equipped with the inherited
min--max--min order.
It remains to relate the actual finite-$n$ domain $\bb\Omega_n$ and the limiting domain $\bb\Omega^*$ to the ambient domain $\bb\Omega$.
For this, define the domain/nullspace event
\begin{equation}
	\label{event:N_n}
	\mathcal N_n(s)
	=
	\left\{
		\abs{
			\operatorname{val}(F_n;\bb\Omega_n)
			-
			\operatorname{val}(F_n;\bb\Omega)
		}
		\leq s
	\right\}
	\cap
	\left\{
		\abs{
			\operatorname{val}(F;\bb\Omega)
			-
			\operatorname{val}(F;\bb\Omega^*)
		}
		\leq s
	\right\}\,.
\end{equation}
The event $\mathcal N_n(s)$ is controlled jointly with
$\mathcal U_n(s)$ in Proposition~\ref{prop:uniform_AO_control} using the
value-Lipschitz property, the exact finite and limiting nullspace identities,
and the uniform-convergence and structural-bias bounds below.
Combining \eqref{eq:AO_as_val_Fn}, \eqref{eq:bar_phi_val}, \eqref{eq:value_lipschitz}, and \eqref{event:N_n}, we obtain the following implication:
\begin{equation}
	\label{eq:AO_value_implication}
	\mathcal R_n
	\cap
	\mathcal U_n(s)
	\cap
	\mathcal N_n(s)
	\Longrightarrow
	\abs{
		\phi(\bg,\bh)-\bar\phi
	}
	\leq 3s\,.
\end{equation}
Indeed, on this event,
$$
\begin{aligned}
	\abs{\phi(\bg,\bh)-\bar\phi}
	&=
	\abs{
		\operatorname{val}(F_n;\bb\Omega_n)
		-
		\operatorname{val}(F;\bb\Omega^*)
	}\\
	&\leq
	\abs{
		\operatorname{val}(F_n;\bb\Omega_n)
		-
		\operatorname{val}(F_n;\bb\Omega)
	}\\
	&+
	\abs{
		\operatorname{val}(F_n;\bb\Omega)
		-
		\operatorname{val}(F;\bb\Omega)
	}\\
	&+
	\abs{
		\operatorname{val}(F;\bb\Omega)
		-
		\operatorname{val}(F;\bb\Omega^*)
	}\\
	&\leq 3s\,.
\end{aligned}
$$
The main probabilistic statement of this section is thus the following.

\begin{proposition}[Uniform convergence and moving-domain control of the scalar AO objective]
\label{prop:uniform_AO_control}
For every $s>0$, there exist
$t_0=t_0(s)\in(0,t_\star]$, constants $C_s,c_s>0$, and
$N=N(s)\in\mathbb N$ such that, for all $n>N$, the following hold.
\begin{enumerate}[label=(\roman*)]
    \item
	\begin{equation}
	\label{eq:uniform_AO_convergence_conditional}
		\mathds{1}_{\mathcal S_n(t_0)}
		\Pr\left(
			\lnot\mathcal U_n(s)
			\mid
			\mathcal G_n
		\right)
		\leq
		C_s\exp\{-c_sn\}\,.
	\end{equation}

	\item
	\begin{equation}
	\label{eq:uniform_AO_domain_inclusion}
		\mathcal S_n(t_0)
		\cap
		\mathcal U_n\left(\frac{s}{4}\right)
		\subseteq
		\mathcal U_n(s)
		\cap
		\mathcal N_n(s)\,.
	\end{equation}
\end{enumerate}
Consequently,
\begin{equation}
\label{eq:uniform_AO_joint_conditional}
	\mathds{1}_{\mathcal S_n(t_0)}
	\Pr\left(
		\left\{\lnot\mathcal U_n(s)\right\}
		\cup
		\left\{\lnot\mathcal N_n(s)\right\}
		\mid
		\mathcal G_n
	\right)
	\leq
	C_s\exp\{-c_sn\}\,.
\end{equation}

The constants depend only on $s$, $\kappa,\alpha,\bb\Gamma,\theta_0^*,\theta_P^*$
and on the constants in Section~\ref{sec:setup}.
\end{proposition}

\begin{proof}
	Fix $s>0$ and set
	$$
		q=\frac{s}{4}\,.
	$$
	\begin{enumerate}[label=(\roman*)]
	\item
	By Lemma~\ref{lemma:F_n_structural_bias} applied with input tolerance
	$q$, there exist $t_B=t_B(q)>0$ and $N_0\in\mathbb N$ such that, for
	all $n>N_0$,
	$$
		\mathds{1}_{\mathcal S_n(t_B)}
		\underset{\bb\omega\in\bb\Omega}{\sup}
		\abs{
			\bar F_n(\bb\omega)-F(\bb\omega)
		}
		\leq
		\frac{q}{3}\,.
	$$
	Set $t_0 = \min\{t_\star,t_B\}\,.$
	Since $\mathcal S_n(t_0) \subseteq \mathcal S_n(t_B)\,,$
	the same structural-bias bound holds on $\mathcal S_n(t_0)$.
	\begin{equation}
		\label{eq:uniform_AO_bias_bound}
		\mathds{1}_{\mathcal S_n(t_0)}
		\underset{\bb\omega\in\bb\Omega}{\sup}
		\abs{
			\bar F_n(\bb\omega)-F(\bb\omega)
		}
		\leq
		\frac{q}{3}\,.
	\end{equation}

	Let $L_n$ be the maximum of the five coordinatewise Lipschitz
	coefficients in Lemma~\ref{lemma:ao_lipschitz_constants}. By that
	lemma, there exist deterministic constants $\bar C<\infty$ and
	$C,c>0$ such that, for all sufficiently large $n$,
	\begin{equation}
		\label{eq:uniform_AO_Lipschitz_event}
		\Pr\left(
			L_n>\bar C
			\mid
			\mathcal G_n
		\right)
		\leq
		C\exp\{-cn\}\,.
	\end{equation}
	On the event $\{L_n\leq\bar C\}$,
	Lemma~\ref{lemma:lipschitz_ao} and the Cauchy--Schwarz inequality \citep[Theorem~1.37(d)]{rudin+etal:1976} give
	$$
		\abs{ F_n(\bb\omega)-F_n(\bb\omega') } \leq \bar C \left( \abs{\sigma-\sigma'} + \abs{r-r'} + \vnorm{\bu-\bu'}_2 + \abs{\theta-\theta'} + \abs{t-t'} \right) \leq \sqrt{5}\,\bar C \vnorm{ \bb\omega-\bb\omega' }_2\,,
	$$
	for every $\bb\omega,\bb\omega'\in\bb\Omega$.

	Let $\bar C_F$ be the deterministic Lipschitz constant from
	Lemma~\ref{lemma:limit_ao_lipschitz}, and define
	$$
		\bar C_0
		=
		\sqrt{5}\,\bar C+\bar C_F,
		\quad
		\varepsilon
		=
		\min\left\{
			1,\frac{q}{3\bar C_0}
		\right\}\,.
	$$
	Identify $\bb\Omega$ with its natural embedding in $\Re^6$.
	Since $\bb\Omega$ is a bounded subset of $\Re^6$, the standard
	volumetric covering bound
	(see, for example, \citealt[][Section~4.2]{vershynin:2018})
	provides a deterministic $\varepsilon$-net
	$\bb\Omega_\varepsilon\subseteq\bb\Omega$ satisfying
	\begin{equation}
		\label{eq:uniform_AO_net_cardinality}
		\abs{\bb\Omega_\varepsilon}
		\leq
		\left(
			1+\frac{K_\Omega}{\varepsilon}
		\right)^6\,,
	\end{equation}
	for a constant $K_\Omega<\infty$ depending only on $C_\sigma,C_\beta,C_\theta,C_t,c_r\,.$

	Suppose that $\mathcal S_n(t_0)$ and $\{L_n\leq\bar C\}$ hold.
	For any $\bb\omega\in\bb\Omega$, choose
	$\bb\omega_j\in\bb\Omega_\varepsilon$ such that
	$$
		\vnorm{
			\bb\omega-\bb\omega_j
		}_2
		\leq
		\varepsilon\,.
	$$
	Using \eqref{eq:uniform_AO_bias_bound} and the Lipschitz bounds for
	$F_n$ and $F$, we obtain
	$$
		\begin{aligned}
		\abs{
			F_n(\bb\omega)-F(\bb\omega)
		}
		&\leq
		\abs{
			F_n(\bb\omega_j)-\bar F_n(\bb\omega_j)
		}
		+
		\abs{
			\bar F_n(\bb\omega_j)-F(\bb\omega_j)
		}\\
		&+
		\abs{
			F_n(\bb\omega)-F_n(\bb\omega_j)
		}
		+
		\abs{
			F(\bb\omega)-F(\bb\omega_j)
		}\\
		&\leq
		\underset{\bb\omega_j\in\bb\Omega_\varepsilon}{\max}
		\abs{
			F_n(\bb\omega_j)-\bar F_n(\bb\omega_j)
		}
		+
		\frac{q}{3}
		+
		\bar C_0\varepsilon\\
		&\leq
		\underset{\bb\omega_j\in\bb\Omega_\varepsilon}{\max}
		\abs{
			F_n(\bb\omega_j)-\bar F_n(\bb\omega_j)
		}
		+
		\frac{2q}{3}\,.
		\end{aligned}
	$$
	Consequently,
	\begin{equation}
		\label{eq:uniform_AO_net_reduction}
		\mathcal S_n(t_0)
		\cap
		\{L_n\leq\bar C\}
		\cap
		\left\{
			\lnot\mathcal U_n(q)
		\right\}
		\subseteq
		\left\{
			\underset{\bb\omega_j\in\bb\Omega_\varepsilon}{\max}
			\abs{
			F_n(\bb\omega_j)-\bar F_n(\bb\omega_j)
			}
			>
			\frac{q}{3}
		\right\}\,.
	\end{equation}

	Since $\mathcal S_n(t_0)\in\mathcal G_n$, conditioning
	\eqref{eq:uniform_AO_net_reduction} on $\mathcal G_n$ and using
	\eqref{eq:uniform_AO_Lipschitz_event} gives
	$$
		\begin{aligned}
		&\mathds{1}_{\mathcal S_n(t_0)}
		\Pr\left(
			\lnot\mathcal U_n(q)
			\mid
			\mathcal G_n
		\right)\leq
		\mathds{1}_{\mathcal S_n(t_0)}
		\Pr\left(
			\underset{\bb\omega_j\in\bb\Omega_\varepsilon}{\max}
			\abs{
			F_n(\bb\omega_j)-\bar F_n(\bb\omega_j)
			}
			>
			\frac{q}{3}
			\mid
			\mathcal G_n
		\right)
		+
		C\exp\{-cn\} \,.
		\end{aligned}
	$$
	The net is deterministic, and the constants in
	Lemma~\ref{lemma:F_n_pointwise_conc} are uniform over
	$\bb\omega\in\bb\Omega$. Therefore, a union bound and
	\eqref{eq:uniform_AO_net_cardinality} yield
	$$
		\begin{aligned}
		&\mathds{1}_{\mathcal S_n(t_0)}
		\Pr\left(
			\underset{\bb\omega_j\in\bb\Omega_\varepsilon}{\max}
			\abs{
			F_n(\bb\omega_j)-\bar F_n(\bb\omega_j)
			}
			>
			\frac{q}{3}
			\mid
			\mathcal G_n
		\right)\leq
		C
		\left(
			1+\frac{K_\Omega}{\varepsilon}
		\right)^6
		\exp\left\{
			-\frac{cnq^2}{9}
		\right\} \,.
		\end{aligned}
	$$
	For fixed $s>0$, $\varepsilon$ and the covering-number factor are
	deterministic constants independent of $n$ and of the realised value
	of $\mathcal G_n$. After adjusting constants, we conclude that
	$$
		\mathds{1}_{\mathcal S_n(t_0)}
		\Pr\left(
			\lnot\mathcal U_n(q)
			\mid
			\mathcal G_n
		\right)
		\leq
		C_s\exp\{-c_sn\}\,.
	$$
	Since $\mathcal U_n(q)\subseteq\mathcal U_n(s)$, it follows that
	$$
		\mathds{1}_{\mathcal S_n(t_0)}
		\Pr\left(
			\lnot\mathcal U_n(s)
			\mid
			\mathcal G_n
		\right)
		\leq
		C_s\exp\{-c_sn\}\,.
	$$
	This proves part~(i), namely
	\eqref{eq:uniform_AO_convergence_conditional}.

	\item
	Let
	$$
		\Delta_n
		=
		\underset{\bb\omega\in\bb\Omega}{\sup}
		\abs{
			F_n(\bb\omega)-\bar F_n(\bb\omega)
		}\,,
	$$
	and
	$$
		\Delta_n'
		=
		\underset{\bb\omega\in\bb\Omega}{\sup}
		\abs{
			\bar F_n(\bb\omega)-F(\bb\omega)
		}\,.
	$$
	On
	$\mathcal S_n(t_0)\cap\mathcal U_n(q)$,
	the structural-bias bound obtained in part~(i) gives
	$
		\Delta_n'
		\leq
		q / 3 
	$.
	Therefore,
	$$
	\begin{aligned}
		\Delta_n
		\leq
		\underset{\bb\omega\in\bb\Omega}{\sup}
		\abs{
			F_n(\bb\omega)-F(\bb\omega)
		}
		+
		\Delta_n'\leq
		q+\frac{q}{3}
		=
		\frac{s}{3}\,.
	\end{aligned}
	$$
	By \eqref{eq:value_lipschitz} and
	Lemma~\ref{lemma:nullspace_no_help_barFn}, the triangle inequality gives
	$$
		\abs{
			\operatorname{val}(F_n;\bb\Omega_n)
		-
		\operatorname{val}(F_n;\bb\Omega)
		}
		\leq
		\Delta_n+0+\Delta_n
		=
		2\Delta_n
		\leq
		\frac{2s}{3}
		\leq
		s\,.
	$$
	Moreover, Lemma~\ref{lemma:nullspace_no_help_F} gives
	$$
		\operatorname{val}(F;\bb\Omega)
		=
		\operatorname{val}(F;\bb\Omega^*)\,.
	$$
	Hence $\mathcal N_n(s)$ holds. Since $q=s/4\leq s$,
	$\mathcal U_n(q)\subseteq\mathcal U_n(s)$, and therefore
	$$
		\mathcal S_n(t_0)
		\cap
		\mathcal U_n\left(\frac{s}{4}\right)
		\subseteq
		\mathcal U_n(s)
		\cap
		\mathcal N_n(s)\,.
	$$
	\end{enumerate}
	Taking complements in
	\eqref{eq:uniform_AO_domain_inclusion} gives
	$$
	\begin{aligned}
		\mathcal S_n(t_0)
		\cap
		\left(
			\left\{\lnot\mathcal U_n(s)\right\}
			\cup
			\left\{\lnot\mathcal N_n(s)\right\}
		\right)\subseteq
		\mathcal S_n(t_0)
		\cap
		\left\{
			\lnot\mathcal U_n(q)
		\right\}\,.
	\end{aligned}
	$$
	Since $\mathcal S_n(t_0)\in\mathcal G_n$, the stronger
	quarter-tolerance bound from part~(i) yields
	$$
	\begin{aligned}
		\mathds{1}_{\mathcal S_n(t_0)}
		\Pr\left(
			\left\{\lnot\mathcal U_n(s)\right\}
			\cup
			\left\{\lnot\mathcal N_n(s)\right\}
			\mid
			\mathcal G_n
		\right)\leq
		\mathds{1}_{\mathcal S_n(t_0)}
		\Pr\left(
			\lnot\mathcal U_n(q)
			\mid
			\mathcal G_n
		\right)
		\leq
		C_s\exp\{-c_sn\}\,.
	\end{aligned}
	$$
	This proves \eqref{eq:uniform_AO_joint_conditional}. 
	Taking $N(s)$ to be
	the maximum of the finite thresholds required above completes the
	proof.
\end{proof}

\subsection{Properties of the limiting AO}
\label{sec:limiting-AO-properties}

Before deriving the value gaps used in the comparison argument, we establish
that the limiting AO is well posed and that its compact, natural-domain, and
ambient formulations are equivalent. For the range-restricted and ambient
inner domains, write
\begin{equation}
\label{eq:bb-k}
    U^*
    =
    U\cap\range(\bb\Gamma),
    \quad
    \bb K^*
    =
    U^*\times[-C_\theta,C_\theta]\times[0,C_t],
\quad
    \bb K
    =
    U\times[-C_\theta,C_\theta]\times[0,C_t] \,.
\end{equation}

The compact and ambient inner problems are profiled as
\begin{equation}
\label{eq:compact_profiles}
\begin{aligned}
    H(\sigma,r)
    &=
    \underset{(\bu,\theta,t)\in\bb K^*}{\min}
    F(\sigma,r,\bu,\theta,t),
    &
    V(\sigma)
    &=
    \underset{r\in[c_r,1]}{\max}
    H(\sigma,r),
    \\
    H_0(\sigma,r)
    &=
    \underset{(\bu,\theta,t)\in\bb K}{\min}
    F(\sigma,r,\bu,\theta,t),
    &
    V_0(\sigma)
    &=
    \underset{r\in[c_r,1]}{\max}
    H_0(\sigma,r) \,.
\end{aligned}
\end{equation}
Thus, in the notation of \eqref{eq:val_operator},
\begin{equation}
\label{eq:ambient_compact_profiles}
    \operatorname{val}(F;\bb\Omega^*)
    =
    \underset{\sigma\in[0,C_\sigma]}{\min}V(\sigma),
    \quad
    \operatorname{val}(F;\bb\Omega)
    =
    \underset{\sigma\in[0,C_\sigma]}{\min}V_0(\sigma) \,.
\end{equation}

\begin{proposition}[Well-posedness and equivalent characterisations of the limiting AO]
\label{prop:limiting_AO_properties}
The compactification constants may be chosen consistently with
Propositions~\ref{prop:mdypl-po-link} and
\ref{prop:scalarisation_events} so that the following statements hold.
\begin{enumerate}[label=(\roman*)]
    \item
    The compact problems on $\bb\Omega^*$ and $\bb\Omega$, and the
    natural-domain problem in \eqref{eq:limit_ao_natural}, are all attained
    at the same unique saddle chain
    $$
        (\sigma^*,r^*,\bu^*,\theta^*,t^*) \,.
    $$
    Their common value is $\bar\phi$. In particular,
    $$
        \operatorname{val}(F;\bb\Omega)
        =
        \operatorname{val}(F;\bb\Omega^*)
        =
        \bar\phi \,.
    $$
    The common saddle satisfies
    $$
        0<\sigma^*<C_\sigma,
        \quad
        c_r<r^*<1,
        \quad
        0<t^*<C_t,
        \quad
        \vnorm{\bu^*}_2<C_\beta,
        \quad
        \abs{\theta^*}<C_\theta \,.
    $$

    \item
    Set
    $$
        \lambda^*
        =
        \frac{t^*}{r^*} \,.
    $$
    Then $(\sigma^*,\lambda^*,\bu^*,\theta^*)$ is the unique solution of
    \eqref{eq:FOCs} in
    $$
        (0,\infty)^2
        \times
        \range(\bb\Gamma)
        \times
        \Re \,.
    $$
    The remaining saddle coordinates are recovered as
    $$
        r^*
        =
        \frac{\sigma^*\sqrt\kappa}{\lambda^*},
        \quad
        t^*
        =
        \sigma^*\sqrt\kappa \,.
    $$
\end{enumerate}
\end{proposition}

\begin{proof}
The compatibility of the compactification constants is the choice made in
the proof of Proposition~\ref{prop:scalarisation_events}, together with the
preceding deterministic enlargements of $C_\sigma$, $C_\beta$, and
$C_\theta$.

Lemma~\ref{lemma:limiting_AO_natural_domain} gives attainment, uniqueness,
interiority, and equality of the compact problem on $\bb\Omega^*$ and the
natural-domain problem. Lemma~\ref{lemma:ambient_limiting_saddle_chain}
shows that enlarging the compact signal domain to obtain $\bb\Omega$ leaves
the value and the unique saddle chain unchanged. This proves part~(i).
Part~(ii) follows from Lemma~\ref{lemma:limiting_AO_focs}.
\end{proof}

\subsection{Restricted value comparison and PO localisation}
\label{sec:cgmt-cost-comparisons}
We now use the saddle and ambient-domain properties of
Proposition~\ref{prop:limiting_AO_properties} to localise the PO through
value comparisons similar to those of \citet{thrampoulidis+etal:2018} by forcing either
the radial variable or the signal and intercept variables away from the
limiting saddle. Deterministic value gaps for the corresponding limiting
AOs are transferred first to the finite AOs and then, through the
conditional CGMT, to the restricted POs. The resulting separation rules
out every PO optimiser outside the prescribed neighbourhood. Throughout
this subsection, a minimisation over an empty set has value $+\infty$.

For $\epsilon>0$, define the radial localising set
$$
    \mathcal S_{\bw}^{\epsilon}
    =
    \left\{
        \bw\in W_p:
        \abs{
            \frac{\vnorm{\bw}_2}{\sqrt p}
            -
            \sigma^*
        }
        <
        \epsilon
    \right\}\,,
$$
and let $\mathcal S_{\bw}^{\epsilon,c}$ be its complement in $W_p$.
Further, define the signal and intercept localising set
$$
    \mathcal S_{\bv,\theta}^{\epsilon}
    =
    \left\{
        (\bv,\theta)\in V_p\times[-C_\theta,C_\theta]:
        \vnorm{\bb\Gamma_p^{1/2}\bv-\bu^*}_2<\epsilon,\,
        \abs{\theta-\theta^*}<\epsilon
    \right\}\,,
$$
and let $\mathcal S_{\bv,\theta}^{\epsilon,c}$ be its complement in
$V_p\times[-C_\theta,C_\theta]$. Thus the target localising set for the PO
variables is
$$
    \mathcal S^{\epsilon}
    =
    \mathcal S_{\bw}^{\epsilon}
    \times
    \mathcal S_{\bv,\theta}^{\epsilon} \,.
$$
To prove localisation, we show that the PO value is strictly larger than
the unrestricted value $\Psi$ under the radial bad-set restriction and
under every signal and intercept restriction used to cover
$\mathcal S_{\bv,\theta}^{\epsilon,c}$. Both types of restriction can only
increase the PO value and thus the purpose of the comparison argument is to show
that this increase is bounded away from zero. The radial variable $\bw$ is
the outer minimisation variable in the PO of \eqref{eq:lagrangian5} and can
therefore be restricted directly.
By contrast, $(\bv,\theta)$ are
minimised inside the profile $\psi$, so restricting them changes the
profile and produces a separate PO--AO pair rather than merely restricting
the outer minimisation set as in
\citet[Theorem~3(ii)]{thrampoulidis+etal:2018}.
Restricting the radial variable gives the value
\begin{equation}
\label{eq:PO_restricted_w}
    \Psi_{\bw}^{\epsilon,c}
    =
    \underset{\bw\in\mathcal S_{\bw}^{\epsilon,c}}{\min}
    \underset{\blambda\in\mathcal B^n_1}{\max}
    \left\{
        -\frac1{n\sqrt p}\blambda^\top\bH_2\bw
        +
        \psi(\blambda)
    \right\} \,.
\end{equation}
To formulate the corresponding restrictions of the signal and intercept
variables, condition on $\mathcal C_n$ and let
$D\subseteq V_p\times[-C_\theta,C_\theta]$ be a nonempty compact convex
set. Recall $L(\bb\zeta;\blambda)$ from
Proposition~\ref{thm:pathwise_cgmt_localisation}, and define the
$D$-restricted profile and PO value by
\begin{equation}
\label{eq:restricted_profile}
     \Psi_D
    =
    \underset{\bw\in W_p}{\min}
    \underset{\blambda\in\mathcal B^n_1}{\max}
    \left\{
        -\frac1{n\sqrt p}\blambda^\top\bH_2\bw
        +
        \psi_D(\blambda)
    \right\}, \quad
     \psi_D(\blambda)
    =
    \underset{\bb\zeta\in Z_p(D)}{\min}
    L(\bb\zeta;\blambda), \quad
    Z_p(D)
    =
    D\times\mathcal B^n_{C_\eta}
    \,.
\end{equation}
Thus $\Psi_D$ is the value of the same conditional PO with only
$(\bv,\theta)$ restricted to $D$.
Moreover,
since $Z_p(D)$ is contained in the unrestricted domain,
$$
    \psi_D(\blambda)\geq\psi(\blambda)
    \quad\text{for every }\blambda\in\mathcal B^n_1,
    \quad
    \Psi_D\geq\Psi \,.
$$
Conditional on $\mathcal C_n$, compactness of $Z_p(D)$ and joint
continuity of $L$ give attainment and continuity of $\psi_D$ by Berge's maximum
theorem (see, for example,
\citealt[][Theorem~17.31.1]{aliprantis+et+al:2006}). Since
$L(\bb\zeta;\cdot)$ is affine for every fixed $\bb\zeta$, its pointwise
infimum $\psi_D$ is concave in $\blambda$. Since $Z_p(D)$ is compact and
convex, Lemma~\ref{lemma:psi_conv_conc} applies with $Z=Z_p(D)$, and
\eqref{eq:restricted_profile} is an admissible conditional CGMT instance
of Proposition~\ref{thm:pathwise_cgmt_localisation} with $D_0=D$. 

Since $\mathcal S_{\bv,\theta}^{\epsilon,c}$ need not be convex, it cannot
be used as a single restriction $D$ while retaining this convex PO--AO
structure. We therefore cover it by finitely many compact convex
pieces, apply \eqref{eq:restricted_profile} to each piece, and combine the
resulting comparisons over the fixed finite index set.

Fix a $1/2$-net $\mathcal A\subset\mathbb S^1$ and enumerate it as
$$
    \mathcal A
    =
    \left\{
        \bb a_i\in\mathbb S^1:
        i\in\mathcal I_{\mathcal A}
    \right\},
    \quad
    \mathcal I_{\mathcal A}
    =
    \{1,\ldots,\abs{\mathcal A}\}\,,
$$
such that
$$
    \forall\,\bb x\in\mathbb S^1
    \quad
    \exists\,i\in\mathcal I_{\mathcal A}:
    \vnorm{\bb x-\bb a_i}_2\leq\frac12 \,.
$$
For instance, one may take
$\abs{\mathcal A}\leq5^2$
\citep[Corollary~4.2.13]{vershynin:2018}. For such a cover, define the index set
$
    \mathcal J
    =
    \mathcal I_{\mathcal A}\cup\{+,-\}
$, and define the intercept intervals
$$
    \Theta_+^\epsilon
    =
    [-C_\theta,C_\theta]
    \cap
    [\theta^*+\epsilon,\infty),
    \quad
    \Theta_-^\epsilon
    =
    [-C_\theta,C_\theta]
    \cap
    (-\infty,\theta^*-\epsilon] \,.
$$
For $j\in\mathcal J$, define
\begin{equation}
\label{eq:finite_bad_pieces}
    C_{p,j}^{\epsilon}
    =
    \begin{cases}
        \left\{
            (\bv,\theta)
            \in
            V_p\times[-C_\theta,C_\theta]:
            \bb a_j^\top
            \left(
                \bb\Gamma_p^{1/2}\bv-\bu^*
            \right)
            \geq
            \dfrac{\epsilon}{2}
        \right\},
        & j\in\mathcal I_{\mathcal A},
        \\[3ex]
        V_p\times\Theta_j^\epsilon,
        & j\in\{+,-\}
    \end{cases}\,.
\end{equation}
Lemma~\ref{lemma:finite_convex_bad_cover} shows that
$$
    \mathcal S_{\bv,\theta}^{\epsilon,c}
    \subseteq
    \bigcup_{j\in\mathcal J}C_{p,j}^{\epsilon}\,,
$$
and that every nonempty piece is compact and convex.
Set
$$
    Z_{p,j}^{\epsilon}
    =
    Z_p(C_{p,j}^{\epsilon})
    =
    C_{p,j}^{\epsilon}\times\mathcal B^n_{C_\eta},
    \quad
    \mathcal J_n^\epsilon
    =
    \left\{
        j\in\mathcal J:
        C_{p,j}^{\epsilon}\neq\emptyset
    \right\}
    \,.
$$
For $j\in\mathcal J_n^\epsilon$, apply
\eqref{eq:restricted_profile} with $D=C_{p,j}^{\epsilon}$ and write
$$
    \psi_j^\epsilon(\blambda)
    =
    \psi_{C_{p,j}^{\epsilon}}(\blambda)
    =
    \underset{\bb\zeta\in Z_{p,j}^{\epsilon}}{\min}
    L(\bb\zeta;\blambda),
    \quad
    \Psi_j^\epsilon
    =
    \Psi_{C_{p,j}^{\epsilon}} \,,
$$
For $j\notin\mathcal J_n^\epsilon$, set
$\Psi_j^\epsilon=+\infty$. For every
$j\in\mathcal J_n^\epsilon$, the argument following
\eqref{eq:restricted_profile} gives a valid conditional PO--AO pair.
Consequently, it is enough to establish the required value separation
for the finitely many values
$\{\Psi_j^\epsilon:j\in\mathcal J_n^\epsilon\}$.

The scalarisation of Section~\ref{sec:ao-scalar} carries the restrictions
above to the scalar variables $\bb\omega=(\sigma,r,\bu,\theta,t)$ through
$$
\frac{\vnorm{\bw}_2}{\sqrt p}\longmapsto\sigma,
\quad
\bv\longmapsto\bu=\bb\Gamma_p^{1/2}\bv,
\quad
\theta\longmapsto\theta \,.
$$
Since $\bw\mapsto\vnorm{\bw}_2/\sqrt p$ maps $W_p$ onto $[0,C_\sigma]$, the
radial restriction becomes
$$
\mathcal S_{\sigma}^{\epsilon,c}
    =
    \left\{
        \sigma\in[0,C_\sigma]:
        \abs{\sigma-\sigma^*}\geq\epsilon
    \right\}\,,
$$
which is nonempty whenever $\mathcal S_{\bw}^{\epsilon,c}$ is. For
$U_0\subseteq U$ and $j\in\mathcal J$, define, in correspondence
with \eqref{eq:finite_bad_pieces},
\begin{equation}
\label{eq:scalar_bad_pieces}
    \bb K_j^\epsilon(U_0)
    =
    \begin{cases}
        \left\{
            (\bu,\theta,t)
            \in
            U_0\times[-C_\theta,C_\theta]\times[0,C_t]:
            \bb a_j^\top(\bu-\bu^*)
            \geq
            \dfrac{\epsilon}{2}
        \right\},
        & j\in\mathcal I_{\mathcal A},
        \\[3ex]
        U_0\times\Theta_j^\epsilon\times[0,C_t],
        & j\in\{+,-\}
    \end{cases}\,.
\end{equation}
By \eqref{eq:finite_signal_domain}, the change of variables
$$
    \bu
    =
    \bb\Gamma_p^{1/2}\bv\,,
$$
maps $V_p$ onto
$U_n=U\cap\range(\bb\Gamma_p)$. For
$j\in\mathcal I_{\mathcal A}$, it transforms the signal half-space in
\eqref{eq:finite_bad_pieces} into
$$
    \bb a_j^\top(\bu-\bu^*)
    \geq
    \frac{\epsilon}{2}\,,
$$
while for $j\in\{+,-\}$ it leaves the restriction
$\theta\in\Theta_j^\epsilon$ unchanged. Hence the scalar image of
$C_{p,j}^\epsilon\times[0,C_t]$ is exactly
$\bb K_j^\epsilon(U_n)$ for every $j\in\mathcal J$.
Recall $\bb K$ from \eqref{eq:bb-k} and set
$$
    \bb K_n
    =
    U_n\times[-C_\theta,C_\theta]\times[0,C_t] \,.
$$
For $j\in\mathcal J$, the finite and ambient restricted scalar domains are
\begin{equation}
\label{eq:restricted_scalar_domains}
\begin{aligned}
    \bb\Omega_{n,\sigma}^{\epsilon,c}
    =
    \mathcal S_{\sigma}^{\epsilon,c}
    \times[c_r,1]\times\bb K_n,
    \quad
    \bb\Omega_{\sigma}^{\epsilon,c}
    =
    \mathcal S_{\sigma}^{\epsilon,c}
    \times[c_r,1]\times\bb K,
    \quad 
    \bb\Omega_{n,j}^{\epsilon}
    =
    [0,C_\sigma]\times[c_r,1]\times\bb K_j^\epsilon(U_n),
    \quad 
    \bb\Omega_j^{\epsilon}
    =
    [0,C_\sigma]\times[c_r,1]\times\bb K_j^\epsilon(U) \,.
\end{aligned}
\end{equation}
With that, finally, define the restricted limiting AO values 
\begin{equation}
\label{eq:limiting_restricted_values}
\bar\phi_{\sigma}^{\epsilon,c}
=
\operatorname{val}
\left(
    F;
    \bb\Omega_{\sigma}^{\epsilon,c}
\right),
\quad
\bar\phi_j^\epsilon
=
\operatorname{val}
\left(
    F;
    \bb\Omega_j^\epsilon
\right),
\quad j\in\mathcal J \,.
\end{equation}
Lemma~\ref{lemma:limiting_AO_gap} proves strict gaps for all restricted
limiting AOs. For the radial restriction, the claim is immediate when
$\mathcal S_{\sigma}^{\epsilon,c}$ is empty. Otherwise, continuity of
$V_0$ and uniqueness of its minimiser $\sigma^*$ give
$$
    \bar\phi_{\sigma}^{\epsilon,c}
    =
    \underset{\sigma\in\mathcal S_{\sigma}^{\epsilon,c}}{\min}
    V_0(\sigma)
    >
    V_0(\sigma^*)
    =
    \bar\phi\,,
$$
since $\sigma^*\notin\mathcal S_{\sigma}^{\epsilon,c}$. 
For a signal or intercept piece, define the corresponding restricted inner
and outer profiles by
$$
    H_j^\epsilon(\sigma,r)
    =
    \underset{(\bu,\theta,t)\in\bb K_j^\epsilon(U)}{\min}
    F(\sigma,r,\bu,\theta,t),
    \quad
    V_j^\epsilon(\sigma)
    =
    \underset{r\in[c_r,1]}{\max}
    H_j^\epsilon(\sigma,r) \,.
$$
If $\bar\phi_j^\epsilon=\bar\phi$, then
$V_j^\epsilon\geq V_0$ and uniqueness of $\sigma^*$ first give
$V_j^\epsilon(\sigma^*)=\bar\phi$. Consequently,
$$
    \bar\phi
    =
    H_0(\sigma^*,r^*)
    \leq
    H_j^\epsilon(\sigma^*,r^*)
    \leq
    V_j^\epsilon(\sigma^*)
    =
    \bar\phi\,.
$$
Thus the restricted inner minimum at $(\sigma^*,r^*)$ equals
$\bar\phi$. Attainment and uniqueness of the unrestricted inner
minimiser would require $(\bu^*,\theta^*,t^*)$ to lie in a bad piece, which
is impossible by construction. 
Lemma~\ref{lemma:limiting_AO_gap} thus shows that there is a deterministic $\eta_\epsilon>0$, independent of
$n$, such that
\begin{equation}
\label{eq:limiting_barriers}
    \bar\phi_{\sigma}^{\epsilon,c}
    \geq
    \bar\phi+3\eta_\epsilon,
    \quad
    \bar\phi_j^\epsilon
    \geq
    \bar\phi+3\eta_\epsilon,
    \quad
    j\in\mathcal J \,.
\end{equation}
All that remains is to propagate the limiting AO value gaps through the CGMT pipeline back to the 
PO. The AO counterparts of \eqref{eq:PO_restricted_w} and $\Psi_j^\epsilon$ are
$$
    \phi_{\bw}^{\epsilon,c}(\bg,\bh)
    =
    \underset{\bw\in\mathcal S_{\bw}^{\epsilon,c}}{\min}
    \underset{\blambda\in\mathcal B^n_1}{\max}
    \left\{
        -\frac1{n\sqrt p}
        \left(
            \vnorm{\bw}_2\bg^\top\blambda
            +
            \vnorm{\blambda}_2\bw^\top\bh
        \right)
        +
        \psi(\blambda)
    \right\}\,,
$$
and, for $j\in\mathcal J_n^\epsilon$,
$$
    \phi_j^\epsilon(\bg,\bh)
    =
    \underset{\bw\in W_p}{\min}
    \underset{\blambda\in\mathcal B^n_1}{\max}
    \left\{
        -\frac1{n\sqrt p}
        \left(
            \vnorm{\bw}_2\bg^\top\blambda
            +
            \vnorm{\blambda}_2\bw^\top\bh
        \right)
        +
        \psi_j^\epsilon(\blambda)
    \right\} \,.
$$
If $\mathcal S_{\bw}^{\epsilon,c}=\emptyset$, set
$\phi_{\bw}^{\epsilon,c}=+\infty$, and set $\phi_j^\epsilon=+\infty$ for
$j\notin\mathcal J_n^\epsilon$. The parametrisations below are
understood on the events that the corresponding restricted domains are
nonempty. In the notation of
Proposition~\ref{thm:pathwise_cgmt_localisation}, the three PO--AO pairs
are parameterised by
\begin{equation}
\label{eq:restricted_po_ao_pairs}
\begin{aligned}
    (\Psi,\phi)
    &: &
    W_0&=W_p,
    &
    D_0&=V_p\times[-C_\theta,C_\theta],
    \\
    (\Psi_{\bw}^{\epsilon,c},\phi_{\bw}^{\epsilon,c})
    &: &
    W_0&=\mathcal S_{\bw}^{\epsilon,c},
    &
    D_0&=V_p\times[-C_\theta,C_\theta],
    \\
    (\Psi_j^\epsilon,\phi_j^\epsilon)
    &: &
    W_0&=W_p,
    &
    D_0&=C_{p,j}^\epsilon,
    \quad j\in\mathcal J_n^\epsilon \,.
\end{aligned}
\end{equation}
Lemma~\ref{lemma:restricted_piece_scalarisation} shows that, on one event of
conditional exponentially high probability, the unrestricted and
restricted AOs are represented by the same objective $F_n$, with only the
feasible domain changed:
$$
    \phi(\bg,\bh)
    =
    \operatorname{val}(F_n;\bb\Omega_n),
    \quad
    \phi_{\bw}^{\epsilon,c}(\bg,\bh)
    =
    \operatorname{val}(F_n;\bb\Omega_{n,\sigma}^{\epsilon,c}), \quad
    \phi_j^\epsilon(\bg,\bh)
    =
    \operatorname{val}(F_n;\bb\Omega_{n,j}^{\epsilon}),
    \quad j\in\mathcal J_n^\epsilon \,.
$$
Lemma~\ref{lemma:restricted_AO_value_bounds} combines this simultaneous
scalarisation with Proposition~\ref{prop:uniform_AO_control}. In each
restricted comparison, enlarging $U_n$ to $U$ enlarges only the innermost
minimisation domain and can therefore only decrease the min--max--min value.
Uniform control of $F_n-F$ on $\bb\Omega$ then gives the required lower
bounds relative to \eqref{eq:limiting_restricted_values}, simultaneously
over $j\in\mathcal J_n^\epsilon$.


Lemma~\ref{lemma:restricted_PO_value_barriers} applies
Proposition~\ref{thm:pathwise_cgmt_localisation} in the required tail
directions. The upper-tail comparison controls the unrestricted PO. The
lower-tail comparison controls the radial restriction and every active
signal or intercept piece. 
Since
$\mathcal J_n^\epsilon\subseteq\mathcal J$ and $\mathcal J$ is
deterministic and finite, a single uniform union bound covers all
active pieces. Define
\begin{equation}
\label{event:P_n_eps}
    \mathcal P_n^\epsilon
    =
    \left\{
        \Psi\leq\bar\phi+\eta_\epsilon
    \right\}
    \cap
    \left\{
        \Psi_{\bw}^{\epsilon,c}
        \geq
        \bar\phi_{\sigma}^{\epsilon,c}-\eta_\epsilon
    \right\}
    \cap
    \bigcap_{j\in\mathcal J_n^\epsilon}
    \left\{
        \Psi_j^\epsilon
        \geq
        \bar\phi_j^\epsilon-\eta_\epsilon
    \right\} \,.
\end{equation}
Lemma~\ref{lemma:restricted_PO_value_barriers} gives
$t_0=t_0(\epsilon)>0$, constants $C_\epsilon,c_\epsilon>0$, and
$N_\epsilon\in\mathbb N$ such that, for all $n>N_\epsilon$,
\begin{equation}
\label{eq:P_n_eps_prob}
\mathds{1}_{\mathcal S_n(t_0)}
\Pr\left(
    \lnot\mathcal P_n^\epsilon
    \mid
    \mathcal G_n
\right)
\leq
C_\epsilon\exp\{-c_\epsilon n\} \,.
\end{equation}
Together with \eqref{eq:limiting_barriers}, this yields the local
separation
\begin{equation}
\label{eq:PO_separation}
    \Psi
    \leq
    \bar\phi+\eta_\epsilon
    <
    \bar\phi+2\eta_\epsilon
    \leq
    \min
    \left\{
        \Psi_{\bw}^{\epsilon,c},
        \underset{j\in\mathcal J_n^\epsilon}{\min}\,
        \Psi_j^\epsilon
    \right\}
    \quad\text{on }\mathcal P_n^\epsilon \,.
\end{equation}

For all sufficiently large $n$, let
$(\hat\bv,\hat\bw,\hat\theta,\hat{\bb\eta})$
be the minimising variables fixed in
Proposition~\ref{prop:mdypl-po-link}, and set
$$
    \hat{\bb\zeta}
    =
    (\hat\bv,\hat\theta,\hat{\bb\eta}) \,.
$$
On $\mathcal E_n$, Lemma~\ref{lemma:compact_PO_uniqueness} shows that
these minimising variables are unique and that $\hat\bw$ is the unique
minimiser of the profiled PO in \eqref{eq:PO}.

\begin{proposition}[PO localisation]
\label{prop:PO_localisation}
For every $\epsilon>0$, there are $t_0=t_0(\epsilon)>0$, constants
$C_\epsilon,c_\epsilon>0$, and $N_\epsilon\in\mathbb N$ such that, for all
$n>N_\epsilon$,
$$
    \mathds{1}_{\mathcal S_n(t_0)}
    \Pr\left(
        (\hat\bw,\hat\bv,\hat\theta)
        \notin
        \mathcal S^\epsilon
        \mid
        \mathcal G_n
    \right)
    \leq
    C_\epsilon\exp\{-c_\epsilon n\} \,.
$$
Consequently, after adjusting the constants,
$$
    \Pr\left(
        (\hat\bw,\hat\bv,\hat\theta)
        \notin
        \mathcal S^\epsilon
    \right)
    \leq
    C_\epsilon\exp\{-c_\epsilon n\} \,.
$$
\end{proposition}

\begin{proof}
Let $\mathcal P_n^\epsilon$ be as in \eqref{event:P_n_eps}. By
\eqref{eq:P_n_eps_prob}, and since $\Pr(\lnot\mathcal E_n)=0$ for $n\geq p+1$, it is enough to prove
$$
    \mathcal P_n^\epsilon\cap\mathcal E_n
    \subseteq
    \left\{
        (\hat\bw,\hat\bv,\hat\theta)
        \in
        \mathcal S^\epsilon
    \right\} \,.
$$
Work on $\mathcal P_n^\epsilon\cap\mathcal E_n$, where the compact PO has a unique primal minimiser by Lemma~\ref{lemma:compact_PO_uniqueness}. If
$\hat\bw\in\mathcal S_{\bw}^{\epsilon,c}$, then, using $\hat\bw$ as a
feasible point in \eqref{eq:PO_restricted_w} and its optimality in the
unrestricted profiled PO,
$$
    \Psi_{\bw}^{\epsilon,c}
    \leq
    \underset{\blambda\in\mathcal B^n_1}{\max}
    \left\{
        -\frac1{n\sqrt p}\blambda^\top\bH_2\hat\bw
        +
        \psi(\blambda)
    \right\}
    =
    \Psi \,.
$$
This contradicts \eqref{eq:PO_separation}. Hence
$\hat\bw\in\mathcal S_{\bw}^{\epsilon}$.
If
$(\hat\bv,\hat\theta)\in\mathcal S_{\bv,\theta}^{\epsilon,c}$,
Lemma~\ref{lemma:finite_convex_bad_cover} gives
$j\in\mathcal J$ such that
$(\hat\bv,\hat\theta)\in C_{p,j}^{\epsilon}$. Then
$j\in\mathcal J_n^\epsilon$ and
$\hat{\bb\zeta}\in Z_{p,j}^{\epsilon}$. By
\eqref{eq:restricted_profile},
$$
    \psi_j^\epsilon(\blambda)
    \leq
    L(\hat{\bb\zeta};\blambda)
    \quad
    \text{for every }
    \blambda\in\mathcal B^n_1 \,.
$$
Using $\hat\bw$ as a feasible point in the restricted PO therefore gives
$$
    \Psi_j^\epsilon
    \leq
    \underset{\blambda\in\mathcal B^n_1}{\max}
    \left\{
        -\frac1{n\sqrt p}\blambda^\top\bH_2\hat\bw
        +
        \psi_j^\epsilon(\blambda)
    \right\}
    \leq
    \underset{\blambda\in\mathcal B^n_1}{\max}
    \left\{
        -\frac1{n\sqrt p}\blambda^\top\bH_2\hat\bw
        +
        L(\hat{\bb\zeta};\blambda)
    \right\}
    =
    \Psi \,.
$$
The final equality follows from the primal optimality of
$(\hat\bv,\hat\bw,\hat\theta,\hat{\bb\eta})$ in
\eqref{eq:lagrangian4} and the value identity in
\eqref{eq:lagrangian5}. This again contradicts
\eqref{eq:PO_separation}. Thus
$(\hat\bv,\hat\theta)\in\mathcal S_{\bv,\theta}^{\epsilon}$, and hence
$(\hat\bw,\hat\bv,\hat\theta)\in\mathcal S^\epsilon$ on
$\mathcal P_n^\epsilon\cap\mathcal E_n$.
Consequently, using $\Pr(\lnot\mathcal E_n)=0$,
$$
\begin{aligned}
    &\mathds{1}_{\mathcal S_n(t_0)}
    \Pr\left(
        (\hat\bw,\hat\bv,\hat\theta)
        \notin
        \mathcal S^\epsilon
        \mid
        \mathcal G_n
    \right)
    \leq
    \mathds{1}_{\mathcal S_n(t_0)}
    \Pr\left(
        \lnot\mathcal P_n^\epsilon
        \mid
        \mathcal G_n
    \right)
    \leq
    C_\epsilon\exp\{-c_\epsilon n\} \,.
\end{aligned}
$$
This proves the conditional bound.
Taking expectations and splitting over $\mathcal S_n(t_0)$ gives
$$
\begin{aligned}
    \Pr\left(
        (\hat\bw,\hat\bv,\hat\theta)
        \notin
        \mathcal S^\epsilon
    \right)
    \leq
    \expect\left[
        \mathds{1}_{\mathcal S_n(t_0)}
        \Pr\left(
            (\hat\bw,\hat\bv,\hat\theta)
            \notin
            \mathcal S^\epsilon
            \mid
            \mathcal G_n
        \right)
    \right]
 +
    \Pr\left(
        \lnot\mathcal S_n(t_0)
    \right) \,.
\end{aligned}
$$
The structural concentration bound for $\mathcal S_n(t_0)$ gives the
unconditional estimate after adjusting constants.
\end{proof}

\subsection{Convergence on test functions}
\label{sec:test-function-convergence}

At last, we translate the localisation of the PO in \eqref{eq:lagrangian5}
into the convergence on test function result
of Theorem~\ref{thm:mdypl-convergence}.
The main argument relies on the conditional spherical representation of the 
signal-orthogonal component of the PO in \eqref{eq:lagrangian5}, which
extends the distributional equivalence of
\citet[Proposition~2.1]{zhao+etal:2022} and which we subsequently replace by a 
Gaussian  
surrogate $\bb z_p$ and propagate through the pseudo-Lipschitz test functions. 

For all sufficiently large $n$, let
$(\hat\bv,\hat\bw,\hat\theta,\hat{\bb\eta})$
be the minimising variables fixed in
Proposition~\ref{prop:mdypl-po-link}. By
Lemma~\ref{lemma:compact_PO_uniqueness}, they are unique on
$\mathcal E_n$. Set
\begin{equation}
\label{eq:PO-slope}
    \hat\bbeta
    =
    \bb B\hat\bv+\bb E\hat\bw  \,.
\end{equation}
Further, let $\hat \sigma = \vnorm{\hat \bw}_2 / \sqrt p$ and introduce the notation
$ \bb{\mathfrak Z}
    =
    (
        \hat\bv,
        \hat\theta,
        \hat{\bb\eta},
        \hat\sigma
    )$.
Since $s=\rank(\bb B)\leq2$ and $p>2$, the space
$\range(\bb B)^\perp$ is nontrivial. Let $\bb\xi_n\sim\mathrm N(\bb0_p,\bb I_p)$ be independent of the original
variables, and define
\begin{equation}
\label{eq:compact_PO_direction}
    \bb d
    =
    \begin{cases}
        \displaystyle
        \frac{\bb E\hat\bw}{\vnorm{\hat\bw}_2},
        &
        \hat\sigma>0,
        \\
        \displaystyle
        \frac{\bb P^\perp\bb\xi_n}
        {\vnorm{\bb P^\perp\bb\xi_n}_2},
        &
        \hat\sigma=0,
        \quad
        \vnorm{\bb P^\perp\bb\xi_n}_2>0,
        \\
        \bb E\bb e_1,
        &
        \hat\sigma=0,
        \quad
        \vnorm{\bb P^\perp\bb\xi_n}_2=0\,,
    \end{cases}
\end{equation}
where $\bb e_1$ is the first standard basis vector of $\Re^{p-s}$. Then
$\bb d\in\range(\bb B)^\perp$, $\vnorm{\bb d}_2=1$, and
\begin{equation}
\label{eq:compact_PO_polar}
    \hat\bbeta
    =
    \bb B\hat\bv
    +
    \hat\sigma\sqrt p \bb d\,.
\end{equation}
Writing
$
    \mathbb S_{\bb B}^\perp
    =
    \left\{
        \bx\in\range(\bb B)^\perp:
        \vnorm{\bx}_2=1
    \right\}
$,
Lemma~\ref{lemma:compact_PO_spherical_representation} shows that,
conditionally on $\mathcal C_n$, the direction $\bb d$ is uniform on
$\mathbb S_{\bb B}^\perp$ and independent of $\bb{\mathfrak Z}$.
Now, define the Gaussian surrogate $\bb z_p$ for the noise
$\sqrt p \bb d$, namely
\begin{equation}
\label{eq:mdypl-gaussian-surrogate}
    \bb z_p
    =
    \vnorm{\bb P^\perp\bb\xi_n}_2\bb d
    +
    \bb P\bb\xi_n\,.
\end{equation}
On $\mathcal E_\sigma=\{\hat\sigma>0\}$,
Lemma~\ref{lemma:compact_PO_spherical_representation} and the independence
of $\bb\xi_n$ imply that, conditionally on
$\sigma(\mathcal C_n,\bb{\mathfrak Z})$, the vector $\bb d$ is uniform on
the unit sphere of $\range(\bb B)^\perp$ and independent of $\bb\xi_n$.
The Gaussian polar decomposition
\citep[Exercise~3.3.7]{vershynin:2018} gives
$\vnorm{\bb P^\perp\bb\xi_n}_2\bb d
\sim\mathrm N(\bb0_p,\bb P^\perp)$ independently of
$\bb P\bb\xi_n\sim\mathrm N(\bb0_p,\bb P)$ on $\mathcal E_\sigma$, while
\eqref{eq:compact_PO_direction} gives
$\bb z_p=\bb\xi_n$ on $\lnot\mathcal E_\sigma$.
Since $\mathcal E_\sigma\in\sigma(\bb{\mathfrak Z})$, both branches have
the same conditional law, so
\begin{equation}
\label{eq:gaussian-surrogate-law}
    \bb z_p
    \mid
    \mathcal C_n,
    \bb{\mathfrak Z}
    \sim
    \mathrm N(\bb0_p,\bb I_p)\,,
\end{equation}
and $\bb z_p$ is independent of
$\sigma(\mathcal C_n,\bb{\mathfrak Z})$. In particular, it is standard
Gaussian conditionally on $\mathcal G_n$.

Moreover, orthogonality in \eqref{eq:mdypl-gaussian-surrogate} gives
$$
    \frac1p\vnorm{\sqrt p\bb d-\bb z_p}_2^2
    =
    \left(
        1-\vnorm{\bb P^\perp\bb\xi_n}_2/\sqrt p
    \right)^2
    +
    \vnorm{\bb P\bb\xi_n}_2^2/p\,.
$$
Conditionally on $\mathcal G_n$, the two norms are Gaussian norms in
dimensions $p-s$ and $s$, respectively, with the second equal to zero when
$s=0$. Since $s\leq2$ and
$\abs{1-\sqrt{(p-s)/p}}\leq2/p$, Gaussian norm concentration
\citep[Theorem~3.1.1]{vershynin:2018} and a union bound give a universal
constant $c>0$ such that, for every $\delta\in(0,1]$ and all sufficiently
large $n$,
\begin{equation}
\label{eq:gaussian-surrogate-direction-bound}
    \Pr\left(
        \vnorm{\sqrt p\bb d-\bb z_p}_2/\sqrt p
        >
        \delta
        \mid
        \mathcal G_n
    \right)
    \leq
    4\exp\{-c\delta^2p\}
    \quad
    \text{almost surely}\,.
\end{equation}
Now, fix a pseudo-Lipschitz function $\psi:\Re^3\to\Re$ of order two with
Lipschitz-constant $L_\psi$. With $\bbeta^*=\bb B\bv^*$ from
\eqref{eq:limiting-alignment}, write
$\bb r_p=\hat\bbeta-\bbeta^*$ and define the empirical averages
\begin{equation}
\label{eq:test_function_Tp}
\begin{aligned}
    T_p^{\textrm{\tiny DY}}
    &=
    \frac1p
    \sum_{j=1}^p
    \psi
    \left(
        \betady_j-\bbeta_j^*,
        \bbeta_{0,j},
        \bbeta_{P,j}
    \right),
    \quad
    &T_p
    &=
    \frac1p
    \sum_{j=1}^p
    \psi
    \left(
        [\bb r_p]_j,
        \bbeta_{0,j},
        \bbeta_{P,j}
    \right),
    \\
    T_p^{\bb z}
    &=
    \frac1p
    \sum_{j=1}^p
    \psi
    \left(
        \sigma^*[\bb z_p]_j,
        \bbeta_{0,j},
        \bbeta_{P,j}
    \right),
    \quad
    &\bar T_p
    &=
    \frac1p
    \sum_{j=1}^p
    \expect_G
    \left[
        \psi
        \left(
            \sigma^*G,
            \bbeta_{0,j},
            \bbeta_{P,j}
        \right)
    \right] \,, 
\end{aligned}
\end{equation}
together with the limit
$
    T
    =
    \expect
    [
        \psi
        (
            \sigma^*G,
            \bar\beta_0,
            \bar\beta_P
        )
    ]
$,
where
$(\bar\beta_0,\bar\beta_P)\sim\pi_{\{\bar\beta_0,\bar\beta_P\}}$ and
$G\sim\mathrm N(0,1)$ is independent of $(\bar\beta_0,\bar\beta_P)$.
$T_p^{\textrm{\tiny DY}}$ is the average from
\eqref{eq:test_fun_conc} that we wish to link to $T$. 
We shall do so by arguing that (i) $T_p^{\textrm{\tiny DY}} = T_p$ on the MDYPL containment event,
(ii) $\bb r_p$ is close to $\sigma^*\bb z_p$ in normalised Euclidean norm, by PO localisation and the
surrogate bound, (iii) this transfers to $\abs{T_p - T_p^{\bb z}}$ by pseudo-Lipschitz stability, while
\eqref{eq:gaussian-surrogate-law} identifies $\bar T_p$ as the
$\mathcal G_n$-conditional expectation of $T_p^{\bb z}$, around which it
concentrates, and (iv) $\bar T_p$
converges to $T$ by the $W_2$ convergence of the signal array assumed in Section~\ref{sec:setup}.

Thus, for $\epsilon>0$, define
\begin{equation}
\label{eq:test-function-events}
\begin{aligned}
    \mathcal E_n^{\mathrm{vec}}(\epsilon)
    &=
    \left\{
        \vnorm{
            \bb r_p-\sigma^*\bb z_p
        }_2 / \sqrt p 
        \leq
        \epsilon
    \right\},
    \\
    \mathcal E_n^{\mathrm{PL}}(\epsilon)
    &=
    \left\{
        \abs{T_p-T_p^{\bb z}}
        \leq
        \epsilon
    \right\},
    \\
    \mathcal E_n^{\mathrm{G}}(\epsilon)
    &=
    \left\{
        \abs{T_p^{\bb z}-\bar T_p}
        \leq
        \epsilon
    \right\}\,.
\end{aligned}
\end{equation}
We bound the four events in turn, conditionally on $\mathcal G_n$ where a
conditional bound is available, and on the structural event
$\mathcal S_n(t)$ of \eqref{event:S_n}, on which all conditional
constants are deterministic. Since $t_\star\leq1$, on
$\mathcal S_n(t)$ with $t\in(0,t_\star]$, we record that 
\begin{equation}
\label{eq:structural-moment-bound}
    \frac1p
    \sum_{j=1}^p
    \left(
        \bbeta_{0,j}^2
        +
        \bbeta_{P,j}^2
    \right)
    =
    \trace(\bb\Gamma_p)
    \leq
    \trace(\bb\Gamma)
    +
    2\mnorm{\bb\Gamma_p-\bb\Gamma}_2
    \leq
    \trace(\bb\Gamma)
    +
    2
    =
    M_B\,,
\end{equation}
which we will use repeatedly in our error expansions. 
By
Proposition~\ref{prop:mdypl-po-link}, on $\mathcal M_n\cap\mathcal E_n$,
$
    \hat\theta=\thetady$, 
    $
    \hat\bbeta=\betady$, 
    $
    \hat{\bb\eta}=\etady$, 
so that $T_p^{\textrm{\tiny DY}}=T_p$ with
$\Pr(\lnot(\mathcal M_n\cap\mathcal E_n))\leq C\exp\{-cn\}$. 

Next, we work with the PO minimiser of \eqref{eq:lagrangian4} in the form of \eqref{eq:compact_PO_polar}, 
where we replace the noise $\bb d$ by the Gaussian surrogate $\bb z_p$ from \eqref{eq:mdypl-gaussian-surrogate} to get 
$$
    \hat\bbeta-\bbeta^*-\sigma^*\bb z_p
=\bb B(\hat\bv-\bv^*)+\sqrt p(\hat\sigma-\sigma^*)\bb d
+\sigma^*(\sqrt p\bb d-\bb z_p) \,. 
$$
Now, using that 
$\bb B^\top\bb B=p\bb\Gamma_p$, $\vnorm{\bb d}_2=1$, and the operator
H\"older continuity of the positive-semidefinite square root
\citep[][Theorem~X.1.1]{bhatia:1997}, we obtain the upper bound
\begin{equation}
\label{eq:compact_PO_gaussian_l2_decomposition}
\begin{aligned}
    \frac1{\sqrt p}
    \vnorm{
        \bb r_p
        -
        \sigma^*\bb z_p
    }_2
    &\leq
    \vnorm{
        \bb\Gamma_p^{1/2}\hat\bv
        -
        \bu^*
    }_2
    +
    \mnorm{
        \bb\Gamma_p
        -
        \bb\Gamma
    }_2^{1/2}
    \vnorm{\bv^*}_2
    +
    \abs{
        \hat\sigma
        -
        \sigma^*
    }
    +
    \frac{\sigma^*}{\sqrt p}
    \vnorm{
        \sqrt p\bb d
        -
        \bb z_p
    }_2\,.
\end{aligned}
\end{equation}
For $\epsilon\in(0,1]$, set
\begin{equation}
\label{eq:test-function-tolerance}
    \epsilon'
    =
    \frac{\epsilon}{3+\sigma^*},
    \quad
    t_\epsilon
    =
    \min
    \left\{
        t_0(\epsilon'),
        \frac{(\epsilon')^2}
        {(1+\vnorm{\bv^*}_2)^2},
        t_\star
    \right\} \,, 
\end{equation}
with $t_0(\epsilon')$ as in Proposition~\ref{prop:PO_localisation}. On
$\mathcal S_n(t_\epsilon)$,
we have $\mnorm{\bb\Gamma_p-\bb\Gamma}_2^{1/2}\vnorm{\bv^*}_2\leq\epsilon'$ and on the 
localisation event $\{(\hat\bw,\hat\bv,\hat\theta)\in\mathcal S^{\epsilon'}\}$ of Proposition~\ref{prop:PO_localisation}, 
the first
and third terms of \eqref{eq:compact_PO_gaussian_l2_decomposition} are
each less than $\epsilon'$. 
Finally, on 
$\{\vnorm{\sqrt p\bb d-\bb z_p}_2 / \sqrt p \leq\epsilon'\}$, the last term
is at most $\sigma^*\epsilon'$. Hence
$$
    \mathcal S_n(t_\epsilon)
    \cap
    \left\{
        (\hat\bw,\hat\bv,\hat\theta)
        \in
        \mathcal S^{\epsilon'}
    \right\}
    \cap
    \left\{
        \vnorm{
            \sqrt p\bb d-\bb z_p
        }_2 /
        \sqrt p
        \leq
        \epsilon'
    \right\}
    \subseteq
    \mathcal E_n^{\mathrm{vec}}(\epsilon)\,.
$$
Since $t_\epsilon\leq t_0(\epsilon')$,
Proposition~\ref{prop:PO_localisation},
\eqref{eq:gaussian-surrogate-direction-bound} with
$\delta=\epsilon'$, and a union bound give for all $n$ large enough, 
\begin{equation}
\label{eq:test_function_vec_bound}
    \mathds{1}_{\mathcal S_n(t_\epsilon)}
    \Pr
    \left(
        \lnot\mathcal E_n^{\mathrm{vec}}(\epsilon)
        \mid
        \mathcal G_n
    \right)
    \leq
    C_\epsilon
    \exp\{-c_\epsilon n\}
    +
    4\exp\{-c(\epsilon')^2p\}\,.
\end{equation}
Now let us translate this localisation and noise replacement into similarity of the averages $T_p$ and $T_p^{\bb z}$. 
Fix a precision level $u>0$ and with $M_B$ from \eqref{eq:structural-moment-bound}, set
$$
    M_z
    =
    2(\sigma^*)^2+M_B,
    \quad
    M_x
    =
    2+4(\sigma^*)^2+M_B,
    \quad
    L_0
    =
    \left\{
        3
        (
            1+M_x+M_z
        )
    \right\}^{1/2}, \quad 
    \delta_u
    =
    \min
    \left\{
        1,
        \frac{u}{2(1+L_\psi)L_0}
    \right\} \,, 
$$
and define $t_u$ as the $t$ from \eqref{eq:test-function-tolerance} at the tolerance level
$\epsilon=\delta_u$. 
By \eqref{eq:gaussian-surrogate-law} and Gaussian norm concentration
\citep[Theorem~3.1.1]{vershynin:2018},
$$
    \Pr
    \left(
        \frac1p
        \vnorm{\bb z_p}_2^2
        >
        2
        \mid
        \mathcal G_n
    \right)
    \leq
    2\exp\{-cp\}\,.
$$
Work on
$\mathcal S_n(t_u)
\cap\mathcal E_n^{\mathrm{vec}}(\delta_u)
\cap\{p^{-1}\vnorm{\bb z_p}_2^2\leq2\}$,
and collect the arguments of the two averages into
$$
    \bb a_j
    =
    (
        [\bb r_p]_j,
        \bbeta_{0,j},
        \bbeta_{P,j}
    ),
    \quad
    \bb b_j
    =
    (
        \sigma^*[\bb z_p]_j,
        \bbeta_{0,j},
        \bbeta_{P,j}
    ),
    \quad
    j=1,\ldots,p\,,
$$
so that $T_p=p^{-1}\sum_{j=1}^p\psi(\bb a_j)$, 
$T_p^{\bb z}=p^{-1}\sum_{j=1}^p\psi(\bb b_j)$ and 
$$
    \frac1p
    \sum_{j=1}^p
    \vnorm{\bb a_j-\bb b_j}_2^2
    =
    \frac1p
    \vnorm{\bb r_p-\sigma^*\bb z_p}_2^2\,,
$$
which is controlled by
$\mathcal E_n^{\mathrm{vec}}(\delta_u)$.
It remains to bound the two empirical second moments. Splitting the squared
norms into their surrogate and signal parts and using
\eqref{eq:structural-moment-bound} together with
$p^{-1}\vnorm{\bb z_p}_2^2\leq2$,
$$
    \frac1p
    \sum_{j=1}^p
    \vnorm{\bb b_j}_2^2
    =
    (\sigma^*)^2
    \frac{\vnorm{\bb z_p}_2^2}{p}
    +
    \frac1p
    \sum_{j=1}^p
    \left(
        \bbeta_{0,j}^2+\bbeta_{P,j}^2
    \right)
    \leq
    2(\sigma^*)^2+M_B
    =
    M_z\,.
$$
The same split gives that, upon $\mathcal E_n^{\mathrm{vec}}(\delta_u) \cap \{p^{-1}\vnorm{\bb z_p}_2^2\leq2\}$, 
\begin{equation}
    \label{eq:split-one} 
     \begin{aligned}
      \frac1p\sum_{j=1}^p\vnorm{\bb a_j}_2^2
    &=
    \frac1p\vnorm{\bb r_p}_2^2
    +
      \frac1p\sum_{j=1}^p(\bbeta_{0,j}^2+\bbeta_{P,j}^2) \\ 
      & \leq \frac2p\vnorm{\bb r_p - \sigma^* \bb z_p}_2^2 + \frac{2(\sigma^*)^2}p \vnorm{\bb z_p}_2^2 +
      \frac1p\sum_{j=1}^p(\bbeta_{0,j}^2+\bbeta_{P,j}^2) \\  
      &\leq 2 \delta_u^2 + 4(\sigma^*)^2 + M_B \\ 
      & \leq 2 + 4(\sigma^*)^2 + M_B \\ 
      & = M_x \,, 
   \end{aligned}
\end{equation}
where we used that 
$\bb r_p=(\bb r_p-\sigma^*\bb z_p)+\sigma^*\bb z_p$ and
$(a+b)^2\leq2a^2+2b^2$. 
Consequently,
\begin{equation}
    \label{eq:pl-expansion} 
    \begin{aligned}
    \abs{
        T_p-T_p^{\bb z}
    }
    & = 
    \abs{ \frac1p \sum_{j = 1}^p  \psi(\bb a_j) - \psi(\bb b_j) } \\
    &\leq \frac1p \sum_{j = 1}^p  \abs{\psi(\bb a_j) - \psi(\bb b_j) } \\ 
    &\leq  
    L_\psi\frac1p \sum_{j = 1}^p (1 + \vnorm{\bb a_j}_2 + \vnorm{\bb b_j}_2) \vnorm{\bb a_j - \bb b_j}_2 \\ 
    &\leq   L_\psi
    \left[
        \frac1p
        \sum_{j=1}^p
        \left(
            1+\vnorm{\bb a_j}_2+\vnorm{\bb b_j}_2
        \right)^2
    \right]^{1/2}
    \left[
        \frac1p
        \sum_{j=1}^p
        \vnorm{\bb a_j-\bb b_j}_2^2
    \right]^{1/2}
    \\
    &\leq
    L_\psi
    \left\{
        3
        \left(
            1+M_x+M_z
        \right)
    \right\}^{1/2}
    \frac1{\sqrt p}
    \vnorm{
        \bb r_p-\sigma^*\bb z_p
    }_2
    \\
    &=
    L_\psi L_0
    \frac1{\sqrt p}
    \vnorm{
        \bb r_p-\sigma^*\bb z_p
    }_2
    \leq
    L_\psi L_0\delta_u
    \\
    &\leq
    L_\psi L_0
    \frac{u}{2(1+L_\psi)L_0}\\ 
    &
    =
    \frac{L_\psi}{1+L_\psi}
    \frac u2 \\ 
    & \leq
    \frac u2\,,
\end{aligned}
\end{equation}
where we used triangle inequality, the pseudo-Lipschitz property of $\psi$ termwise, followed
by Cauchy--Schwarz
\citep[Theorem~1.37(d)]{rudin+etal:1976}, then
$(1+a+b)^2\leq3(1+a^2+b^2)$ with the two preceding moment bounds, then
the definition of $L_0$ and $\mathcal E_n^{\mathrm{vec}}(\delta_u)$, and
finally the definition of $\delta_u$. 
Therefore, 
$$
    \mathcal S_n(t_u)
    \cap
    \mathcal E_n^{\mathrm{vec}}(\delta_u)
    \cap
    \left\{
        \vnorm{\bb z_p}_2^2 / p
        \leq
        2
    \right\}
    \subseteq
    \mathcal E_n^{\mathrm{PL}}(u/2)\,.
$$
A union bound with
\eqref{eq:test_function_vec_bound} and $p/n\to\kappa\in(0,1)$ gives,
after adjusting constants, for all $n$ large enough,
\begin{equation}
\label{eq:test_function_pl_bound}
    \mathds{1}_{\mathcal S_n(t_u)}
    \Pr
    \left(
        \lnot\mathcal E_n^{\mathrm{PL}}(u/2)
        \mid
        \mathcal G_n
    \right)
    \leq
    C_u\exp\{-c_un\}\,.
\end{equation}
We next prove the concentration of $T_p^{\bb z}$ around
$\bar T_p$. By \eqref{eq:gaussian-surrogate-law}, conditionally on
$\mathcal G_n$ the coordinates of $\bb z_p$ are independent standard
Gaussian variables, so the summands
$$
    Z_{p,j}
    =
    \psi
    \left(
        \sigma^*[\bb z_p]_j,
        \bbeta_{0,j},
        \bbeta_{P,j}
    \right)
    -
    \expect_G
    \left[
        \psi
        \left(
            \sigma^*G,
            \bbeta_{0,j},
            \bbeta_{P,j}
        \right)
    \right],
    \quad
    j=1,\ldots,p\,,
$$
are independent with mean zero. 
The part of $\psi$ that does not depend on the
surrogate cancels under centring, i.e. for
$$
    D_{p,j}(z)
    =
    \psi
    \left(
        \sigma^*z,
        \bbeta_{0,j},
        \bbeta_{P,j}
    \right)
    -
    \psi
    \left(
        0,
        \bbeta_{0,j},
        \bbeta_{P,j}
    \right) \,, 
$$
we have
$Z_{p,j}=D_{p,j}([\bb z_p]_j)-\expect_G[D_{p,j}(G)]$, and
pseudo-Lipschitzness applied to
$(\sigma^*z,\bbeta_{0,j},\bbeta_{P,j})$ and
$(0,\bbeta_{0,j},\bbeta_{P,j})$
gives
$$
    \abs{D_{p,j}(z)}
    \leq
    L_\psi  \sigma^* \max\{\sigma^*, 2\}
    \left(
        \abs z
        +
        z^2
        +
        \vnorm{(\bbeta_{0,j},\bbeta_{P,j})}_2
        \abs z
    \right)\,.
$$
Sub-exponential centring
\citep[Exercise~2.7.10]{vershynin:2018} and the subexponential bounds
$\vnorm{G}_{\psi_1}<\infty$ and $\vnorm{G^2}_{\psi_1}<\infty$ (cf. \citealt[Lemmas~2.7.6\&2.7.7]{vershynin:2018}), give, conditional on $\mathcal G_n$, 
$$
    \begin{aligned}
        \vnorm{Z_{p,j}}_{\psi_1} &= \vnorm{D_{p,j}([\bb z_p]_j) - \expect_G[D_{p,j}(G)]}_{\psi_1} \\ 
        & \leq  C\vnorm{D_{p,j}([\bb z_p]_j) }_{\psi_1} \\ 
        &\leq C L_\psi \sigma^* \max\{2, \sigma^*\} \left(\vnorm{[\bb z_p]_j}_{\psi_1} + \vnorm{[\bb z_p]_j^2}_{\psi_1} + \vnorm{(\bbeta_{0,j},\bbeta_{P,j})}_2\vnorm{[\bb z_p]_j}_{\psi_1} \right) \\
        &\leq C_{\psi} (1 + \vnorm{(\bbeta_{0,j},\bbeta_{P,j})}_2) \\ 
        &= K_{p,j} \,, 
    \end{aligned}
$$
where $C_\psi$ depends only on $L_\psi$, $\sigma^*$ and the subexponential norms of $\abs G, G^2$. 

Now, on $\mathcal S_n(t_u)$, the structural moment bound in \eqref{eq:structural-moment-bound} gives
$\sum_{j=1}^pK_{p,j}^2\leq C_{\psi,M_B}p$ and, since
$\max_j\vnorm{(\bbeta_{0,j},\bbeta_{P,j})}_2^2\leq M_Bp$, also
$\max_jK_{p,j}\leq C_{\psi,M_B}\sqrt p$. 
Thus, Bernstein's inequality for
independent centred sub-exponential variables
\citep[Theorem~2.8.1]{vershynin:2018}, applied conditionally on
$\mathcal G_n$, gives
\begin{equation}
\label{eq:test_function_lln_bound}
\begin{aligned}
    \mathds{1}_{\mathcal S_n(t_u)}
    \Pr
    \left(
        \lnot\mathcal E_n^{\mathrm{G}}(u/2)
        \mid
        \mathcal G_n
    \right)
    &=
    \mathds{1}_{\mathcal S_n(t_u)}
    \Pr
    \left(
        \abs{
            T_p^{\bb z}-\bar T_p
        }
        >
        \frac u2
        \mid
        \mathcal G_n
    \right)
    \\
    &=
    \mathds{1}_{\mathcal S_n(t_u)}
    \Pr
    \left(
        \abs{
            \sum_{j=1}^pZ_{p,j}
        }
        >
        \frac{pu}2
        \mid
        \mathcal G_n
    \right)
    \\
    &\leq
    2\exp
    \left\{
        -c
        \min
        \left(
            \frac{(pu/2)^2}
            {\sum_{j=1}^pK_{p,j}^2},
            \frac{pu/2}
            {\underset{1\leq j\leq p}{\max}K_{p,j}}
        \right)
    \right\}
    \\
    &\leq
    2\exp
    \left\{
        -c
        \min
        \left(
            C_{\psi,M_B}^{-1}pu^2,
            C_{\psi,M_B}^{-1}\sqrt p\,u
        \right)
    \right\}
    \\
    &\leq
    C_u\exp\{-c_u\sqrt n\}\,.
\end{aligned}
\end{equation}


The remaining step involves no event: both $\bar T_p$ and $T$ are
$\mathcal G_n$-measurable, so conditioning on $\mathcal G_n$ is vacuous
for it. By
Assumption~\ref{ass:betas_dist}, the empirical measures
$\pi_p=p^{-1}\sum_{j=1}^p\delta_{(\bbeta_{0,j},\bbeta_{P,j})}$ converge
almost surely in $W_2$ to $\pi_{\{\bar\beta_0,\bar\beta_P\}}$. Writing
$\varphi(\bb b)=\expect_G[\psi(\sigma^*G,\bb b)]$ and
$A_\psi=\abs{\psi(\bb 0)}+2L_\psi$, pseudo-Lipschitzness at
$(0,\bb 0)$ with $t+t^2\leq1+2t^2$ gives
$\abs{\psi(z,\bb b)}\leq A_\psi(1+z^2+\vnorm{\bb b}_2^2)$ and hence,
since $\expect_G[(\sigma^*G)^2]=(\sigma^*)^2$,
$\abs{\varphi(\bb b)}\leq A_\psi(1+(\sigma^*)^2)(1+\vnorm{\bb b}_2^2)$,
while dominated convergence gives continuity of $\varphi$. Continuous functions of at
most quadratic growth are admissible test functions for $W_2$ 
convergence \citep[][Definition~6.8(iv) \& Theorem~6.9]{villani:2009}, so that 
\begin{equation}
\label{eq:test_function_signal_replacement}
    \bar T_p
    =
    \int\varphi\,\mathrm d\pi_p
    \overset{\mathrm{a.s.}}{\longrightarrow}
    \int\varphi\,\mathrm d\pi_{\{\bar\beta_0,\bar\beta_P\}}
    =
    T\,.
\end{equation}

The four bounds and \eqref{eq:test_function_signal_replacement} combine
as follows.
\begin{proposition}[Convergence on test functions]
\label{prop:test_function_convergence}
Assume the conditions of Section~\ref{sec:setup}, and fix a
pseudo-Lipschitz function $\psi:\Re^3\to\Re$ of order two with constant
$L_\psi$.

\begin{enumerate}[label=(\roman*)]
    \item For every $u>0$ there exist $N\in\mathbb N$ and constants $C_u,c_u>0$ such that, for
    $n > N$,
    $$
        \Pr
        \left(
            \abs{
                T_p-\bar T_p
            }
            >
            u
        \right)
        \leq
        C_u
        \exp\{-c_u\sqrt n\}\,,
    $$
    where the constants $c_u, C_u$ and $N$ depend only on $u$, $L_\psi$,
    $\kappa,\alpha,\bb\Gamma,\theta_0^*,\theta_P^*$
    and on the constants in Section~\ref{sec:setup}.
    \item
    $$
        T_p^{\textrm{\tiny DY}}
        \overset{\mathrm{a.s.}}{\longrightarrow}
        T\,.
    $$
\end{enumerate}

\end{proposition}

\begin{proof}
    We prove each claim in turn. 
    \begin{enumerate}[label=(\roman*)]
        \item For part~(i), on
$\mathcal E_n^{\mathrm{PL}}(u/2)\cap\mathcal E_n^{\mathrm{G}}(u/2)$,
$$
    \abs{
        T_p-\bar T_p
    }
    \leq
    \abs{T_p-T_p^{\bb z}}
    +
    \abs{T_p^{\bb z}-\bar T_p}
    \leq
    u\,,
$$
the two summands being bounded by $u/2$ on
$\mathcal E_n^{\mathrm{PL}}(u/2)$ and $\mathcal E_n^{\mathrm{G}}(u/2)$
respectively. Since $\mathcal S_n(t_u)$ is
$\mathcal G_n$-measurable, taking expectations in
\eqref{eq:test_function_pl_bound} and
\eqref{eq:test_function_lln_bound} and splitting over
$\mathcal S_n(t_u)$ gives
$$
    \Pr
    \left(
        \abs{
            T_p-\bar T_p
        }
        >
        u
    \right)
    \leq
    C_u\exp\{-c_un\}
    +
    C_u'\exp\{-c_u'\sqrt n\}
    +
    \Pr
    \left(
        \lnot\mathcal S_n(t_u)
    \right)\,,
$$
where the first term is the pseudo-Lipschitz bound
\eqref{eq:test_function_pl_bound} and the second the Gaussian bound
\eqref{eq:test_function_lln_bound}. The structural concentration bound
following \eqref{event:S_n} gives
$\Pr(\lnot\mathcal S_n(t_u))\leq C\exp\{-cn\}$. Since
$\exp\{-cn\}\leq\exp\{-c\sqrt n\}$, this proves part~(i) after
adjusting constants.
        \item Part~(ii) follows from two applications of the first Borel--Cantelli
lemma (cf. \citealt[p.~308]{shiryaev:2016}). First, the bound in
part~(i) is summable, so applying the lemma along $u=1/m$ and
intersecting over $m\in\mathbb N$ gives
$T_p-\bar T_p\to0$ almost surely. Second,
$\Pr(\lnot(\mathcal M_n\cap\mathcal E_n))\leq C\exp\{-cn\}$ is summable,
so $\mathcal M_n\cap\mathcal E_n$ holds eventually almost surely and
hence $T_p^{\textrm{\tiny DY}}=T_p$ eventually almost surely. Combining
the two with \eqref{eq:test_function_signal_replacement} proves
part~(ii).
\end{enumerate}
\end{proof}
We can finally complete the proof of Theorem~\ref{thm:mdypl-convergence}:
\begin{proof}[Proof of Theorem~\ref{thm:mdypl-convergence}]
For each $m\in\mathbb N$, the unconditional bound in
Proposition~\ref{prop:PO_localisation}, applied with $\epsilon=1/m$,
is summable. The first Borel--Cantelli lemma
(cf. \citealt[p.~308]{shiryaev:2016}), followed by intersection over
$m\in\mathbb N$, therefore gives
\begin{equation}
\label{eq:compact_PO_coordinates_as}
    \hat\sigma
    \overset{\mathrm{a.s.}}{\longrightarrow}
    \sigma^*,
    \quad
    \hat\bu
    =
    \bb\Gamma_p^{1/2}\hat\bv
    \overset{\mathrm{a.s.}}{\longrightarrow}
    \bu^*,
    \quad
    \hat\theta
    \overset{\mathrm{a.s.}}{\longrightarrow}
    \theta^*\,.
\end{equation}
Assumption~\ref{ass:betas_conv} and the first Borel--Cantelli lemma give
$\bb\Gamma_p\to\bb\Gamma$ almost surely. Since
$\bb\Gamma^{1/2}\bv^*=\bu^*$, the operator H\"older continuity of the
positive-semidefinite square root
\citep[][Theorem~X.1.1]{bhatia:1997} gives
\begin{equation}
\label{eq:fixed-alignment-square-root-limit}
    \bb\Gamma_p^{1/2}\bv^*
    \overset{\mathrm{a.s.}}{\longrightarrow}
    \bu^*\,.
\end{equation}
By Proposition~\ref{prop:mdypl-po-link},
$\Pr(\lnot(\mathcal M_n\cap\mathcal E_n))\leq C\exp\{-cn\}$ is summable, so
the first Borel--Cantelli lemma gives that
$\mathcal M_n\cap\mathcal E_n$ holds eventually almost surely, and hence
\begin{equation}
\label{eq:mdypl-identification-eventually}
    \hat\theta
    =
    \thetady,
    \quad
    \hat\bbeta
    =
    \betady
    \quad
    \text{eventually almost surely}\,.
\end{equation}
Part~(i) is the intercept component of
\eqref{eq:compact_PO_coordinates_as} combined with
\eqref{eq:mdypl-identification-eventually}. Part~(ii) is
Proposition~\ref{prop:test_function_convergence}(ii), because
$T_p^{\textrm{\tiny DY}}$ is the empirical average in
\eqref{eq:test_fun_conc}.
\end{proof}

The same argument records the following vector-level consequence.

\begin{lemma}[Gaussian $\ell_2$-coupling for the MDYPL slope]
\label{lemma:mdypl-gaussian-l2-coupling}
On the product extension of
Section~\ref{sec:test-function-convergence}, the vector $\bb z_p$ in
\eqref{eq:mdypl-gaussian-surrogate} satisfies
$\bb z_p\sim\mathrm N(\bb0_p,\bb I_p)$, is independent of
$\sigma(\mathcal C_n,\bb{\mathfrak Z})$, and
\begin{equation}
\label{eq:mdypl-gaussian-l2-coupling}
    \frac1{\sqrt p}
    \vnorm{
        \betady
        -
        \bb B\bv^*
        -
        \sigma^*\bb z_p
    }_2
    \overset{\mathrm{a.s.}}{\longrightarrow}
    0\,.
\end{equation}
\end{lemma}

\begin{proof}
The law and independence of $\bb z_p$ are
\eqref{eq:gaussian-surrogate-law}. For every fixed
$\epsilon\in(0,1]$, with $\epsilon'$ as in
\eqref{eq:test-function-tolerance}, taking expectations in
\eqref{eq:test_function_vec_bound}, splitting over
$\mathcal S_n(t_\epsilon)$, and using the structural concentration bound
following \eqref{event:S_n} gives
$$
    \Pr
    \left(
        \lnot\mathcal E_n^{\mathrm{vec}}(\epsilon)
    \right)
    \leq
    C_\epsilon\exp\{-c_\epsilon n\}
    +
    4\exp\{-c(\epsilon')^2p\}
    +
    \Pr
    \left(
        \lnot\mathcal S_n(t_\epsilon)
    \right)\,.
$$
Since $p/n\to\kappa\in(0,1)$, the right-hand side is summable. The first
Borel--Cantelli lemma, applied along $\epsilon=1/m$ and followed by
intersection over $m\in\mathbb N$, therefore gives
$$
    \frac1{\sqrt p}
    \vnorm{
        \hat\bbeta
        -
        \bb B\bv^*
        -
        \sigma^*\bb z_p
    }_2
    \overset{\mathrm{a.s.}}{\longrightarrow}
    0\,.
$$
Finally, \eqref{eq:mdypl-identification-eventually} gives
$\hat\bbeta=\betady$ eventually almost surely, proving
\eqref{eq:mdypl-gaussian-l2-coupling}.
\end{proof}

\section{Proofs for Section~\ref{sec:pred}}
\label{app:pred}

\subsection{Proof of Theorem~\ref{thm:prediction-classification}}

\begingroup
\renewcommand{\thetheorem}{\arabic{theorem}}
\setcounter{theorem}{1}
\begin{theorem}
\label{thm:prediction-classification}
Assume the conditions of Section~\ref{sec:setup}, and fix
$\alpha\in(0,1)$. The statements below hold for $j=\mathrm{DY}$ and, when
$v_{1,\alpha}^*\neq0$, also for $j=\mathrm{adj}$.
\begin{enumerate}[label=(\roman*)]
    \item Conditionally on $\mathcal F_n$, almost surely,
    \begin{equation}
    \label{eq:prediction-joint-limit}
        \left(
            \eta_0(\bx_{\rm new}),
            \hat\eta_\alpha^j(\bx_{\rm new}),
            Y_{\rm new}
        \right)
        \overset{\mathrm d}{\longrightarrow}
        \left(
            S_0,
            S_\alpha^j,
            Y
        \right) \,.
    \end{equation}
    The convergence also holds unconditionally.
    \item Suppose that $\var(S_\alpha^j)>0$. Then, for every fixed
    $h\in\Re$ and $a\in\{-1,1\}$,
    \begin{equation}
    \label{eq:prediction-risk-convergence}
        \mathcal E_{n,\alpha}^{j,a}(h)
        \overset{\mathrm{a.s.}}{\longrightarrow}
        \mathcal E_\alpha^{j,a}(h) \,,
    \end{equation}
    where
    \begin{equation}
    \label{eq:prediction-limiting-risk}
        \mathcal E_\alpha^{j,a}(h)
        =
        \expect\left[
            \rho'(S_0)
            \Phi_{\rm N}\left\{
                \frac{
                    a\{h-\overline S_\alpha^j(\bb Q)\}
                }{
                    \tau_\alpha^j
                }
            \right\}
            +
            \{1-\rho'(S_0)\}
            \Phi_{\rm N}\left\{
                \frac{
                    a\{\overline S_\alpha^j(\bb Q)-h\}
                }{
                    \tau_\alpha^j
                }
            \right\}
        \right] \,,
    \end{equation}
    and which is understood by continuity if $\tau_\alpha = 0$. Moreover,
    \begin{equation}
    \label{eq:prediction-unconditional-risk-convergence}
        \Pr\left(
            Y_{\rm new}
            \neq
            \widehat Y_\alpha^{j,a}(\bx_{\rm new};h)
        \right)
        =
        \expect\left\{
            \mathcal E_{n,\alpha}^{j,a}(h)
        \right\}
        \longrightarrow
        \mathcal E_\alpha^{j,a}(h) \,.
    \end{equation}
\end{enumerate}
\end{theorem}
\endgroup

\begin{proof}
Fix $\alpha\in(0,1)$ and assume throughout that $v_{1,\alpha}^*\neq0$, so
that $\check{\bbeta}_\alpha$ and $\hat\eta_\alpha^{\mathrm{adj}}$ in
Section~\ref{subsec:mer} are defined. If $v_{1,\alpha}^*=0$, the statements for
$j=\mathrm{adj}$ are void and the argument below applies after
deleting every $\mathrm{adj}$ coordinate. Write
\begin{equation}
\label{eq:prediction-proof-notation}
    \bbeta_\alpha^*
    =
    v_{1,\alpha}^*\bbeta_0+v_{2,\alpha}^*\bbeta_P
    =
    \bb B\bv_\alpha^*,
    \quad
    \bb\Delta_\alpha
    =
    \betady_\alpha-\bbeta_\alpha^*,
    \quad
    \bb M_\alpha
    =
    \begin{bmatrix}
        \bbeta_0 & \betady_\alpha & \check{\bbeta}_\alpha
    \end{bmatrix}
    \in\Re^{p\times3}\,,
\end{equation}
so that, by definition of the limiting logits in Section~\ref{subsec:mer},
$
    \betady_\alpha
    =
    \bb B\bv_\alpha^*+\bb\Delta_\alpha
$
and
$
    \check{\bbeta}_\alpha
    =
    \bbeta_0+\bb\Delta_\alpha/v_{1,\alpha}^*
$.

We shall prove each claim of the Theorem in turn. 
\begin{enumerate}[label=(\roman*)]
    \item Conditional on $\mathcal{F}_n$, the scores $\bb M_\alpha^\top \bx_{\textrm{new}}$ are multivariate Gaussian with 
    variance-covariance matrix $\bb M_\alpha^\top \bb M_\alpha / p$, whose limit $\bSigma_\alpha$ we determine below.
    Theorem~\ref{thm:mdypl-convergence}(ii), applied to the test functions
    $\psi(d,b_0,b_P)\in\{db_0,db_P,d^2\}$, gives
    \begin{equation}
    \label{eq:prediction-proof-residual-aggregates}
    \begin{aligned}
        \frac{\bbeta_0^\top\bb\Delta_\alpha}{p}
        &\overset{\mathrm{a.s.}}{\longrightarrow}0,
        &
        \frac{\bbeta_P^\top\bb\Delta_\alpha}{p}
        &\overset{\mathrm{a.s.}}{\longrightarrow}0,
        &
        \frac{\vnorm{\bb\Delta_\alpha}_2^2}{p}
        &\overset{\mathrm{a.s.}}{\longrightarrow}
        (\sigma_\alpha^*)^2\,,
    \end{aligned}
    \end{equation}
    while Assumption~\ref{ass:betas_conv} and the first Borel--Cantelli
    lemma (cf. \citealt[p.~308]{shiryaev:2016}) give
    $\bb\Gamma_p\overset{\mathrm{a.s.}}{\longrightarrow}\bb\Gamma$. Since
    $\bv_\alpha^{*\top}\bb\Gamma\bv_\alpha^*=\vnorm{\bu_\alpha^*}_2^2$,
    expanding $\bb M_\alpha^\top\bb M_\alpha$ through
    \eqref{eq:prediction-proof-notation} yields
    \begin{equation}
    \label{eq:prediction-proof-gram-limit}
        \frac{\bb M_\alpha^\top\bb M_\alpha}{p}
        \overset{\mathrm{a.s.}}{\longrightarrow}
        \bSigma_\alpha
        =
        \begin{bmatrix}
            \gamma^2
            & c_\alpha
            & \gamma^2
            \\
            c_\alpha
            & \vnorm{\bu_\alpha^*}_2^2+(\sigma_\alpha^*)^2
            & c_\alpha+\dfrac{(\sigma_\alpha^*)^2}{v_{1,\alpha}^*}
            \\
            \gamma^2
            & c_\alpha+\dfrac{(\sigma_\alpha^*)^2}{v_{1,\alpha}^*}
            & \gamma^2+\dfrac{(\sigma_\alpha^*)^2}{(v_{1,\alpha}^*)^2}
        \end{bmatrix},
        \quad
        c_\alpha
        =
        \gamma^2v_{1,\alpha}^*+\varphi v_{2,\alpha}^*\,.
    \end{equation}
    By \eqref{eq:prediction-limiting-logits},
    $
        \bSigma_\alpha
        =
        \cov\left(
            S_0,S_\alpha^{\textrm{\tiny DY}},S_\alpha^{\textrm{adj}}
        \right)
    $.
    Conditionally on $\mathcal F_n$, by defintion of the limiting logits in Sectin~\ref{subsec:mer} and
    $\bx_{\rm new}\sim\mathrm N(\bb0_p,p^{-1}\bb I_p)$ independent of
    $\mathcal F_n$, the score triple
    $$
        \bb\eta_{n,\alpha}
        =
        \left(
            \eta_0(\bx_{\rm new}),
            \hat\eta_\alpha^{\textrm{\tiny DY}}(\bx_{\rm new}),
            \hat\eta_\alpha^{\textrm{adj}}(\bx_{\rm new})
        \right)^\top\,,
    $$
    is Gaussian with mean $(\theta_0,\thetady_\alpha,\theta_0)^\top$ and
    covariance matrix $\bb M_\alpha^\top\bb M_\alpha/p$.
    Assumption~\ref{ass:theta-conv} with the first Borel--Cantelli lemma and
    Theorem~\ref{thm:mdypl-convergence}(i) give
    $
        (\theta_0,\thetady_\alpha,\theta_0)^\top
        \overset{\mathrm{a.s.}}{\longrightarrow}
        (\theta_0^*,\theta_\alpha^*,\theta_0^*)^\top
    $.
    Gaussian characteristic functions are continuous in the mean vector and
    the covariance matrix, so \eqref{eq:prediction-proof-gram-limit} and thus \citep[Chapter~3.3,~Theorem~1]{shiryaev:2016}, conditionally on
    $\mathcal F_n$, almost surely,
    \begin{equation}
    \label{eq:prediction-proof-score-triple-limit}
        \bb\eta_{n,\alpha}
        \overset{\mathrm d}{\longrightarrow}
        \left(
            S_0,
            S_\alpha^{\textrm{\tiny DY}},
            S_\alpha^{\textrm{adj}}
        \right)^\top\,.
    \end{equation}
    Moreover, $\varepsilon_{\rm new}\sim\mathrm{Unif}([0,1])$ is independent of
    $(\mathcal F_n,\bx_{\rm new})$ and $\varepsilon$ is independent of
    $(\bb Q,G)$, so \eqref{eq:prediction-proof-score-triple-limit} holds
    jointly with $\varepsilon_{\rm new}$ appended in the left-hand side and
    $\varepsilon$ in the right-hand side.
    The map
    \begin{equation}
        \label{eq:step-fun} 
          (t,s_1,s_2,e)
        \mapsto
        \left(
            t,s_1,s_2,\mathds{1}\{e<\rho'(t)\}
        \right)\,,
    \end{equation}
    does not depend on $n$ and is continuous except on $\{e=\rho'(t)\}$,
    which has probability zero under the limiting law because $\varepsilon$
    is uniform and independent of $S_0$. The continuous mapping theorem for almost everywhere continuous functions \citep[Chapter~3.8,~Theorem~2]{shiryaev:2016}, applied to \eqref{eq:step-fun}, therefore
    gives, conditionally on $\mathcal F_n$, 
    \begin{equation}
    \label{eq:prediction-proof-full-limit}
        \left(
            \eta_0(\bx_{\rm new}),
            \hat\eta_\alpha^{\textrm{\tiny DY}}(\bx_{\rm new}),
            \hat\eta_\alpha^{\textrm{adj}}(\bx_{\rm new}),
            Y_{\rm new}
        \right)
        \overset{\mathrm d}{\longrightarrow}
        \left(
            S_0,
            S_\alpha^{\textrm{\tiny DY}},
            S_\alpha^{\textrm{adj}},
            Y
        \right)\,.
    \end{equation}
    Retaining coordinates $(1,2,4)$ and $(1,3,4)$ in
    \eqref{eq:prediction-proof-full-limit} proves
    \eqref{eq:prediction-joint-limit} for $j=\mathrm{DY}$ and
    $j=\mathrm{adj}$, respectively. For the unconditional statement, let $f$
    be bounded and continuous on $\Re^3$. Then, by the convergence in distribution \citep[Chapter~2.10,~Definition~4]{shiryaev:2016}, almost surely, 
    $$
        \expect\left[
            f\left(
                \eta_0(\bx_{\rm new}),
                \hat\eta_\alpha^j(\bx_{\rm new}),
                Y_{\rm new}
            \right)
            \mid
            \mathcal F_n
        \right]
        \longrightarrow
        \expect\left[
            f(S_0,S_\alpha^j,Y)
        \right]\,,
    $$
    and these conditional expectations are bounded by $\vnorm{f}_\infty$, so
    dominated convergence \citep[Chapter~2.6,~Theorem~3]{shiryaev:2016} and the tower property \citep[Chapter~3.8~Equation~(4)]{shiryaev:2016} give
    \eqref{eq:prediction-joint-limit} unconditionally.

    \item Fix $h\in\Re$,
    $a\in\{-1,1\}$ and $j\in\{\mathrm{DY},\mathrm{adj}\}$. By
    \eqref{eq:prediction-threshold-classifier} and
    \eqref{eq:prediction-finite-risk},
    \begin{equation}
            \label{eq:prediction-proof-risk-decomposition}
    \begin{aligned}
        \mathcal E_{n,\alpha}^{j,a}(h)
        &=
        \Pr
            \left(
                \{a(\hat\eta_\alpha^j(\bx_{\rm new}) - h) \leq 0 \}   \cap \{Y_{\rm new} = 1\} 
            \mid
            \mathcal F_n
            \right) \\ 
            &+
              \Pr
            \left(
                \{a(\hat\eta_\alpha^j(\bx_{\rm new}) - h) > 0 \}    \cap \{Y_{\rm new} = 0 \} 
            \mid
            \mathcal F_n
            \right) \\
            &= \expect \left[Y_{\rm new} \mathds{1}\{a(\hat\eta_\alpha^j(\bx_{\rm new}) - h) \leq 0\} \mid \mathcal{F}_n \right] \\ 
            &+ \expect \left[ (1-Y_{\rm new}) \mathds{1}\{a(\hat\eta_\alpha^j(\bx_{\rm new}) - h) > 0\} \mid \mathcal{F}_n \right] \\  
            &= \expect \left[\rho'(\eta_0(\bx_{\rm new})) \mathds{1}\{a(\hat\eta_\alpha^j(\bx_{\rm new}) - h) \leq 0\} \mid \mathcal{F}_n \right] \\ 
            &+ \expect \left[ \rho'(-\eta_0(\bx_{\rm new})) \mathds{1}\{a(\hat\eta_\alpha^j(\bx_{\rm new}) - h) > 0\} \mid \mathcal{F}_n \right] \\  
            &= \expect \left[\rho'(\eta_0(\bx_{\rm new})) \Pr \left(a(\hat\eta_\alpha^j(\bx_{\rm new}) - h) \leq 0 \mid \eta_0(\bx_{\rm new}), \mathcal F_n \right) \mid \mathcal{F}_n \right] \\ 
            &+ \expect \left[ \rho'(-\eta_0(\bx_{\rm new}))  \Pr \left(a(\hat\eta_\alpha^j(\bx_{\rm new}) - h) > 0 \mid \eta_0(\bx_{\rm new}), \mathcal F_n \right) \mid \mathcal{F}_n \right] \,,
    \end{aligned}
    \end{equation}
    where the third equality follows from the law of total expectation with respect to
        $\sigma(\mathcal F_n,\bx_{\rm new})\supseteq\mathcal F_n$ and using that $
            \expect[Y_{\rm new}\mid\mathcal F_n,\bx_{\rm new}]
            =\rho'(\eta_0(\bx_{\rm new}))
        $ and the last equality follows from the law of total expectation with respect to
        $
            \mathcal F_n
            \subseteq
            \sigma(\mathcal F_n,\eta_0(\bx_{\rm new}))
            \subseteq
            \sigma(\mathcal F_n,\bx_{\rm new})
        $. 
    As argued in (i), the pair 
    $(\eta_0(\bx_{\rm new}),\hat\eta_\alpha^j(\bx_{\rm new}))$ is bivariate
    Gaussian conditionally on $\mathcal F_n$, so the inner conditional
    probabilities in \eqref{eq:prediction-proof-risk-decomposition} are
    normal distribution functions, and \eqref{eq:prediction-limiting-risk}
    is the limiting form of exactly this expression.
    Observe that the third equality in
    \eqref{eq:prediction-proof-risk-decomposition} is of the form 
    $$
    \begin{aligned}
        \mathcal E_{n,\alpha}^{j,a}(h)
        =
        \expect\left[
            g\left(
                \eta_0(\bx_{\rm new}),
                \hat\eta_\alpha^j(\bx_{\rm new})
            ; h, a\right)
            \mid
            \mathcal F_n
        \right],\quad 
        g(t,s ; h, a)
        =
        \rho'(t)\mathds{1}\{a(s-h)\leq0\}
        +\rho'(-t)\mathds{1}\{a(s-h)>0\}\,,
    \end{aligned}
    $$
   The function $g(\cdot, \cdot; h, a)$ is bounded by one and
continuous except when $s = h$. Since $\var(S_\alpha^j)>0$, this is a set of measure zero for 
 the law of 
    $(S_0,S_\alpha^j)$. The continuous mapping
    theorem \citep[Chapter~3.8,~Theorem~2]{shiryaev:2016}, followed by dominated convergence (cf.
    \citep[Chapter~2.6,~Theorem~3]{shiryaev:2016}), therefore
    gives, almost surely,
    \begin{equation}
    \label{eq:prediction-proof-risk-indicator}
        \mathcal E_{n,\alpha}^{j,a}(h)
        \longrightarrow
        \expect\left[
            g(S_0,S_\alpha^j; h, a )
        \right]\,.
    \end{equation}
    Repeating the last equality of
    \eqref{eq:prediction-proof-risk-decomposition} in the limit, now
    conditioning on $\bb Q$ and using the conditional score moments of Section~\ref{subsec:mer},
    gives, when
    $\tau_\alpha^j>0$,
    $$
    \begin{aligned}
        \Pr\left(
            a(S_\alpha^j-h)\leq0
            \mid
            \bb Q
        \right)
        =
        \Phi_{\rm N}\left\{
            \frac{
                a\{h-\overline S_\alpha^j(\bb Q)\}
            }{
                \tau_\alpha^j
            }
        \right\},
        \quad 
        \Pr\left(
            a(S_\alpha^j-h)>0
            \mid
            \bb Q
        \right)
        =
        \Phi_{\rm N}\left\{
            \frac{
                a\{\overline S_\alpha^j(\bb Q)-h\}
            }{
                \tau_\alpha^j
            }
        \right\}\,,
    \end{aligned}
    $$
    and substituting into \eqref{eq:prediction-proof-risk-indicator} proves
    \eqref{eq:prediction-risk-convergence} with limit
    \eqref{eq:prediction-limiting-risk}. If $\tau_\alpha^j=0$, the same
    formula is understood as the corresponding pointwise limit of the normal
    distribution function.
    Finally, by \eqref{eq:prediction-finite-risk} and the law of total
    expectation \citep[Chapter~3.8~Equation~(4)]{shiryaev:2016},
    $$
        \Pr\left(
            Y_{\rm new}
            \neq
            \widehat Y_\alpha^{j,a}(\bx_{\rm new};h)
        \right)
        =
        \expect\left[
            \mathcal E_{n,\alpha}^{j,a}(h)
        \right]\,.
    $$
    The integrand satisfies $0\leq\mathcal E_{n,\alpha}^{j,a}(h)\leq1$ and,
    by \eqref{eq:prediction-risk-convergence}, converges almost surely to
    $\mathcal E_\alpha^{j,a}(h)$ and dominated convergence \citep[Chapter~2.6,~Theorem~3]{shiryaev:2016} gives
    \eqref{eq:prediction-unconditional-risk-convergence}. 
\end{enumerate}
This concludes the proof.
\end{proof}

\subsection{Proof of Proposition~\ref{prop:gaussian-score-prediction}}

\begingroup
\renewcommand{\thetheorem}{\arabic{theorem}}
\renewcommand{\theproposition}{\arabic{theorem}}
\setcounter{theorem}{0}
\begin{proposition}
\label{prop:gaussian-score-prediction}
Under the Gaussian-score setup in
\eqref{eq:prediction-generic-parameters},
suppose that $\chi\neq0$. Then
\begin{equation}
\label{eq:prediction-optimal-score-rule}
    \mathds{1}\left\{
        \chi(S-h_S)>0
    \right\}
	\, ,
\end{equation}
is the unique minimiser, up to almost-sure equivalence, of
$
    \Pr(Y\neq g(S))
$
over all $\sigma(S)$-measurable functions $g:\Re\to\{0,1\}$. 
\end{proposition}
\endgroup

\begin{proof}
We recall that throughout $\chi\neq0$. By \eqref{eq:prediction-generic-level-set}, the prediction rule in
\eqref{eq:prediction-optimal-score-rule} is
$$
    g^*(s)
    =
    \mathds{1}\{\pi_S(s)>1/2\}
    =
    \mathds{1}\{\chi(s-h_S)>0\}\,.
$$
Let $g:\Re\to\{0,1\}$ be Borel measurable. Conditioning on $S$ and using that 
$\pi_S(S)=\Pr(Y=1\mid S)$, we have that 
$$
    \Pr(Y\neq g(S)\mid S)
    =
    \pi_S(S)+\{1-2\pi_S(S)\}g(S)\,.
$$
Subtracting the same identity for $g^*$ and taking expectations,
$$
\begin{aligned}
    \Pr(Y\neq g(S))
    -
    \Pr(Y\neq g^*(S))
    &=
    \expect\left[
        \{1-2\pi_S(S)\}\{g(S)-g^*(S)\}
    \right]
    \\
    &=
    \expect\left[
        \abs{2\pi_S(S)-1}
        \mathds{1}\{g(S)\neq g^*(S)\}
    \right]
    \geq
    0\,.
\end{aligned}
$$
The second equality holds because the integrand is zero when
$g(S)=g^*(S)$, and otherwise either $\pi_S(S)>1/2$ which then requires
$g(S)-g^*(S)=-1$, or $\pi_S(S)\leq1/2$ and $g(S)-g^*(S)=1$, so that in both
cases the integrand is $\abs{2\pi_S(S)-1}$. This proves optimality of $g^*$. 

To see uniqueness, note that 
the strict monotonicity used to obtain
\eqref{eq:prediction-generic-level-set} also gives
$$
    \{\pi_S(S)=1/2\}
    =
    \{S=h_S\}\,.
$$
Since $\nu>0$ in \eqref{eq:prediction-generic-parameters}, $S$ is a
nondegenerate Gaussian random variable, and thus the event $\{S=h_S\}$  has
probability zero and $\abs{2\pi_S(S)-1}>0$ almost surely. Hence
$$
    \Pr(Y\neq g(S))
    =
    \Pr(Y\neq g^*(S))
$$
if and only if $g(S)=g^*(S)$ almost surely.

This concludes the proof. 
\end{proof}

\subsection{Proof of Proposition~\ref{prop:prediction-comparison}}

\begingroup
\renewcommand{\thetheorem}{\arabic{theorem}}
\renewcommand{\theproposition}{\arabic{theorem}}
\setcounter{theorem}{1}
\begin{proposition}
\label{prop:prediction-comparison}
Assume the conditions of Theorem~\ref{thm:prediction-classification} with
$\gamma^2>0$. For $j\in\{\mathrm{DY},\mathrm{adj}\}$ and every
$\alpha\in(0,1)$ for which $S_\alpha^j$ is defined and informative, consider
\begin{equation}
\label{eq:prediction-mer-logits}
    \underset{a \in \{-1,1\}}{\min}\underset{h\in\Re}{\inf} \, \mathcal E_\alpha^{j,a}(h)
    =
    \mathcal E^*((R_\alpha^j)^2),
    \, 
    \mathcal E^*(r)
    =
    \expect\left[
        \rho'(S_0)\,\Phi_{\rm N}\left(\frac{h_r-S_0}{\varsigma_r}\right)
        +\left\{1-\rho'(S_0)\right\}
        \Phi_{\rm N}\left(\frac{S_0-h_r}{\varsigma_r}\right)
    \right] \,,
\end{equation}
for $\varsigma_r^2=\gamma^2(1-r)/r$ and
$h_r=-\theta_0^*(1-r)/r$, $r \in (0,1)$ and where for $r=1$, we set
$\mathcal E^*(1)=\expect[\min\{\rho'(S_0),1-\rho'(S_0)\}]$.

Then $\mathcal E^*$ is strictly decreasing on $(0,1]$, and when the respective scores are defined and informative, 
\begin{enumerate}[label=(\roman*)]
    \item up to ties, and if attained, 
    \begin{equation}
    \label{eq:prediction-optimal-alpha}
        \underset{\alpha \in (0,1)}{\arg\min}\, \mathcal E^*((R_\alpha^j)^2)
        =
        \underset{\alpha \in (0,1)}{\arg\max}\,\left(R_\alpha^j\right)^2 \,, 
    \end{equation}
    which for the adjusted logit, is the $\alpha$ that minimises slope error, i.e. 
    $$
        \underset{\alpha \in (0,1)}{\arg\min}\,
        \mathcal E^*((R_\alpha^{\rm adj})^2)
        =
        \underset{\alpha \in (0,1)}{\arg\min}\,
        \frac{(\sigma_\alpha^*)^2}{(v_{1,\alpha}^*)^2}\,.
    $$
    \item oracle
    adjustment lowers the optimised error if and only if
    $(R_\alpha^{\mathrm{adj}})^2>(R_\alpha^{\textrm{\tiny DY}})^2$, which
    holds whenever $\varphi=0$, $v_{1,\alpha}^*\neq0$ and
    $\delta^2(v_{2,\alpha}^*)^2>0$.
\end{enumerate}
\end{proposition}
\endgroup

\begin{proof}
Let $S$ satisfy \eqref{eq:prediction-generic-parameters} with $\gamma^2>0$
and $\chi\neq0$, write $r=R^2$, and set
$$
    W_r
    =
    \theta_0^*+\frac{\gamma^2}{\chi}(S-\mu)\,.
$$
Since $\cov(S_0,W_r)=\gamma^2=\var(S_0)$ and
$\var(W_r)=\gamma^2/r$, joint Gaussianity gives
\begin{equation}
\label{eq:prediction-proof-channel}
    W_r=S_0+\epsilon_r,
    \quad
    \epsilon_r\perp S_0,
    \quad
    \epsilon_r\sim\mathrm N(0,\varsigma_r^2)\,.
\end{equation}
Moreover, the generic prediction threshold $h_S = \mu - (\nu^2 / \chi) \theta_0^*$ gives
$$
    \theta_0^*+\frac{\gamma^2}{\chi}(h_S-\mu)=h_r,
    \quad
    \frac{\chi^2}{\gamma^2}(W_r-h_r)=\chi(S-h_S)\,.
$$
Thus the optimal rule in
Proposition~\ref{prop:gaussian-score-prediction} is
$\mathds 1(W_r>h_r)$, and
\begin{equation}
\label{eq:prediction-proof-attained}
    \underset{a\in\{-1,1\}}{\min}\ 
    \underset{h\in\Re}{\inf}\,
    \Pr\left[Y\neq\mathds 1\{a(S-h)>0\}\right]
    =
    \Pr\left\{Y\neq\mathds 1(W_r>h_r)\right\}\,.
\end{equation}
For $0<r<1$, 
$$
\begin{aligned}
    \Pr\left\{Y\neq\mathds 1(W_r>h_r)\right\}
    =
    \expect\left[
        \rho'(S_0)\Phi_{\rm N}\left(\frac{h_r-S_0}{\varsigma_r}\right)
        +
        \{1-\rho'(S_0)\}
        \Phi_{\rm N}\left(\frac{S_0-h_r}{\varsigma_r}\right)
    \right]
    =
    \mathcal E^*(r)\,.
\end{aligned}
$$
At $r=1$, $W_1=S_0$ and $h_1=0$, so the left-hand side of
\eqref{eq:prediction-proof-attained} equals
$\expect[\min\{\rho'(S_0),1-\rho'(S_0)\}]$, because $\rho'(S_0)<1/2$ if and
only if $S_0<0$. 
Substituting the means, variances and covariances of
$S_\alpha^{\textrm{\tiny DY}}$ and $S_\alpha^{\mathrm{adj}}$ from
\eqref{eq:prediction-limiting-logits} into
$R^2=\chi^2/(\gamma^2\nu^2)$ gives
\eqref{eq:prediction-information-logits}, and hence
\eqref{eq:prediction-mer-logits}.

It remains to prove strict monotonicity. Fix $0<r_1<r_2\leq1$ and couple
the corresponding channels as
$$
    W_2=S_0+\epsilon_2,
    \quad
    W_1=W_2+\xi\,,
$$
where $\epsilon_2$ and $\xi$ are independent of $S_0$ and of each other,
with
$$
    \var(\epsilon_2)
    =
    \gamma^2\left(\frac1{r_2}-1\right),
    \quad
    \var(\xi)
    =
    \gamma^2\left(\frac1{r_1}-\frac1{r_2}\right)>0\,,
$$
so that $(S_0,W_k)$ has the law of \eqref{eq:prediction-proof-channel} with
$r=r_k$. For $\pi_k=\expect[\rho'(S_0)\mid W_k]$,
$k\in\{1,2\}$, we have by construction 
$\pi_1=\expect[\pi_2\mid W_1]$, and, since
$\{W_k>h_{r_k}\}=\{\pi_k>1/2\}$ by
\eqref{eq:prediction-generic-level-set} applied to $W_k$, the error
$\mathcal E^*(r_k)$ of $\mathds 1(W_k>h_{r_k})$ is
$\expect[\min\{\pi_k,1-\pi_k\}]$. By the strict monotonicity used to
obtain \eqref{eq:prediction-generic-level-set}, $\pi_2-1/2$ has the sign
of $W_2-h_{r_2}$, while
$$
    \var(W_2\mid W_1)
    =
    \frac{\gamma^2(r_2-r_1)}{r_2^2}>0\,,
$$
so that, conditionally on $W_1$, $W_2$ has full support and
$\pi_2-1/2$ takes both signs with positive probability. The strict
triangle inequality therefore gives, almost surely,
$$
\begin{aligned}
    \min\{\pi_1,1-\pi_1\}
    &=
    \frac12-\abs{\expect[\pi_2-1/2\mid W_1]}
    \\
    &>
    \frac12-\expect[\abs{\pi_2-1/2}\mid W_1]
    =
    \expect\left[\min\{\pi_2,1-\pi_2\}\mid W_1\right] \,.
\end{aligned}
$$
Taking expectations yields
$\mathcal E^*(r_1)>\mathcal E^*(r_2)$, so $\mathcal E^*$ is strictly
decreasing on $(0,1]$.

Part~(i) follows immediately. For the adjusted logit,
$$
    (R_\alpha^{\mathrm{adj}})^2
    =
    \frac{\gamma^2}{\gamma^2+\tau_\alpha^2},
    \quad
    \tau_\alpha^2
    =
    \frac{(\sigma_\alpha^*)^2}{(v_{1,\alpha}^*)^2}\,,
$$
so maximising $(R_\alpha^{\mathrm{adj}})^2$ is equivalent to minimising
$(\sigma_\alpha^*)^2/(v_{1,\alpha}^*)^2$.
For part~(ii), strict monotonicity gives the stated comparison. If
$\varphi=0$, then
$\vnorm{\bu_\alpha^*}_2^2
=\gamma^2(v_{1,\alpha}^*)^2+\delta^2(v_{2,\alpha}^*)^2$, and therefore
$$
    \left(R_\alpha^{\textrm{\tiny DY}}\right)^2
    =
    \frac{\gamma^2(v_{1,\alpha}^*)^2}{
        \gamma^2(v_{1,\alpha}^*)^2
        +\delta^2(v_{2,\alpha}^*)^2
        +(\sigma_\alpha^*)^2
    },
    \quad
    \left(R_\alpha^{\mathrm{adj}}\right)^2
    =
    \frac{\gamma^2(v_{1,\alpha}^*)^2}{
        \gamma^2(v_{1,\alpha}^*)^2
        +(\sigma_\alpha^*)^2
    }\,,
$$
so $(R_\alpha^{\mathrm{adj}})^2\geq(R_\alpha^{\textrm{\tiny DY}})^2$,
strictly if and only if $\delta^2(v_{2,\alpha}^*)^2>0$.
\end{proof}

\section{Proofs for Section~\ref{sec:testing}}
\label{sec:proofs-testing}

\subsection{Finite-contrast Gaussian limit}

\begin{lemma}[Finite-contrast Gaussian limit]
\label{lemma:finite-contrast-replacement}
Assume the conditions of Section~\ref{sec:setup}. Let
$\bb D_p\in\Re^{k\times p}$ be $\mathcal G_n$-measurable, with $k$ fixed,
and suppose that
\begin{equation}
\label{eq:fixed-contrast-proof-normalisation}
    \bb D_p\bb D_p^\top
    =
    \bb I_k,
    \quad
    \mnorm{
        \bb D_p
        \bb B
        \bb\Gamma_p^{+1/2}
    }_2
    =
    \mathcal O_{\mathrm p}(1)\,.
\end{equation}
Let
$(\hat\bv,\hat\bw,\hat\theta,\hat{\bb\eta})$
be the minimising variables fixed in
Proposition~\ref{prop:mdypl-po-link}, and let $\bb d$ be as in
\eqref{eq:compact_PO_direction}. Then
\begin{equation}
\label{eq:finite-contrast-alignment-replacement}
\begin{aligned}
    \bb D_p\bb B(\hat\bv-\bv^*)
    &=
    o_{\mathrm p}(1),
    \\
    \mnorm{
        \bb D_p\bb P\bb D_p^\top
    }_2
    &\overset{\mathrm p}{\longrightarrow}
    0\,.
\end{aligned}
\end{equation}
Moreover,
\begin{equation}
\label{eq:finite-contrast-projection-replacement}
    \sqrt p\,\bb D_p\bb d
    \overset{\mathrm d}{\longrightarrow}
    \mathrm N(\bb0_k,\bb I_k)\,,
\end{equation}
and
\begin{equation}
\label{eq:fixed-contrast-gaussian-limit}
    \frac1{\sigma^*}
    \bb D_p
    \left(
        \betady-\bb B\bv^*
    \right)
    \overset{\mathrm d}{\longrightarrow}
    \mathrm N(\bb0_k,\bb I_k)\,.
\end{equation}
\end{lemma}

\begin{proof}
By \eqref{eq:compact_PO_coordinates_as} and
\eqref{eq:fixed-alignment-square-root-limit},
$$
    \vnorm{
        \bb\Gamma_p^{1/2}
        (
            \hat\bv-\bv^*
        )
    }_2
    \overset{\mathrm{a.s.}}{\longrightarrow}
    0\,.
$$
Since
$\operatorname{null}(\bb\Gamma_p)=\operatorname{null}(\bb B)$,
$$
    \bb B
    =
    \bb B
    \bb\Gamma_p^{+1/2}
    \bb\Gamma_p^{1/2}\,.
$$
Consequently,
$$
\begin{aligned}
    \vnorm{
        \bb D_p\bb B
        (
            \hat\bv-\bv^*
        )
    }_2
    &\leq
    \mnorm{
        \bb D_p
        \bb B
        \bb\Gamma_p^{+1/2}
    }_2
    \vnorm{
        \bb\Gamma_p^{1/2}
        (
            \hat\bv-\bv^*
        )
    }_2
    =
    o_{\mathrm p}(1)\,.
\end{aligned}
$$
Since
$\bb P=\bb B\bb\Gamma_p^+\bb B^\top/p$,
$$
\begin{aligned}
    \mnorm{
        \bb D_p
        \bb P
        \bb D_p^\top
    }_2
    &\leq
    \frac1p
    \mnorm{
        \bb D_p
        \bb B
        \bb\Gamma_p^{+1/2}
    }_2^2
    \overset{\mathrm p}{\longrightarrow}
    0\,.
\end{aligned}
$$
This proves
\eqref{eq:finite-contrast-alignment-replacement}.

Conditionally on $\mathcal C_n$,
Lemma~\ref{lemma:compact_PO_spherical_representation} gives
$$
    \sqrt p\,\bb D_p\bb d
    \overset{\mathrm d}{=}
    \frac{
        \bb D_p\bb P^\perp\bb g_p
    }{
        \vnorm{\bb P^\perp\bb g_p}_2/\sqrt p
    },
    \quad
    \bb g_p
    \sim
    \mathrm N(\bb0_p,\bb I_p)\,,
$$
where $\bb g_p$ may be taken independent of
$\sigma(\mathcal C_n,\bb{\mathfrak Z})$.
Conditionally on $\mathcal G_n$,
$$
    \bb D_p\bb P^\perp\bb g_p
    \sim
    \mathrm N
    \left(
        \bb0_k,
        \bb D_p\bb P^\perp\bb D_p^\top
    \right)\,,
$$
and, by \eqref{eq:finite-contrast-alignment-replacement},
$$
    \bb D_p\bb P^\perp\bb D_p^\top
    =
    \bb I_k
    -
    \bb D_p\bb P\bb D_p^\top
    \overset{\mathrm p}{\longrightarrow}
    \bb I_k\,.
$$
Moreover, with $s=\rank(\bb P)\leq2$,
$$
    \vnorm{\bb P^\perp\bb g_p}_2^2
    \mid
    \mathcal G_n
    \sim
    \chi_{p-s}^2,
    \quad
    \frac{
        \vnorm{\bb P^\perp\bb g_p}_2
    }{
        \sqrt p
    }
    \overset{\mathrm p}{\longrightarrow}
    1\,.
$$
Cram\'er--Wold and Slutsky's theorem
\citep[Sections~1.5.2 and~1.5.4]{serfling:1980} prove
\eqref{eq:finite-contrast-projection-replacement}.

Set
$$
    \hat\bbeta
    =
    \bb B\hat\bv+\bb E\hat\bw,
    \quad
    \hat\sigma
    =
    \frac{\vnorm{\hat\bw}_2}{\sqrt p}\,.
$$
The exact PO decomposition is
\begin{equation}
\label{eq:fixed-contrast-compact-decomposition}
\begin{aligned}
    \frac1{\sigma^*}
    \bb D_p
    \left(
        \hat\bbeta-\bb B\bv^*
    \right)
    =
    \frac1{\sigma^*}
    \bb D_p\bb B
    (
        \hat\bv-\bv^*
    )
    +
    \frac{\hat\sigma}{\sigma^*}
    \sqrt p\,\bb D_p\bb d\,.
\end{aligned}
\end{equation}
Proposition~\ref{prop:limiting_AO_properties} gives $\sigma^*>0$, and
\eqref{eq:compact_PO_coordinates_as} gives
$\hat\sigma\overset{\mathrm p}{\longrightarrow}\sigma^*$.
Combining
\eqref{eq:finite-contrast-alignment-replacement},
\eqref{eq:finite-contrast-projection-replacement}, and
\eqref{eq:fixed-contrast-compact-decomposition} yields
$$
    \frac1{\sigma^*}
    \bb D_p
    \left(
        \hat\bbeta-\bb B\bv^*
    \right)
    \overset{\mathrm d}{\longrightarrow}
    \mathrm N(\bb0_k,\bb I_k)\,.
$$
Proposition~\ref{prop:mdypl-po-link} gives
$\hat\bbeta=\betady$ with probability tending to one. This proves
\eqref{eq:fixed-contrast-gaussian-limit}.
\end{proof}

\subsection{Proof of Theorem~\ref{thm:isotropic-adjusted-z}}
\label{sec:proof-isotropic-adjusted-z}

\begingroup
\renewcommand{\thetheorem}{\arabic{theorem}}
\setcounter{theorem}{2}
\begin{theorem}
\label{thm:isotropic-adjusted-z}
Assume the conditions of Section~\ref{sec:setup}. Let
$I\subseteq\{1,\ldots,p\}$ be a deterministic sequence of index sets
with $\abs{I}=k$ fixed, write
$
    \bb B_I
    =
    \bb J_I^\top\bb B
$,
assume the leverage condition
\begin{equation}
\label{eq:isotropic-adjusted-z-block-leverage}
    \mnorms{
        \bb B_I
        \bb\Gamma_p^{+1/2}
    }_2
    =
    \mathcal O_{\mathrm p}(1) \,,
\end{equation}
and let
$\bv^*=[v_1^*,v_2^*]^\top$
be defined in \eqref{eq:limiting-alignment}. Then
\begin{equation}
\label{eq:isotropic-fixed-block-limit}
    \frac1{\sigma^*}
    \left(
        \betady_I-\bb B_I\bv^*
    \right)
    \overset{\mathrm d}{\longrightarrow}
    \mathrm N(\bb0_k,\bb I_k) \,.
\end{equation}
Consequently, for $v_1^* \neq 0$ and a deterministic $b_I^0\in\Re^k$, under
$
    H_0:
    \bbeta_{0,I}=b_I^0 
$, 
we have
\begin{equation}
\label{eq:isotropic-adjusted-z-limit}
    \bb Z_I^{\mathrm{adj}}(b_I^0)
    =
    \frac{\check{\bbeta}_I-b_I^0}{(\sigma^* / v_1^*)}
    \overset{\mathrm d}{\longrightarrow}
    \mathrm N(\bb0_k,\bb I_k) \,.
\end{equation}
\end{theorem}
\endgroup

\begin{proof}
Take $\bb D_p=\bb J_I^\top$. Then
$\bb D_p\bb D_p^\top=\bb I_k$, and
\eqref{eq:isotropic-adjusted-z-block-leverage} is exactly the leverage
assumption of Lemma~\ref{lemma:finite-contrast-replacement}.
Hence \eqref{eq:fixed-contrast-gaussian-limit} gives
\eqref{eq:isotropic-fixed-block-limit}.
If $v_1^*\neq0$, then under $H_0$,
$$
\begin{aligned}
    \betady_I-\bb B_I\bv^*
    &=
    \betady_I
    -
    v_1^*\bbeta_{0,I}
    -
    v_2^*\bbeta_{P,I}
    \\
    &=
    v_1^*
    \left(
        \check{\bbeta}_I-b_I^0
    \right)\,,
\end{aligned}
$$
which gives \eqref{eq:isotropic-adjusted-z-limit}.
\end{proof}

\section{Proofs for Section~\ref{sec:beyond-isotropic}}
\label{sec:proofs-extensions}

\subsection{Proof of
Corollary~\ref{cor:gaussian-covariance-empirical-law}}
\label{sec:proof-gaussian-covariance-empirical-law}

\begingroup
\renewcommand{\thetheorem}{\arabic{theorem}}
\renewcommand{\thecorollary}{\arabic{theorem}}
\setcounter{theorem}{0}
\begin{corollary}
\label{cor:gaussian-covariance-empirical-law}
Assume that the transformed parameters
$\widetilde\theta_0,\widetilde\theta_P,
\widetilde\bbeta_0,\widetilde\bbeta_P$ satisfy the
conditions of Section~\ref{sec:setup}, and let
$\bv^*=[v_1^*,v_2^*]^\top$,
$\widetilde\theta^*$ and $\sigma^*$ denote the corresponding limiting state
parameters of \eqref{eq:FOCs}.
\begin{enumerate}[label=(\roman*)]
    \item Suppose that, almost surely,
    \begin{equation}
    \label{eq:original-intercept-shift-conditions}
        \bb m_p\longrightarrow\bb m,
        \quad
        \limsup_{p\to\infty} \, 
        \bb m_p^\top\widetilde{\bb\Gamma}_p^+\bb m_p<\infty,
        \quad
        q_p\log p\longrightarrow0\,,
    \end{equation}
    for a deterministic
    $\bb m=(m_0,m_P)^\top\in\Re^2$. Then
    \begin{equation}
    \label{eq:gaussian-covariance-original-intercept-limit}
        \thetady\overset{\mathrm{a.s.}}{\longrightarrow}
        \widetilde\theta^*-\bb m^\top\bv^*
        =\widetilde\theta^*-v_1^*m_0-v_2^*m_P\,.
    \end{equation}

    \item Suppose that
    $\sup_p \, \mnorms{\bb R_p}_2<\infty$, 
    and, almost surely,
    \begin{equation}
    \label{eq:gaussian-covariance-standardised-signal-law}
        \frac1p\sum_{j=1}^p
        \delta_{(\sqrt p\,\tau_j\bbeta_{0,j},
        \sqrt p \tau_j\bbeta_{P,j})}
        \overset{W_2}{\longrightarrow}\pi_\Sigma\,,
    \end{equation}
    for a deterministic probability law
    $\pi_\Sigma$. Then, for any pseudo-Lipschitz function $\psi$
    of order two,
    \begin{equation}
    \label{eq:gaussian-covariance-test-function-law}
    \begin{aligned}
        \frac1p\sum_{j=1}^p
        \psi\left(
            \sqrt p\tau_j
            \{\betady_j-v_1^*\bbeta_{0,j}-v_2^*\bbeta_{P,j}\},
            \sqrt p\tau_j\bbeta_{0,j},
            \sqrt p\tau_j\bbeta_{P,j}
        \right)
       \overset{\mathrm{a.s.}}{\longrightarrow}
        \expect\left[
            \psi(\sigma^*G,\bar\beta_{0,\Sigma},\bar\beta_{P,\Sigma})
        \right]\,,
    \end{aligned}
    \end{equation}
    where
    $
        (\bar\beta_{0,\Sigma},\bar\beta_{P,\Sigma})\sim\pi_\Sigma$, 
        $G\sim\mathrm N(0,1)
    $,
    independent from each other.
\end{enumerate}
\end{corollary}
\endgroup

\begin{proof}
Work throughout with the centred and whitened problem in
\eqref{eq:affine-reparameterisation}. Let
$(\hat\bbeta,\hat\bv,\hat\sigma,\bb d,\bb z_p)$ denote its PO from \eqref{eq:lagrangian5} and
Gaussian-surrogate quantities, and write $\widetilde{\mathcal C}_n$ for the
field $\mathcal C_n$ in \eqref{eq:sigma-fields} formed from the transformed
problem.

We shall prove each claim in turn. 
\begin{enumerate}[label=(\roman{*})]
    \item 
By \eqref{eq:invariance}, write 
$$
    \thetady = \widetilde{\theta}^{\textrm{\tiny DY}} - \left(p^{-1/2} \bb L_p^{-1} \bmu_p\right)^\top \widetilde{\bbeta}^{\textrm{\tiny DY}} \,,
$$  
where $\widetilde{\theta}^{\textrm{\tiny DY}}, \widetilde{\bbeta}^{\textrm{\tiny DY}}$ are the whitened, centred MDYPL estimates. By Theorem~\ref{thm:mdypl-convergence}(i), 
$$
\widetilde{\theta}^{\textrm{\tiny DY}} \overset{\textrm{a.s.}}{\longrightarrow} \widetilde \theta^* \,. 
$$
Now let $\hat \bbeta$ be the primal minimiser from the PO of \eqref{eq:lagrangian5} for the transformed parameters $\widetilde{\bb B}_p, \widetilde{\theta}_0, \widetilde{\theta}_P$.
Then on the containment event $\mathcal{M}_n \cap \mathcal{E}_n$ from \eqref{event:E_n} and \eqref{event:M_n}, we have that 
\begin{equation}
    \label{eq:po-decomp} 
    \widetilde{\bbeta}^{\textrm{\tiny DY}} = \hat \bbeta = \widetilde{\bb B}_p \hat \bv + \hat \sigma \sqrt{p} \bb d \,, 
\end{equation}
where the second equality follows from \eqref{eq:PO-slope} with  $\bb d\in\range(\widetilde{\bb B}_p)^\perp$, $\vnorm{\bb d}_2=1$. 
Define the whitened covariate mean 
$$
\bb a_p = p^{-1/2} \bb L_p^{-1} \bmu_p
$$
and its component that is orthogonal to $\widetilde{\bb B}_p$, i.e. 
$$
    \bb a_p^\perp = \bb P^\perp_{{\widetilde{\bb B}_p}} \bb a_p = p^{-1/2} \bb L_p^{-1} \bmu_p^\perp, \quad \bb \mu^\perp_p = \bb \mu_p - \bb \Sigma_p \bb B \widetilde{\bb \Gamma}_p^{+} \bb m_p \,. 
$$
Then, by \eqref{eq:po-decomp} and since $\bb d \in \range(\widetilde{\bb B}_p)^\perp$, 
$$
    \begin{aligned}
        \bb \mu_p^\top \betady &= \bb a_p^\top \left(\widetilde{\bb B}_p \hat \bv + \hat \sigma \sqrt p \bb d \right) 
        = \bb m_p^\top \hat \bv + \hat \sigma  \langle \bb a_p^\perp,\sqrt p \bb d \rangle \,. 
    \end{aligned}
$$
To control the first term, note that 
$$
    \begin{aligned}
        \abs{\bb m_p^\top (\hat \bv - \bv^*)} &= \abs{\bb m_p^\top \widetilde{\bb \Gamma}_p^{+1/2} \widetilde{\bb \Gamma}_p^{1/2} (\hat \bv - \bv^*)} \\ 
        &\leq \vnorm{\bb m_p^\top \widetilde{\bb \Gamma}_p^{+1/2}}_2 \left\{ \vnorm{ \widetilde{\bb \Gamma}_p^{1/2} \hat \bv - \bu^*}_2  + \vnorm{(\widetilde{\bb \Gamma}_p^{1/2} -  \widetilde{\bb \Gamma}^{1/2}) \bv^* }_2\right\} \\ 
        &\leq  \vnorm{\bb m_p^\top \widetilde{\bb \Gamma}_p^{+1/2}}_2 \left\{ \vnorm{ \widetilde{\bb \Gamma}_p^{1/2} \hat \bv - \bu^*}_2  + \mnorm{\widetilde{\bb \Gamma}_p  - \widetilde{\bb \Gamma}}_2^{1/2} \vnorm{\bv^*}_2\right\} \,,
    \end{aligned}
$$
where we used the Cauchy--Schwarz inequality, the definition of $\bu^*$ from Section~\ref{sec:ao-limit}, the triangle inequality and operator H\"older continuity of the positive-semidefinite
square root \citep[Theorem~X.1.1]{bhatia:1997}. Now, the term $\vnorm{\bb m_p^\top \widetilde{ \bb \Gamma}_p^{+1/2}}_2$ is $\mathcal{O}(1)$ almost surely 
by the assumption of \eqref{eq:original-intercept-shift-conditions}. The second term $\vnorm{\widetilde{\bb \Gamma}_p^{1/2} \hat \bv - \bu^*}_2 \to 0$ almost surely by the transformed versions of
\eqref{eq:compact_PO_coordinates_as} and
\eqref{eq:fixed-alignment-square-root-limit} with limiting signal matrix $\widetilde{\bb \Gamma}$. The transformed version of Assumption~\ref{ass:betas_conv} and the first Borel--Cantelli lemma give
$\widetilde{\bb\Gamma}_p\to\widetilde{\bb\Gamma}$ almost surely, and hence
$\mnorms{\widetilde{\bb\Gamma}_p-\widetilde{\bb\Gamma}}_2^{1/2}\to0$.
Since $\bv^*$ is a fixed finite vector, this controls the remaining term.

It thus remains to control $\hat \sigma \langle \bb a_p^\perp, \sqrt p \bb d \rangle$. Since by \eqref{eq:compact_PO_coordinates_as}, $\hat{\sigma} \to \sigma^*$ almost surely, it is sufficient to establish that $\langle \bb a_p^\perp, \sqrt p \bb d \rangle \to 0$ almost surely. 
Conditionally on $\widetilde{\mathcal C}_n$,
Lemma~\ref{lemma:compact_PO_spherical_representation} gives that $\bb d$ is
uniform on the unit sphere of $\range(\widetilde{\bb B}_p)^\perp$, whose
dimension is $p-s$ for $s=\rank(\widetilde{\bb B}_p)\leq2$. Since, conditionally on $\widetilde{\mathcal C}_n$, $\sqrt p \bb d$ is subgaussian
\citep[Theorem~3.4.6]{vershynin:2018} with constant subgaussian norm, and 
$\bb a_p^\perp \in \range(\widetilde{\bb B}_p)^\perp$ is $\widetilde{\mathcal C}_n$-measurable, with squared norm $q_p$, a Hoeffding bound \citep[Theorem~2.6.3]{vershynin:2018}, therefore gives absolute constants
$C,c>0$ such that, for every $\epsilon>0$ and all sufficiently large $p$,
\begin{equation}
    \label{eq:inner-prod-bound} 
        \mathds{1}\{q_p > 0 \} \Pr
    \left(
        \abs{
            \langle
                \bb a_p^\perp,
                \sqrt p \bb d
            \rangle
        }
        >
        \epsilon
        \mid
        \widetilde{\mathcal C}_n
    \right)
    \leq
    C
    \exp
    \left\{
        -
        \frac{c\epsilon^2}{q_p}
    \right\}\,. 
\end{equation}
Now, on the
$\widetilde{\mathcal C}_n$-measurable event
$\{q_p\log p\leq c\epsilon^2/3\}$, the right-hand side of \eqref{eq:inner-prod-bound} is at most
$Cp^{-3}$. 
On $\{q_p=0\}$, trivially, $\langle \bb a_p^\perp,\sqrt p\,\bb d \rangle = 0$. 
Taking expectations, applying the first Borel--Cantelli lemma
(cf. \citealt[p.~308]{shiryaev:2016}), and using
$q_p\log p\to0$ almost surely gives
$$
    \sqrt p
    \abs{
        \langle
            \bb a_p^\perp,
            \bb d
        \rangle
    }
    \leq
    \epsilon
$$
eventually almost surely. Intersecting over $\epsilon=1/m$,
$m\in\mathbb N$, proves $\langle \bb a_p^\perp , \sqrt p \bb d \rangle \to 0$ almost surely.
Finally, the summability of the compact-containment failure probabilities in
Lemma~\ref{lemma:containment} and the first Borel--Cantelli lemma give
$$
    \bb \mu_p^\top \betady \overset{\textrm{a.s.}}{\longrightarrow} \bb m^\top \bv^* \,,
$$
and we conclude that
$$
    \lim\limits_{n \to \infty} \thetady = \widetilde\theta^* - \bb m^\top \bv^* \,, 
$$
almost surely as required.

\item We start with the identity of \eqref{eq:invariance}. Let $\widetilde{\bb r}_p = \widetilde{\bb \beta}^{\textrm{\tiny DY}} - \widetilde{\bb B}_p \bv^*$, $\bb z_p$ be the Gaussian surrogate of \eqref{eq:mdypl-gaussian-surrogate} 
from the transformed whitened problem, 
and for $\bb A_p = \bb \Delta_p^{-1/2} \bb L_p^{-\top}$, write 
\begin{equation}
    \label{eq:betady-reverse}
    \sqrt{p} \bb \Delta_p^{-1/2} \left(\betady - \bb B \bv^* \right) = \bb A_p \widetilde{\bb r}_p = \sigma^* \bb A_p  \bb z_p + \bb A_p \left(\widetilde{\bb r}_p - \sigma^* \bb z_p \right) \,. 
\end{equation}
Now, recall that $[\bb \Delta_p^{-1/2}]_{jj}=\tau_j$, and set
$\bb b_p=\sqrt p\,\bb \Delta_p^{-1/2}\bb B$,
$\bb r_p=\sqrt p\,\bb \Delta_p^{-1/2}(\betady-\bb B\bv^*)$,
$\bb x_p=[\bb r_p,\bb b_p]$, and
$\bb y_p=[\sigma^*\bb A_p\bb z_p,\bb b_p]$.
Then \eqref{eq:betady-reverse}, together with manipulations analogous to
\eqref{eq:split-one} and \eqref{eq:pl-expansion}, yields
\begin{equation}
    \label{eq:first-tau-expansion} 
    \begin{aligned}
    \abs{\frac{1}{p} \sum_{j =1}^{p} \psi(\bb x_{p,j}) - \psi(\bb y_{p,j})} 
    &\leq L_{\psi} \frac{1}{p} \sum_{j = 1}^{p} (1 + \vnorm{\bb x_{p,j}}_2 + \vnorm{\bb y_{p,j}}_2) \abs{[\bb A_p(\widetilde{\bb r}_p - \sigma^* \bb z_p)]_j} \\ 
    & \leq C \left(1 + \vnorm{\bb \Delta_p^{-1/2} \betady}_2 + \vnorm{\bb \Delta_p^{-1/2} \bbeta_0}_2 +  \vnorm{\bb \Delta_p^{-1/2} \bbeta_P}_2 + \frac{\vnorm{\bb A_p \bb z_p}_2}{\sqrt{p}} \right) \frac{\vnorm{\bb A_p \left(\widetilde{\bb r}_p - \sigma^* \bb z_p\right)}_2}{\sqrt p} \,,
    \end{aligned}
\end{equation}
for some constant $C$ that depends on $L_\psi, \bv^*, \sigma^*$. Now note that by \eqref{eq:invariance}, 
\begin{equation}
    \label{eq:bb-bounds} 
    \begin{aligned}
        \frac{\vnorm{\sqrt{p} \bb \Delta_p^{-1/2} \betady}_2^2}{p} &= \frac{\vnorm{ \bb \Delta_p^{-1/2} \bb L_p^{-\top} \widetilde{\bb \beta}^{\textrm{\tiny DY}}}_2^2}{p} 
        \leq \mnorm{\bb R_p} \frac{\vnorm{ \widetilde{\bb \beta}^{\textrm{\tiny DY}}}_2^2}{p} 
        = \mathcal{O}(1) \,, 
    \end{aligned}
\end{equation}
almost surely by the assumption on $\bb R_p$, the containment event of Lemma~\ref{lemma:containment} applied to $\widetilde{\bb \beta}^{\textrm{\tiny DY}}$ and the first Borel--Cantelli lemma. 
Similarly, 
$$
     \vnorm{\bb \Delta_p^{-1/2} \bbeta_0}_2 = \mathcal{O}(1), \quad   \vnorm{\bb \Delta_p^{-1/2} \bbeta_P}_2 = \mathcal{O}(1) \,, 
$$
almost surely by the $W2$-convergence assumption of their joint empirical distribution function \citep[][Definition~6.8(iv) \& Theorem~6.9]{villani:2009}. Finally, by standard Gaussian concentration \citep[Theorem~3.1.1]{vershynin:2018} and the assumption in Corollary~\ref{cor:gaussian-covariance-empirical-law}(ii), 
$$
    \frac{\vnorm{\bb A_p \bb z_p}_2}{\sqrt{p}} \leq  \mnorm{\bb R_p}_2 \frac{\vnorm{\bb z_p}_2}{\sqrt p} = \mathcal{O}(1) \,, 
$$
almost surely. Hence, by Lemma~\ref{lemma:mdypl-gaussian-l2-coupling} applied to the whitened problem,
$$
    \frac{\vnorm{\bb A_p\left(\widetilde{\bb r}_p-\sigma^*\bb z_p\right)}_2}{\sqrt p}
    \leq
    \mnorm{\bb R_p}_2^{1/2}
    \frac{\vnorm{\widetilde{\bb r}_p-\sigma^*\bb z_p}_2}{\sqrt p}
    \overset{\mathrm{a.s.}}{\longrightarrow} 0\,.
$$
Consequently, since an almost surely bounded sequence times an almost surely null sequence is almost surely null, the last line in \eqref{eq:first-tau-expansion} goes to zero almost surely. 
It remains to control 
$$
    T_p^{\bb z}
    =
    \frac1p
    \sum_{j=1}^p
    \psi(\bb y_{p,j})\,.
$$
For that, condition on $\mathcal G_n=\sigma(\bb B,\theta_0,\theta_P)$ from \eqref{eq:sigma-fields}. The matrix $\bb A_p$ is deterministic, the $\bb b_{p,j}$ are
$\mathcal G_n$-measurable, and $\mathcal G_n\subseteq\widetilde{\mathcal C}_n$ (since $\bmu_p$ and $\bb L_p$ are deterministic and $\bb L_p$ is invertible, $(\bb B,\theta_0,\theta_P)$ are deterministic functions of the transformed signal and intercepts),
so $\bb z_p$ is independent of $\mathcal G_n$ by
Lemma~\ref{lemma:mdypl-gaussian-l2-coupling}. By
assumption (ii) of Corollary~\ref{cor:gaussian-covariance-empirical-law},
$\mnorm{\bb R_p}_2\leq K$ for a deterministic $K<\infty$.
To get Lipschitz-control over the function $\psi$, for 
$M>0$, introduce the clamp 
$c_M(x)=\sign(x)\min\{\abs{x},M\}$ and define the clamped average 
$$
    T_{p,M}
    =
    T_{p,M}(\bb z_p),
    \quad
    T_{p,M}(\bb x)
    =
    \frac1p
    \sum_{j=1}^p
    \psi
    \left(
        \sigma^* c_M\{[\bb A_p\bb x]_j\},
        \bb b_{p,j}
    \right)\,.
$$
Pseudo-Lipschitzness, Cauchy--Schwarz and
$\mnorm{\bb A_p}_2^2=\mnorm{\bb R_p}_2\leq K$ give
$$
\begin{aligned}
    \abs{T_{p,M}(\bb x)-T_{p,M}(\bb y)}
    &\leq
    \frac{C}{p}
    \sum_{j=1}^p
    \left(
        1+M+\vnorm{\bb b_{p,j}}_2
    \right)
    \abs{[\bb A_p(\bb x-\bb y)]_j}
    \\
    &\leq
    C
    \left(
        \frac{1+M^2+B_p}{p}
    \right)^{1/2}
    \vnorm{\bb x-\bb y}_2\,,
\end{aligned}
$$
where $C$ depends only on $K$, $\sigma^*$ and the pseudo-Lipschitz constant of
$\psi$ and 
$$
    B_p 
    = \frac{1}{p} \sum_{j=1}^{p} \vnorm{\bb b_{p,j}}_2^2 \,.
$$
Thus, on $\{B_p\leq B\}$, the clamped average $T_{p,M}$ is Lipschitz with Lipschitz constant upper bounded by $\sqrt{C_B(1 + M^2) / p}$ for some constant $C_B$ large enough. 
Hence, conditionally on $\mathcal G_n$, Gaussian concentration
\citep[Theorem~5.2.2]{vershynin:2018} yields for every
$\epsilon>0$,
\begin{equation}
    \label{eq:prob-bound-lips} 
        \mathds{1}\{B_p\leq B\}
    \Pr
    \left(
        \abs{
            T_{p,M}
            -
            \expect[T_{p,M}\mid\mathcal G_n]
        }
        >
        \epsilon
        \mid
        \mathcal G_n
    \right)
    \leq
    2
    \exp
    \left\{
        -
        \frac{c_Bp\epsilon^2}{1+M^2}
    \right\}\,.
\end{equation}
For $p \geq 2$, set
$
    M =
    2\sqrt{\log p}
$, 
so that
\begin{equation}
    \label{eq:lips-bound} 
   \begin{aligned}
     \mathds{1}\{B_p \leq B \} \Pr \left( \abs{
            T_{p,M}
            -
            \expect[T_{p,M}\mid\mathcal G_n]
        }
        >
        \epsilon
        \mid
        \mathcal G_n\right) 
        \leq 2
    \exp
    \left\{
        -
        \frac{
            c_Bp\epsilon^2
        }{
            1+4\log p
        }
    \right\} \,. 
   \end{aligned}
\end{equation}
Thus, taking expectations, one gets 
$$
    \Pr \left( \abs{T_{p,M} - \expect[T_{p,M} \mid \mathcal G_n]} > \epsilon \cap B_p \leq B \right) \leq  2
    \exp
    \left\{
        -
        \frac{
            c_Bp\epsilon^2
        }{
            1+4\log p
        }
    \right\} \,, 
$$
which is summable in $n$ and thus by the first Borel--Cantelli lemma, $\{\abs{T_{p,M} - \expect[T_{p,M} \mid \mathcal G_n]} > \epsilon \} \cap \{ B_p \leq B\}$ occurs only finitely often almost surely. 
Separately, for $B$ sufficiently large, by \eqref{eq:bb-bounds}, $\{ B_p > B \}$ finitely often almost surely. Hence, since 
$$
    \{\abs{T_{p,M} - \expect[T_{p,M} \mid \mathcal G_n] }> \epsilon  \} \subseteq [\{\abs{T_{p,M} - \expect[T_{p,M} \mid \mathcal G_n]} > \epsilon\} \cap  \{B_p \leq B\} ] \cup \{ B_p > B \}\,, 
$$
the first Borel--Cantelli lemma, applied along
$\epsilon=1/m$, $m\in\mathbb N$, gives
$$
    T_{p,M}
    -
    \expect[T_{p,M}\mid\mathcal G_n]
    \overset{\mathrm{a.s.}}{\longrightarrow}
    0\,.
$$
Thus, finally, we need to control the truncation bias 
$
     T_p^{\bb z} - T_{p,M}  
$. Pseudo-Lipschitzness gives
\begin{equation}
    \begin{aligned}
        \label{eq:T_p-bound} 
         \abs{T_p^{\bb z}-T_{p,M}} &\leq \sigma^* L_\psi \frac1p \sum_{j=1}^{p} (1 + 2 \sigma^* \abs{[\bb A_p\bb z_p]_j} + 2\vnorm{\bb b_{p,j}}_2) (\abs{[\bb A_p \bb z_p]_j} - M)_+ \,. 
    \end{aligned}
\end{equation}
Now, consider the expectations of the $j$th summand 
$$
    \begin{aligned}
        \expect[(1 + 2 \sigma^* \abs{[\bb A_p\bb z_p]_j} + 2\vnorm{\bb b_{p,j}}_2) (\abs{[\bb A_p \bb z_p]_j} - M)_+ \mid \mathcal G_n] &= (1 + 2 \vnorm{\bb b_{p,j}}_2) \expect[(\abs{[\bb A_p \bb z_p]_j} - M)_+ \mid \mathcal G_n]  \\ 
        &+ 2 \sigma^* \expect[ \abs{[\bb A_p \bb z_p]_j} (\abs{[\bb A_p \bb z_p]_j} - M)_+\mid \mathcal G_n] \,.
    \end{aligned}
$$
Now, since $\diag(\bb R_p)=\bb1$, conditionally on $\mathcal G_n$ each
$[\bb A_p \bb z_p]_j$ is standard normal. Thus,  
\begin{equation}
    \label{eq:first-exp} 
        \begin{aligned}
        \expect[(\abs{[\bb A_p \bb z_p]_j} - M)_+ \mid \mathcal G_n] &= 2 \int_{M}^\infty (x - M) \phi_{\rm N}(x) \mathrm{d}x \\ 
        &= 2 \left(\int_M^\infty x \phi_{\rm N}(x) \mathrm{d}x - M \Phi_{\rm N}(-M) \right) \\ 
        &= 2 (\phi_{\rm N}(M) - M \Phi_{\rm N}(-M)) \\  
        &\leq 2 \phi_{\rm N}(M) \,, 
    \end{aligned}
\end{equation}
where $\phi_{\rm N}, \Phi_{\rm N}$ denote the standard normal pdf and cdf function, respectively, and where the last equality follows from integration by parts. 
Similarly, 
\begin{equation}
    \begin{aligned}
        \expect[ \abs{[\bb A_p \bb z_p]_j} (\abs{[\bb A_p \bb z_p]_j} - M)_+\mid \mathcal G_n] &= 2 \int_{M}^\infty x (x - M) \phi_{\rm N}(x) \mathrm{d}x \\ 
        &= 2 \left(\int_M^\infty x^2 \phi_{\rm N}(x) \mathrm{d}x - M \phi_{\rm N}(M)  \right) \\ 
        &= 2 \left(\left[-x \phi_{\rm N}(x)\right]_M^\infty + \int_M^\infty  \phi_{\rm N}(x)  \mathrm{d}x  - M \phi_{\rm N}(M) \right) \\ 
        &= 2\int_M^\infty \phi_{\rm N}(x) \mathrm{d} x \\ 
        &\leq 2\int_M^\infty \frac{x}{M} \phi_{\rm N}(x) \mathrm{d} x \\ 
        &= 2\frac{\phi_{\rm N}(M)}{M} \\ 
        &\leq 2 \phi_{\rm N}(M) \,, 
    \end{aligned}
\end{equation}
where we used \eqref{eq:first-exp} and integration parts again and the last line holds for $M \geq 1$. Thus, 
\begin{equation}
        \label{eq:T_p-exp-conv}
        \begin{aligned}
        \expect\left[\abs{T_p^{\bb z}-T_{p,M}} \mid \mathcal G_n \right] &\leq \sigma^* L_\psi \frac{1}{p} \sum_{j = 1}^{p} (1 + 2 \vnorm{\bb b_{p,j}}_2 + 2 \sigma^*)2 \phi_{\rm N}(M) \\ 
        &= \exp\{-M^2/2 \} \sqrt{\frac{2}{\pi}} \sigma^* L_\psi \frac{1}{p} \sum_{j = 1}^{p} (1 + 2 \vnorm{\bb b_{p,j}}_2 + 2 \sigma^*)  \\  
        &\leq \exp\{-M^2/2 \} \sqrt{\frac{2}{\pi}} \sigma^* L_\psi (1 + 2 \sigma^* + 2 \sqrt B_p)  \,. 
    \end{aligned}
\end{equation}
Therefore, by Markov's inequality \citep[Proposition~1.2.4]{vershynin:2018}, for any $\epsilon > 0$, 
$$
    \Pr \left(\abs{T_p^{\bb z}-T_{p,M}} > \epsilon \mid \mathcal G_n \right) \leq \frac{1}{\epsilon} \exp\{-M^2/2 \} \sqrt{\frac{2}{\pi}} \sigma^* L_\psi (1 + 2 \sigma^* + 2 \sqrt B_p) \,.  
$$
Hence, on $\{B_p\leq B\}$, for $M = 2 \sqrt{\log p}$, this
conditional expectation bound gives
$$
    \Pr\left(\abs{T_p^{\bb z}-T_{p,M}} > \epsilon \cap \{B_p \leq B \} \right) \leq C_B \epsilon^{-1} p^{-2} \,, 
$$
for some large enough constant $C_B$ depending on $\sigma^*, L_\psi$. This upper bound is summable, and thus by the first Borel--Cantelli lemma, 
$\abs{T_p^{\bb z}-T_{p,M}} > \epsilon \cap \{B_p \leq B \}$ finitely often almost surely. Now since also $\{ B_p \leq B \}$ eventually almost surely for $B$ large enough, 
we conclude that 
$$
    \abs{T_p^{\bb z} - T_{p,M}} \overset{\mathrm{a.s.}}{\longrightarrow} 0 \,. 
$$
Again, since $B_p\leq B$ eventually almost surely, it holds by \eqref{eq:T_p-exp-conv}, 
$$
    \abs{\expect[T_p^{\bb z}\mid\mathcal G_n]
    -
    \expect[T_{p,M}\mid\mathcal G_n]} \leq \expect\left[\abs{T_p^{\bb z} - T_{p,M}} \mid \mathcal G_n\right] \overset{\mathrm{a.s.}}{\longrightarrow} 0 \,, 
$$
eventually almost surely. 
Moreover, $\diag(\bb R_p)=\bb1$ implies
$$
    \bar T_p
    :=
    \expect[T_p^{\bb z}\mid\mathcal G_n]
    =
    \frac1p
    \sum_{j=1}^p
    \expect_G
    \left[
        \psi(\sigma^*G,\bb b_{p,j})
    \right]\,.
$$
Consequently,
$$
    T_p^{\bb z}
    -
    \bar T_p
    \overset{\mathrm{a.s.}}{\longrightarrow}
    0\,.
$$

Finally,
$\expect_G[\psi(\sigma^*G,\bb b)]$ is continuous and has at
most quadratic growth. The $W_2$ convergence
\eqref{eq:gaussian-covariance-standardised-signal-law}
and \citet[Theorem~6.9]{villani:2009}, therefore gives
$$
    \bar T_p
    \overset{\mathrm{a.s.}}{\longrightarrow}
    \int
    \expect_G[\psi(\sigma^*G,\bb b)]
    \,\pi_\Sigma(\mathrm d\bb b)
    =
    \expect
    \left[
        \psi(\sigma^*G,\bar\beta_{0,\Sigma},\bar\beta_{P,\Sigma})
    \right]\,.
$$
Therefore
$$
    T_p^{\bb z}
    \overset{\mathrm{a.s.}}{\longrightarrow}
    \expect
    \left[
        \psi
        \left(
            \sigma^*G,
            \bar\beta_{0,\Sigma},
            \bar\beta_{P,\Sigma}
        \right)
    \right] \,.
$$
Combining the last displays with \eqref{eq:first-tau-expansion} proves the claim.
\end{enumerate}
\end{proof}

\subsection{Proof of Corollary~\ref{cor:gaussian-covariance-adjusted-z}}
\label{sec:proof-gaussian-covariance-adjusted-z}

\begingroup
\renewcommand{\thetheorem}{\arabic{theorem}}
\renewcommand{\thecorollary}{\arabic{theorem}}
\setcounter{theorem}{1}
\begin{corollary}

\label{cor:gaussian-covariance-adjusted-z}
Assume that the centred and whitened triangular array satisfies the
conditions of Section~\ref{sec:setup}. Let
$I\subseteq\{1,\ldots,p\}$ be a deterministic sequence of index sets
with $\abs{I}=k$ fixed, and suppose that
\begin{equation}
\label{eq:covariance-adjusted-z-leverage}
    \mnorms{
        \sqrt p
        \bb C_{I,p}^{1/2}
        \bb J_I^\top
        \bb B
        \widetilde{\bb\Gamma}_p^{+1/2}
    }_2
    =
    \mathcal O_{\mathrm p}(1) \,.
\end{equation}
Let
$\bv^*=[v_1^*,v_2^*]^\top$
be the canonical alignment vector of the centred and whitened problem, and
suppose that $v_1^*\neq0$. For deterministic
$b_I^0\in\Re^k$, under
$
    H_0:
    \bbeta_{0,I}=b_I^0
$,
\begin{equation}
\label{eq:covariance-adjusted-z-limit}
    \bb Z_I^{\mathrm{adj}}(b_I^0)
    =
    \sqrt p  \bb C_{I,p}^{1/2} \frac{\check{\bbeta}_I-b_I^0}{(\sigma^* / v_1^*)}
    \overset{\mathrm d}{\longrightarrow}
    \mathrm N(\bb0_k,\bb I_k)\,.
\end{equation}

\end{corollary}
\endgroup

\begin{proof}
Let
\begin{equation}
\label{eq:covariance-adjusted-z-whitened-contrast}
    \bb D_{I,p}
    =
    \bb C_{I,p}^{1/2}
    \bb J_I^\top
    \bb L_p^{-\top}\,.
\end{equation}
Since
$\bTheta_p=\bb L_p^{-\top}\bb L_p^{-1}$,
\begin{equation}
\begin{aligned}
    \bb D_{I,p}\bb D_{I,p}^\top
    &=
    \bb C_{I,p}^{1/2}
    \bb J_I^\top
    \bTheta_p
    \bb J_I
    \bb C_{I,p}^{1/2}
    =
    \bb I_k\,.
\end{aligned}
\end{equation}
Moreover,
\begin{equation}
\begin{aligned}
    \bb D_{I,p}
    \widetilde{\bb B}_p
    \widetilde{\bb\Gamma}_p^{+1/2}
    =
    \sqrt p\,
    \bb C_{I,p}^{1/2}
    \bb J_I^\top
    \bb B
    \widetilde{\bb\Gamma}_p^{+1/2}\,,
\end{aligned}
\end{equation}
so \eqref{eq:covariance-adjusted-z-leverage} is exactly the leverage
assumption of Lemma~\ref{lemma:finite-contrast-replacement} in the centred
and whitened coordinates. That lemma gives
\begin{equation}
\begin{aligned}
    &\frac1{\sigma^*}
    \bb D_{I,p}
    \left(
        \widetilde{\bbeta}^{\textrm{\tiny DY}}
        -
        \widetilde{\bb B}_p\bv^*
    \right)
    \overset{\mathrm d}{\longrightarrow}
    \mathrm N(\bb0_k,\bb I_k)\,.
\end{aligned}
\end{equation}
Using
$$
    \widetilde{\bbeta}^{\textrm{\tiny DY}}
    =
    \sqrt p\,\bb L_p^\top\betady,
    \quad
    \widetilde{\bb B}_p
    =
    \sqrt p\,\bb L_p^\top\bb B\,,
$$
the left-hand side equals
$$
    \frac{\sqrt p}{\sigma^*}
    \bb C_{I,p}^{1/2}
    \left(
        \betady_I-\bb J_I^\top\bb B\bv^*
    \right)\,.
$$
Under $H_0$,
\begin{equation}
\begin{aligned}
    \bb Z_I^{\mathrm{adj}}(b_I^0)
    &=
    \frac{\sqrt p\,v_1^*}{\sigma^*}
    \bb C_{I,p}^{1/2}
    \left(
        \check{\bbeta}_I-b_I^0
    \right)
    \\
    &=
    \frac{\sqrt p}{\sigma^*}
    \bb C_{I,p}^{1/2}
    \left(
        \betady_I-\bb J_I^\top\bb B\bv^*
    \right)\,.
\end{aligned}
\end{equation}
This proves \eqref{eq:covariance-adjusted-z-limit}.
\end{proof}

\subsection{Proof of Proposition~
\ref{prop:unknowns-o-identification-consistency}}
\label{sec:proof-unknowns-o-identification-consistency}

\begingroup
\renewcommand{\thetheorem}{\arabic{theorem}}
\renewcommand{\theproposition}{\arabic{theorem}}
\setcounter{theorem}{2}
\begin{proposition}
\label{prop:unknowns-o-identification-consistency}
Assume the conditions of Section~\ref{sec:setup}. Then
\begin{equation}
\label{eq:unknowns-o-consistency}
    (\widehat\theta_0,\widehat\gamma^2,\widehat\varphi)
    \overset{\mathrm p}{\longrightarrow}
    (\theta_0^*,\gamma^2,\varphi).
\end{equation}
\end{proposition}
\endgroup

\begin{proof}
Let
$
    \pi_p
    =
    \pi(\theta_0,\gamma_p)$, $
    \zeta_p
    =
    \zeta(\theta_0,\gamma_p)$, $
    c_p
    =
    c(\theta_0,\gamma_p)$, $
    b_p
    =
    \varphi_pc_p
$. 
By \eqref{eq:unknowns-response-population-moments} and
\eqref{eq:cov-y-eta_p}, the conditional population targets of
$\widehat\pi$, $\widehat\zeta^2$, and $\widehat b$ are
$\pi_p$, $\zeta_p^2$, and $b_p$, respectively.
By the assumptions of Section~\ref{sec:setup},
$
    \expect\left[
        (\widehat\pi-\pi_p)^2
        \mid
        \mathcal G_n
    \right]
    =
    \pi_p(1-\pi_p)/n
    \leq
    1/(4n)
$.
The untruncated statistic $\widehat\zeta^2$ is the order-two
$U$-statistic with kernel $pY_iY_j\bx_i^\top\bx_j$
\citep[Section~5.1]{serfling:1980}, whose conditional expectation is
$
    p\vnorm{\expect[\bx_iY_i\mid\mathcal G_n]}_2^2
    =
    \zeta_p^2
$
by \eqref{eq:unknowns-response-population-moments}.
The order-two Hoeffding variance decomposition
\citep[Section~5.2.1, Lemma~A]{serfling:1980}
gives
\begin{equation}
\label{eq:unknowns-response-u-variance-exact}
\begin{aligned}
    \var(\widehat\zeta^2\mid\mathcal G_n)
    =
    \frac{4p^2}{n}
    \expect[\bx_iY_i\mid\mathcal G_n]^\top
    \var(\bx_iY_i\mid\mathcal G_n)
    \expect[\bx_iY_i\mid\mathcal G_n]
    +
    \frac{2p^2}{n(n-1)}
    \trace\left\{
        \var(\bx_iY_i\mid\mathcal G_n)^2
    \right\}\,.
\end{aligned}
\end{equation}
Since $Y_i\in\{0,1\}$,
$$
    \bb0
    \preceq
    \var(\bx_iY_i\mid\mathcal G_n)
    \preceq
    \expect(\bx_i\bx_i^\top\mid\mathcal G_n)
    =
    \frac{1}{p}\bb I_p\,,
$$
and hence
\begin{equation}
\label{eq:unknowns-response-u-variance-bound}
    \var(\widehat\zeta^2\mid\mathcal G_n)
    \leq
    \frac{4\zeta_p^2}{n}
    +
    \frac{2p}{n(n-1)}\,.
\end{equation}
For the aligned moment, the sample covariance $\widehat b$ is the
order-two $U$-statistic with kernel
$
    (\bx_i-\bx_j)^\top\bbeta_P
    (Y_i-Y_j) /2
$.
Its conditional expectation is $b_p$, while, since $Y_i\in\{0,1\}$ and
$(\bx_1-\bx_2)^\top\bbeta_P\mid\mathcal G_n\sim\mathrm N(0,2\delta_p^2)$,
the conditional second moment of the kernel is at most
$$
    \frac{1}{4}
    \expect\left[
        \{(\bx_1-\bx_2)^\top\bbeta_P\}^2
        \mid
        \mathcal G_n
    \right]
    =
    \frac{\delta_p^2}{2}\,.
$$
The same variance decomposition therefore gives
\begin{equation}
\label{eq:unknowns-response-b-variance-bound}
    \expect\left[
        (\widehat b-b_p)^2
        \mid
        \mathcal G_n
    \right]
    \leq
    \frac{3\delta_p^2}{n}\,.
\end{equation}
By the assumptions of Section~\ref{sec:setup},
$\gamma_p^2$ and $\delta_p^2$ are bounded in probability and
$p/n\to\kappa$. Since
$
    \zeta_p=\gamma_pc_p
$
and
$
    0<c_p\leq1/4
$,
the right-hand sides of
\eqref{eq:unknowns-response-u-variance-bound} and
\eqref{eq:unknowns-response-b-variance-bound} converge to zero in
probability. Chebyshev's inequality \citep[Corollary~1.2.5]{vershynin:2018} therefore gives 
$$
    (
        \widehat\pi-\pi_p,
        \widehat\zeta^2-\zeta_p^2,
        \widehat b-b_p
    )
    \overset{\mathrm p}{\longrightarrow}
    \bb0_3\,.
$$
Since $\zeta_p\geq0$ and $\widehat\zeta=\{\max(\widehat\zeta^2,0)\}^{1/2}$,
we have
$
    \abs{\widehat\zeta-\zeta_p}^2
    \leq
    \abs{\widehat\zeta^2-\zeta_p^2}
$,
and hence
\begin{equation}
\label{eq:unknowns-response-moment-consistency}
    (
        \widehat\pi-\pi_p,
        \widehat\zeta-\zeta_p,
        \widehat b-b_p
    )
    \overset{\mathrm p}{\longrightarrow}
    \bb0_3\,.
\end{equation}
Now, to invert the relation of $\widehat{\pi}, \widehat{\zeta}$ to $\widehat{\gamma}, \widehat{\theta}_0$, note that 
$(\theta,g)\mapsto\{\pi(\theta,g),\zeta(\theta,g)\}$ is
continuously differentiable on $\Re^2$, with Jacobian determinant
$$
    \expect\{\rho''(\theta+gZ)\}
    \expect\{Z^2\rho''(\theta+gZ)\}
    -
    \expect\{Z\rho''(\theta+gZ)\}^2
    >0\,,
$$
where the inequality is the weighted Cauchy--Schwarz inequality, and is strict
because $\rho''>0$ and $Z$ is nondegenerate.
The inverse function theorem
\citep[Theorem~9.24]{rudin+etal:1976}
therefore applies at $(\theta_0^*,\gamma)$ and the assumptions of Section~\ref{sec:setup}, continuity of $\pi$ and
$\zeta$, and \eqref{eq:unknowns-response-moment-consistency} give
$$
    (\widehat\pi,\widehat\zeta)
    \overset{\mathrm p}{\longrightarrow}
    \{\pi(\theta_0^*,\gamma),\zeta(\theta_0^*,\gamma)\}\,.
$$
Hence, with probability tending to one,
\eqref{eq:unknowns-response-inversion} has a local solution converging to
$(\theta_0^*,\gamma)$. Since
$
    \zeta(\theta,g)=gc(\theta,g)
$
has the same sign as $g$ and $\widehat\zeta\geq0$, this solution satisfies
$g\geq0$, and the identification argument in the main text makes it the
unique solution in $\Re\times[0,\infty)$. Thus
\begin{equation}
\label{eq:unknowns-response-channel-consistency}
    (
        \widehat\theta_0,
        \widehat\gamma
    )
    \overset{\mathrm p}{\longrightarrow}
    (
        \theta_0^*,
        \gamma
    )\,.
\end{equation}
By continuity and positivity of $c$,
$$
    c(\widehat\theta_0,\widehat\gamma)
    \overset{\mathrm p}{\longrightarrow}
    c(\theta_0^*,\gamma)
    >0\,.
$$
Moreover,
$$
    b_p
    =
    \varphi_pc(\theta_0,\gamma_p)
    \overset{\mathrm p}{\longrightarrow}
    \varphi c(\theta_0^*,\gamma)\,.
$$
It follows from
\eqref{eq:unknowns-response-moment-consistency} that
$$
    \widetilde\varphi
    =
    \frac{
        \widehat b
    }{
        c(\widehat\theta_0,\widehat\gamma)
    }
    \overset{\mathrm p}{\longrightarrow}
    \varphi\,.
$$
Finally,
$
    \abs{\widehat\varphi-\widetilde\varphi}
    =
    (\abs{\widetilde\varphi}-\widehat\gamma\delta_p)_+
    \overset{\mathrm p}{\longrightarrow}
    (\abs{\varphi}-\gamma\delta)_+
    =
    0
$,
because $\widehat\gamma\delta_p\overset{\mathrm p}{\longrightarrow}\gamma\delta$ and
$\abs{\varphi}\leq\gamma\delta$ by positive semidefiniteness of
$\bb\Gamma$. Hence also $\widehat\varphi\overset{\mathrm p}{\longrightarrow}\varphi$, which concludes the proof. 
\end{proof}

\section{Auxiliary Lemmas}
\label{sec:aux-lemmas}

\subsection{Lemmas for Section~\ref{sec:po}}

\subsubsection{MDYPL containment}

\begin{lemma}[Joint MDYPL containment]
	\label{lemma:containment}
	Let $\thetady, \betady$ be the MDYPL estimates of \eqref{eq:og}, $\etady = \thetady \bb 1 + \bX\betady \,,$
	and let $\hat{\blambda}^{\textrm{\tiny DY}}=\bYtil-\bb\rho'(\etady)$.
	Then there exist constants
	$C_{\beta,0},C_{\theta,0},C_{\eta,0},C,c>0$ and $N\in\mathbb N$
	such that the following holds. For any constants
	$C_\beta\geq C_{\beta,0}$, $C_\theta\geq C_{\theta,0}$ and
	$C_\eta\geq C_{\eta,0}$, define
	$$
		\mathcal{M}_n
		=
		\left\{
			\frac{\vnorm{\betady}_2}{\sqrt p} \leq C_\beta,\,
			\abs{\thetady} \leq C_\theta,\,
			\frac{\vnorm{\etady}_2}{\sqrt n} \leq C_\eta,\,
			\frac{\vnorm{\hat{\blambda}^{\textrm{\tiny DY}}}_2}{\sqrt n} \leq 1
		\right\} \,.
	$$
	Then, for all $n>N$,
	$$
		\Pr\left(\lnot\mathcal M_n\right)
		\leq
		C\exp\{-cn\} \,.
	$$
	The constants $C_{\beta,0},C_{\theta,0},C_{\eta,0},C,c$ depend only on
	the fixed quantities $\kappa,\alpha,\bb\Gamma,\theta_0^*,\theta_P^*$ and
	the constants defined in Section~\ref{sec:setup}.
\end{lemma}

\begin{proof}
Let $K$ be the constant of Lemma~\ref{lemma:boundedness}. Take
$C_{\beta,0}=K$, $C_{\theta,0}=K$, and $C_{\eta,0}=3K$. Define
$$
    \widetilde{\bH}
    =
    [\bb1,\bH],
    \quad
    \mathcal E_{\beta,\theta}
    =
    \left\{
        \frac1{\sqrt p}
        \vnorm{
            \begin{bmatrix}
                \sqrt p\,\thetady\\
                \betady
            \end{bmatrix}
        }_2
        \leq K
    \right\},
    \quad
    \mathcal E_{\widetilde{\bH}}
    =
    \left\{
        \frac{\mnorm{\widetilde{\bH}}_2}{\sqrt n}
        \leq3
    \right\}\,.
$$
We have
$$
    \etady
    =
    \frac1{\sqrt p}
    \widetilde{\bH}
    \begin{bmatrix}
        \sqrt p\,\thetady\\
        \betady
    \end{bmatrix}\,.
$$
On
$\mathcal E_{\beta,\theta}\cap\mathcal E_{\widetilde{\bH}}$,
$$
    \frac{\vnorm{\etady}_2}{\sqrt n}
    \leq
    \frac{\mnorm{\widetilde{\bH}}_2}{\sqrt n}
    \frac1{\sqrt p}
    \vnorm{
        \begin{bmatrix}
            \sqrt p\,\thetady\\
            \betady
        \end{bmatrix}
    }_2
    \leq3K\,.
$$
For these baseline constants,
$$
    \lnot\mathcal M_n
    \subseteq
    \lnot\mathcal E_{\beta,\theta}
    \cup
    \lnot\mathcal E_{\widetilde{\bH}}
    \cup
    \left\{
        \vnorm{\hat{\blambda}^{\textrm{\tiny DY}}}_2
        >
        \sqrt n
    \right\}\,.
$$
Lemma~\ref{lemma:boundedness} controls the first event. Since
$p/n\to\kappa\in(0,1)$, Lemma~\ref{lemma:sv_ub}, applied for all
sufficiently large $n$ with $t=\sqrt n-\sqrt p$, controls the second
by $C\exp\{-cn\}$. Since $\Ytil_i\in(0,1)$ and
$\rho'(\etady_i)\in(0,1)$, one has
$\abs{\hat\lambda_i^{\textrm{\tiny DY}}}<1$ for every $i$, so the
third event is empty. The claim follows by a union bound after adjusting
$N,C,c$. Enlarging $C_\beta,C_\theta,C_\eta$ enlarges $\mathcal M_n$,
so the same bound remains valid.
\end{proof}

\begin{lemma}[MDYPL boundedness] \label{lemma:boundedness}
	Let $\thetady, \betady$ be the MDYPL estimates of \eqref{eq:og}.
	Then, in the setup of Section~\ref{sec:setup}, there exists an $N \in \mathbb{N}$ and constants $K, C, c$ such that
	for all $n > N$,
	\begin{equation}
		\label{eq:mdypl-containment}
		\Pr \left(\frac{1}{\sqrt{p}} \vnorm{\begin{bmatrix}\sqrt{p}\thetady \\ \betady \end{bmatrix}}_2 > K  \right) \leq C \exp\{-c n\} \,,
	\end{equation}
	where $K,C,c$ depend only on the constants $\kappa, \alpha, \bb \Gamma, \theta^*_0, \theta_P^*$ and on the constants in
Assumptions~\ref{ass:betas_conv} and~\ref{ass:theta-conv}.
\end{lemma}

\begin{proof}
For
$c=\alpha y+(1-\alpha)\rho'(z)$, if $x\geq0$, then
$c\leq1-(1-\alpha)\rho'(-z)$ and $\rho(x)>x$, so
\[
    \rho(x)-cx
    >
    (1-\alpha)\rho'(-z)x.
\]
If $x\leq0$, then $c\geq(1-\alpha)\rho'(z)$ and
$\rho(x)>0$, so
\[
    \rho(x)-cx
    >
    (1-\alpha)\rho'(z)\abs{x}.
\]
Hence, in both cases,
\[
    \rho(x)
    -
    \{\alpha y+(1-\alpha)\rho'(z)\}x
    >
    (1-\alpha)
    \min\{\rho'(z),\rho'(-z)\}
    \abs{x}.
\]
Applying this bound coordinatewise with
$x=\thetady+\bx_i^\top\betady$,
$z=\theta_P+\bx_i^\top\bbeta_P$, and $y=Y_i$, and using optimality of
$\thetady,\betady$, it holds that
	$$
	\begin{aligned}
		n \log(2) &=-\ell(0, \bb 0_p; \bYtil, \bX) \\
		& \geq -\ell(\thetady, \betady; \bYtil, \bX) \\
		& \geq (1- \alpha)\sum_{i = 1}^{n} \omega(\theta_P + \bx_i^\top \bbeta_P) \abs{\thetady + \bx_i^\top \betady} \,,
	\end{aligned}
	$$
	for $\omega(x) = \min\{\rho'(x), \rho'(-x) \} = \omega(\abs{x})$.
	Define
	$$
	\mathcal{I} = \left\{i \in [n]: \abs{\theta_P + \bx_i^\top \bbeta_P} \leq \frac{3}{\sqrt{\epsilon p }} \vnorm{\begin{bmatrix}\sqrt{p} \theta_P \\ \bbeta_P \end{bmatrix}}_2 \right\} \,.
	$$
	Let
	$$
	\omega_\epsilon = \omega\left( \frac{3}{\sqrt{\epsilon p }} \vnorm{\begin{bmatrix}\sqrt{p} \theta_P \\ \bbeta_P \end{bmatrix}}_2\right) \,.
	$$
	Recall that $\bh_i^\top$ is the $i$th row of $\bH=\sqrt p\bX$, so that $\bx_i=\bh_i/\sqrt p$ with $\bh_i\sim\mathrm N(\bb0_p,\bb I_p)$, and that $\widetilde{\bH}=[\bb1,\bH]$ has rows $\widetilde{\bh}_i^\top=(1,\bh_i^\top)$. Then, we get
	\begin{equation}
	\label{eq:bound_chain}
	\begin{aligned}
		\frac{\log(2)}{1 - \alpha} & \geq \frac{1}{n}\sum_{i = 1}^{n} \omega(\theta_P + \bx_i^\top \bbeta_P) \abs{ \thetady + \bx_i^\top \betady }  \\
		&\geq \frac{1}{n} \sum_{i \in \mathcal{I}} \omega(\theta_P + \bx_i^\top \bbeta_P) \abs{ \thetady + \bx_i^\top \betady }  \\
		&\geq  \frac{1}{n} \sum_{i \in \mathcal I} \omega_\epsilon \abs{ \thetady + \bx_i^\top \betady } \\
		& = \omega_\epsilon\frac{1}{n} \sum_{i \in \mathcal I} \frac{1}{\sqrt{p}} \abs{ [1, \bh_i^\top] \begin{bmatrix} \sqrt{p} \thetady \\ \betady \end{bmatrix}} \\
		&\geq  \frac{\omega_\epsilon}{\sqrt{p}} \vnorm{\begin{bmatrix} \sqrt{p} \thetady \\ \betady \end{bmatrix}}_2   \underset{\substack{\bu \in \Re^{p+1}: \\ \vnorm{\bu}_2  = 1}}{\inf} \, \frac{1}{n} \sum_{i \in \mathcal I} \abs{ \widetilde{\bh}_i^\top \bu } \\
		&=  \frac{\omega_\epsilon}{\sqrt{p}} \vnorm{\begin{bmatrix} \sqrt{p} \thetady \\ \betady \end{bmatrix}}_2   \underset{\substack{\bu \in \Re^{p+1}: \\ \vnorm{\bu}_2  = 1}}{\inf} \, \left( \frac{1}{n}\sum_{i = 1}^n \abs{ \widetilde{\bh}_i^\top \bu }  - \frac{1}{n}\sum_{i \notin \mathcal{I}} \abs{ \widetilde{\bh}_i^\top \bu }\right)\\
		&\geq \frac{\omega_\epsilon}{\sqrt{p}} \vnorm{\begin{bmatrix} \sqrt{p} \thetady \\ \betady \end{bmatrix}}_2   \underset{\substack{\bu \in \Re^{p+1}: \\ \vnorm{\bu}_2  = 1}}{\inf} \, \left( \frac{1}{n}\sum_{i = 1}^n \abs{ \widetilde{\bh}_i^\top \bu }  - \frac{\sqrt{n - \abs{\mathcal I}}}{n} \sqrt{ \sum_{i \notin \mathcal{I}}  \left(\widetilde{\bh}_i^\top \bu \right)^2 }\right)\\
		&\geq  \frac{\omega_\epsilon}{\sqrt{p}} \vnorm{\begin{bmatrix} \sqrt{p} \thetady \\ \betady \end{bmatrix}}_2   \underset{\substack{\bu \in \Re^{p+1}: \\ \vnorm{\bu}_2  = 1}}{\inf} \, \left( \frac{1}{n}\sum_{i = 1}^n \abs{ \widetilde{\bh}_i^\top \bu }  - \frac{\sqrt{n - \abs{\mathcal I}}}{n} \sqrt{ \sum_{i = 1}^n  \left(\widetilde{\bh}_i^\top\bu \right)^2 }\right)\\
		&\geq  \frac{\omega_\epsilon}{\sqrt{p}} \vnorm{\begin{bmatrix} \sqrt{p} \thetady \\ \betady \end{bmatrix}}_2  \left\{ \underset{\substack{\bu \in \Re^{p+1}: \\ \vnorm{\bu}_2  = 1}}{\inf} \, \left( \frac{1}{n}\sum_{i = 1}^n \abs{\widetilde{\bh}_i^\top \bu } \right) -  \sqrt{\frac{n - \abs{\mathcal I}}{n}} \underset{\substack{\bu \in \Re^{p + 1}:\\ \vnorm{\bu}_2 = 1}}{\sup}\, \left(\frac{1}{\sqrt{n}} \sqrt{\sum_{i = 1}^n  \left(\widetilde{\bh}_i^\top \bu \right)^2 }\right)\right\} \\
		&=  \frac{\omega_\epsilon}{\sqrt{p}} \vnorm{\begin{bmatrix} \sqrt{p} \thetady \\ \betady \end{bmatrix}}_2  \left\{ \underset{\substack{\bu \in \Re^{p+1}: \\ \vnorm{\bu}_2  = 1}}{\inf} \, \frac{1}{n} \vnorm{\widetilde{\bH} \bu}_1 -  \sqrt{\frac{n - \abs{\mathcal I}}{n}} \frac{1}{\sqrt{n}}\mnorm{\widetilde{\bH}}_2 \right\} \,,
	\end{aligned}
	\end{equation}
	The third line follows by the definition of $\mathcal{I}$ and since $\omega(x)$ is decreasing in $\abs{x}$ and the seventh line follows
	by the Cauchy--Schwarz inequality \citep[Theorem~1.37(d)]{rudin+etal:1976}.

	Now let $c_{\ell_1}(\kappa)$ be the lower bound of Lemma~\ref{lemma:ell_1_lb}, $\epsilon = (c_{\ell_1}(\kappa) / 12)^2$ and $C_{\theta_P} = \max \{1,\abs{\theta_P^*}\}\,, C_{\beta_P} = \max\{1, \delta\}$, and define the events
	$$
	\begin{aligned}
		\mathcal{E}_1 &= \left\{ \underset{\substack{\bu \in \Re^{p+1}: \\
		\vnorm{\bu}_2 = 1}}{\inf} \,\frac{1}{n} \vnorm{\widetilde{\bH} \bu}_1 \geq \frac{c_{\ell_1}(\kappa)}{2}\right\} \\
		\mathcal{E}_2 &= \left\{ \frac{\mnorm{\widetilde{\bH}}_2}{\sqrt{n}} \leq 3 \right\} \\
		\mathcal{E}_3 &= \left\{ \frac{1}{\sqrt{p}} \vnorm{\begin{bmatrix}\sqrt{p}\theta_P \\ \bbeta_P \end{bmatrix}}_2 \leq 2 \sqrt{C_{\theta_P}^2 + C_{\beta_P}^2} \right\}\,.
	\end{aligned}
	$$
	By Lemma~\ref{lemma:ell_1_lb}, applied with the fixed gap
	$c_{\ell_1}(\kappa)/2$, there exist
	$n_1\in\mathbb N$ and $c_1=c_1(\kappa)>0$ such that, for all
	$n>n_1$,
	$$
		\Pr(\mathcal E_1)
		\geq
		1-4\exp\{-c_1n\}\,.
	$$
	Since
	$
		1-\sqrt{p/n}
		\to
		1-\sqrt\kappa
		>0
	$,
	there exists $n_2\in\mathbb N$ such that, for all $n>n_2$,
	$
		1-\sqrt{p/n}
		\geq
		(1-\sqrt\kappa)/2
	$.
	Taking $t=\sqrt n-\sqrt p$ in
	Lemma~\ref{lemma:sv_ub} therefore gives
	$$
		\Pr(\mathcal E_2)
		\geq
		1
		-
		2\exp\left\{
			-\frac{c_{\mathrm{sv}}}{4}
			(1-\sqrt\kappa)^2n
		\right\}\,.
	$$
	Finally, by the modelling assumptions \ref{ass:betas_conv} and \ref{ass:theta-conv}, there exists an $n_3 \in \mathbb{N}$ such that for all $n > n_3$,
	$$
		\Pr\left(\mathcal{E}_3\right) \geq 1 - C_\Theta \exp \{ -c_\Theta n\} - C_\Gamma \exp \{-c_\Gamma n \} \,.
	$$
	Lemma~\ref{lemma:set_bound} implies that, conditional on $\mathcal{E}_2$,
	$$
	\sqrt{\frac{n - \abs{\mathcal{I}}}{n}} \frac{1}{\sqrt{n}} \mnorm{\widetilde{\bH}}_2 \leq 3\sqrt{\epsilon} \leq c_{\ell_1}(\kappa) / 4 \,.
	$$
	Conditionally on the event $\mathcal{E}_1 \cap \mathcal{E}_2 \cap \mathcal{E}_3$, \eqref{eq:bound_chain} yields that
	$$
	\frac{\log(2)}{1 - \alpha} \geq \omega\left(\frac{72 \sqrt{C_{\theta_P}^2 + C_{\beta_P}^2}}{c_{\ell_1}(\kappa)} \right) \left\{ \frac{c_{\ell_1}(\kappa)}{2} - \frac{c_{\ell_1}(\kappa)}{4} \right\} \frac{1}{\sqrt{p}} \vnorm{\begin{bmatrix} \sqrt{p} \thetady \\ \betady \end{bmatrix}}_2 \,.
	$$
	From this we conclude that
	$$
	\frac{1}{\sqrt{p}} \vnorm{\begin{bmatrix} \sqrt{p} \thetady \\ \betady \end{bmatrix}}_2 \leq \frac{4 \log(2)}{1 - \alpha} \left\{\omega\left(\frac{72 \sqrt{C_{\theta_P}^2 +  C_{\beta_P}^2}}{c_{\ell_1}(\kappa)} \right)c_{\ell_1}(\kappa) \right\}^{-1} \,.
	$$
	Thus, for $n' = \max\{n_1, n_2, n_3 \}$, it follows that for all $n > n'$,
	$$
	\begin{aligned}
		\Pr \left(\mathcal{E}_1 \cap \mathcal{E}_2 \cap \mathcal{E}_3 \right) &\geq 1 - \Pr(\lnot \mathcal{E}_1) - \Pr(\lnot \mathcal{E}_2) - \Pr(\lnot \mathcal{E}_3) \\ &\geq 1 - 4 \exp\{-c_1n\} - 2 \exp\{-c_2 n\} - 2C_3 \exp \{-c_3 n \} \,,
	\end{aligned}
	$$
	for $c_2=c_{\mathrm{sv}}(1-\sqrt\kappa)^2/4$, $C_3 = \max\{ C_\Theta, C_\Gamma\}\,, c_3 = \min\{c_\Theta, c_\Gamma \}$.
	The claim follows by taking $N = \max\{n_1, n_2, n_3 \}$ and $c = \min \{c_1, c_2, c_3\}$, $C = 6 + 2 C_3$.
\end{proof}

\begin{lemma} \label{lemma:ell_1_lb}
	Let
	\begin{equation}
	\label{eq:c_ell1_kappa}
		c_{\ell_1}(\kappa) = \underset{r \in [0,1]}{\inf} \, \mathrm{E}[(W(r) - \tau(r))_{+}] \,,
	\end{equation}
	with $(x)_{+} = \max\{0,x\}$, $W(r) = \abs{r Z - \sqrt{1 - r^2}}$, $Z \sim \mathrm{N}(0,1)$ and where $\tau(r)$ is the unique nonnegative root of
	$$
		\mathrm{E}[\min\{W(r)^2, \tau(r)^2 \}] = r^2 \kappa \,.
	$$
	Then  $c_{\ell_1}(\kappa) > 0$ for any $\kappa \in (0,1)$.

	Furthermore, let $\widetilde{\bH} = \begin{bmatrix} \bb 1 & \bb G \end{bmatrix}$, where $\bb G$ is an $n \times p$ matrix with i.i.d. $\mathrm{N}(0,1)$ entries and $p/n \to \kappa \in (0,1)$ as $n \to \infty$.
	Then, for every fixed $\epsilon>0$, there exist $N\in\mathbb{N}$ and $c=c(\kappa)>0$ such that, for all $n>N$,
	$$
	    \Pr\left(
	    \underset{\substack{\bu \in \Re^{p + 1}:\\ \vnorm{\bu}_2 = 1}}{\inf}
	    \frac{1}{n}\vnorm{\widetilde{\bH}  \bu}_1
	    < c_{\ell_1}(\kappa) - \epsilon
	    \right)
	    \leq
	    4\exp\{-cn\epsilon^2/4\} \,,
	$$
	where $c_{\ell_1}(\kappa)$ is as defined in \eqref{eq:c_ell1_kappa}.
\end{lemma}
\begin{proof}
	Consider the minimisation problem:
	$$
	\underset{\substack{\bu \in \Re^{p + 1}:\\ \vnorm{\bu}_2 = 1}}{\inf}\, \frac{1}{n}\vnorm{\widetilde{\bH}  \bu}_1\,.
	$$
	We shall translate it into a min--max problem and employ the lower-tail comparison of Theorem~\ref{thm:cgmt}(i) to analyse its behaviour.
	\paragraph*{Constraint set}
	For $\bu\in\Re^{p+1}$ with $\vnorm{\bu}_2=1$, write
	$\bu=[\theta,\bw^\top]^\top$, where $\vnorm{\bw}_2\leq1$ and
	$\abs{\theta}=\sqrt{1-\vnorm{\bw}_2^2}$. The minimisation may be
	restricted to vectors with nonpositive first coordinate. Indeed, define
	$$
	\tilde{\bu}
	=
	\begin{cases}
	\bu, & \theta\leq0,\\
	-\bu, & \theta>0\,.
	\end{cases}
	$$
	Then $\vnorm{\tilde{\bu}}_2=1$, the first coordinate of
	$\tilde{\bu}$ is $-\abs{\theta}$, and
	$$
	\vnorm{\widetilde{\bH}\tilde{\bu}}_1
	=
	\vnorm{\widetilde{\bH}\bu}_1 \,.
	$$
	Conversely, any such vector lies in the original constraint set. So it is equivalent to minimise over all unit vectors in $\Re^{p+1}$ or over all vectors
	$$
	\tilde{\bu} = \begin{bmatrix}
		-\sqrt{1 - \vnorm{\bw}_2^2} \\
		\bw
	\end{bmatrix}\,,
	$$
	with $\bw \in \Re^p: \vnorm{\bw}_2 \leq 1$.
	Hence, rewrite the minimisation as
	$$
	\frac{1}{n}\underset{\substack{\bu \in \Re^{p + 1}: \\ \vnorm{\bu}_2 = 1}}{\inf}\, \vnorm{\widetilde{\bH} \bu}_1 = \frac{1}{n}\underset{\substack{\bu \in \Re^{p}: \\ \vnorm{\bu}_2 \leq 1}}{\inf}\, \vnorm{\bb G \bu - \sqrt{1 - \vnorm{\bu}_2^2} \bb 1}_1 \,.
	$$
	\paragraph*{Dual norm}
	Now note that
	$$
	\vnorm{\bb a}_1 = \underset{\vnorm{\bv}_\infty \leq 1}{\sup} \, \bb a^\top \bv \,.
	$$
	Therefore,
	\begin{equation}
		\begin{aligned} \label{eq:PO1}
	\frac{1}{n}\underset{\substack{\bu \in \Re^{p + 1}: \\ \vnorm{\bu}_2 = 1}}{\inf}\, \vnorm{\widetilde{\bH} \bu}_1
	&= \frac{1}{n}\underset{\substack{\bu \in \Re^{p}: \\ \vnorm{\bu}_2 \leq 1}}{\inf}\, \vnorm{\bb G \bu - \sqrt{1 - \vnorm{\bu}_2^2} \bb 1}_1 \\
	&= \frac{1}{n}\underset{\substack{\bu \in \Re^{p }: \\ \vnorm{\bu}_2 \leq 1}}{\inf} \, \underset{\substack{\bv \in \Re^{n}:\\ \vnorm{\bv}_\infty \leq 1}}{\sup} \, \bv^\top  \left( \bb G \bu - \sqrt{1 - \vnorm{\bu}_2^2} \bb 1\right) \\
		&= \frac{1}{n}\underset{\substack{\bu \in \Re^{p }: \\ \vnorm{\bu}_2 \leq 1}}{\inf}\, \underset{\substack{\bv \in \Re^{n}:\\ \vnorm{\bv}_\infty \leq 1}}{\sup} \, \bv^\top   \bb G \bu - \sqrt{1 - \vnorm{\bu}_2^2} \bv^\top \bb 1 \,.
	\end{aligned}
	\end{equation}
	The constraint sets of \eqref{eq:PO1} are compact, and the term
	$$
	\psi(\bu, \bv) = - \sqrt{1 - \vnorm{\bu}_2^2} \bv^\top \bb 1\,,
	$$
	is continuous on the product of these constraint sets and thus the lower-tail part of Theorem~\ref{thm:cgmt} (i) applies. 
	The outer factor $n^{-1}$ is handled by applying the comparison to the 
	unnormalised PO/AO values and rescaling the threshold, as recorded 
	below.
	\paragraph*{PO \& AO}
	Define the PO
	$$
	\Phi(\widetilde{\bH}) =  \frac{1}{n}\underset{\substack{\bu \in \Re^{p }: \\ \vnorm{\bu}_2 \leq 1}}{\inf} \, \underset{\substack{\bv \in \Re^{n}:\\ \vnorm{\bv}_\infty \leq 1}}{\sup} \, \bv^\top \bb G \bu - \sqrt{1 - \vnorm{\bu}_2^2} \bv^\top \bb 1 \,,
	$$
	and for two independent vectors $\bb g \in \Re^n$, $\bb h \in \Re^p$ with i.i.d. $\mathrm{N}(0,1)$,
	define the corresponding AO as
	$$
	\phi(\bb g, \bb h) = \frac{1}{n}\underset{\substack{\bu \in \Re^{p }: \\ \vnorm{\bu}_2 \leq 1}}{\inf} \, \underset{\substack{\bv \in \Re^{n}:\\ \vnorm{\bv}_\infty \leq 1}}{\sup} \, \vnorm{\bu}_2 \bb g^\top \bv + \vnorm{\bv}_2 \bb h^\top \bu  - \sqrt{1 - \vnorm{\bu}_2^2} \bv^\top \bb 1 \,.
	$$
	Let $\overline{\Phi}(\widetilde{\bH}) = n\Phi(\widetilde{\bH})$ and 
	$\overline{\phi}(\bb g,\bb h) = n\phi(\bb g,\bb h)$ denote the 
	corresponding unnormalised values. By the lower-tail part of 
	Theorem~\ref{thm:cgmt} (i), for every $t \in \Re$,
	$$
	\Pr\left(\Phi(\widetilde{\bH}) < t\right)
	=
	\Pr\left(\overline{\Phi}(\widetilde{\bH}) < nt\right)
	\leq
	2\Pr\left(\overline{\phi}(\bb g,\bb h) \leq nt\right)
	=
	2\Pr\left(\phi(\bb g,\bb h) \leq t\right) \,.
	$$
	\paragraph*{Radial decomposition}
	Any $\bu \in \Re^{p}$ can be decomposed into a length $r$ and a unit vector along its direction. Hence the AO is equivalent to
	\begin{equation}
	\label{eq:rad_AO}
	\phi(\bb g, \bb h) =	\frac{1}{n}\underset{\substack{r \in [0,1] \\ \bu \in \Re^{p}: \\ \vnorm{\bu}_2 = 1}}{\inf}\,  \underset{\substack{\bv \in \Re^{n}:\\ \vnorm{\bv}_\infty \leq 1}}{\sup} \, r \bb g^\top \bv + r \vnorm{\bb v}_2 \bb h^\top \bu  - \sqrt{1 - r^2} \bv^\top \bb 1 \,.
	\end{equation}
	For every $\bu\in\Re^p$ with $\vnorm{\bu}_2=1$, the
	Cauchy--Schwarz inequality gives
	$$
	\vnorm{\bv}_2\bb h^\top\bu
	\geq
	-\vnorm{\bv}_2\vnorm{\bb h}_2 \,.
	$$
	Equality is attained at $\bu=-\bb h/\vnorm{\bb h}_2$ when
	$\bb h\neq\bb0_p$, and by every feasible $\bu$ when
	$\bb h=\bb0_p$. Hence, \eqref{eq:rad_AO} becomes
	\begin{equation}
	\label{eq:rad_AO1}
	\phi(\bb g, \bb h) = \frac{1}{n}\underset{r \in [0,1]}{\inf} \, \underset{\substack{\bv \in \Re^{n}:\\ \vnorm{\bv}_\infty \leq 1}}{\sup} \,  \bv^\top \left(r \bb g - \sqrt{1 - r^2} \bb 1 \right) - r \vnorm{\bv}_2 \vnorm{\bb h}_2  \,.
	\end{equation}
	\paragraph*{Reparameterisation}
	Define
	$$
	\bb b(r)=r\bb g-\sqrt{1-r^2}\bb1,
	\quad
	w_i(r)=\abs{b_i(r)},
	\quad
	\lambda(r)=r\frac{\vnorm{\bb h}_2}{\sqrt n},
	\quad
	s(\bb a)=\left(\frac1n\sum_{i=1}^na_i^2\right)^{1/2}\,.
	$$
	For every $\bv$ with $\vnorm{\bv}_\infty\leq1$, setting
	$a_i=\abs{v_i}$ gives
	$$
	\bv^\top\bb b(r)
	\leq
	\sum_{i=1}^na_iw_i(r),
	\quad
	\vnorm{\bv}_2=\vnorm{\bb a}_2 \,. 
	$$
	Conversely, for every $\bb a\in[0,1]^n$, equality is attained by
	choosing $v_i=a_i\sign\{b_i(r)\}$ when $b_i(r)\neq0$ and
	$v_i=a_i$ when $b_i(r)=0$. Therefore \eqref{eq:rad_AO1} is
	equivalent to
	$$
	\phi(\bb g,\bb h)
	=
	\underset{r\in[0,1]}{\inf}
	\underset{\bb a\in[0,1]^n}{\sup} \, 
	\frac1n\sum_{i=1}^na_iw_i(r)-\lambda(r)s(\bb a)
    \,. 
	$$
	For any fixed $r \in [0,1]$ consider the inner maximisation, which is given by
	$$
		\underset{\bb a \in \Re^n}{\sup} \, \frac{1}{n} \sum_{i = 1}^{n} a_i w_i(r) - \lambda(r) s(\bb a), \quad \text{such that } a_i \in [0,1], \, \text{for } i = 1,\ldots,n \,.
	$$
	Equivalently, define the box-constrained minimisation
	\begin{equation}
		\label{eq:inner_min}
		m(r)
		=
		-\underset{\bb a \in [0,1]^n}{\inf}\, 
		-\frac1n\sum_{i=1}^na_iw_i(r)
		+\lambda(r)s(\bb a)
		\,. 
	\end{equation}
	The Lagrangian associated with \eqref{eq:inner_min} is given by
	$$
		\mathcal{L}(\bb a, \bb \mu, \bb \nu) = -\frac{1}{n} \sum_{i = 1}^{n} a_i w_i(r) + \lambda(r) s(\bb a) + \sum_{i = 1}^{n} \mu_i (a_i - 1) - \sum_{i =1}^{n} \nu_i a_i \,.
	$$
	For every fixed $r\in[0,1]$, the objective in \eqref{eq:inner_min}
	is continuous and convex, and the feasible box $[0,1]^n$ is compact.
	Hence an optimal solution exists. Slater's condition holds, for example,
	at every $\bb a\in(0,1)^n$, so the KKT conditions are necessary and
	sufficient for optimality (see, for example,
	\citealt[][Theorem~3.78]{beck:2017}).
	\paragraph*{KKT conditions}
	Thus, consider the KKT conditions.
	\begin{enumerate}
		\item \textbf{Stationarity}: For $\bb a \neq \bb 0_n$, $a_i$ must satisfy:
		$$
		\frac{\partial \mathcal{L}(\bb a, \bb{\mu}, \bb{\nu})}{\partial a_i} = -\frac{w_i(r)}{n} + \lambda(r) \frac{a_i}{n \cdot s(\bb a)} + \mu_i - \nu_i = 0 \,.
		$$
		If $\bb a= \bb 0_n$, then $\bb 0_n$ must be in the subdifferential of $\mathcal{L}(\bb 0_n, \bb \mu, \bb \nu)$, denote it by $\partial \mathcal{L}(\bb 0_n, \bb \mu, \bb \nu)$.
		Now
		$$
			\begin{aligned}
				\partial \mathcal{L}(\bb 0_n, \bb \mu, \bb \nu) &= \left\{\bb s \in \Re^n: \mathcal{L}(\bb a, \bb \mu, \bb \nu) \geq \mathcal{L}(\bb 0_n, \bb \mu, \bb \nu) + \bb s^\top \bb a, \, \forall \bb a \in \Re^n \right\} \\
				&= \left\{\bb s \in \Re^n: \frac{\lambda(r)}{\sqrt{n}} \vnorm{\bb a}_2 \geq \bb a^\top \left( \bb s + \frac{1}{n} \bb w(r) + \bb \nu - \bb \mu \right),  \, \forall \bb a \in \Re^n \right\} \\
				&= \left\{\bb s \in \Re^n: \frac{\lambda(r)}{\sqrt{n}} \geq \vnorm{\bb s + \frac{1}{n} \bb w(r) + \bb \nu - \bb \mu}_2 \right\} \\
				&= \left\{\bb s \in \Re^n: \bb s = \lambda(r) \bu - \frac{1}{n} \bb w(r) - \bb \nu + \bb \mu, \, \bb{u} \in \Re^n: \vnorm{\bu}_2 \leq \frac{1}{\sqrt{n}}   \right\}\,.
			\end{aligned}
		$$
		\item \textbf{Primal feasibility}: $a_i \in [0,1]$ ($i = 1, \ldots, n)$
		\item \textbf{Dual Feasibility}: $\mu_i, \nu_i \geq 0$ ($i = 1, \ldots, n)$
		\item \textbf{Complementary slackness}: $\mu_i (a_i - 1) = 0$, $\nu_i a_i = 0$
	\end{enumerate}
	Consider the following cases for candidate solutions:
	\begin{enumerate}[label=(\roman{*})]
\item 
		$a_i\in(0,1)$: By complementary slackness,
		$\mu_i=\nu_i=0$. Define $\tau_n(r)=\lambda(r)/s(\bb a)$.
		Stationarity gives
		$$
		\tau_n(r)a_i=w_i(r) \,.
		$$
		If $\lambda(r)=0$, stationarity therefore requires $w_i(r)=0$.
		In particular, when $r=0$, $w_i(0)=1$, so no interior coordinate
		can occur. If $\lambda(r)>0$, then $\tau_n(r)>0$ and
		$$
		a_i=\frac{w_i(r)}{\tau_n(r)}<1
		\quad\Longleftrightarrow\quad
		w_i(r)<\tau_n(r) \,.
		$$
		\item $a_i = 1$: By complementary slackness, $\nu_i = 0$. By stationarity: $\tau_n(r) + n \mu_i = w_i(r) \,.$
		If $r = 0$, $\tau_n(r) = 0$, $w_i(r) = 1$ so that $\mu_i = 1/n$ is feasible. If $r \neq 0$, then stationarity requires that $\mu_i = (w_i(r) - \tau_n(r))/n$.
		This is only dual feasible if $w_i(r) \geq \tau_n(r)$.
		\item $a_i = 0$, $\bb a \neq \bb 0_n$: By complementary slackness, $\mu_i = 0$. By stationarity, $- n \nu_i = w_i(r) \,.$
		If $r = 0$, $w_i(r) = 1$ which requires that $\nu_i = -1/n$, which violates dual feasibility. If $r \neq 0$, then stationarity and dual feasibility can only hold
		simultaneously if $w_i(r) = 0$, which requires that $r g_i = \sqrt{1 - r^2} \,,$
		By dual feasibility, this requires that $\nu_i = 0$.

\item 
		$\bb a=\bb0_n$: For every $\bb a\in[0,1]^n$, the
		Cauchy--Schwarz inequality gives
		$$
		\frac1n\sum_{i=1}^na_iw_i(r)-\lambda(r)s(\bb a)
		\leq
		s(\bb a)
		\left\{
		\frac{\vnorm{\bb w(r)}_2}{\sqrt n}
		-
		r\frac{\vnorm{\bb h}_2}{\sqrt n}
		\right\}\,.
		$$
		Hence $\bb a=\bb0_n$ is optimal if
		$\vnorm{\bb w(r)}_2\leq r\vnorm{\bb h}_2$. Conversely, if
		$\vnorm{\bb w(r)}_2>r\vnorm{\bb h}_2$, choose $t>0$ sufficiently
		small that $t\bb w(r)\in[0,1]^n$. At $\bb a=t\bb w(r)$, the
		fixed-radius objective equals
		$$
		\frac{t}{n}\vnorm{\bb w(r)}_2
		\left\{
		\vnorm{\bb w(r)}_2-r\vnorm{\bb h}_2
		\right\}
		>0\,.
		$$
		Therefore $\bb a=\bb0_n$ is optimal if and only if
		$$
		\frac{\vnorm{\bb w(r)}_2}{\sqrt n}
		-
		r\frac{\vnorm{\bb h}_2}{\sqrt n}
		\leq0\,.
		$$
		We now replace this possibility by a deterministic high-probability
		exclusion event.
			Fix $\kappa_+\in(\kappa,1)$. Choose
			$\eta_\kappa\in(0,1-\kappa_+)$, and set $\Delta_\kappa = \sqrt{1-\eta_\kappa}-\sqrt{\kappa_+}>0\,.$
			Define
			$$
			\mathcal E_n^{\mathrm{nz}}
			=
			\left\{
			\underset{r\in[0,1]}{\sup}
			\left|
			\frac{1}{n}\vnorm{\bb w(r)}_2^2-1
			\right|
			\leq \eta_\kappa
			\right\}
			\cap
			\left\{
			\frac{1}{n}\vnorm{\bh}_2^2\leq \kappa_+
			\right\} \,. 
			$$
			For all sufficiently large $n$, there exist constants $C,c>0$, depending
			only on $\kappa$, such that
			$$
			    \Pr\left(\lnot\mathcal E_n^{\mathrm{nz}}\right)
			    \leq
			    C\exp\{-cn\} \,. 
			$$
			Indeed,
			$$
			\begin{aligned}
			\underset{r \in [0,1]}{\sup}
			\left|
			\frac{1}{n}\vnorm{\bb w(r)}_2^2 -1
			\right|
			=
			\underset{r \in [0,1]}{\sup}
			\left|
			r^2\left(\frac{1}{n}\sum_{i=1}^n g_i^2-1\right)
			-
			2r\sqrt{1-r^2}\frac{1}{n}\sum_{i=1}^n g_i
			\right| \leq
			\left|
			\frac{1}{n}\sum_{i=1}^n g_i^2-1
			\right|
			+
			\left|
			\frac{1}{n}\sum_{i=1}^n g_i
			\right| \,,
			\end{aligned}
			$$
			The scalar Gaussian average is controlled by the standard one-dimensional
			Gaussian tail bound, while the two quadratic terms are controlled by
			chi-square concentration. (see, for example,
			\citealt[][Exercise~3.3.7 and Theorem~3.1.1]{vershynin:2018}). Hence the
			right-hand side has an exponentially small upper tail. The second event also
			has an exponentially small complement by the same chi-square concentration,
			since $p/n\to\kappa<\kappa_+$ and $\vnorm{\bh}_2^2\sim\chi_p^2$.
			
			On $\mathcal E_n^{\mathrm{nz}}$,
			$$
			\begin{aligned}
			\underset{r\in[0,1]}{\inf}
			\left\{
			\frac{1}{\sqrt n}\vnorm{\bb w(r)}_2
			-
			r\frac{1}{\sqrt n}\vnorm{\bh}_2
			\right\}
			&\geq
			\sqrt{1-\eta_\kappa}
			-
			\sqrt{\kappa_+}
			=
			\Delta_\kappa
			>
			0 \,.
			\end{aligned}
			$$
			Consequently, on $\mathcal E_n^{\mathrm{nz}}$, the case
			$\bb a=\bb 0_n$ is not optimal for any $r\in[0,1]$.
		\end{enumerate}
		Set
		$$
		A_{\ell_1,n}(r,\tau)
		=
		\frac{1}{n}\sum_{i=1}^{n}\min\{\tau^2,w_i(r)^2\},
		\quad
		\lambda_n^2(r)
		=
		r^2\frac{\vnorm{\bh}_2^2}{n}\,,
		$$
		and
		$$
		A_{\ell_1,n}^{\max}(r)
		=
		\frac{1}{n}\sum_{i=1}^{n}w_i(r)^2\,.
		$$
		For each $r\in[0,1]$, define $\tau_n(r)$ as follows. If
		$\lambda_n^2(r)<A_{\ell_1,n}^{\max}(r)$, let $\tau_n(r)$ be the unique
		nonnegative solution of
		$$
		A_{\ell_1,n}(r,\tau_n(r))
		=
		\lambda_n^2(r)\,.
		$$
		If $\lambda_n^2(r)\geq A_{\ell_1,n}^{\max}(r)$, set
		$\tau_n(r)=0$ as a harmless convention.
		Now define the extended fixed-radius value
		$$
		\widehat\phi_n(r;\bg,\bh)
		=
		\begin{cases}
		\displaystyle
		\frac{1}{n}\sum_{i=1}^n
			\left(w_i(r)-\tau_n(r)\right)_+,
		&
		\lambda_n^2(r)<A_{\ell_1,n}^{\max}(r),
		\\
		0,
		&
		\lambda_n^2(r)\geq A_{\ell_1,n}^{\max}(r)\,.
		\end{cases}
		$$
		We now verify the fixed-radius value. If
		$\lambda_n^2(r)\geq A_{\ell_1,n}^{\max}(r)$, the characterisation
		above shows that $\bb a=\bb0_n$ is optimal and the value is zero.
		Suppose that
		$$
		\lambda_n^2(r)<A_{\ell_1,n}^{\max}(r) \,.
		$$
		The map $\tau\mapsto A_{\ell_1,n}(r,\tau)$ is continuous and
		strictly increasing on
		$[0,\max_{1\leq i\leq n}w_i(r)]$, with endpoint values zero and
		$A_{\ell_1,n}^{\max}(r)$. Hence the defining equation for
		$\tau_n(r)$ has a unique solution on this interval. If
		$\lambda(r)=0$, then $\tau_n(r)=0$ and the fixed-radius objective
		is maximised at $\bb a^*=\bb1_n$, with value
		$$
		\frac1n\sum_{i=1}^nw_i(r)
		=
		\frac1n\sum_{i=1}^n\{w_i(r)-\tau_n(r)\}_+\,.
		$$
		Suppose finally that $\lambda(r)>0$. Then $\tau_n(r)>0$. Define
		$$
		a_i^*
		=
		\min\left\{\frac{w_i(r)}{\tau_n(r)},1\right\},
		\quad
		\mu_i^*
		=
		\frac{\{w_i(r)-\tau_n(r)\}_+}{n},
		\quad
		\nu_i^*=0\,.
		$$
		The threshold equation yields
		$$
		\begin{aligned}
		s(\bb a^*)^2
		=
		\frac1n\sum_{i=1}^n
		\min\left\{\frac{w_i(r)^2}{\tau_n(r)^2},1\right\}
		=
		\frac{A_{\ell_1,n}(r,\tau_n(r))}{\tau_n(r)^2}
		=
		\frac{\lambda(r)^2}{\tau_n(r)^2} 
		\end{aligned}\,.
		$$
		Thus $\lambda(r)/s(\bb a^*)=\tau_n(r)$. The displayed
		$\bb a^*,\bb\mu^*,\bb\nu^*$ satisfy primal feasibility, dual
		feasibility, complementary slackness, and, for every $i$,
		$$
		-\frac{w_i(r)}n
		+
		\lambda(r)\frac{a_i^*}{n s(\bb a^*)}
		+
		\mu_i^*
		-
		\nu_i^*
		=0  \,.
		$$
		They therefore satisfy the KKT conditions, so $\bb a^*$ is optimal.
		Moreover,
		$$
		\begin{aligned}
		\frac1n\sum_{i=1}^na_i^*w_i(r)-\lambda(r)s(\bb a^*) =
		\frac1n\sum_{i=1}^n\{w_i(r)-\tau_n(r)\}_+ \,.
		\end{aligned}
		$$
		Thus $\widehat\phi_n(r;\bg,\bh)$ is exactly the value of the
		fixed-$r$ inner maximisation for every $r\in[0,1]$.
		Therefore
		\begin{equation}
		    \label{eq:ao-r}
		    \phi(\bg,\bh)
		    =
		    \underset{r\in[0,1]}{\inf}
		    \widehat\phi_n(r;\bg,\bh) \,.
		\end{equation}
		
		On $\mathcal E_n^{\mathrm{nz}}$, the second branch in the definition of
		$\widehat\phi_n$ never occurs. Hence, on $\mathcal E_n^{\mathrm{nz}}$,
		$$
		\widehat\phi_n(r;\bg,\bh)
		=
		\frac{1}{n}\sum_{i=1}^{n}
			\left(w_i(r)-\tau_n(r)\right)_+,
		\quad r\in[0,1]\,.
		$$
		For later use, write
		$$
		\phi(r,\tau;\bg,\bh)
		=
			\frac{1}{n}\sum_{i=1}^{n}\left(w_i(r)-\tau\right)_+ \,.
		$$

\paragraph*{Limiting AO}

We now concern ourselves with the limiting value of the AO.
For this, first define the limiting counterparts to the AO
quantities.

Let $Z \sim \mathrm{N}(0,1)$ and define
$$
    W(r) = \abs{r Z - \sqrt{1-r^2}} \,.
$$
For each $r \in [0,1]$, define
$$
    A_{\ell_1}(r,\tau) = \expect\left[\min\{\tau^2, W(r)^2\}\right]\,,
$$
and let $\tau(r)\geq0$ be the unique solution of
\begin{equation}
    \label{eq:tau-pop-eq}
    A_{\ell_1}(r,\tau(r)) =r^2 \kappa\,.
\end{equation}
In particular, $\tau(0)=0$, whereas $\tau(r)>0$ for $r>0$.
Existence, uniqueness, continuity, and uniform boundedness are
established in Lemma~\ref{lemma:A_ell1_props}(iv).

Now let
$$
\varphi(r, \tau) = \expect \left[ (W(r)-\tau)_{+} \right] \,.
$$
Then define the limiting AO integrand at radius $r$ as
$$
    \varphi(r) = \varphi(r, \tau(r)) = \expect\left[(W(r) - \tau(r))_{+}\right] \,.
$$
Then, the limiting AO value is
$$
    \phi^* = \underset{r \in [0,1]}{\inf} \varphi(r) \,.
$$

\paragraph*{Convergence of the AO to the limiting AO}

	We will show that
	\begin{equation}
	\label{eq:ao_conv}
	    \phi(\bg,\bh)\overset{\mathrm p}{\longrightarrow} \phi^* \,.
	\end{equation}
	The threshold representation is used only on the event
	$\mathcal E_n^{\mathrm{nz}}$, whose complement is exponentially small.
	Thus it is enough to prove uniform convergence of the extended fixed-radius
	values on this event. More precisely, for every fixed $\delta>0$,
	$$
	\begin{aligned}
	\Pr\left(
	\left|\phi(\bg,\bh)-\phi^*\right|>\delta
	\right)
	&\leq
	\Pr\left(\lnot\mathcal E_n^{\mathrm{nz}}\right) +
	\Pr\left(
	\underset{r\in[0,1]}{\sup}
	\left|
	\widehat\phi_n(r;\bg,\bh)-\varphi(r)
	\right|>\delta,\ 
	\mathcal E_n^{\mathrm{nz}}
	\right) \,. 
	\end{aligned}
	$$
	The first term is exponentially small, and the second term tends to zero by
	the uniform threshold argument below.
	
	To establish \eqref{eq:ao_conv}, we consider the following steps:
	\begin{enumerate}
	\item For $T > 0$,
	\begin{equation}
	\label{eq:unif-Aellone-phi}
	\begin{aligned}
		\underset{\substack{r \in [0,1] \\ \tau \in [0, T]}}{\sup} \, \abs{\phi(r, \tau; \bg, \bh) - \varphi(r,\tau)} &\overset{\mathrm p}{\longrightarrow} 0 \\
		\underset{\substack{r \in [0,1] \\ \tau \in [0, T]}}{\sup} \, \abs{A_{\ell_1,n}(r, \tau) - A_{\ell_1}(r,\tau)} &\overset{\mathrm p}{\longrightarrow} 0\,.
	\end{aligned}
	\end{equation}
	\item For any $r \in [0,1]$,
		$\underset{r \in [0,1]}{\sup} \, \abs{\tau_n(r)- \tau(r)}  \overset{\mathrm p}{\longrightarrow} 0$
		\item $\underset{r \in [0,1]}{\sup} \, \abs{\widehat\phi_n(r;\bg,\bh) - \varphi(r)} \overset{\mathrm p}{\longrightarrow} 0$
		\item $\abs{\phi(\bg,\bh) - \phi^*} \overset{\mathrm p}{\longrightarrow} 0$
		\item $\expect\left[\phi(\bg,\bh)\right] \to \phi^*$
		\item Tail bound based on Lipschitz Gaussian concentration
	\end{enumerate}
As a first step, we show that $\phi(r,\tau;\bg,\bh)$ and
$A_{\ell_1,n}(r,\tau)$ converge uniformly in probability to
$\varphi(r,\tau)$ and $A_{\ell_1}(r,\tau)$, respectively, on
$[0,1]\times[0,T]$. For every fixed $z\in\Re$, the maps
\begin{equation}
	\label{eq:uniform_maps}
	\begin{aligned}
		(r,\tau)
		&\mapsto
		\left(\abs{rz-\sqrt{1-r^2}}-\tau\right)_+,\\
		(r,\tau)
		&\mapsto
		\min\left\{\tau^2,(rz-\sqrt{1-r^2})^2\right\}
	\end{aligned}
\end{equation}
are continuous on $[0,1]\times[0,T]$. Uniformly over this parameter
set,
$$
\begin{aligned}
\left(\abs{rz-\sqrt{1-r^2}}-\tau\right)_+
&\leq \abs{z}+1,\\
\min\left\{\tau^2,(rz-\sqrt{1-r^2})^2\right\}
&\leq 2z^2+2 \,. 
\end{aligned}
$$
Both envelopes are integrable under the standard normal law. Dominated
convergence gives continuity of $\varphi$ and $A_{\ell_1}$ on the compact
parameter set, and the uniform law of large numbers
(see, for example, \citealt[][Lemma~2.4]{newey+mcfadden:1994}) gives
$$
\begin{aligned}
\underset{\substack{r\in[0,1]\\\tau\in[0,T]}}{\sup}
\abs{\phi(r,\tau;\bg,\bh)-\varphi(r,\tau)}
&\overset{\mathrm p}{\longrightarrow}0,\\
\underset{\substack{r\in[0,1]\\\tau\in[0,T]}}{\sup}
\abs{A_{\ell_1,n}(r,\tau)-A_{\ell_1}(r,\tau)}
&\overset{\mathrm p}{\longrightarrow}0\,.
\end{aligned}
$$
	On $\mathcal E_n^{\mathrm{nz}}$, for each $r\in[0,1]$, the second branch
	in the definition of $\widehat\phi_n$ never occurs, and $\tau_n(r)$ is the
	unique solution of
	$$
	A_{\ell_1,n}(r,\tau_n(r))
	=
	\lambda_n^2(r)\,.
	$$
	Similarly, $\tau(r)$ is the unique solution of
	\eqref{eq:tau-pop-eq} (see Lemma~\ref{lemma:A_ell1_props}, (iv)). We now show that $\tau_n(r)$ converges uniformly over $r$ to $\tau(r)$.
Fix $\epsilon>0$. Define
$$
    m_\epsilon^+
    =
    \underset{r\in[0,1]}{\min}
    \left\{
        A_{\ell_1}(r,\tau(r)+\epsilon)-r^2\kappa
    \right\}>0 \,.
$$
The minimum exists by continuity and compactness, and is strictly positive by
strict monotonicity where applicable and by the extreme value theorem
\citep[Theorem~4.16]{rudin+etal:1976}. Let $R_\epsilon=\{r\in[0,1]:\tau(r)\geq\epsilon\} \,.$
If $R_\epsilon\neq\emptyset$, define
$$
    m_\epsilon^-
    =
    \underset{r\in R_\epsilon}{\min}
    \left\{
        r^2\kappa-A_{\ell_1}(r,\tau(r)-\epsilon)
    \right\}>0 \,.
$$
This minimum exists because $R_\epsilon$ is compact, and it is strictly positive
by strict monotonicity for $r>0$ and the extreme value theorem
\citep[Theorem~4.16]{rudin+etal:1976}. Set
$$
    m_\epsilon
    =
    \begin{cases}
        m_\epsilon^+,
        & R_\epsilon=\emptyset,\\
        \min\{m_\epsilon^+,m_\epsilon^-\},
        & R_\epsilon\neq\emptyset\,.
    \end{cases}
$$
Since $\vnorm{\bh}_2^2\sim\chi_p^2$ and $p/n\to\kappa$, chi-square
concentration gives \citep[Equation~(3.1)]{vershynin:2018}
$$
\frac{\vnorm{\bh}_2^2}{n}
=
\frac pn\frac{\vnorm{\bh}_2^2}{p}
\overset{\mathrm p}{\longrightarrow}
\kappa\,.
$$
Consequently,
$$
\underset{r\in[0,1]}{\sup}
\abs{\lambda_n^2(r)-r^2\kappa}
=
\abs{\frac{\vnorm{\bh}_2^2}{n}-\kappa}
\overset{\mathrm p}{\longrightarrow}0\,.
$$
Choose $T>\sup_{r\in[0,1]}\tau(r)+\epsilon$, which is possible by
Lemma~\ref{lemma:A_ell1_props}~(iv). Combining the preceding convergence
with \eqref{eq:unif-Aellone-phi}, the intersection of the events
$$
\underset{r\in[0,1]}{\sup}
\abs{\lambda_n^2(r)-r^2\kappa}
\leq\frac{m_\epsilon}{4}
$$
and
$$
\underset{\substack{r\in[0,1]\\\tau\in[0,T]}}{\sup}
\abs{A_{\ell_1,n}(r,\tau)-A_{\ell_1}(r,\tau)}
\leq\frac{m_\epsilon}{4}
$$
has probability tending to one. On this intersection, for every
$r\in[0,1]$,
$$
\begin{aligned}
	    A_{\ell_1,n}(r,\tau(r)+\epsilon)
	    &\ge A_{\ell_1}(r,\tau(r)+\epsilon) - \frac{m_\epsilon}{4}
	    \ge r^2\kappa + m_\epsilon^+ - \frac{m_\epsilon}{4}
	    \ge r^2\kappa + \frac{3m_\epsilon}{4}
	    \ge \lambda_n^2(r) + \frac{m_\epsilon}{2} \,.
	\end{aligned}
	$$
	Since $\tau\mapsto A_{\ell_1,n}(r,\tau)$ is nondecreasing and
	$\tau_n(r)$ solves $A_{\ell_1,n}(r,\tau_n(r))=\lambda_n^2(r)$, this gives
	$\tau_n(r)\leq\tau(r)+\epsilon$.
If $r\in R_\epsilon$, then
$$
\begin{aligned}
	    A_{\ell_1,n}(r,\tau(r)-\epsilon)
	    &\leq A_{\ell_1}(r,\tau(r)-\epsilon) + \frac{m_\epsilon}{4}
	    \leq r^2\kappa - m_\epsilon^- + \frac{m_\epsilon}{4}
	    \leq r^2\kappa - \frac{3m_\epsilon}{4}
	    \leq \lambda_n^2(r) - \frac{m_\epsilon}{2} \,.
	\end{aligned}
	$$
	and hence $\tau_n(r)\geq\tau(r)-\epsilon$.
	For $r\notin R_\epsilon$, the lower bound is automatic because
	$\tau_n(r)\geq0>\tau(r)-\epsilon$. Thus, on
	$\mathcal E_n^{\mathrm{nz}}$, simultaneously for all
	$r\in[0,1]$, $\tau_n(r)\in[\tau(r)-\epsilon,\tau(r)+\epsilon]$ on an event
	whose probability tends to one. Since
	$\Pr(\lnot\mathcal E_n^{\mathrm{nz}})\to0$, we obtain
	\begin{equation}
	    \label{eq:unif-tau}
	    \underset{r\in[0,1]}{\sup}
	    \abs{\tau_n(r)-\tau(r)}
	    \overset{\mathrm p}{\longrightarrow}0 \,.
	\end{equation}
	On $\mathcal E_n^{\mathrm{nz}}$, recall that
	$$
	   \widehat\phi_n(r;\bg,\bh)
	   =
	   \phi(r,\tau_n(r);\bg,\bh),
	   \quad
	   \varphi(r)
	   =
	   \varphi(r,\tau(r)) \,.
	$$
	Hence, decompose
	\begin{equation}
	\label{eq:phi_varphi}
	    \begin{aligned}
			\abs{\widehat\phi_n(r;\bg,\bh) - \varphi(r)} &= \abs{\phi(r, \tau_n(r); \bg, \bh) - \phi(r, \tau(r); \bg, \bh) + \phi(r, \tau(r); \bg, \bh) - \varphi(r, \tau(r))} \\
			&\leq \abs{\phi(r, \tau_n(r); \bg, \bh) - \phi(r, \tau(r); \bg, \bh) } + \abs{\phi(r, \tau(r); \bg, \bh) - \varphi(r, \tau(r))} \,.
		\end{aligned}
	\end{equation}
Now, note that
$$
	\phi(r, \tau_n(r); \bg, \bh) - \phi(r, \tau(r); \bg, \bh) = \frac{1}{n} \sum_{i = 1}^n \left\{(w_i(r)- \tau_n(r))_{+} - (w_i(r) - \tau(r))_{+}\right\} \,.
$$
For fixed $w$, the map $\tau \mapsto (w - \tau)_{+}$ is $1$-Lipschitz.
Indeed, without loss of generality, let $\tau > \tau'$ and consider
\begin{enumerate}
	\item $w > \tau > \tau'$: $\abs{(w - \tau)_{+} - (w - \tau')_{+}} = \abs{\tau - \tau'} \,.$
	\item $\tau \geq w > \tau'$: $\abs{(w - \tau)_{+} - (w - \tau')_{+}} = \abs{w - \tau'} = w - \tau' \leq \tau - \tau' = \abs{\tau - \tau'} \,.$
	\item $\tau > \tau' \geq w$: $\abs{(w - \tau)_{+} - (w - \tau')_{+}} = 0 \leq \abs{\tau - \tau'} \,.$
\end{enumerate}
Thus
$$
	\begin{aligned}
		\abs{\phi(r, \tau_n(r); \bg, \bh) - \phi(r, \tau(r); \bg, \bh)} & \leq \frac{1}{n} \sum_{i = 1}^n \abs{(w_i(r)- \tau_n(r))_{+} - (w_i(r) - \tau(r))_{+}}  \\
		&\leq \frac{1}{n} \sum_{i = 1}^n \abs{\tau_n(r) - \tau(r)} \\
		&=  \abs{\tau_n(r) - \tau(r)} \,.
	\end{aligned}
$$
Hence, by \eqref{eq:unif-Aellone-phi}, \eqref{eq:unif-tau}, equation \eqref{eq:phi_varphi} becomes
	$$
		\begin{aligned}
			 \underset{r \in [0,1]}{\sup}\, \abs{\widehat\phi_n(r;\bg,\bh) - \varphi(r)} &\leq  \underset{r \in [0,1]}{\sup}\,\abs{\phi(r, \tau_n(r); \bg, \bh) - \phi(r, \tau(r); \bg, \bh) } +  \underset{r \in [0,1]}{\sup}\,\abs{\phi(r, \tau(r); \bg, \bh) - \varphi(r, \tau(r))} \\
			&\leq \underset{r \in [0,1]}{\sup}\,  \abs{\tau_n(r) - \tau(r)} +  \underset{\substack{r \in [0,1] \\ \tau \in [0, T]}}{\sup} \,\abs{\phi(r, \tau; \bg, \bh) - \varphi(r, \tau)} \\
			&\overset{\mathrm p}{\longrightarrow} 0 \,.
		\end{aligned}
	$$
	The convergence above holds unconditionally because
	$\Pr(\lnot\mathcal E_n^{\mathrm{nz}})\to0$.
	Finally, recall that the AO value is
	$$
	    \phi(\bg,\bh) = \underset{r \in [0,1]}{\inf} \, \widehat\phi_n(r;\bg,\bh) \,,
	$$
and that the limiting AO value is
$$
    \phi^* = \underset{r \in [0,1]}{\inf} \,  \varphi(r) \,.
$$
	Let
	$$
	\delta_n =
	\underset{r \in [0,1]}{\sup}
	\abs{\widehat\phi_n(r;\bg,\bh)-\varphi(r)} \,.
	$$
	Then
	$$
	\varphi(r) = \widehat\phi_n(r;\bg,\bh) + (\varphi(r) - \widehat\phi_n(r;\bg,\bh)) \leq \widehat\phi_n(r;\bg,\bh) + \delta_n \,,
	$$
	and consequently $\phi^* \leq  \phi(\bg,\bh) + \delta_n \,.$
	Similarly, $\phi^* \geq  \phi(\bg,\bh) - \delta_n \,,$
	and in conclusion
	\begin{equation}
		\label{eq:ao_conv_bound}
		\abs{\phi(\bg,\bh)-\phi^*}\leq \delta_n \overset{\mathrm p}{\longrightarrow}0 \,.
	\end{equation}
	We now show that $\expect[\phi(\bg,\bh)]\to\phi^*$. The raw AO values are
	uniformly integrable. Indeed, for every feasible radius, the inner supremum is
	nonnegative because $\bv=\bb 0_n$ is feasible. Hence
	$\phi(\bg,\bh)\geq0$. On the other hand, taking $r=0$ in
	\eqref{eq:rad_AO1} gives the value one, and therefore $0\leq \phi(\bg,\bh)\leq 1\,.$
	Thus $\{\phi(\bg,\bh)\}_{n\geq1}$ is uniformly integrable. Together with
	\eqref{eq:ao_conv_bound}, this implies $\expect[\phi(\bg,\bh)]\to \phi^* \,.$
	Finally, we establish a concentration inequality for $\phi(\bg,\bh)$ around
	its mean. This step uses the original AO value, before the KKT threshold
	representation. Thus no exceptional event is needed here. Namely,
	$$
	\phi(\bg,\bh)
	=
	\frac{1}{n}
	\underset{\substack{\bu \in \Re^p:\\ \vnorm{\bu}_2 \leq 1}}{\inf}
	\,
	\underset{\substack{\bv \in \Re^n:\\ \vnorm{\bv}_\infty \leq 1}}{\sup}
	\left\{
	\vnorm{\bu}_2 \bg^\top\bv
	+
	\vnorm{\bv}_2 \bh^\top\bu
	-
	\sqrt{1-\vnorm{\bu}_2^2}\,\bv^\top\bb 1
	\right\} \,. 
	$$
	Only the first two terms depend on $(\bg,\bh)$. For fixed
	$(\bu,\bv)$, define
	$$
	F_n(\bg,\bh;\bu,\bv)
	=
	\frac{1}{n}
	\left\{
	\vnorm{\bu}_2 \bg^\top\bv
	+
	\vnorm{\bv}_2 \bh^\top\bu
	-
	\sqrt{1-\vnorm{\bu}_2^2}\,\bv^\top\bb 1
	\right\} \,. 
	$$
		For fixed $(\bu,\bv)$, the gradients with respect to $\bg$ and $\bh$
		are $\nabla_{\bg}F_n = \vnorm{\bu}_2\bv/n$, $\nabla_{\bh}F_n = \vnorm{\bv}_2\bu/n\,.$
		Hence $\vnorm{\nabla_{\bg}F_n}_2 \leq 1/\sqrt n$, $\vnorm{\nabla_{\bh}F_n}_2 \leq 1/\sqrt n\,.$
	Thus the joint gradient with respect to $(\bg,\bh)$ is bounded by
	$\sqrt{2/n}$. Hence, for each fixed $(\bu,\bv)$, the map
	$(\bg,\bh)\mapsto F_n(\bg,\bh;\bu,\bv)$ is $\sqrt{2/n}$-Lipschitz for the
	Euclidean norm on $\Re^{n+p}$. Taking a supremum over $\bv$ and an
	infimum over $\bu$ cannot increase the Lipschitz constant. Indeed, if every
	$f_\alpha$ is $L$-Lipschitz, then
	$$
	\abs{
	\underset{\alpha}{\sup}\, f_\alpha(x)
	-
	\underset{\alpha}{\sup}\, f_\alpha(y)
	}
	\leq
	\underset{\alpha}{\sup}\,\abs{f_\alpha(x)-f_\alpha(y)}
	\leq
	L\vnorm{x-y}_2\,,
	$$
	and the same argument applies to infima. Therefore the AO value map $(\bg,\bh)\mapsto\phi(\bg,\bh)$
	is also $\sqrt{2/n}$-Lipschitz.
	
	By Gaussian concentration for Lipschitz functions
	(see, e.g., \citealt[][Theorem~5.2.2]{vershynin:2018}), there exists a
	constant $c>0$ such that, for all $t>0$,
	$$
	    \Pr\left(
	    \abs{\phi(\bg,\bh)-\expect[\phi(\bg,\bh)]}>t
	    \right)
	    \leq
	    2\exp\{-cnt^2\} \,. 
	$$
	Since $\expect[\phi(\bg,\bh)]\to\phi^*$, for every fixed $\epsilon>0$ and
	all sufficiently large $n$,
	$\abs{\expect[\phi(\bg,\bh)]-\phi^*}\leq\epsilon/2$. Hence
	$$
		\Pr\left(\phi(\bg,\bh)\leq\phi^*-\epsilon\right) \leq \Pr\left( \phi(\bg,\bh)-\expect[\phi(\bg,\bh)]\leq-\epsilon/2 \right) \leq 2\exp\{-cn\epsilon^2/4\} \,.
	$$

	\paragraph*{Tail bound for the PO}
	Combining the lower-tail comparison from Theorem~\ref{thm:cgmt} with the AO
	concentration bound gives, for every fixed $\epsilon>0$ and all sufficiently
	large $n$,
	$$
		\Pr\left( \Phi(\widetilde{\bH}) < \phi^*-\epsilon \right) \leq 2\Pr\left( \phi(\bg,\bh) \leq \phi^*-\epsilon \right) \leq 4\exp\{-cn\epsilon^2/4\} \,.
	$$
	Now let $c_{\ell_1}(\kappa)=\phi^*$. Strict positivity is shown in
	Lemma~\ref{lemma:A_ell1_props}, (v). This concludes the proof.
\end{proof}

\begin{lemma}\label{lemma:A_ell1_props}
For $r \in [0,1]$, $\tau \in \Re_{\geq 0}$, let
$$
A_{\ell_1}(r, \tau) = \mathrm{E}\left[\min\{\tau^2, W(r)^2 \} \right] \,,
$$
where
$$
W(r) = \abs{rZ - \sqrt{1 - r^2}} \,,
$$
and $Z \sim \mathrm{N}(0,1)$. Then
\begin{enumerate}[label = (\roman{*})]
	\item For any $r\in(0,1]$, the map
	$\tau\mapsto A_{\ell_1}(r,\tau)$ is strictly increasing on
	$\Re_{\geq0}$. For $r=0$,
	$A_{\ell_1}(0,\tau)=\min\{\tau^2,1\}$, so it is strictly increasing
	on $[0,1]$ and constant on $[1,\infty)$.
	\item $\lim\limits_{\tau \to 0} A_{\ell_1}(r, \tau) = 0$, $\lim\limits_{\tau \to \infty} A_{\ell_1}(r, \tau) = 1$
	\item For any $\tau >0 $, $r \mapsto A_{\ell_1}(r,\tau)$ is strictly decreasing in $r$ on $[0,1]$
	\item For $\kappa \in (0,1)$, $\tau(r)$, defined as the solution to
	\begin{equation}
		\label{eq:taur}
		A_{\ell_1}(r, \tau) = r^2 \kappa \,,
	\end{equation}
	exists, is unique, is continuous in $r$, and bounded uniformly over $r \in [0,1]$.
	\item Let
	$$
	c_{\ell_1}(\kappa) = \underset{r \in [0,1]}{\inf} \, \expect[(W(r)-\tau(r))_{+}] \,,
	$$
	with $(x)_+ = \max \{0, x\}$ and $\tau(r)$ defined by \eqref{eq:taur}. Then for any $\kappa \in (0,1)$,
	$$
	c_{\ell_1}(\kappa) > 0 \,.
	$$
\end{enumerate}
\end{lemma}
\begin{proof}
	We shall prove each claim in turn.
	\begin{enumerate}[label=(\roman{*})]
		\item For $r=0$, $W(0)=1$, so
		$A_{\ell_1}(0,\tau)=\min\{\tau^2,1\}$, which is strictly
		increasing on $[0,1]$ and constant on $[1,\infty)$.
		For $r>0$, fix $0\leq\tau'<\tau$. The difference $\min\{W(r)^2,\tau^2\}-\min\{W(r)^2,(\tau')^2\}$
		is nonnegative and is strictly positive on $\{W(r)>\tau'\}$.
		Since $W(r)$ is the absolute value of a nondegenerate normal
		random variable, $\Pr(W(r)>\tau')>0$. Hence $A_{\ell_1}(r,\tau)>A_{\ell_1}(r,\tau') \,,$
		as required.
		\item  By dominated convergence,
        $$
            \begin{aligned}
				\lim\limits_{\tau \to \infty} \mathrm{E}[\min\{W(r)^2, \tau^2\}] &= \mathrm{E}[\lim\limits_{\tau \to \infty} \min\{W(r)^2, \tau^2\}] = \mathrm{E}[W(r)^2] = 1 \\
				\lim\limits_{\tau \to 0} \mathrm{E}[\min\{W(r)^2, \tau^2\}] &= \mathrm{E}[\lim\limits_{\tau \to 0} \min\{W(r)^2, \tau^2\}] = 0\,.
			\end{aligned}
        $$
\item  		Fix $\tau>0$. For $r\in(0,1)$, set
		$
			s=\sqrt{1-r^2}$, $
			Y=rZ-s
		$. 
		Fix $r_1\in(0,1)$ and a compact interval
		$I\subset(0,1)$ containing $r_1$. For every fixed $Z$, the map
		$$
			u\mapsto
			\min\left\{
				\tau^2,
				\left(uZ-\sqrt{1-u^2}\right)^2
			\right\}
		$$
		is absolutely continuous on $I$, and, wherever its derivative
		exists,
		$$
			\left|
				\frac{\partial}{\partial u}
				\min\left\{
					\tau^2,
					\left(uZ-\sqrt{1-u^2}\right)^2
				\right\}
			\right|
			\leq
			C_I(1+Z^2)
		$$
		for a finite deterministic constant $C_I$. Since
		$\Pr(\abs{Y}=\tau)=0$, differentiation under the expectation gives
		\begin{equation}
		\label{eq:A-ell1-r-derivative}
			\frac{\partial}{\partial r}
			A_{\ell_1}(r,\tau)
			=
			2\expect\left[
				\mathds{1}\{\abs{Y}<\tau\}
				Y
				\left(
					Z+\frac{r}{s}
				\right)
			\right]\,.
		\end{equation}
		Let
		$$
			f_r(y)
			=
			\frac1r
			\phi\left(
				\frac{y+s}{r}
			\right)
		$$
		be the density of $Y$. Since
		$
			f_r'(y)
			=
			-\frac{y+s}{r^2}f_r(y)
		$, 
		integration by parts in \eqref{eq:A-ell1-r-derivative} gives
		\begin{equation}
			\frac{\partial}{\partial r}
			A_{\ell_1}(r,\tau)
			=
			\frac2r
			\int_{-\tau}^{\tau}
			y(y+s)f_r(y)\mathrm dy
			+
			\frac{2r}{s}
			\int_{-\tau}^{\tau}
			yf_r(y) \mathrm dy
			=
			-2r
			\left[
				\tau\{f_r(\tau)+f_r(-\tau)\}
				+
				\frac{r^2}{s}
				\{f_r(\tau)-f_r(-\tau)\}
			\right] \,. 
		\end{equation}
		Writing
		$
			k=\frac{\tau s}{r^2}>0 
		$ 
		we have
		$$
			\frac{
				f_r(-\tau)-f_r(\tau)
			}{
				f_r(-\tau)+f_r(\tau)
			}
			=
			\tanh(k)
			<
			k
			=
			\frac{\tau s}{r^2} \,, 
		$$
		where $\tanh(k)<k$ follows because
		$k-\tanh(k)$ vanishes at zero and has derivative
		$\tanh^2(k)>0$ for $k>0$. Hence
		$$
			\tau\{f_r(\tau)+f_r(-\tau)\}
			-
			\frac{r^2}{s}
			\{f_r(-\tau)-f_r(\tau)\}
			>
			0\,, 
		$$
		and therefore
		$$
			\frac{\partial}{\partial r}
			A_{\ell_1}(r,\tau)
			<
			0,
			\qquad
			r\in(0,1) \,. 
		$$

		Finally, $(r,\tau)\mapsto A_{\ell_1}(r,\tau)$ is jointly continuous
		on $[0,1]\times\Re_{\geq0}$ by dominated convergence, since
		$$
			\min\{\tau^2,W(r)^2\}
			\leq
			W(r)^2
			\leq
			Z^2+1 \,. 
		$$
		The strict decrease on $(0,1)$ therefore extends to $[0,1]$ by
		continuity at the endpoints.
        \item For $r=0$, the equation is $A_{\ell_1}(0,\tau)=0$, whose
	   unique solution is $\tau(0)=0$. For $r>0$, continuity of
	   $\tau\mapsto A_{\ell_1}(r,\tau)$ follows from the joint continuity
	   established in part~(iii). Hence the map
	   $\tau \mapsto A_{\ell_1}(r, \tau) - r^2 \kappa$ is continuous.
	   It is negative at $\tau=0$ and has limit $1-r^2\kappa\geq 1-\kappa>0$
	   as $\tau\to\infty$. By the intermediate value theorem \citep[Theorem~4.23]{rudin+etal:1976}, there must be
	   a value $\tau \geq0$ such that
	   $A_{\ell_1}(r, \tau) - r^2\kappa = 0$. Uniqueness follows
	   by strict monotonicity of $A_{\ell_1}$ in $\tau$ for $r>0$.
	   To establish uniform boundedness of $\tau(r)$, for any $r \in [0,1]$, let $\tau(r)$ be the solution to $A_{\ell_1}(r, \tau) - r^2 \kappa = 0 \,.$
	   Then, since $r \mapsto A_{\ell_1}(r, \tau)$ is decreasing in $r$,
	   $$
		A_{\ell_1}(1, \tau(r)) - \kappa \leq A_{\ell_1}(1, \tau(r)) - r^2 \kappa  \leq  A_{\ell_1}(r, \tau(r))- r^2 \kappa = 0 \,.
	   $$
	   Now since $\tau \mapsto A_{\ell_1}(1,\tau)$ is strictly increasing in $\tau$, it follows that $\tau(1) \geq \tau(r)$ for any $r \in [0,1]$.

	   Finally to show that $\tau(r)$ is continuous, consider any convergent sequence $\{r_n\}_{n \in \mathbb{N}}, r_n \in [0,1]$ with limit $r^* \in [0,1]$.
	   Concurrently, let $\tau_n$ be the unique solution to $A_{\ell_1}(r_n, \tau) - r_n^2\kappa = 0$. Since $\tau_{n} \in [0, \tau(1)]$, by Bolzano--Weierstrass \citep[Theorem~3.6(b)]{rudin+etal:1976}, there
	   must exist a convergent subsequence, call it $\tau_{n_k}$ with limit say $\bar{\tau}$.
	   Then by joint continuity of $A_{\ell_1}(r,\tau)$, it follows that
	   $$
		0 = \lim\limits_{k \to \infty} A_{\ell_1}(r_{n_k}, \tau_{n_k}) - r^2_{n_k} \kappa = A_{\ell_1}(r^*, \bar \tau) - (r^*)^2\kappa \,.
	   $$
	   By uniqueness of $\tau(r)$, $\bar{\tau} = \tau(r^*)$ and thus every convergent subsequence must converge to $\tau(r^*)$. But then also
	   $\tau_{n} \to \tau(r^*)$. If this were not the case, then $\tau_n$ would have to have a non-convergent subsequence for which $\abs{\tau_{n_k} - \bar \tau} > \epsilon$.
	   For $\tau_{n_k}$, which is still bounded, another convergent subsequence must exist and thus converge to $\bar \tau$, which contradicts $\abs{\tau_{n_k} - \bar \tau} > \epsilon$.
	   \item For existence, note that the map $(w,\tau) \mapsto (w - \tau)_{+}$ is continuous and that
	   $(W(r) - \tau(r))_{+} \leq W(r) \leq \abs{Z} + 1$ is integrable. Hence, by dominated convergence $(w, \tau) \mapsto \expect[(w - \tau)_{+}]$ is
	   continuous. Further $r \mapsto \tau(r)$ is continuous. Thus,
	   $r \mapsto \expect[(W(r)-\tau(r))_{+}]$ is continuous and by the extreme value theorem, there exists a $r^* \in [0,1]$ such that
	   $$
		\expect[(W(r^*) - \tau(r^*))_{+}] = \underset{r \in [0,1]}{\inf} \, \expect[(W(r) - \tau(r))_{+}] \,.
 	   $$
	   Now assume that for this $r^*$, $c_{\ell_1}(\kappa) = 0$. Then it must hold that $\tau(r^*) \geq W(r^*)$ $Z$-almost everywhere.
	   In that case $A_{\ell_1}(r^*, \tau(r^*)) = \expect[W(r^*)^2] = 1$. But then $A_{\ell_1}(r^*, \tau(r^*)) - (r^*)^2 \kappa > 0$, so that
	   $\tau(r^*)$ is no solution to \eqref{eq:taur}, a contradiction. Thus, $c_{\ell_1}(\kappa) >0$ for all $\kappa \in (0,1)$.
	\end{enumerate}
	This concludes the proof.
\end{proof}

\subsubsection{Product-domain constraints}

\begin{lemma}[Orthogonal coordinate decomposition]
\label{lemma:orth_decomp}
Let $\bB\in\Re^{p\times 2}$, $ \rank(\bb B) = s$ and
$$
\bb P = \bB (\bB^\top \bB)^{+} \bB^\top\,, 
$$
and $\bb E \in\Re^{p\times (p - s)}$ has orthonormal columns and satisfies
$$
\bb E \bb E^\top = \bb P^{\perp} = \bb{I}_p - \bb P \,. 
$$
Then, for every $\bbeta\in\Re^p$ there exist unique $\bv \in \operatorname{range}(\bB^\top)$,
$\bw \in \Re^{p-s}: \bb E \bw \in \operatorname{ker}(\bb B^\top)$ such that
$$
\bbeta = \bB\bv + \bb E\bw \,.
$$
\end{lemma}

\begin{proof}
Since $\bb P$ is the orthogonal projector onto
$\range(\bb B)$ and
$\bb E\bb E^\top=\bb P^\perp$ with orthonormal columns,
$\range(\bb E)=\ker(\bb B^\top)$ and hence
$\bb B^\top\bb E=\bb0$. For any $\bbeta\in\Re^p$, take
$\bv=\bb B^+\bbeta\in\range(\bb B^\top)$ and
$\bw=\bb E^\top\bbeta$. Then
$\bb B\bv=\bb P\bbeta$ and
$\bb E\bw=\bb P^\perp\bbeta$, so
$\bbeta=\bb B\bv+\bb E\bw$. If two such representations exist,
orthogonality gives $\bb E(\bw_1-\bw_2)=0$ and
$\bb B(\bv_1-\bv_2)=0$. Since $\bb E$ has full column rank and
$\range(\bb B^\top)\cap\ker(\bb B)=\{\bb0\}$, both differences
vanish.
\end{proof}

\begin{lemma}[Product-domain relaxation preserves the compact MDYPL saddle]
\label{lemma:product_domain_relaxation}
Let $\bb B\in\Re^{p\times 2}$, set $s=\rank(\bb B)$, and let
$\bE\in\Re^{p\times(p-s)}$ satisfy
$$
	\bE^\top\bE=\bb I_{p-s},
	\quad
	\range(\bE)=\range(\bb B)^\perp \,.
$$
Let $\bH_1=\bH\bb B$ and $\bH_2=\bH\bE$.
For constants $C_\sigma \geq C_\beta > 0$, consider the sets
$$
V_p=\{\bv\in\operatorname{range}(\bb B^\top):\vnorm{\bb B \bv}_2\le C_\beta\sqrt p\},
\quad
W_p=\{\bw\in\Re^{p-s}:\vnorm{\bw}_2\le C_\sigma\sqrt p\} \,,
$$
and
$$
    \mathcal K^p_{C_\beta}
    =
    \left\{
        (\bv,\bw):
        \bv\in\operatorname{range}(\bB^\top),\,
        \bw\in\Re^{p-s},\,
        \vnorm{\bB\bv+\bE\bw}_2\le C_\beta\sqrt p
    \right\} \,.
$$
Let
$$
    d(\bv,\bw,\theta,\bb{\eta})
    =
    \bb{\eta}-\theta\bb{1}
    -
    \frac{1}{\sqrt p}\bH_1\bv
    -
    \frac{1}{\sqrt p}\bH_2\bw, \quad L(\bv, \bw, \theta, \bb \eta, \blambda) = \frac{1}{n}
\left\{
    \bb{1}^\top\bb \rho(\bb{\eta})-\bYtil^\top\bb{\eta}
    +
    \bb{\lambda}^\top d(\bv,\bw,\theta,\bb{\eta})
\right\} \,, 
$$
and 
define the coupled-domain and product-domain saddle problems, respectively, as
$$
\Psi_{\mathcal K}
=
\underset{\substack{(\bv,\bw)\in\mathcal K^p_{C_\beta}\\
\abs{\theta}\leq C_\theta \\
\bb{\eta}\in\mathcal B^n_{C_\eta}}}{\min}
\underset{\bb{\lambda}\in\mathcal B^n_1}{\max}
\, L(\bv, \bw, \theta, \bb \eta, \blambda) , \quad 
\Psi_{V \times W}
=
\underset{\substack{\bv\in V_p,\ \bw\in W_p\\
\abs{\theta}\le C_\theta\\ 
\bb{\eta}\in\mathcal B^n_{C_\eta}}}{\min}
\underset{\bb{\lambda}\in\mathcal B^n_1}{\max} \,
L(\bv, \bw, \theta, \bb \eta, \blambda)  \,.
$$
On the event $\mathcal{E}_n \cap \mathcal M_n$ from \eqref{event:E_n} and \eqref{event:M_n} respectively, 
the two values are equal, i.e
$
    \Psi_{\mathcal K}=\Psi_{V \times W}
$. 
Moreover, every optimiser of the product-domain problem satisfies
$
    d(\bv,\bw,\theta,\bb{\eta})=\bb{0}
$,
and consequently, every product-domain optimiser has the same
$
    (\theta,\bbeta,\bb{\eta})
$
coordinates as the original MDYPL optimiser of \eqref{eq:og}, where $\bbeta=\bB\bv+\bE\bw$.
\end{lemma}

\begin{proof}
Work on $\mathcal M_n\cap\mathcal E_n$. On $\mathcal{E}_n$, the intercept augmented design matrix $\widetilde{\bb H} = [\bb1,\bH]$ has full column rank and
the minimiser of \eqref{eq:og} exists and is unique \citep[Theorem~1]{rigon+aliverti:2023}.
Let $\etady = \thetady\bb{1} + \bH\betady/\sqrt p$
and write the orthogonal coordinate representation of the original MDYPL solution as $\betady = \bB\hat\bv+\bE\hat\bw\,.$
On $\mathcal M_n\cap\mathcal E_n$, $\vnorm{\betady}_2\le C_\beta\sqrt p, \quad \abs{\thetady}\le C_\theta, \quad \etady\in\mathcal B^n_{C_\eta}\,.$
Hence $(\hat\bv,\hat\bw)\in\mathcal K^p_{C_\beta} \,.$ Since
$\bB^\top\bE=\bb0$, $\bE^\top\bE=\bb I_{p-s}$, and
$C_\sigma\geq C_\beta$, this also gives
$(\hat\bv,\hat\bw)\in V_p\times W_p \,.$
Thus the original MDYPL solution is feasible for both the coupled-domain
and product-domain saddle problems.

Let $\hat{\bb{\lambda}}^{\textrm{\tiny DY}} = \bYtil-\rho'(\etady) \,.$
Since $\Ytil_i\in(0,1)$ and $\rho'((\etady)_i)\in(0,1)$ for every $i$, $\abs{\hat{\bb \lambda}_i^{\textrm{\tiny DY}}}<1, \quad i=1,\ldots,n \,.$
Therefore $\vnorm{\hat{\bb{\lambda}}^{\textrm{\tiny DY}}}_2<\sqrt n \,,$
so $\hat{\bb{\lambda}}^{\textrm{\tiny DY}}$ lies in the relative interior of $\mathcal B^n_1$.

Recall the empirical loss $G_n$ from \eqref{eq:ao_scalarisation_psi_def}, so that the Lagrangian of the statement is
$$
    L(\bv,\bw,\theta,\bb{\eta},\bb{\lambda})
    =
    G_n(\bb{\eta})
    +
    \frac{1}{n}\bb{\lambda}^\top d(\bv,\bw,\theta,\bb{\eta}) \,.
$$
Because $\etady$ is the fitted predictor of the original MDYPL solution, $d(\hat\bv,\hat\bw,\thetady,\etady)=\bb{0} \,.$
On $\mathcal E_n$, $(\thetady,\betady)$ is the genuine MDYPL
optimiser. The KKT equations for \eqref{eq:og-constr} identify
$\hat{\blambda}^{\textrm{\tiny DY}}
=\bYtil-\bb\rho'(\etady)$ as its equality-constraint multiplier and give
$$
\bb1^\top\hat{\blambda}^{\textrm{\tiny DY}}=0,
\quad
\bH^\top\hat{\blambda}^{\textrm{\tiny DY}}=\bb0_p\,.
$$
Since $\bH_1=\bH\bB$ and $\bH_2=\bH\bE$, it follows that $\bH_1^\top\hat{\bb{\lambda}}^{\textrm{\tiny DY}}=\bb{0}, \quad \bH_2^\top\hat{\bb{\lambda}}^{\textrm{\tiny DY}}=\bb{0}\,.$
Consequently, for every $(\bv,\bw,\theta,\bb{\eta})$,
$$
    L(\bv,\bw,\theta,\bb{\eta},\hat{\bb{\lambda}}^{\textrm{\tiny DY}})
    =
    \frac{1}{n}
    \left\{
        \bb{1}^\top\rho(\bb{\eta})
        -
        \bYtil^\top\bb{\eta}
        +
        (\hat{\bb{\lambda}}^{\textrm{\tiny DY}})^\top\bb{\eta}
    \right\} \,.
$$
Using $\hat{\bb{\lambda}}^{\textrm{\tiny DY}}=\bYtil-\rho'(\etady)$, this becomes
$$
    L(\bv,\bw,\theta,\bb{\eta},\hat{\bb{\lambda}}^{\textrm{\tiny DY}})
    =
    \frac{1}{n}
    \sum_{i=1}^n
    \left\{
        \rho(\bb{\eta}_i)
        -
        \rho'((\etady)_i)\bb{\eta}_i
    \right\} \,.
$$
The right-hand side is convex in $\bb{\eta}$ and is minimised at $\bb{\eta}=\etady$. Hence
$$
    L(\bv,\bw,\theta,\bb{\eta},\hat{\bb{\lambda}}^{\textrm{\tiny DY}})
    \ge
    L(\hat\bv,\hat\bw,\thetady,\etady,\hat{\bb{\lambda}}^{\textrm{\tiny DY}})
    =
    G_n(\etady)\,,
$$
for every $(\bv,\bw,\theta,\bb{\eta})$.

Now fix any point in the product domain. Since $\mathcal B^n_1 = \{\bb{\lambda}\in\Re^n:\vnorm{\bb{\lambda}}_2\le\sqrt n\} \,,$
we have
$$
	\underset{\bb{\lambda}\in\mathcal B^n_1}{\max} L(\bv,\bw,\theta,\bb{\eta},\bb{\lambda}) = G_n(\bb{\eta}) + \frac{1}{n} \underset{\bb{\lambda}\in\mathcal B^n_1}{\max} \bb{\lambda}^\top d(\bv,\bw,\theta,\bb{\eta}) = G_n(\bb{\eta}) + \frac{1}{\sqrt n} \vnorm{d(\bv,\bw,\theta,\bb{\eta})}_2\,.
$$
In particular,
$$
    \underset{\bb{\lambda}\in\mathcal B^n_1}{\max}
    L(\bv,\bw,\theta,\bb{\eta},\bb{\lambda})
    \ge
    L(\bv,\bw,\theta,\bb{\eta},\hat{\bb{\lambda}}^{\textrm{\tiny DY}})
    \ge
    G_n(\etady) \,.
$$
At the original MDYPL solution, because the equality constraint is satisfied,
$$
    \underset{\bb{\lambda}\in\mathcal B^n_1}{\max}
    L(\hat\bv,\hat\bw,\thetady,\etady,\bb{\lambda})
    =
    G_n(\etady) \,.
$$
Therefore the product-domain value is $\Psi_{V \times W}=G_n(\etady)\,.$
The same argument applies to the coupled-domain problem, since the original MDYPL solution belongs to $\mathcal K^p_{C_\beta}$, and hence $\Psi_{\mathcal K}=G_n(\etady)=\Psi_{V \times W}\,.$

It remains to show that the product-domain problem introduces no infeasible optimiser. Let $(\tilde\bv,\tilde\bw,\tilde\theta,\tilde{\bb{\eta}})$ be a product-domain optimiser. Then
$$
    \underset{\bb{\lambda}\in\mathcal B^n_1}{\max}
    L(\tilde\bv,\tilde\bw,\tilde\theta,\tilde{\bb{\eta}},\bb{\lambda})
    =
    G_n(\etady) \,.
$$
Towards a contradiction, assume that $d(\tilde\bv,\tilde\bw,\tilde\theta,\tilde{\bb{\eta}})\ne \bb{0} \,.$
Then
$$
\begin{aligned}
    \underset{\bb{\lambda}\in\mathcal B^n_1}{\max}
    L(\tilde\bv,\tilde\bw,\tilde\theta,\tilde{\bb{\eta}},\bb{\lambda})
    -
    L(\tilde\bv,\tilde\bw,\tilde\theta,\tilde{\bb{\eta}},\hat{\bb{\lambda}}^{\textrm{\tiny DY}}) &=
    \frac{1}{\sqrt n}
    \vnorm{d(\tilde\bv,\tilde\bw,\tilde\theta,\tilde{\bb{\eta}})}_2
    -
    \frac{1}{n}
    (\hat{\bb{\lambda}}^{\textrm{\tiny DY}})^\top
    d(\tilde\bv,\tilde\bw,\tilde\theta,\tilde{\bb{\eta}}) \\
    &\geq 
    \left(
        \frac{1}{\sqrt n}
        -
        \frac{\vnorm{\hat{\bb{\lambda}}^{\textrm{\tiny DY}}}_2}{n}
    \right)
    \vnorm{d(\tilde\bv,\tilde\bw,\tilde\theta,\tilde{\bb{\eta}})}_2
    >
    0 \,.
\end{aligned}
$$
The final inequality uses $\vnorm{\hat{\bb{\lambda}}^{\textrm{\tiny DY}}}_2<\sqrt n$. Since
$
    L(\tilde\bv,\tilde\bw,\tilde\theta,\tilde{\bb{\eta}},\hat{\bb{\lambda}}^{\textrm{\tiny DY}})
    \ge
    G_n(\etady) 
$, 
we obtain
$$
    \underset{\bb{\lambda}\in\mathcal B^n_1}{\max}
    L(\tilde\bv,\tilde\bw,\tilde\theta,\tilde{\bb{\eta}},\bb{\lambda})
    >
    G_n(\etady) \,,
$$
contradicting optimality of $(\tilde\bv,\tilde\bw,\tilde\theta,\tilde{\bb{\eta}})$. Therefore $d(\tilde\bv,\tilde\bw,\tilde\theta,\tilde{\bb{\eta}})=\bb{0} \,.$
Thus every product-domain optimiser is feasible for the equality-constrained MDYPL problem and attains the original MDYPL objective value. Since, on $\mathcal M_n\cap\mathcal E_n$, the original MDYPL optimiser is unique,
$\tilde\theta=\thetady, \quad \tilde{\bb{\eta}}=\etady, \quad \bB\tilde\bv+\bE\tilde\bw=\betady\,.$
This proves the claim.
\end{proof}

\begin{lemma}[Uniqueness of the compact PO]
\label{lemma:compact_PO_uniqueness}
On $\mathcal E_n$, the compact problem in
\eqref{eq:lagrangian4} has a saddle point, and all its saddle points have
the same minimising variables, denoted by
$(\hat\bv,\hat\bw,\hat\theta,\hat{\bb\eta})$. Consequently,
$\hat\bw$ is the unique minimiser of the profiled PO in the final line of
\eqref{eq:lagrangian5}, and the inner minimising variables at every saddle
point of the corresponding fixed-$\hat\bw$ problem are
$(\hat\bv,\hat\theta,\hat{\bb\eta})$.
\end{lemma}

\begin{proof}
Maximising over $\blambda\in\mathcal B_1^n$ reduces
\eqref{eq:lagrangian4} to minimising, over its compact convex primal domain,
the continuous convex function
$$
    \frac1n
    \left\{
        \bb1^\top\bb\rho(\bb\eta)-\bYtil^\top\bb\eta
    \right\}
    +
    \frac1{\sqrt n}
    \vnorm{
        \bb\eta-\theta\bb1
        -\frac1{\sqrt p}\bH_1\bv
        -\frac1{\sqrt p}\bH_2\bw
    }_2\,.
$$
Hence a minimiser exists. Let
$(\bv,\bw,\theta,\bb\eta)$ and
$(\bv',\bw',\theta',\bb\eta')$
be minimisers, with residuals $\br$ and $\br'$, respectively. Their
midpoint is feasible and, by convexity, is also a minimiser. Strict
convexity of
$\bb1^\top\bb\rho(\bb\eta)-\bYtil^\top\bb\eta$
in $\bb\eta$ then forces $\bb\eta=\bb\eta'$. The endpoint objective
values give
$\vnorm{\br}_2=\vnorm{\br'}_2$, while midpoint optimality gives
$$
    \vnorm{
        \frac{\br+\br'}2
    }_2
    =
    \frac{\vnorm{\br}_2+\vnorm{\br'}_2}{2}\,.
$$
Equality in the Euclidean triangle inequality, together with equality of
the norms, gives $\br=\br'$. Hence
$$
    \theta\bb1
    +
    \frac1{\sqrt p}\bH_1\bv
    +
    \frac1{\sqrt p}\bH_2\bw
    =
    \theta'\bb1
    +
    \frac1{\sqrt p}\bH_1\bv'
    +
    \frac1{\sqrt p}\bH_2\bw'\,.
$$
On $\mathcal E_n$, injectivity of
$$
    (\theta,\bbeta)
    \longmapsto
    \theta\bb1+\frac1{\sqrt p}\bH\bbeta
$$
followed by Lemma~\ref{lemma:orth_decomp} gives
$$
    (\bv,\bw,\theta,\bb\eta)
    =
    (\bv',\bw',\theta',\bb\eta')\,.
$$

The Lagrangian in \eqref{eq:lagrangian4} is continuous, convex in its
minimising variables and affine in $\blambda$ on compact convex domains.
Sion's theorem and attainment therefore give a saddle point. The
minimising variables of every saddle point minimise the primal objective
above and hence equal
$(\hat\bv,\hat\bw,\hat\theta,\hat{\bb\eta})$.

For each fixed $\bw$, the Sion identity in \eqref{eq:lagrangian5}
identifies the profiled value with the attained minimum of the same primal
objective over $(\bv,\theta,\bb\eta)$. Thus the profiled minimisers are
exactly the $\bw$-coordinates of minimisers of the primal objective, so
$\hat\bw$ is unique. At $\hat\bw$, the inner minimising variables at every
saddle point attain this fixed-$\hat\bw$ primal minimum and therefore equal
$(\hat\bv,\hat\theta,\hat{\bb\eta})$.
\end{proof}

\subsection{Lemmas for Section~\ref{sec:po-ao-pair}}



\begin{lemma}[Concave $\psi$ and ball-to-radius minimax reduction]
	\label{lemma:psi_conv_conc}
	Let
	$$
		Z
		\subseteq
		\Re^d\times\Re\times\Re^n\,,
	$$
	be a nonempty compact convex set, with coordinates
	$\bb\zeta=(\bv,\theta,\bb\eta)$. Let
	$S_{\blambda}\subseteq\Re^n$ be nonempty, compact, and convex, and let
	$S_{\bw}$ be a nonempty compact subset of a finite-dimensional Euclidean
	space. Define
	\begin{equation}
		\label{eq:psi_conv_conc_def}
		\psi(\blambda)
		=
		\underset{\bb\zeta=(\bv,\theta,\bb\eta)\in Z}{\min} \,
		\left\{
			\frac{1}{n}
			\left\{
				\bb 1^\top\bb\rho(\bb\eta)
				-
				\bYtil^\top\bb\eta
			\right\}
			+
			\frac{1}{n}
			\blambda^\top\bb\eta
			-
			\frac{1}{n}
			\theta\bb 1^\top\blambda
			-
			\frac{1}{n\sqrt p}
			\blambda^\top\bH_1\bv
		\right\} \,.
	\end{equation}
	Then the following hold.
	\begin{enumerate}[label=(\roman*),ref=\thelemma(\roman*)]
		\item \label{lemma:psi_conv_conc_concavity} The function $\psi$ is concave in $\blambda$ and continuous on $S_{\blambda}$. Moreover, the minimisation defining $\psi$ is over a compact convex set, its objective is convex in $\bb\zeta$ for fixed $\blambda$, affine in $\blambda$ for fixed $\bb\zeta$, and continuous in all variables. Consequently, the map
		\begin{equation}
			\label{eq:psi_conv_conc_po_map}
			(\bw,\blambda)
			\mapsto
			-
			\frac{1}{n\sqrt p}
			\blambda^\top\bH_2\bw
			+
			\psi(\blambda)\,,
		\end{equation}
		is continuous on $S_{\bw}\times S_{\blambda}$, affine in $\bw$, and
		concave in $\blambda$. In particular, it is convex--concave whenever
		$S_{\bw}$ is convex.

		\item \label{lemma:psi_conv_conc_ball_to_radius} For
		$\bb\zeta=(\bv,\theta,\bb\eta)\in Z$ and $\bb b\in\Re^n$, recall $G_n(\bb\eta)$ from \eqref{eq:ao_scalarisation_psi_def} and define
		\begin{equation}
			\label{eq:psi_conv_conc_ell_A_def}
			\bb a_{\bb b}(\bb\zeta)
			=
			\bb\eta
			-
			\theta\bb 1
			-
			\frac{1}{\sqrt p}
			\bH_1\bv
			-
			\bb b \,.
		\end{equation}
		Then, for every $R>0$, every $c\geq0$, and every $\bb b\in\Re^n$,
		\begin{equation}
			\label{eq:psi_conv_conc_ball_to_radius}
			\begin{aligned}
			\underset{\blambda\in\Re^n:\vnorm{\blambda}_2\leq R}{\max}
				\underset{\bb\zeta\in Z}{\min} \,
				\left\{
					G_n(\bb\eta)
					+
					\frac{1}{n}
					\blambda^\top
					\bb a_{\bb b}(\bb\zeta)
					-
					\frac{c}{n}
					\vnorm{\blambda}_2
				\right\} 
				=
				\underset{\tilde r\in[0,R]}{\max}
				\underset{\bb\zeta\in Z}{\min} \,
				\left\{
					G_n(\bb\eta)
					+
					\frac{\tilde r}{n}
					\vnorm{
						\bb a_{\bb b}(\bb\zeta)
					}_2
					-
					\frac{c\tilde r}{n}
				\right\} \,.
			\end{aligned}
		\end{equation}
	\end{enumerate}
\end{lemma}

\begin{proof}
We prove each part in turn. 
\begin{enumerate}[label=(\roman*)]
\item	Define
	\begin{equation}
		\label{eq:psi_conv_conc_L_def}
		L(\bb\zeta;\blambda)
		=
		\frac{1}{n}
		\left\{
			\bb 1^\top\bb\rho(\bb\eta)
			-
			\bYtil^\top\bb\eta
		\right\}
		+
		\frac{1}{n}
		\blambda^\top\bb\eta
		-
		\frac{1}{n}
		\theta\bb 1^\top\blambda
		-
		\frac{1}{n\sqrt p}
		\blambda^\top\bH_1\bv \,.
	\end{equation}
	For fixed $\blambda$, the map $\bb\zeta\mapsto L(\bb\zeta;\blambda)$ is
	convex and continuous on $Z$, because $\bb 1^\top\bb\rho(\bb\eta)$
	is convex in $\bb\eta$ and all remaining dependence on $\bb\zeta$ is
	affine. For fixed $\bb\zeta$, the map
	$\blambda\mapsto L(\bb\zeta;\blambda)$ is affine. Hence $\psi$ is the
	pointwise infimum of affine functions of $\blambda$, and is therefore
	concave.

	Since $Z$ is compact and $L$ is continuous, the minimum defining
	$\psi$ is attained. By Berge's maximum theorem (see, for example,
	\citealt[][Theorem~17.31.1]{aliprantis+et+al:2006}), applied to the
	constant compact correspondence $\blambda\mapsto Z$,
	the function $\psi$ is continuous on $S_{\blambda}$.

	The map in \eqref{eq:psi_conv_conc_po_map} is affine in $\bw$ for
	fixed $\blambda$. For fixed $\bw$, it is the sum of an affine function
	of $\blambda$ and the concave function $\psi(\blambda)$, and hence is
	concave in $\blambda$. Continuity follows from continuity of the
	bilinear term and continuity of $\psi$. If $S_{\bw}$ is convex, its
	affine dependence on $\bw$ gives the stated convex--concave property.

	\item For $\bb\zeta\in Z$ and $\blambda\in\Re^n$, define
	\begin{equation}
		\label{eq:psi_conv_conc_K_def}
		K(\bb\zeta,\blambda)
		=
		G_n(\bb\eta)
		+
		\frac{1}{n}
		\blambda^\top
		\bb a_{\bb b}(\bb\zeta)
		-
		\frac{c}{n}
		\vnorm{\blambda}_2 \,.
	\end{equation}
	For fixed $\blambda$, the map $\bb\zeta\mapsto K(\bb\zeta,\blambda)$ is convex and continuous on $Z$. For fixed $\bb\zeta$, the map $\blambda\mapsto K(\bb\zeta,\blambda)$ is concave and continuous on the ball $\{\blambda\in\Re^n:\vnorm{\blambda}_2\leq R\}$, because it is the sum of an affine function of $\blambda$ and $-c\vnorm{\blambda}_2/n$, with $c\geq0$. Since both constraint sets are compact and convex, Sion's min--max theorem (see, for example, \citealt[][Theorem~3]{simons:1995}) gives
	\begin{equation}
		\label{eq:psi_conv_conc_first_sion}
		\underset{\blambda\in\Re^n:\vnorm{\blambda}_2\leq R}{\max}
		\underset{\bb\zeta\in Z}{\min} \,
		K(\bb\zeta,\blambda)
		=
		\underset{\bb\zeta\in Z}{\min}
		\underset{\blambda\in\Re^n:\vnorm{\blambda}_2\leq R}{\max} \,
		K(\bb\zeta,\blambda) \,.
	\end{equation}
	For fixed $\bb\zeta\in Z$, write $\blambda=\tilde r\bu$, where $\tilde r\in[0,R]$ and $\bu\in\mathbb S^{n-1}$. Then
	\begin{equation}
		\label{eq:psi_conv_conc_support}
		\begin{aligned}
			&\underset{\blambda\in\Re^n:\vnorm{\blambda}_2\leq R}{\max} \,
			\left\{
				\frac{1}{n}
				\blambda^\top
				\bb a_{\bb b}(\bb\zeta)
				-
				\frac{c}{n}
				\vnorm{\blambda}_2
			\right\} 
			=
			\underset{\tilde r\in[0,R]}{\max} \,
			\left\{
				\frac{\tilde r}{n}
				\vnorm{
					\bb a_{\bb b}(\bb\zeta)
				}_2
				-
				\frac{c\tilde r}{n}
			\right\} \,,
		\end{aligned}
	\end{equation}
	where we used 
	\begin{equation}
		\label{eq:psi_conv_conc_support_identity}
		\underset{\bu\in\mathbb S^{n-1}}{\max} \,
		\bu^\top
		\bb a_{\bb b}(\bb\zeta)
		=
		\vnorm{
			\bb a_{\bb b}(\bb\zeta)
		}_2 \,.
	\end{equation}
	Thus
	\begin{equation}
		\label{eq:psi_conv_conc_before_second_sion}
		\begin{aligned}
			&\underset{\blambda\in\Re^n:\vnorm{\blambda}_2\leq R}{\max}
			\underset{\bb\zeta\in Z}{\min} \,
			K(\bb\zeta,\blambda)
			=
			\underset{\bb\zeta\in Z}{\min}
			\underset{\tilde r\in[0,R]}{\max} \,
			\left\{
				G_n(\bb\eta)
				+
				\frac{\tilde r}{n}
				\vnorm{
					\bb a_{\bb b}(\bb\zeta)
				}_2
				-
				\frac{c\tilde r}{n}
			\right\} \,.
		\end{aligned}
	\end{equation}
	For fixed $\tilde r\in[0,R]$, the map
	$$
		\bb\zeta
		\mapsto
		G_n(\bb\eta)
		+
		\frac{\tilde r}{n}
		\vnorm{
			\bb a_{\bb b}(\bb\zeta)
		}_2
		-
		\frac{c\tilde r}{n}\,,
	$$
	is convex and continuous on $Z$, because $\tilde r\geq0$ and $\bb a_{\bb b}$ is affine. For fixed $\bb\zeta$, the same expression is affine, hence concave and continuous, in $\tilde r$. Applying the same theorem on $Z\times[0,R]$ gives
	\begin{equation}
		\label{eq:psi_conv_conc_second_sion}
		\begin{aligned}
			&\underset{\bb\zeta\in Z}{\min}
			\underset{\tilde r\in[0,R]}{\max} \,
			\left\{
				G_n(\bb\eta)
				+
				\frac{\tilde r}{n}
				\vnorm{
					\bb a_{\bb b}(\bb\zeta)
				}_2
				-
				\frac{c\tilde r}{n}
			\right\} 
			=
			\underset{\tilde r\in[0,R]}{\max}
			\underset{\bb\zeta\in Z}{\min} \,
			\left\{
				G_n(\bb\eta)
				+
				\frac{\tilde r}{n}
				\vnorm{
					\bb a_{\bb b}(\bb\zeta)
				}_2
				-
				\frac{c\tilde r}{n}
			\right\} \,.
		\end{aligned}
	\end{equation}
	Combining \eqref{eq:psi_conv_conc_first_sion}--\eqref{eq:psi_conv_conc_second_sion} proves \eqref{eq:psi_conv_conc_ball_to_radius}.
\end{enumerate}
\end{proof}

\subsection{Lemmas for Section~\ref{sec:ao-scalar}}

\subsubsection{Small-$r$ environment and exclusion}
\label{subsec:small_r_exclusion}

\begin{lemma}[Conditional small-$r$ environment]
\label{lemma:small_r_environment}
Fix $t_0>0$. There exist deterministic constants
$$
	d_r>0,
	\qquad
	a\in(0,1/2),
	\qquad
	K,C_L\in(0,\infty)\,,
$$
and constants $C,c>0$ and $N\in\mathbb N$ such that, for all $n>N$,
\begin{equation}
\label{eq:small_r_environment_probability}
	\mathds{1}_{\mathcal S_n(t_0)}
	\Pr\left(
		\lnot\mathcal E_n^{\mathrm{sr}}(d_r,a,K,C_L)
		\mid
		\mathcal G_n
	\right)
	\leq
	C\exp\{-cn\}\,.
\end{equation}
The constants depend only on
$$
	\kappa,
	\alpha,
	\bb\Gamma,
	\theta_0^*,
	\theta_P^*,
	t_0,
	C_\beta,
	C_\theta,
	C_\sigma\,,
$$
and on the constants in the probability bounds of
Section~\ref{sec:setup}. In particular, they do not depend on $C_\eta$.
Moreover, deterministically, for every $\sigma\in[0,C_\sigma]$ and every
$\bb b\in\sigma\bg+\mathcal R_0$,
\begin{equation}
\label{eq:small_r_deterministic_envelope}
	\abs{b_i-\eta_{\Ytil,i}}
	\leq
	L_i,
	\qquad
	i=1,\ldots,n\,.
\end{equation}
\end{lemma}

\begin{proof}
Let $s=\rank(\bb B)\in\{0,1,2\}$ and choose, by a fixed
$\mathcal G_n$-measurable orthonormalisation rule, an orthonormal basis
$\{\bb b_1,\ldots,\bb b_s\}$ of $\operatorname{col}(\bb B)$. With
$\bh_i^\top$ denoting the $i$th row of $\bH$, set
$$
	\bb r_i
	=
	(1,\bh_i^\top\bb b_1,\ldots,\bh_i^\top\bb b_s)^\top
	\in
	\Re^{s+1}\,,
$$
with $\bb r_i=1$ when $s=0$. Conditional on $\mathcal G_n$, the last
$s$ coordinates of $\bb r_i$ are independent standard Gaussian variables,
so
$$
	\expect[
		\bb r_i\bb r_i^\top
		\mid
		\mathcal G_n
	]
	=
	\bb I_{s+1}\,.
$$
Define
$$
	\bSigma_n
	=
	\frac1n\sum_{i=1}^n\bb r_i\bb r_i^\top,
	\qquad
	\bmu_n
	=
	\frac1n\sum_{i=1}^n\bb r_i\eta_{\Ytil,i},
	\qquad
	s_n
	=
	\frac1n\sum_{i=1}^n\eta_{\Ytil,i}^2\,,
$$
and
$$
	\bar{\bmu}_n
	=
	\expect[
		\bb r_i\eta_{\Ytil,i}
		\mid
		\mathcal G_n
	],
	\qquad
	\bar s_n
	=
	\expect[
		\eta_{\Ytil,i}^2
		\mid
		\mathcal G_n
	]\,.
$$

On $\mathcal S_n(t_0)$, both $\abs{\theta_P}$ and
$\vnorm{\bbeta_P}_2/\sqrt p$ are bounded by deterministic constants. Thus
$$
	\eta_{P,i}
	=
	\theta_P
	+
	\frac1{\sqrt p}\bh_i^\top\bbeta_P
$$
is conditionally subgaussian with a deterministic norm bound. If
$q=\rho'(x)$, the two possible values
$(1-\alpha)q$ and $\alpha+(1-\alpha)q$ of the pseudo-response satisfy
$$
	\abs{\logit{(1-\alpha)q}}
	\vee
	\abs{\logit{\{\alpha+(1-\alpha)q\}}}
	\leq
	\abs{x}+C_\alpha
$$
for a finite constant $C_\alpha$ depending only on $\alpha$. Consequently,
conditional on $\mathcal G_n$ and on $\mathcal S_n(t_0)$,
\begin{equation}
\label{eq:small_r_eta_subgaussian}
	\vnorm{\eta_{\Ytil,i}}_{\psi_2}
	\leq
	K_\eta
\end{equation}
for a deterministic $K_\eta<\infty$, uniformly in $i$ and $n$.

By \citealt[][Lemmas~2.7.6 and~2.7.7]{vershynin:2018}, squares and
products of subgaussian variables are sub-exponential. Since
$s+1\leq3$, there is a deterministic constant $K_0<\infty$ such that,
conditionally on $\mathcal G_n$ and on $\mathcal S_n(t_0)$,
$$
	\vnorm{\eta_{\Ytil,i}^2-\bar s_n}_{\psi_1}
	\leq
	K_0,
	\qquad
	\vnorm{\bb r_i\eta_{\Ytil,i}-\bar{\bmu}_n}_{\psi_1}
	\leq
	K_0
$$
coordinatewise, and
$$
	\vnorm{\bb r_i\bb r_i^\top-\bb I_{s+1}}_{\psi_1}
	\leq
	K_0
$$
entrywise. Bernstein's inequality for sums of independent sub-exponential variables
\citep[Theorem~2.8.1]{vershynin:2018}, applied conditionally on
$\mathcal G_n$, and a union bound over the fixed number of coordinates give,
for every fixed $u>0$,
\begin{equation}
\label{eq:small_r_moment_concentration}
\begin{aligned}
	\mathds{1}_{\mathcal S_n(t_0)}
	\Pr\Bigg(
		\max\Big\{
			\abs{s_n-\bar s_n},
			\vnorm{\bmu_n-\bar{\bmu}_n}_2,
			\mnorm{\bSigma_n-\bb I_{s+1}}_2
		\Big\}
		>u
		\mid
		\mathcal G_n
	\Bigg)
	\leq
	C\exp\{-cn\min(u^2,u)\}\,.
\end{aligned}
\end{equation}

We next obtain a population residual that is bounded away from zero. Let
$$
	q_{0,i}
	=
	\rho'\left(
		\theta_0
		+
		\frac1{\sqrt p}\bh_i^\top\bbeta_0
	\right),
	\qquad
	q_{P,i}
	=
	\rho'\left(
		\theta_P
		+
		\frac1{\sqrt p}\bh_i^\top\bbeta_P
	\right)\,.
$$
Because $\bbeta_0,\bbeta_P\in\operatorname{col}(\bb B)$, both
$q_{0,i}$ and $q_{P,i}$ are measurable with respect to
$\sigma(\bb r_i,\mathcal G_n)$. Conditional on this field,
$$
	B_i
	=
	\mathds{1}
	\left\{
		\varepsilon_i<q_{0,i}
	\right\}
	\sim
	\operatorname{Bernoulli}(q_{0,i}),
	\qquad
	\Ytil_i
	=
	(1-\alpha)q_{P,i}+\alpha B_i\,.
$$
The two corresponding values of $\eta_{\Ytil,i}$ differ by
$$
	\int_{(1-\alpha)q_{P,i}}^{(1-\alpha)q_{P,i}+\alpha}
	\frac{du}{u(1-u)}
	\geq
	4\alpha\,.
$$
It follows that
$$
	\var(
		\eta_{\Ytil,i}
		\mid
		\bb r_i,\mathcal G_n
	)
	\geq
	16\alpha^2q_{0,i}(1-q_{0,i})\,.
$$
On $\mathcal S_n(t_0)$, the pair
$$
	\left(
		\theta_0,
		\frac{\vnorm{\bbeta_0}_2^2}{p}
	\right)
$$
ranges over a deterministic compact subset of $\Re\times[0,\infty)$.
The map
$$
	(\theta,v)
	\longmapsto
	\expect\left[
		\rho'(\theta+\sqrt v Z)
		\{1-\rho'(\theta+\sqrt v Z)\}
	\right],
	\qquad
	Z\sim\mathrm N(0,1)\,,
$$
is continuous and strictly positive on this compact set. Hence there is a
deterministic $\underline\mu_0>0$ such that, on
$\mathcal S_n(t_0)$,
$$
	\expect[
		q_{0,i}(1-q_{0,i})
		\mid
		\mathcal G_n
	]
	\geq
	\underline\mu_0\,.
$$
With
$$
	\tau_0
	=
	4\alpha\sqrt{\underline\mu_0}\,,
$$
for every $\bb a\in\Re^{s+1}$,
$$
\begin{aligned}
	\expect\left[
		(\eta_{\Ytil,i}-\bb r_i^\top\bb a)^2
		\mid
		\mathcal G_n
	\right]
	&\geq
	\expect\left[
		\var(
			\eta_{\Ytil,i}
			\mid
			\bb r_i,\mathcal G_n
		)
		\mid
		\mathcal G_n
	\right]
	\\
	&\geq
	\tau_0^2\,.
\end{aligned}
$$
Since $\expect[\bb r_i\bb r_i^\top\mid\mathcal G_n]=\bb I_{s+1}$,
this is equivalent to
\begin{equation}
\label{eq:small_r_population_residual}
	\bar s_n
	-
	\vnorm{\bar{\bmu}_n}_2^2
	\geq
	\tau_0^2 \,.
\end{equation}

Let $\bb U=[\bb b_1,\ldots,\bb b_s]$. Since
$\operatorname{col}(\bb U)=\operatorname{col}(\bb B)$,
$$
	\operatorname{col}(\bH\bb U)
	=
	\operatorname{col}(\bH\bb B)
	=
	\operatorname{col}(\bH_1) \,.
$$
Thus the matrix with rows $\bb r_i^\top$ has column space
$\vspan\{\bb1,\operatorname{col}(\bH_1)\}$. Let
$\bb\Pi_n^\perp$ denote the orthogonal projector onto the complement of
this space. Then
\begin{equation}
\label{eq:small_r_projection_residual_identity}
	\frac1n\vnorm{\bb\Pi_n^\perp\bEtaYtil}_2^2
	=
	\underset{\bb a\in\Re^{s+1}}{\inf}
	\frac1n\sum_{i=1}^n
	(\eta_{\Ytil,i}-\bb r_i^\top\bb a)^2 \,.
\end{equation}

By \eqref{eq:small_r_eta_subgaussian}, there is a deterministic
$B<\infty$ such that $\bar s_n\leq B^2$ on
$\mathcal S_n(t_0)$. Conditional Cauchy--Schwarz gives
$$
	\vnorm{\bar{\bmu}_n}_2
	\leq
	\expect[
		\vnorm{\bb r_i}_2^2
		\mid
		\mathcal G_n
	]^{1/2}
	\bar s_n^{1/2}
	\leq
	\sqrt3B \,.
$$
Fix $\delta_0\in(0,1/2]$ and suppose that the maximum in
\eqref{eq:small_r_moment_concentration} is at most $\delta_0$. Then, for
every $\bb a\in\Re^{s+1}$,
$$
\begin{aligned}
	\frac1n\sum_{i=1}^n
	(\eta_{\Ytil,i}-\bb r_i^\top\bb a)^2
	&=
	s_n-2\bb a^\top\bmu_n+\bb a^\top\bSigma_n\bb a
	\\
	&\geq
	s_n-2\bb a^\top\bmu_n+(1-\delta_0)\vnorm{\bb a}_2^2
	\\
	&\geq
	s_n-
	\frac{\vnorm{\bmu_n}_2^2}{1-\delta_0}\,.
\end{aligned}
$$
Writing $m=\vnorm{\bar{\bmu}_n}_2\leq\sqrt3B$ and using
$$
	\frac{(m+\delta_0)^2-(1-\delta_0)m^2}{1-\delta_0}
	\leq
	2\delta_0(m+1)^2\,,
$$
we obtain from \eqref{eq:small_r_population_residual} that
$$
	\frac1n\vnorm{\bb\Pi_n^\perp\bEtaYtil}_2^2
	\geq
	\tau_0^2
	-
	\delta_0
	\left\{
		1+2(\sqrt3B+1)^2
	\right\}\,.
$$
Choose $\delta_0$ so that the subtracted term is at most
$\tau_0^2/2$, and set
$$
	a_0
	=
	\frac{\tau_0}{2}\,.
$$
By \eqref{eq:small_r_moment_concentration},
\begin{equation}
\label{eq:small_r_empirical_residual_probability}
	\mathds{1}_{\mathcal S_n(t_0)}
	\Pr\left(
		\frac{\vnorm{\bb\Pi_n^\perp\bEtaYtil}_2}{\sqrt n}
		<a_0
		\mid
		\mathcal G_n
	\right)
	\leq
	C\exp\{-cn\} \,.
\end{equation}

Choose $\delta\in(0,1/2)$ sufficiently small that
$$
	q_\delta
	=
	\frac{\sqrt\kappa+\delta}{1-\delta}: \quad 
	q_\delta^2<1-\delta \,,
$$
which is possible because $\kappa<1$, and define
\begin{equation}
\label{eq:small_r_gap_constant}
	d_r
	=
	a_0
	\sqrt{1-\delta-q_\delta^2}
	>
	0\,.
\end{equation}

We now establish the envelope bounds. By \eqref{eq:ao_H_repar}, on the
common completed AO probability space,
$$
	\frac1{\sqrt p}\bH_1
	=
	\bb Z\bb\Gamma_p^{1/2}\,,
$$
where the rows $\bz_i^\top$ of $\bb Z$ are conditionally independent
standard Gaussian vectors. If $\bv\in V_p$ and
$\bu=\bb\Gamma_p^{1/2}\bv$, then
$$
	\vnorm{\bu}_2^2
	=
	\bv^\top\bb\Gamma_p\bv
	=
	\frac1p\vnorm{\bb B\bv}_2^2
	\leq
	C_\beta^2\,.
$$
Consequently,
\begin{equation}
\label{eq:small_r_signal_envelope}
	\underset{\bv\in V_p}{\sup}
	\abs{
		\frac1{\sqrt p}(\bH_1\bv)_i
	}
	\leq
	C_\beta\vnorm{\bz_i}_2 \,.
\end{equation}
Together with \eqref{eq:small_r_eta_subgaussian} and the Gaussianity of
$g_i$, this shows that $L_i$ is conditionally subgaussian and $L_i^2$ is
conditionally sub-exponential, with deterministic norm bounds on
$\mathcal S_n(t_0)$. Choose $C_L$ larger than twice a uniform upper bound
for $\expect[L_i^2\mid\mathcal G_n]$. Bernstein's inequality \citep[Theorem~2.8.1]{vershynin:2018} then gives
\begin{equation}
\label{eq:small_r_envelope_energy_probability}
	\mathds{1}_{\mathcal S_n(t_0)}
	\Pr\left(
		\frac1n\sum_{i=1}^nL_i^2>C_L
		\mid
		\mathcal G_n
	\right)
	\leq
	C\exp\{-cn\} \,.
\end{equation}

Let
$$
	Q_{P,i}
	=
	\frac1{\sqrt p}(\bH\bbeta_P)_i\,.
$$
On $\mathcal S_n(t_0)$, $Q_{P,i}$ is conditionally Gaussian with uniformly
bounded variance. Fix $M_0<\infty$. On
$\{\abs{Q_{P,i}}\leq M_0\}\cap\mathcal S_n(t_0)$,
$$
	\Ytil_i
	\geq
	(1-\alpha)\rho'(\theta_P^*-t_0-M_0) \,,
$$
and
$$
	1-\Ytil_i
	\geq
	(1-\alpha)
	\{1-\rho'(\theta_P^*+t_0+M_0)\}\,.
$$
Set
$$
	a
	=
	\min\left\{
		(1-\alpha)\rho'(\theta_P^*-t_0-M_0),
		(1-\alpha)
		\{1-\rho'(\theta_P^*+t_0+M_0)\},
		\frac14
	\right\}\,.
$$
Then $a\in(0,1/2)$ and
\begin{equation}
\label{eq:small_r_boundary_inclusion}
	\{\Ytil_i\notin[a,1-a]\}
	\subseteq
	\{\abs{Q_{P,i}}>M_0\}\,.
\end{equation}
The uniform conditional subgaussian bound for $L_i$ gives a uniform fourth
moment, and conditional Cauchy--Schwarz therefore yields
$$
\begin{aligned}
	\expect\left[
		L_i^2\mathds{1}\{\abs{Q_{P,i}}>M_0\}
		\mid
		\mathcal G_n
	\right]
	&\leq
	\expect[
		L_i^4
		\mid
		\mathcal G_n
	]^{1/2}
	\\
	&\quad\times
	\Pr\left(
		\abs{Q_{P,i}}>M_0
		\mid
		\mathcal G_n
	\right)^{1/2}\,.
\end{aligned}
$$
The right-hand side tends to zero uniformly on $\mathcal S_n(t_0)$ as
$M_0\to\infty$. The same subgaussian bound gives
$$
	\underset{i,n}{\sup}
	\mathds{1}_{\mathcal S_n(t_0)}
	\expect\left[
		L_i^2\mathds{1}\{L_i>K\}
		\mid
		\mathcal G_n
	\right]
	\longrightarrow
	0
$$
as $K\to\infty$. Choose first $M_0$ and then $K$ so that
$$
	\mathds{1}_{\mathcal S_n(t_0)}
	\expect\left[
		L_i^2
		\mathds{1}
		\{\abs{Q_{P,i}}>M_0\ \text{or}\ L_i>K\}
		\mid
		\mathcal G_n
	\right]
	\leq
	\frac{d_r^2}{128}\,.
$$
The variables in this expectation are conditionally independent across
$i$, and their centred versions have uniformly bounded conditional
sub-exponential norms because they are dominated by $L_i^2$. Another
application of the Bernstein inequality \citep[Theorem~2.8.1]{vershynin:2018}, together with
\eqref{eq:small_r_boundary_inclusion}, gives
\begin{equation}
\label{eq:small_r_exceptional_energy_probability}
\begin{aligned}
	\mathds{1}_{\mathcal S_n(t_0)}
	\Pr\left(
		\frac1n
		\sum_{i\notin\mathcal I_n(a,K)}L_i^2
		>
		\frac{d_r^2}{64}
		\mid
		\mathcal G_n
	\right)
	\leq
	C\exp\{-cn\}\,.
\end{aligned}
\end{equation}

It remains to turn the projected residual into the uniform gap. Conditional
on $\mathcal C_n$, the vector $\bEtaYtil$ and the projector
$\bb\Pi_n^\perp$ are fixed, whereas $\bg$ and $\bh$ are independent
standard Gaussian vectors. Set
$$
	a_n
	=
	\frac{\vnorm{\bb\Pi_n^\perp\bEtaYtil}_2}{\sqrt n},
	\qquad
	b_n
	=
	\frac{\vnorm{\bb\Pi_n^\perp\bg}_2}{\sqrt n},
	\qquad
	c_n
	=
	\frac{\vnorm{\bh}_2}{\sqrt n}\,,
$$
and
$$
	d_n
	=
	\frac1n
	\left\langle
		\bb\Pi_n^\perp\bEtaYtil,
		\bb\Pi_n^\perp\bg
	\right\rangle\,.
$$
The rank $m$ of $\bb\Pi_n^\perp$ satisfies $n-3\leq m\leq n$. Since
$\vnorm{\bb\Pi_n^\perp\bg}_2^2\sim\chi_m^2$, Gaussian norm
concentration \citep[Theorem~3.1.1]{vershynin:2018} and
$$
	0
	\leq
	1-\sqrt{m/n}
	\leq
	\frac3n
$$
give, for all sufficiently large $n$,
\begin{equation}
\label{eq:small_r_projected_gaussian_norm}
	\Pr\left(
		\abs{b_n-1}>\delta
		\mid
		\mathcal C_n
	\right)
	\leq
	C\exp\{-cn\}\,.
\end{equation}
Similarly, since $p/n\to\kappa$ and
$\bh\sim\mathrm N(\bb0_{p-s},\bb I_{p-s})$ conditionally on
$\mathcal C_n$,
\begin{equation}
\label{eq:small_r_h_norm}
	\Pr\left(
		c_n>\sqrt\kappa+\delta
		\mid
		\mathcal C_n
	\right)
	\leq
	C\exp\{-cn\}\,.
\end{equation}
Moreover,
$$
	d_n
	=
	\frac1n
	\left\langle
		\bb\Pi_n^\perp\bEtaYtil,
		\bg
	\right\rangle
	\sim
	\mathrm N\left(
		0,
		\frac{a_n^2}{n}
	\right)
$$
conditionally on $\mathcal C_n$. Hence, on $\{a_n\geq a_0\}$,
\begin{equation}
\label{eq:small_r_cross_term}
\begin{aligned}
	\Pr\left(
		\abs{d_n}>\delta(1-\delta)a_n
		\mid
		\mathcal C_n
	\right)
	&\leq
	2\exp\left\{
		-\frac{n\delta^2(1-\delta)^2}{2}
	\right\}\,.
\end{aligned}
\end{equation}
On $\{b_n\geq1-\delta\}$, the complementary event in
\eqref{eq:small_r_cross_term} implies
$\abs{d_n}\leq\delta a_nb_n$.

For every $\sigma\geq0$, projection onto
$\vspan\{\bb1,\operatorname{col}(\bH_1)\}^\perp$ gives
$$
\begin{aligned}
	M_n(\sigma,\bEtaYtil)
	-
	\sigma\frac{\vnorm{\bh}_2}{\sqrt n}
	&\geq
	\frac1{\sqrt n}
	\vnorm{
		\bb\Pi_n^\perp(\bEtaYtil-\sigma\bg)
	}_2
	-
	\sigma c_n
	\\
	&=
	\sqrt{
		a_n^2+\sigma^2b_n^2-2\sigma d_n
	}
	-
	\sigma c_n\,.
\end{aligned}
$$
Suppose that
$$
	a_n\geq a_0,
	\qquad
	b_n\geq1-\delta,
	\qquad
	c_n\leq\sqrt\kappa+\delta,
	\qquad
	\abs{d_n}\leq\delta a_nb_n\,.
$$
With $t=\sigma b_n$ and $q_n=c_n/b_n\leq q_\delta$, one has
$$
	0\leq q_n^2\leq q_\delta^2<1-\delta\,,
$$
so the square root below is real. Moreover,
$$
	\sqrt{
		a_n^2+\sigma^2b_n^2-2\sigma d_n
	}
	-
	\sigma c_n
	\geq
	\sqrt{1-\delta}
	\sqrt{a_n^2+t^2}
	-
	q_nt\,.
$$
By the Cauchy--Schwarz inequality,
$$
	q_nt
	+
	a_n\sqrt{1-\delta-q_n^2}
	\leq
	\sqrt{1-\delta}
	\sqrt{a_n^2+t^2}\,.
$$
Consequently, uniformly over $\sigma\geq0$,
$$
	M_n(\sigma,\bEtaYtil)
	-
	\sigma\frac{\vnorm{\bh}_2}{\sqrt n}
	\geq
	a_n\sqrt{1-\delta-q_n^2}
	\geq
	d_r\,,
$$
where the last step uses \eqref{eq:small_r_gap_constant}.
Combining \eqref{eq:small_r_empirical_residual_probability} and
\eqref{eq:small_r_projected_gaussian_norm}--\eqref{eq:small_r_cross_term}
by the tower property shows that the gap inequality fails with
$\mathcal G_n$-conditional probability at most $C\exp\{-cn\}$ on
$\mathcal S_n(t_0)$.

Finally, if $\bb b\in\sigma\bg+\mathcal R_0$, then
$$
	b_i
	=
	\sigma g_i
	+
	\theta
	+
	\frac1{\sqrt p}(\bH_1\bv)_i
$$
for some $\abs{\theta}\leq C_\theta$ and $\bv\in V_p$, so
\eqref{eq:small_r_deterministic_envelope} follows directly from the
definition of $L_i$. The gap bound,
\eqref{eq:small_r_envelope_energy_probability}, and
\eqref{eq:small_r_exceptional_energy_probability}, combined by a union
bound, prove \eqref{eq:small_r_environment_probability}.
\end{proof}

\begin{proposition}[Deterministic small-$r$ exclusion]
\label{prop:small_r_exclusion}
Assume that $\mathcal E_n^{\mathrm{sr}}(d_r,a,K,C_L)$ holds and that
$C_\eta\geq\sqrt{C_L}$. Set
\begin{equation}
\label{eq:small_r_radius}
	c_{a,K}
	=
	\rho''(\logit{1-a}+K),
	\qquad
	r_0
	=
	\min\left\{
		\frac12,
		\frac{c_{a,K} d_r}{16}
	\right\}\,.
\end{equation}
Let $\mathcal R\subseteq\mathcal R_0$ be nonempty, compact, and convex. For
$\sigma\in[0,C_\sigma]$, $r\geq0$, and $\bb\eta\in\Re^n$, define
$$
	M_{\mathcal R}(\sigma,\bb\eta)
	=
	\frac1{\sqrt n}
	\operatorname{dist}_2
	\left(
		\bb\eta-\sigma\bg,
		\mathcal R
	\right)\,,
$$
$$
	\widetilde{\bb\eta}_{\mathcal R,\sigma,r}
	=
	\underset{\bb\eta\in\Re^n}{\arg\min}
	\left\{
		G_n(\bb\eta)
		+
		rM_{\mathcal R}(\sigma,\bb\eta)
	\right\}\,,
$$
and
$$
	A_{\mathcal R}(\sigma,r)
	=
	-\sigma r\frac{\vnorm{\bh}_2}{\sqrt n}
	+
	\underset{\bb\eta\in\mathcal B^n_{C_\eta}}{\min}
	\left\{
		G_n(\bb\eta)
		+
		rM_{\mathcal R}(\sigma,\bb\eta)
	\right\}\,.
$$
Then, uniformly over $\sigma\in[0,C_\sigma]$ and over all such
$\mathcal R$, the following hold.
\begin{enumerate}[label=(\roman*)]
	\item For every $r\geq0$, the unconstrained minimiser exists uniquely and
	$$
		\frac{
			\vnorm{\widetilde{\bb\eta}_{\mathcal R,\sigma,r}}_2
		}{
			\sqrt n
		}
		\leq
		\sqrt{C_L}
		\leq
		C_\eta\,.
	$$
	Hence it is also the unique minimiser over
	$\mathcal B^n_{C_\eta}$.

	\item For every $r\geq0$,
	\begin{equation}
	\label{eq:small_r_displacement}
		\frac1{\sqrt n}
		\vnorm{
			\widetilde{\bb\eta}_{\mathcal R,\sigma,r}
			-
			\bEtaYtil
		}_2
		\leq
		\frac{2r}{c_{a,K}}
		+
		\frac{d_r}{8}\,.
	\end{equation}
	In particular, the left-hand side is at most $d_r/4$ whenever
	$r\in[0,r_0]$.

	\item For $0\leq r<r'\leq r_0$,
	\begin{equation}
	\label{eq:small_r_profile_increment}
		A_{\mathcal R}(\sigma,r')
		-
		A_{\mathcal R}(\sigma,r)
		\geq
		\frac{3d_r}{4}(r'-r)\,.
	\end{equation}
	Consequently, for every $c\in(0,r_0]$,
	\begin{equation}
	\label{eq:small_r_profile_exclusion}
		\underset{r\in[0,1]}{\max}
		A_{\mathcal R}(\sigma,r)
		=
		\underset{r\in[c,1]}{\max}
		A_{\mathcal R}(\sigma,r)\,.
	\end{equation}
\end{enumerate}
All conclusions remain valid after any further deterministic enlargement of
$C_\eta$.
\end{proposition}

\begin{proof}
Fix $\sigma$ and $\mathcal R$, and write
$$
	M(\bb\eta)
	=
	M_{\mathcal R}(\sigma,\bb\eta),
	\qquad
	\widetilde{\bb\eta}_r
	=
	\widetilde{\bb\eta}_{\mathcal R,\sigma,r}\,.
$$
For $i=1,\ldots,n$, set
$$
	f_i(x)
	=
	\rho(x)-\Ytil_i x\,.
$$
Since $\Ytil_i\in(0,1)$, $f_i$ is strictly convex and coercive, with
unique minimiser $\eta_{\Ytil,i}$. Hence $G_n$ is strictly convex and
coercive. The function $M$ is convex, continuous, nonnegative, and
$1/\sqrt n$-Lipschitz. Therefore
$G_n+rM$ is strictly convex and coercive for every $r\geq0$, so
$\widetilde{\bb\eta}_r$ exists uniquely.

Let
$$
	\bb b
	=
	\sigma\bg
	+
	P_{\mathcal R}
	(\widetilde{\bb\eta}_r-\sigma\bg)\,.
$$
Then $\bb b\in\sigma\bg+\mathcal R\subseteq\sigma\bg+\mathcal R_0$ and
$$
	M(\widetilde{\bb\eta}_r)
	=
	\frac1{\sqrt n}
	\vnorm{\widetilde{\bb\eta}_r-\bb b}_2\,.
$$
For each $i$, let $J_i$ be the closed interval with endpoints
$\eta_{\Ytil,i}$ and $b_i$, and define
$$
	\eta_i'
	=
	P_{J_i}(\widetilde\eta_{r,i})\,.
$$
If $\widetilde\eta_{r,i}\notin J_i$, then, for some
$t_i\in[0,1)$,
$$
	\eta_i'
	=
	(1-t_i)\eta_{\Ytil,i}
	+
	t_i\widetilde\eta_{r,i}\,.
$$
Since $\widetilde\eta_{r,i}\neq\eta_{\Ytil,i}$ and
$\eta_{\Ytil,i}$ is the unique minimiser of $f_i$, convexity, with strict
inequality when $t_i\in(0,1)$, gives
$$
\begin{aligned}
	f_i(\eta_i')
	&\leq
	(1-t_i)f_i(\eta_{\Ytil,i})
	+
	t_i f_i(\widetilde\eta_{r,i})
	\\
	&<
	f_i(\widetilde\eta_{r,i})\,.
\end{aligned}
$$
Moreover, since $b_i\in J_i$ and projection onto an interval is
nonexpansive,
$$
	\abs{\eta_i'-b_i}
	\leq
	\abs{\widetilde\eta_{r,i}-b_i}\,.
$$
Thus, if the vector $\bb\eta'$ differs from
$\widetilde{\bb\eta}_r$, then
$$
	G_n(\bb\eta')
	<
	G_n(\widetilde{\bb\eta}_r)
$$
and
$$
	M(\bb\eta')
	\leq
	\frac1{\sqrt n}\vnorm{\bb\eta'-\bb b}_2
	\leq
	\frac1{\sqrt n}\vnorm{\widetilde{\bb\eta}_r-\bb b}_2
	=
	M(\widetilde{\bb\eta}_r)\,,
$$
contradicting optimality. Hence, for every $i$,
$\widetilde\eta_{r,i}$ lies between $\eta_{\Ytil,i}$ and $b_i$.
By \eqref{eq:small_r_deterministic_envelope},
$$
	\abs{b_i-\eta_{\Ytil,i}}
	\leq
	L_i\,.
$$
The definition of $L_i$ also gives
$\abs{\eta_{\Ytil,i}}\leq L_i$ and $\abs{b_i}\leq L_i$. Therefore
\begin{equation}
\label{eq:small_r_coordinate_envelope}
	\abs{
		\widetilde\eta_{r,i}-\eta_{\Ytil,i}
	}
	\leq
	L_i,
	\qquad
	\abs{\widetilde\eta_{r,i}}
	\leq
	L_i\,.
\end{equation}
The environment event now gives
$$
	\frac{\vnorm{\widetilde{\bb\eta}_r}_2}{\sqrt n}
	\leq
	\left(
		\frac1n\sum_{i=1}^nL_i^2
	\right)^{1/2}
	\leq
	\sqrt{C_L}
	\leq
	C_\eta\,.
$$
Thus the unconstrained minimiser is feasible for
$\mathcal B^n_{C_\eta}$ and is the unique constrained minimiser. This proves
(i).

Set
$$
	\Delta_r
	=
	\frac1{\sqrt n}
	\vnorm{
		\widetilde{\bb\eta}_r-\bEtaYtil
	}_2\,.
$$
Optimality against $\bEtaYtil$ and the Lipschitz property of $M$ give
\begin{equation}
\label{eq:small_r_excess_upper}
\begin{aligned}
	G_n(\widetilde{\bb\eta}_r)
	-
	G_n(\bEtaYtil)
	&\leq
	r
	\left\{
		M(\bEtaYtil)-M(\widetilde{\bb\eta}_r)
	\right\}
	\\
	&\leq
	r\Delta_r\,.
\end{aligned}
\end{equation}
Since $f_i'(\eta_{\Ytil,i})=0$ and $f_i''=\rho''$,
$$
	G_n(\widetilde{\bb\eta}_r)
	-
	G_n(\bEtaYtil)
	=
	\frac1n\sum_{i=1}^nD_i\,,
$$
where
$$
	D_i
	=
	\int_{\eta_{\Ytil,i}}^{\widetilde\eta_{r,i}}
	(\widetilde\eta_{r,i}-u)\rho''(u)\,du
	\geq
	0\,.
$$
If $i\in\mathcal I_n(a,K)$, then
$\abs{\eta_{\Ytil,i}}\leq\logit{1-a}$ and
\eqref{eq:small_r_coordinate_envelope} gives
$\abs{\widetilde\eta_{r,i}-\eta_{\Ytil,i}}\leq K$. Every point between
$\eta_{\Ytil,i}$ and $\widetilde\eta_{r,i}$ therefore belongs to
$[-\logit{1-a}-K,\logit{1-a}+K]$. Since $\rho''$ is even and nonincreasing on
$[0,\infty)$,
$$
	D_i
	\geq
	\frac{c_{a,K}}2
	(\widetilde\eta_{r,i}-\eta_{\Ytil,i})^2,
	\qquad
	i\in\mathcal I_n(a,K)\,.
$$
Using this inequality on $\mathcal I_n(a,K)$,
\eqref{eq:small_r_coordinate_envelope} on its complement, and
\eqref{eq:small_r_excess_upper}, we obtain
$$
\begin{aligned}
	\Delta_r^2
	&\leq
	\frac{2}{c_{a,K}}
	\frac1n\sum_{i=1}^nD_i
	+
	\frac1n
	\sum_{i\notin\mathcal I_n(a,K)}L_i^2
	\\
	&\leq
	\frac{2r}{c_{a,K}}\Delta_r
	+
	\frac{d_r^2}{64}\,.
\end{aligned}
$$
If $x^2\leq ux+v^2$ with $x,u,v\geq0$, then $x\leq u+v$. Hence
\eqref{eq:small_r_displacement} follows. If
$r\leq r_0$, then
$2r/c_{a,K}\leq d_r/8$, and therefore $\Delta_r\leq d_r/4$. This proves
(ii).

Let $0\leq r<r'\leq r_0$ and write
$$
	p(r)
	=
	G_n(\widetilde{\bb\eta}_r)
	+
	rM(\widetilde{\bb\eta}_r)\,.
$$
Optimality of $\widetilde{\bb\eta}_r$, tested at
$\widetilde{\bb\eta}_{r'}$, gives
$$
	p(r')-p(r)
	\geq
	(r'-r)M(\widetilde{\bb\eta}_{r'})\,.
$$
By (i),
$$
	A_{\mathcal R}(\sigma,r)
	=
	-\sigma r\frac{\vnorm{\bh}_2}{\sqrt n}
	+
	p(r)\,.
$$
Using the Lipschitz property of $M$, part (ii), the inclusion
$\mathcal R\subseteq\mathcal R_0$, and the gap component of
$\mathcal E_n^{\mathrm{sr}}(d_r,a,K,C_L)$, we obtain
$$
\begin{aligned}
	A_{\mathcal R}(\sigma,r')
	-
	A_{\mathcal R}(\sigma,r)
	&\geq
	(r'-r)
	\left\{
		M(\widetilde{\bb\eta}_{r'})
		-
		\sigma\frac{\vnorm{\bh}_2}{\sqrt n}
	\right\}
	\\
	&\geq
	(r'-r)
	\left\{
		M(\bEtaYtil)
		-
		\Delta_{r'}
		-
		\sigma\frac{\vnorm{\bh}_2}{\sqrt n}
	\right\}
	\\
	&\geq
	(r'-r)
	\left\{
		M_n(\sigma,\bEtaYtil)
		-
		\sigma\frac{\vnorm{\bh}_2}{\sqrt n}
		-
		\frac{d_r}{4}
	\right\}
	\\
	&\geq
	\frac{3d_r}{4}(r'-r)\,,
\end{aligned}
$$
which is \eqref{eq:small_r_profile_increment}.

For fixed $\sigma$ and $\mathcal R$, continuity of $M$ and compactness of
$\mathcal B^n_{C_\eta}$ give
$$
	B_{\sigma,\mathcal R}
	=
	\underset{\bb\eta\in\mathcal B^n_{C_\eta}}{\sup}
	M_{\mathcal R}(\sigma,\bb\eta)
	<
	\infty\,.
$$
Value comparison in the constrained minimum yields, for
$r,r'\in[0,1]$,
$$
	\abs{
		A_{\mathcal R}(\sigma,r')
		-
		A_{\mathcal R}(\sigma,r)
	}
	\leq
	\left\{
		\sigma\frac{\vnorm{\bh}_2}{\sqrt n}
		+
		B_{\sigma,\mathcal R}
	\right\}
	\abs{r'-r}\,.
$$
Thus the maxima in \eqref{eq:small_r_profile_exclusion} are attained. By
\eqref{eq:small_r_profile_increment}, every $r<c$ satisfies
$A_{\mathcal R}(\sigma,r)\leq A_{\mathcal R}(\sigma,c)$, which proves
\eqref{eq:small_r_profile_exclusion}. The proof uses $C_\eta$ only through
$C_\eta\geq\sqrt{C_L}$, so every conclusion is preserved by enlargement of
$C_\eta$.
\end{proof}

\subsubsection{Upper bounds and Moreau-envelope feasibility}

\begin{lemma}[$t$-upper bound]
	\label{lemma:t_ub}
	Fix $c_r\in(0,1]$. Define the event
	$$
		\mathcal I_n^t
		=
		\left\{
			\vnorm{\bg}_2\leq2\sqrt n
		\right\}\,.
	$$
	Set
	$$
		K_t
		=
		\frac12\left(C_\theta+3C_\beta+2C_\sigma\right)^2
		+
		\frac12
		+
		\log 2\,,
	$$
	and choose
	$$
		C_t
		>
		\frac{2(K_t+\log 2)}{c_r}\,.
	$$
	Then, on $\mathcal I_n^t$, $\mathcal E_n^t(c_r,C_t)$
	holds, where
	$$
		\mathcal E_n^t(c_r,C_t)
		=
		\left\{
			\phi^{r,t}(c_r,C_t;\bg,\bh)
			=
			\phi^{r,t}(c_r,\infty;\bg,\bh)
			=
			\phi^r(c_r;\bg,\bh)
		\right\}\,.
	$$
More precisely, on $\mathcal I_n^t$, for every
$\sigma\in[0,C_\sigma]$ and every $r\in[c_r,1]$, the inner infimum in
\eqref{eq:psi-t-def}, before the outer $\sigma$-minimum and $r$-maximum,
is unchanged when $t\in(0,\infty)$ is replaced by
$t\in(0,C_t]$.
	Moreover, there exist constants $C,c>0$ and a $N\in\mathbb N$ such that, for all $n>N$,
	$$
		\Pr\left(
			\lnot\mathcal I_n^t
			\mid
			\mathcal G_n
		\right)
		\leq
		C\exp\{-cn\}\,.
	$$
	Consequently,
	$$
		\Pr\left(
			\lnot\mathcal E_n^t(c_r,C_t)
			\mid
			\mathcal G_n
		\right)
		\leq
		C\exp\{-cn\}\,.
	$$
	The constants $C,c,N$ are universal (they come from the Gaussian norm concentration for $\bg$ alone); $K_t$ depends only on $C_\theta,C_\beta,C_\sigma$.
\end{lemma}

\begin{proof}
	For fixed $\sigma\in[0,C_\sigma]$ and $r\in[c_r,1]$, consider the inner optimisation after the norm identity \eqref{eq:norm_variational_identity}:
	$$
		\underset{\substack{\bv\in V_p\\
		\abs{\theta}\leq C_\theta\\
		\bb\eta\in\mathcal B^n_{C_\eta}}}{\min}
		\underset{t>0}{\inf}
		\left[
			\frac{r}{2t}
			\frac1n
			\vnorm{
				\bb\eta-\theta\bb1-\frac1{\sqrt p}\bH_1\bv-\sigma\bg
			}_2^2
			+
			\frac{rt}{2}
			+
			\frac1n
			\left\{
				\bb1^\top\bb\rho(\bb\eta)-\bYtil^\top\bb\eta
			\right\}
		\right]\,.
	$$
	We first give a uniform upper bound on this value. Since $0\in V_p$, $0\in[-C_\theta,C_\theta]$, and $0\in\mathcal B^n_{C_\eta}$, we may test the feasible point
	$$
		\bv=0,\quad \theta=0,\quad \bb\eta=0,\quad t=1\,.
	$$
	At this point,
	$$
		\frac{r}{2t}
		\frac1n
		\vnorm{
			\bb\eta-\theta\bb1-\frac1{\sqrt p}\bH_1\bv-\sigma\bg
		}_2^2
		+
		\frac{rt}{2}
		+
		\frac1n
		\left\{
			\bb1^\top\bb\rho(\bb\eta)-\bYtil^\top\bb\eta
		\right\}
		=
		\frac r2\sigma^2\frac{\vnorm{\bg}_2^2}{n}
		+
		\frac r2
		+
		\log 2\,.
	$$
	On $\mathcal I_n^t$, this is bounded above by
	$$
		2C_\sigma^2
		+
		\frac12
		+
		\log 2
		\leq
		K_t\,,
	$$
	since $2C_\sigma^2=\tfrac12(2C_\sigma)^2\leq\tfrac12(C_\theta+3C_\beta+2C_\sigma)^2$.
	Thus the inner infimum is at most $K_t$, uniformly over $\sigma\in[0,C_\sigma]$ and $r\in[c_r,1]$.

	Next, for every $y\in(0,1)$ and every $x\in\Re$,
	$$
		\rho(x)-yx\geq
		\underset{z\in\Re}{\inf}\{\rho(z)-yz\}
		=
		-y\log y-(1-y)\log(1-y)
		\geq
		0
		\geq
		-\log 2\,.
	$$
	Therefore, for every $\bb\eta\in\Re^n$,
	$$
		\frac1n
		\left\{
			\bb1^\top\bb\rho(\bb\eta)-\bYtil^\top\bb\eta
		\right\}
		\geq
		-\log 2\,.
	$$
	Since the square term is nonnegative and $r\geq c_r$, any feasible point with $t>C_t$ has objective value at least
	$$
		\frac{c_rt}{2}-\log 2
		>
		K_t\,,
	$$
	by the choice of $C_t$. Such a point cannot attain the inner infimum, because the feasible point above has value at most $K_t$. Hence, for every $\sigma\in[0,C_\sigma]$ and $r\in[c_r,1]$, the infimum over $t>0$ is unchanged when restricted to $t\in(0,C_t]$.

	It follows that
	$$
		\phi^{r,t}(c_r,C_t;\bg,\bh)
		=
		\phi^{r,t}(c_r,\infty;\bg,\bh)
		=
		\phi^r(c_r;\bg,\bh)\,,
	$$
	on $\mathcal I_n^t$, which proves the deterministic implication.

	It remains to bound the probability of $\mathcal I_n^t$. Conditional on $\mathcal G_n$, the vector $\bg$ is independent of $\mathcal G_n$ and has i.i.d. standard normal entries. Hence $\vnorm{\bg}_2^2\sim\chi_n^2$, and the standard chi-square concentration bound
	(see, for example, \citealt[][Theorem~3.1.1]{vershynin:2018}) gives
	$$
		\Pr\left(
			\vnorm{\bg}_2>2\sqrt n
			\mid
			\mathcal G_n
		\right)
		\leq
		C\exp\{-cn\}\,,
	$$
	for all sufficiently large $n$. This proves the conditional bound.
\end{proof}

\begin{lemma}[Moreau feasibility]
	\label{lemma:moreau_feasibility}
	Fix $c_r\in(0,1]$ and $C_t\in(0,\infty)$, and set $L=C_t/c_r$. Let $C_M<\infty$ and define the centre-control event
	$$
		\mathcal E_n^{\mathrm{ctr}}(C_M)
		=
			\underset{\substack{
			\abs{\theta}\leq C_\theta\\
			\bv\in V_p\\
			\sigma\in[0,C_\sigma]\\
			r\in[c_r,1]\\
			t\in(0,C_t]
			}}{\sup} \, 
			\frac1{\sqrt n}
			\vnorm{
				\theta\bb1
				+
				\frac1{\sqrt p}\bH_1\bv
				+
				\sigma\bg
				+
				\frac tr\bYtil
			}_2
			\leq C_M
		\,.
	$$
If $C_\eta\geq C_M+L$, then, on
$\mathcal E_n^{\mathrm{ctr}}(C_M)$, for every admissible
$(\theta,\bv,\sigma,r,t)$, set
$$
    \lambda=\frac tr,
    \quad
    \bb b
    =
    \theta\bb1
    +
    \frac1{\sqrt p}\bH_1\bv
    +
    \sigma\bg
    +
    \lambda\bYtil \,.
$$
The unconstrained coordinatewise Moreau minimiser belongs to
$\mathcal B^n_{C_\eta}$ and
$$
    M_\rho^{\mathrm{const}}(\bb b,\lambda)
    =
    M_\rho(\bb b,\lambda) \,.
$$
Consequently,
$$
    \mathcal E_n^{\mathrm{ctr}}(C_M)
    \Longrightarrow
    \mathcal E_n^{\mathrm M}(c_r,C_t,C_\eta) \,.
$$
\end{lemma}
\begin{proof}
	Fix $\abs{\theta}\leq C_\theta$, $\bv\in V_p$, $\sigma\in[0,C_\sigma]$, $r\in[c_r,1]$, and $t\in(0,C_t]$,
	and set
	$$
		\lambda=\frac tr,
		\quad
		\bb b=
		\theta\bb1
		+
		\frac1{\sqrt p}\bH_1\bv
		+
		\sigma\bg
		+
		\lambda\bYtil \,.
	$$
	Then $\lambda\in(0,L]$. For each coordinate, let
	$$
		\eta_i^*
		=
		\underset{\eta\in\Re}{\arg\min}
		\left\{
			\rho(\eta)
			+
			\frac{1}{2\lambda}(\eta-b_i)^2
		\right\} \,.
	$$
	Since the objective is strictly convex and differentiable, $\eta_i^*$ is unique and satisfies
	$$
		\rho'(\eta_i^*)
		+
		\frac1\lambda(\eta_i^*-b_i)
		=
		0 \,.
	$$
	Therefore,
	$$
		\eta_i^*
		=
		b_i-\lambda\rho'(\eta_i^*)\,.
	$$
	Because $0<\rho'(\eta_i^*)<1$, we have
	$$
		\abs{\eta_i^*}
		\leq
		\abs{b_i}+\lambda \,.
	$$
	Thus, for $\bb\eta^*=(\eta_1^*,\ldots,\eta_n^*)^\top$,
	$$
		\frac{\vnorm{\bb\eta^*}_2}{\sqrt n}
		\leq
		\frac{\vnorm{\bb b}_2}{\sqrt n}
		+
		\lambda
		\leq
		C_M+L
		\leq
		C_\eta\,,
	$$
	on $\mathcal E_n^{\mathrm{ctr}}(C_M)$. Hence the unconstrained coordinatewise Moreau minimiser is feasible for the constraint $\bb\eta\in\mathcal B^n_{C_\eta}\,.$
	Therefore,
	$$
		M_\rho^{\mathrm{const}}(\bb b,\lambda)
		=
		M_\rho(\bb b,\lambda)\,,
	$$
	for every point in the scalar AO domain. Substituting this equality into the definitions of
	$\phi^{M,\mathrm{const}}(c_r,C_t;\bg,\bh)$ and $\phi^M(c_r,C_t;\bg,\bh)$ gives
	$$
		\phi^{M,\mathrm{const}}(c_r,C_t;\bg,\bh)
		=
		\phi^M(c_r,C_t;\bg,\bh)\,.
	$$
	This is exactly $\mathcal E_n^{\mathrm M}(c_r,C_t,C_\eta)$.
\end{proof}

\begin{lemma}[Probability of the Moreau centre-control event]
	\label{lemma:moreau_center_probability}
	Fix $c_r\in(0,1]$ and $C_t\in(0,\infty)$, and set
	$$
		L=\frac{C_t}{c_r}\,.
	$$
	There exist constants $C_M<\infty$, $C,c>0$ and $N\in\mathbb N$ such that, for all $n>N$,
	$$
		\Pr\left(
			\lnot\mathcal E_n^{\mathrm{ctr}}(C_M)
			\mid
			\mathcal G_n
		\right)
		\leq
		C\exp\{-cn\}\,.
	$$
	The constants depend only on $C_\beta,C_\theta,C_\sigma,c_r,C_t$
	and universal Gaussian concentration constants. In particular, they do not depend on the realised values of $\mathcal G_n$.
\end{lemma}

\begin{proof}
	For every admissible $(\theta,\bv,\sigma,r,t)$, we have
	$$
		\begin{aligned}
		&\frac1{\sqrt n}
		\vnorm{
			\theta\bb1
			+
			\frac1{\sqrt p}\bH_1\bv
			+
			\sigma\bg
			+
			\frac tr\bYtil
		}_2\leq
		\abs{\theta}
		+
		\frac1{\sqrt n}
		\vnorm{
			\frac1{\sqrt p}\bH_1\bv
		}_2
		+
		\sigma\frac{\vnorm{\bg}_2}{\sqrt n}
		+
		\frac tr\frac{\vnorm{\bYtil}_2}{\sqrt n}\,.
		\end{aligned}
	$$
		Since $0<\Ytil_i<1$, it holds that $\vnorm{\bYtil}_2/\sqrt n\leq 1\,.$
		Moreover, $\abs{\theta}\leq C_\theta,\quad \sigma\leq C_\sigma,\quad t/r\leq L\,.$
	It remains to control the signal term and the fresh Gaussian term.

	By \eqref{eq:ao_H_repar}, on the common extension,
	$\bH_1/\sqrt p=\bb Z\bb\Gamma_p^{1/2}$ almost surely, where
	$\bb Z\in\Re^{n\times2}$ has i.i.d.\ $\mathrm N(0,1)$ entries
	conditional on $\mathcal G_n$. For $\bv\in V_p$, set
	$\bu=\bb\Gamma_p^{1/2}\bv\,.$
		Then $\vnorm{\bu}_2^2 = \bv^\top\bb\Gamma_p\bv = \vnorm{\bb B\bv}_2^2/p \leq C_\beta^2\,.$
	Therefore, on the common extension,
	$$
		\underset{\bv\in V_p}{\sup}
		\frac1{\sqrt n}
		\vnorm{
			\frac1{\sqrt p}\bH_1\bv
		}_2
		=
		\underset{\bu\in U\cap\range(\bb\Gamma_p)}{\sup}
		\frac1{\sqrt n}
		\vnorm{\bb Z\bu}_2
		\leq
		C_\beta\frac{\mnorm{\bb Z}_2}{\sqrt n}\,.
	$$
	By the standard Gaussian operator-norm bound for the $n\times2$ matrix $\bb Z$
	(see, for example, \citealt[][Theorem~2.6]{rudelson+vershynin:2010}),
	there exist constants $C,c>0$ and $N_1\in\mathbb N$ such that, for all $n>N_1$,
	$$
		\Pr\left(
			\frac{\mnorm{\bb Z}_2}{\sqrt n}>3
			\mid
			\mathcal G_n
		\right)
		\leq
		C\exp\{-cn\}\,.
	$$
	Similarly, since $\bg\sim\mathrm N(\bb0_n,\bb I_n)$ independently of $\mathcal G_n$, there exist constants $C,c>0$ and $N_2\in\mathbb N$ such that, for all $n>N_2$,
	$$
		\Pr\left(
			\frac{\vnorm{\bg}_2}{\sqrt n}>2
			\mid
			\mathcal G_n
		\right)
		\leq
		C\exp\{-cn\}\,.
	$$
	On the intersection of the two events
	$$
		\left\{
			\frac{\mnorm{\bb Z}_2}{\sqrt n}\leq 3
		\right\}
		\cap
		\left\{
			\frac{\vnorm{\bg}_2}{\sqrt n}\leq 2
		\right\}\,,
	$$
	we have
	$$
		\underset{\substack{
			\abs{\theta}\leq C_\theta\\
			\bv\in V_p\\
			\sigma\in[0,C_\sigma]\\
			r\in[c_r,1]\\
			t\in(0,C_t]
		}}{\sup}
		\frac1{\sqrt n}
		\vnorm{
			\theta\bb1
			+
			\frac1{\sqrt p}\bH_1\bv
			+
			\sigma\bg
			+
			\frac tr\bYtil
		}_2
		\leq
		C_\theta+3C_\beta+2C_\sigma+L\,.
	$$
	Thus take $C_M=C_\theta+3C_\beta+2C_\sigma+L\,.$
	A union bound gives
	$$
		\Pr\left(
			\lnot\mathcal E_n^{\mathrm{ctr}}(C_M)
			\mid
			\mathcal G_n
		\right)
		\leq
		C\exp\{-cn\}\,,
	$$
	after adjusting $C,c,N$.
\end{proof}

\subsection{Lemmas for Section~\ref{sec:ao-limit}}

\subsubsection{Nullspace and moving-domain control}
\label{sec:nullspace-moving-domains}


	For $y\in(0,1)$ and $\lambda\geq0$ write
	$$
		e_\lambda(x,y)
		=
		m_\rho(x+\lambda y,\lambda)-yx-\frac{\lambda}{2}y^2\,,
		\qquad x\in\Re\,,
	$$
	with $m_\rho(x,0)=\rho(x)$; this is the function used again in Lemma~\ref{lemma:limiting_AO_geometry}, and for $\lambda=t/r$ the expectation term of $\bar F_n$ (conditionally on $\mathcal G_n$) and of $F$ is $\expect[e_\lambda(\theta+\bz^\top\bu+\sigma G,\Ytil^{(n)})\mid\mathcal G_n]$ and $\expect[e_\lambda(\theta+\bz^\top\bu+\sigma G,\Ytil)]$, respectively.
	For each $y$, $x\mapsto e_\lambda(x,y)$ is strictly convex: for $\lambda=0$, $e_0(x,y)=\rho(x)-yx$ and $\partial_x^2e_0(x,y)=\rho''(x)>0$; for $\lambda>0$, the Moreau-envelope derivative formula of \citet[Lemma~2, equation~(29)]{salehi+et+al:2019} and the proximal-map derivative formula of \citet[Proposition~6.3]{donoho+montanari:2016} give $\partial_x^2e_\lambda(x,y)=\rho''(\eta)/\{1+\lambda\rho''(\eta)\}>0$ with $\eta=\prox{\lambda\rho}{x+\lambda y}$. In particular $e_\lambda(\cdot,y)$ is convex (see also \citealt[][Chapter~6]{beck:2017}).

	\begin{lemma}[Nullspace never helps for the conditional mean objective]
		\label{lemma:nullspace_no_help_barFn}
		For every $\bb\omega=(\sigma,r,\bu,\theta,t)\in\bb\Omega$, let $\bb\omega_p = (\sigma,r,\bb P_{\bb \Gamma_p}\bu,\theta,t)\,.$
		Then, conditional on $\mathcal G_n$, $\bar F_n(\bb\omega_p) \leq \bar F_n(\bb\omega)\,.$
		Consequently, $\operatorname{val}(\bar F_n;\bb\Omega_n) = \operatorname{val}(\bar F_n;\bb\Omega)\,.$
	\end{lemma}

	\begin{proof}
		Fix $\sigma,r,\theta,t$ and write $\lambda=t/r\,.$
	Let $\bu_R=\bb P_{\bb \Gamma_p}\bu, \quad \bu_N=\bu-\bu_R\,.$
	Then $\bu_R\in\range(\bb\Gamma_p)$ and $\bu_N\in\ker(\bb\Gamma_p)$. Since $\bb P_{\bb \Gamma_p}$ is an orthogonal projector, $\vnorm{\bu_R}_2\leq\vnorm{\bu}_2\,,$
	so if $\bu\in U$ then $\bu_R\in U\cap\range(\bb\Gamma_p)$.

	Conditional on $\mathcal G_n$, let $\bz\sim\mathrm N(\bb0_2,\bb I_2)$.
	Because
	$\range(\bb\Gamma_p)=\range(\bb\Gamma_p^{1/2})$, there exists
	$\bb a\in\Re^2$ such that
	$\bu_R=\bb\Gamma_p^{1/2}\bb a$. Hence
	$\bz^\top\bu_R=\bb a^\top\bb\Gamma_p^{1/2}\bz$ is measurable with respect to
	$\bb\Gamma_p^{1/2}\bz$.
	The pseudo-response variable $\Ytil^{(n)}$ in $\bar F_n$ is a measurable function of $(\bb\Gamma_p^{1/2}\bz,\varepsilon)$ and of the fixed parameters $(\theta_0,\theta_P)$. Since $\bb\Gamma_p^{1/2}\bu_N=\bb0$, the scalar $\bz^\top\bu_N$ is independent of $(\bb\Gamma_p^{1/2}\bz,\Ytil^{(n)},G)$ and has mean zero conditional on $\mathcal G_n$.

		Write $\xi_R=\theta+\bz^\top\bu_R+\sigma G$. Conditioning on
	$(\bb\Gamma_p^{1/2}\bz,\Ytil^{(n)},G)$, under which $\xi_R$ and $\Ytil^{(n)}$ are fixed, convexity of $e_\lambda(\cdot,\Ytil^{(n)})$ and Jensen's inequality \citep[Chapter~2.7]{shiryaev:2016} give
	$$
		\expect\left[ e_\lambda \left( \xi_R+\bz^\top\bu_N,\Ytil^{(n)} \right) \mid \bb\Gamma_p^{1/2}\bz,\Ytil^{(n)},G,\mathcal G_n \right] \quad\geq e_\lambda \left( \xi_R,\Ytil^{(n)} \right)\,.
	$$
		Taking expectations gives $\bar F_n(\sigma,r,\bu_R,\theta,t) \leq \bar F_n(\sigma,r,\bu,\theta,t)\,.$
		The terms $-\sigma r\mu_{\bh}/\sqrt n+rt/2$, with $\mu_{\bh} = \expect[\vnorm{\bh}_2 \mid \mathcal G_n]$,
		do not depend on $\bu$, so they are unaffected.

		It remains to pass from the pointwise inequality to the value. Since $\bb\Omega_n\subseteq\bb\Omega$, we always have $\operatorname{val}(\bar F_n;\bb\Omega) \leq \operatorname{val}(\bar F_n;\bb\Omega_n)\,,$
		because the innermost operation over $\bu$ is a minimisation and $\bb\Omega$ has the larger $\bu$-domain. Conversely, the projection argument above shows that every $\bu\in U$ can be replaced by $\bb P_{\bb \Gamma_p}\bu\in U\cap\range(\bb\Gamma_p)$ without increasing $\bar F_n$. Hence enlarging the $\bu$-domain from $\bb\Omega_n$ to $\bb\Omega$ cannot lower the value. Thus,
		$\operatorname{val}(\bar F_n;\bb\Omega_n) = \operatorname{val}(\bar F_n;\bb\Omega)\,.$
	\end{proof}

	\begin{lemma}[Strict nullspace penalty for the limiting objective]
	\label{lemma:nullspace_no_help_F}
		Let $\bb P_{\bb \Gamma}$ be the orthogonal projector onto $\range(\bb\Gamma)$.
		For $\bu\in U$, write $\bu_R=\bb P_{\bb \Gamma}\bu, \quad \bu_N=\bu-\bu_R\,.$
		Then, for every $(\sigma,r,\bu,\theta,t)\in\bb\Omega\,,$
		we have $F(\sigma,r,\bu_R,\theta,t) \leq F(\sigma,r,\bu,\theta,t)\,.$
		Moreover, if $\bu_N\neq0$, then the inequality is strict: $F(\sigma,r,\bu_R,\theta,t) < F(\sigma,r,\bu,\theta,t)\,.$
		Consequently, $\operatorname{val}(F;\bb\Omega) = \operatorname{val}(F;\bb\Omega^*)\,.$
		In particular, for every fixed $(\sigma,r)$, every minimiser of $(\bu,\theta,t) \mapsto F(\sigma,r,\bu,\theta,t)$
		over $U\times[-C_\theta,C_\theta]\times[0,C_t]$
		has $\bu\in\range(\bb\Gamma)$.
	\end{lemma}

	\begin{proof}
		Fix $(\sigma,r,\bu,\theta,t)\in\bb\Omega$ and set $\lambda=t/r\,.$
		Since $r\in[c_r,1]$, we have $\lambda\geq0$. Also,
		$\bu_R\in\range(\bb\Gamma), \quad \bu_N\in\ker(\bb\Gamma), \quad \vnorm{\bu_R}_2\leq\vnorm{\bu}_2\,.$
	Thus $\bu_R\in U\cap\range(\bb\Gamma)$.

	Let
	$$
		\bb Q
		=
		\bb\Gamma^{1/2}\bz,
		\quad
		\bz
		\sim
		\mathrm N(\bb0_2,\bb I_2)\,.
	$$
	The limiting pseudo-response $\Ytil$ is a measurable function of
	$(\bb Q,\varepsilon)$ and of the fixed limiting parameters
	$(\theta_0^*,\theta_P^*)$. Since $\bu_N\in\ker(\bb\Gamma)$,
	$\bb\Gamma^{1/2}\bu_N=\bb0_2\,.$ Hence
	$$
		\operatorname{Cov}
		\left(
			\bz^\top\bu_N,
			\bb Q
		\right)
		=
		\bb\Gamma^{1/2}\bu_N
		=
		\bb0_2\,.
	$$
	Because all variables are jointly Gaussian before applying the response
	map, $\bz^\top\bu_N$ is independent of $\bb Q$. It is also independent of
	$\varepsilon$ and $G$. Therefore
	$$
		\bz^\top\bu_N
		\quad\text{is independent of}\quad
		(\bb Q,\Ytil,G)\,,
	$$
	and
	$$
		\expect[\bz^\top\bu_N]
		=
		0\,.
	$$
	If $\bu_N\neq\bb0$, then
	$$
		\bz^\top\bu_N
		\sim
		\mathrm N
		\left(
			0,
			\vnorm{\bu_N}_2^2
		\right)\,,
	$$
	is nondegenerate.

		Recall that $x\mapsto e_\lambda(x,y)$ is strictly convex for every $y\in(0,1)$ and $\lambda\geq0$ (see the paragraph preceding Lemma~\ref{lemma:nullspace_no_help_barFn}).

	Because $\bu_R\in\range(\bb\Gamma)=\range(\bb\Gamma^{1/2})$, there
	exists $a\in\mathbb R^2$ such that $\bu_R=\bb\Gamma^{1/2}a\,.$
	Hence
	$
		\bz^\top\bu_R
		=
		a^\top\bb\Gamma^{1/2}\bz
		=
		a^\top\bb Q
	$,
	so $\bz^\top\bu_R$ is measurable with respect to $\bb Q$. Write $\xi_R=\theta+\bz^\top\bu_R+\sigma G$. Conditionally on
	$(\bb Q,\Ytil,G)$, the variable $\xi_R$ is fixed, while $\bz^\top\bu_N$ is independent of
	$(\bb Q,\Ytil,G)$ and has mean zero. Jensen's inequality gives
	$$
		\expect
		\left[
			e_\lambda(\xi_R+\bz^\top\bu_N,\Ytil)
			\mid \bb Q,\Ytil,G
		\right]
		\geq
		e_\lambda(\xi_R,\Ytil)\,.
	$$
	If $\bu_N\neq0$, then $\bz^\top\bu_N$ is nondegenerate and
	$e_\lambda(\cdot,\Ytil)$ is strictly convex, so the inequality is strict almost
	surely. Taking expectations gives
	$$
		\expect
		\left[
			e_\lambda
			\left(
				\theta+\bz^\top\bu+\sigma G,\Ytil
			\right)
		\right]
		\geq
		\expect
		\left[
			e_\lambda
			\left(
				\theta+\bz^\top\bu_R+\sigma G,\Ytil
			\right)
		\right]\,,
	$$
	with strict inequality if $\bu_N\neq0$.

		By the definition of $e_\lambda$, this is exactly the comparison of
		the expectation terms in $F$. The remaining terms $-\sigma r\sqrt\kappa+rt/2$
		do not depend on $\bu$. Therefore $F(\sigma,r,\bu_R,\theta,t) \leq F(\sigma,r,\bu,\theta,t)\,,$
	with strict inequality whenever $\bu_N\neq0$.

	It remains to pass from the pointwise comparison to the value identity.
	For every fixed $(\sigma,r)$, the inner minimum over the full signal ball
	is no larger than the inner minimum over its intersection with
	$\range(\bb\Gamma)$. Conversely, projecting every feasible $\bu$ onto
	$\range(\bb\Gamma)$ preserves feasibility and cannot increase $F$, so
	the two inner minima are equal. Taking the remaining operations in the
	profiled order gives
	$\operatorname{val}(F;\bb\Omega)
	=
	\operatorname{val}(F;\bb\Omega^*)$.

	Finally, suppose that $(\tilde\bu,\tilde\theta,\tilde t)$ minimises
	$F(\sigma,r,\cdot)$ over the ambient inner domain and that
	$\tilde\bu\notin\range(\bb\Gamma)$. Then $\tilde\bu_N=\tilde\bu-\bb P_{\bb \Gamma}\tilde\bu\neq0\,,$
		so the strict part gives $F(\sigma,r,\bb P_{\bb \Gamma}\tilde\bu,\tilde\theta,\tilde t) < F(\sigma,r,\tilde\bu,\tilde\theta,\tilde t)\,,$
	contradicting minimality. Hence every ambient inner minimiser has
	$\tilde\bu\in\range(\bb\Gamma)$.
\end{proof}

These deterministic comparisons are combined with the
uniform-convergence ingredients below in
Proposition~\ref{prop:uniform_AO_control}.

\subsubsection{Uniform-convergence ingredients}

\begin{lemma}[Conditional pointwise concentration of $F_n$]
	\label{lemma:F_n_pointwise_conc}
	Fix $\bb\omega=(\sigma,r,\bu,\theta,t)\in\bb\Omega$, and let $F_n(\bb\omega)$ and $\bar F_n(\bb\omega)$ be defined in \eqref{eq:F_n_def} and \eqref{eq:F_bar_n_def}. Then there exist constants $C,c>0$ such that, for all $s>0$ and all $n$ sufficiently large,
	$$
		\Pr\left(
			\abs{
				F_n(\bb\omega)-\bar F_n(\bb\omega)
			}>s
			\mid
			\mathcal G_n
		\right)
		\leq
		C\exp\{-cns^2\}\,.
	$$
	The constants depend only on $C_\sigma,C_\beta,C_\theta,C_t,c_r$
	and on universal Gaussian concentration constants. In particular, they do not depend on the realised value of $\mathcal G_n$.
\end{lemma}
	\begin{proof}
		Let $\lambda=t/r\,.$
	Since $\bb\omega\in\bb\Omega$, we have $0\leq \sigma\leq C_\sigma$, $c_r\leq r\leq 1$, $\vnorm{\bu}_2\leq C_\beta$, $\abs{\theta}\leq C_\theta$, and $0\leq t\leq C_t\,.$
		Thus, $0\leq \lambda\leq L, \quad L=C_t/c_r\,.$

		Conditional on $\mathcal G_n$, the rows $(\bz_i,g_i,\varepsilon_i)$ are independent and identically distributed, and are independent of $\bh$. Define $\xi_i=\theta+\bz_i^\top\bu+\sigma g_i$
		and $f_i(\bb\omega) = m_\rho(\xi_i+\lambda \Ytil_i,\lambda) - \Ytil_i\xi_i - (\lambda/2)\Ytil_i^2\,.$
	Then
	$$
		F_n(\bb\omega)
		=
		\frac1n\sum_{i=1}^n f_i(\bb\omega)
		-
		\sigma r\frac{\vnorm{\bh}_2}{\sqrt n}
		+
		\frac{rt}{2}\,,
	$$
	and $\bar F_n(\bb\omega) = \expect[F_n(\bb\omega)\mid\mathcal G_n]\,.$
	Hence
	$$
		F_n(\bb\omega)-\bar F_n(\bb\omega)
		=
		\frac1n\sum_{i=1}^n
		\left\{
			f_i(\bb\omega)
			-
			\expect[f_i(\bb\omega)\mid\mathcal G_n]
		\right\}
		-
		\sigma r
		\frac{
			\vnorm{\bh}_2-\mu_{\bh}
		}{\sqrt n}\,,
	$$
	where $\mu_{\bh} = \expect[\vnorm{\bh}_2 \mid \mathcal G_n]$

	We first bound the empirical average. Since $0<\Ytil_i<1$ and $0\leq\lambda\leq L$, and since
	$$
		0\leq m_\rho(b,\lambda)\leq \rho(b)\leq \log 2+\abs b\,,
	$$
	for all $b\in\Re$ and $\lambda\geq0$, we have
	$$
		\begin{aligned}
		\abs{f_i(\bb\omega)}
		&\leq
		\log 2+\abs{\xi_i+\lambda \Ytil_i}
		+
		\abs{\Ytil_i\xi_i}
		+
		\frac{\lambda}{2}\Ytil_i^2\leq
		\log 2
		+
		2\abs{\xi_i}
		+
		\frac{3L}{2}\leq
		K_f
		\left(
			1+\vnorm{\bz_i}_2+\abs{g_i}
		\right)\,,
		\end{aligned}
	$$
	where $K_f<\infty$ depends only on $C_\sigma,C_\beta,C_\theta,C_t,c_r\,.$
		Thus, conditional on $\mathcal G_n$,
	$$
		\vnorm{f_i(\bb\omega)}_{\psi_2}\leq K\,,
	$$
	for a deterministic constant $K<\infty$ with the same parameter dependence. Centering changes the subgaussian norm by at most an absolute factor, so
	$$
		\vnorm{
			f_i(\bb\omega)
			-
			\expect[f_i(\bb\omega)\mid\mathcal G_n]
		}_{\psi_2}
		\leq K\,,
	$$
	after increasing $K$. Since the summands are independent conditional on $\mathcal G_n$, the general Hoeffding inequality \citep[Theorem~2.6.3]{vershynin:2018} gives
	$$
		\Pr\left(
			\abs{
				\frac1n\sum_{i=1}^n
				\left\{
					f_i(\bb\omega)
					-
					\expect[f_i(\bb\omega)\mid\mathcal G_n]
				\right\}
			}
			>\frac{s}{2}
			\mid
			\mathcal G_n
		\right)
		\leq
		2\exp\{-c_1ns^2\}\,,
	$$
	for a constant $c_1>0$ depending only on $K$.

	It remains to control the $\bh$ term. Since $\sigma r\leq C_\sigma$, we have
	$$
		\Pr\left(
			\sigma r
			\frac{
				\abs{\vnorm{\bh}_2-\mu_{\bh}}
			}{\sqrt n}
			>
			\frac{s}{2}
			\mid
			\mathcal G_n
		\right)
		\leq
		\Pr\left(
			\abs{\vnorm{\bh}_2-\mu_{\bh}}
			>
			\frac{s\sqrt n}{2C_\sigma}
			\mid
			\mathcal G_n
		\right)\,.
	$$
	Conditional on $\mathcal G_n$, Gaussian concentration for Lipschitz
	functions \citep[Theorem~5.2.2]{vershynin:2018} gives
	$$
		\Pr\left(
			\abs{\vnorm{\bh}_2-\mu_{\bh}}
			>
			u
			\mid
			\mathcal G_n
		\right)
		\leq
		2\exp\{-cu^2\}\,,
	$$
		for a universal constant $c>0$. Taking $u=s\sqrt n/(2C_\sigma)$
	gives
	$$
		\Pr\left(
			\sigma r
			\frac{
				\abs{\vnorm{\bh}_2-\mu_{\bh}}
			}{\sqrt n}
			>
			\frac{s}{2}
			\mid
			\mathcal G_n
		\right)
		\leq
		2\exp\{-c_2ns^2\}\,.
	$$
	A union bound yields
	$$
		\Pr\left(
			\abs{
				F_n(\bb\omega)-\bar F_n(\bb\omega)
			}>s
			\mid
			\mathcal G_n
		\right)
		\leq
		C\exp\{-cns^2\}\,,
	$$
	after adjusting $C,c>0$.
\end{proof}

\begin{lemma}[Structural bias of the conditional AO mean]
	\label{lemma:F_n_structural_bias}
	Let $\bar F_n(\bb\omega)$ and $F(\bb\omega)$ be defined in \eqref{eq:F_bar_n_def} and \eqref{eq:F_limit_def}.
		There exists a constant $K<\infty$, depending only on
		$\kappa,\alpha,\bb\Gamma,\theta_0^*,\theta_P^*, C_\sigma,C_\beta,C_\theta,C_t,c_r$
	such that, for all $t_0\in(0,1]$ and all sufficiently large $n$,
	$$
		\mathds{1}_{\mathcal S_n(t_0)}
		\underset{\bb\omega\in\bb\Omega}{\sup}
		\abs{
			\bar F_n(\bb\omega)-F(\bb\omega)
		}
		\leq
		K
		\left(
			t_0+\sqrt{t_0}
			+
			\abs{
				\frac{\mu_{\bh}}{\sqrt n}
				-
				\sqrt\kappa
			}
		\right), \quad \mu_{\bh} = \expect[\vnorm{\bh}_2 \mid \mathcal G_n] \,.
	$$
	Consequently, for every $s>0$, there exist $t_0=t_0(s)>0$ and $N_s\in\mathbb N$ such that, for all $n>N_s$,
	$$
		\mathds{1}_{\mathcal S_n(t_0)}
		\underset{\bb\omega\in\bb\Omega}{\sup}
		\abs{
			\bar F_n(\bb\omega)-F(\bb\omega)
		}
		\leq
		\frac{s}{3} \,.
	$$
\end{lemma}

\begin{proof}
		Fix $\bb\omega=(\sigma,r,\bu,\theta,t)\in\bb\Omega$
		and set $\lambda=t/r\,.$
		Since $\bb\omega\in\bb\Omega$, we have $0\leq \lambda\leq L, \quad L=C_t/c_r\,.$
		Let $\bb q_n=\bb\Gamma_p^{1/2}\bz, \quad \bb Q=\bb\Gamma^{1/2}\bz, \quad \bz\sim\mathrm N(\bb0_2,\bb I_2)\,,$
	and define
	$$
		\Ytil^{(n)}
		=
		\alpha\mathds{1}
		\left\{
			\varepsilon<\rho'(\theta_0+q_{n,1})
		\right\}
		+
		(1-\alpha)\rho'(\theta_P+q_{n,2})\,,
	$$
	and
	$$
		\Ytil
		=
		\alpha\mathds{1}
		\left\{
			\varepsilon<\rho'(\theta_0^*+Q_1)
		\right\}
		+
		(1-\alpha)\rho'(\theta_P^*+Q_2) \,.
	$$
		Also write $\xi=\theta+\bz^\top\bu+\sigma G, \quad G\sim\mathrm N(0,1)\,,$
	where $(\bz,G,\varepsilon)$ are mutually independent. Then
	$$
		\begin{aligned}
		\bar F_n(\bb\omega)-F(\bb\omega)
		=
		\expect\left[
			m_\rho(\xi+\lambda \Ytil^{(n)},\lambda)
			-
			m_\rho(\xi+\lambda \Ytil,\lambda)
		\right]
		-
		\expect\left[
			(\Ytil^{(n)}-\Ytil)\xi
		\right]
		-
		\frac{\lambda}{2}
		\expect\left[
			(\Ytil^{(n)})^2-\Ytil^2
		\right]
		-
		\sigma r
		\left(
			\frac{\mu_{\bh}}{\sqrt n}
			-
			\sqrt\kappa
		\right) \,.
		\end{aligned}
	$$
	We control the difference between $\Ytil^{(n)}$ and $\Ytil$. On $\mathcal S_n(t_0)$, $\abs{\theta_0-\theta_0^*}\leq t_0$, $\abs{\theta_P-\theta_P^*}\leq t_0$, and $\mnorm{\bb\Gamma_p-\bb\Gamma}_2\leq t_0\,.$
	By the H\"older continuity of the positive semidefinite square-root map
		\citep[][Theorem~X.1.1]{bhatia:1997},
	$$
		\mnorm{\bb\Gamma_p^{1/2}-\bb\Gamma^{1/2}}_2
		\leq
		\mnorm{\bb\Gamma_p-\bb\Gamma}_2^{1/2}
		\leq
		\sqrt{t_0} \,.
	$$
	Using $\underset{x}{\sup}\rho''(x)\leq 1/4$, we have
	$$
		\abs{
			\rho'(\theta_P+q_{n,2})
			-
			\rho'(\theta_P^*+Q_2)
		}
		\leq
		\frac14
		\left(
			t_0
			+
			\sqrt{t_0}\vnorm{\bz}_2
		\right) \,.
	$$
	For the Bernoulli part, conditional on $\bz$, the uniform variable $\varepsilon$ gives
	$$
		\begin{aligned}
		&\expect_\varepsilon
		\left[
			\abs{
				\mathds{1}\{
					\varepsilon<\rho'(\theta_0+q_{n,1})
				\}
				-
				\mathds{1}\{
					\varepsilon<\rho'(\theta_0^*+Q_1)
				\}
			}
			\mid \bz
		\right]=
		\abs{
			\rho'(\theta_0+q_{n,1})
			-
			\rho'(\theta_0^*+Q_1)
		}\leq
		\frac14
		\left(
			t_0
			+
			\sqrt{t_0}\vnorm{\bz}_2
		\right) \,.
		\end{aligned}
	$$
	Combining the two bounds, there is a constant $K_1<\infty$, depending only on $\alpha$, such that
	$$
		\expect_\varepsilon \left[
			\abs{\Ytil^{(n)}-\Ytil}
			\mid
			\bz
		\right]
		\leq
		K_1
		\left(
			t_0+\sqrt{t_0}\vnorm{\bz}_2
		\right) \,.
	$$
	Therefore,
	$$
		\expect\left[
			\abs{\Ytil^{(n)}-\Ytil}
		\right]
		\leq
		K_2(t_0+\sqrt{t_0}) \,.
	$$
	Moreover, since
	$$
		\abs{\xi}
		\leq
		C_\theta+C_\beta\vnorm{\bz}_2+C_\sigma\abs{G}\,,
	$$
	we also get
	$$
		\expect\left[
			(1+\abs{\xi})\abs{\Ytil^{(n)}-\Ytil}
		\right]
		\leq
		K_3(t_0+\sqrt{t_0}) \,.
	$$
	Next, by Lemma~\ref{lemma:moreau_lipschitz}, the map $b\mapsto m_\rho(b,\lambda)$ is $1$-Lipschitz for every $\lambda\geq0$. Hence
	$$
		\abs{
			m_\rho(\xi+\lambda \Ytil^{(n)},\lambda)
			-
			m_\rho(\xi+\lambda \Ytil,\lambda)
		}
		\leq
		\lambda\abs{\Ytil^{(n)}-\Ytil}
		\leq
		L\abs{\Ytil^{(n)}-\Ytil} \,.
	$$
	Also, since $0<\Ytil^{(n)},\Ytil<1$, $\abs{(\Ytil^{(n)})^2-\Ytil^2} \leq 2\abs{\Ytil^{(n)}-\Ytil} \,.$
	Thus,
	$$
		\begin{aligned}
		\abs{\bar F_n(\bb\omega)-F(\bb\omega)}
		&\leq
		L\expect\left[\abs{\Ytil^{(n)}-\Ytil}\right]
		+
		\expect\left[\abs{\xi}\abs{\Ytil^{(n)}-\Ytil}\right]
		+
		L\expect\left[\abs{\Ytil^{(n)}-\Ytil}\right]+
		C_\sigma
		\abs{
			\frac{\mu_{\bh}}{\sqrt n}
			-
			\sqrt\kappa
		}\\
		&\leq
		K
		\left(
			t_0+\sqrt{t_0}
			+
			\abs{
				\frac{\mu_{\bh}}{\sqrt n}
				-
				\sqrt\kappa
			}
		\right) \,.
		\end{aligned}
	$$
	The constant $K$ is uniform over $\bb\omega\in\bb\Omega$, because all coordinates of $\bb\omega$ are restricted to compact sets with deterministic radii. Taking the supremum over $\bb\Omega$ proves the first claim.

	Finally, for each $d\geq1$, let $\widetilde\bh_d \sim \mathrm N(\bb0_d,\bb I_d)\,.$
	Since $\rank(\bb B)\in\{0,1,2\}$ and $p/n\to\kappa$,
	$$
		\underset{j\in\{0,1,2\}}{\sup}
		\abs{
			\frac{
				\expect\left[
					\vnorm{\widetilde\bh_{p-j}}_2
				\right]
			}{\sqrt n}
			-
			\sqrt\kappa
		}
		\to0\,.
	$$
	Consequently,
	$$
		\frac{\mu_{\bh}}{\sqrt n}
		\to
		\sqrt\kappa\,,
	$$
	uniformly over the possible values of $\rank(\bb B)$.
	Given $s>0$, choose $t_0=t_0(s)>0$ small enough so that
	$$
		K(t_0+\sqrt{t_0})\leq \frac{s}{6}\,,
	$$
	and then choose $N_s$ large enough so that, for all $n>N_s$,
	$$
		K
		\abs{
			\frac{\mu_{\bh}}{\sqrt n}
			-
			\sqrt\kappa
		}
		\leq
		\frac{s}{6} \,.
	$$
	This gives
	$$
		\mathds{1}_{\mathcal S_n(t_0)}
		\underset{\bb\omega\in\bb\Omega}{\sup}
		\abs{
			\bar F_n(\bb\omega)-F(\bb\omega)
		}
		\leq
		\frac{s}{3}\,,
	$$
	for all $n>N_s$.
\end{proof}

	\begin{lemma}[Lipschitz bound for the finite AO objective]
		\label{lemma:lipschitz_ao}
		Let $\bb\omega=(\sigma,r,\bu,\theta,t), \quad \bb\omega'=(\sigma',r',\bu',\theta',t')$
		be two points in $\bb\Omega$. Let $\lambda=t/r, \quad \lambda'=t'/r'\,.$
	Then
	$$
		\abs{F_n(\bb\omega)-F_n(\bb\omega')}
		\leq
		L_{\sigma,n}\abs{\sigma-\sigma'}
		+
		L_{r,n}\abs{r-r'}
		+
		L_{u,n}\vnorm{\bu-\bu'}_2
		+
		L_{\theta,n}\abs{\theta-\theta'}
		+
		L_{t,n}\abs{t-t'}\,,
	$$
	where
	$$
		L_{\theta,n}
		=
		1+\frac{\abs{\bYtil^\top\bb1}}{n}, \quad 
		L_{u,n}
		=
		\frac{\mnorm{\bb Z}_2}{\sqrt n}
		+
		\frac{\vnorm{\bb Z^\top\bYtil}_2}{n}, \quad 
		L_{\sigma,n}
		=
		\frac{\vnorm{\bh}_2}{\sqrt n}
		+
		\frac{\vnorm{\bg}_2}{\sqrt n}
		+
		\frac{\abs{\bYtil^\top\bg}}{n}, \quad 
		L_{t,n}
		=
		\frac12
		+
		\frac1{c_r}
		\left(
			\frac{\vnorm{\bYtil}_2}{\sqrt n}
			+
			\frac12
		\right)
		+
		\frac1{2c_r}
		\frac{\vnorm{\bYtil}_2^2}{n}\,,
	$$
	and
	$$
		L_{r,n}
		=
		C_\sigma\frac{\vnorm{\bh}_2}{\sqrt n}
		+
		\frac{C_t}{2}
		+
		\frac{C_t}{c_r^2}
		\left(
			\frac{\vnorm{\bYtil}_2}{\sqrt n}
			+
			\frac12
		\right)
		+
		\frac{C_t}{2c_r^2}
		\frac{\vnorm{\bYtil}_2^2}{n}\,.
	$$
\end{lemma}
\begin{proof}
	Write
	$$
		\Delta_\sigma=\sigma-\sigma',
		\quad
		\Delta_r=r-r',
		\quad
		\Delta_{\bu}=\bu-\bu',
		\quad
		\Delta_\theta=\theta-\theta',
		\quad
		\Delta_t=t-t'\,.
	$$
	Since $r,r'\in[c_r,1]$ and $t,t'\in[0,C_t]$,
	$$
		\abs{
			\frac tr-\frac{t'}{r'}
		}
		=
		\abs{
			\frac{(t-t')r'+t'(r'-r)}{rr'}
		}
		\leq
		\frac{\abs{\Delta_t}}{c_r}
		+
		\frac{C_t}{c_r^2}\abs{\Delta_r}\,.
	$$
	Similarly,
	$$
		\abs{
			\frac{t}{2r}
			-
			\frac{t'}{2r'}
		}
		\leq
		\frac{\abs{\Delta_t}}{2c_r}
		+
		\frac{C_t}{2c_r^2}\abs{\Delta_r}\,.
	$$
	Define
	$$
		\bb x
		=
		\theta\bb1+\bb Z\bu+\sigma\bg+\lambda\bYtil,
		\quad
		\bb x'
		=
		\theta'\bb1+\bb Z\bu'+\sigma'\bg+\lambda'\bYtil\,.
	$$
	By Lemma~\ref{lemma:moreau_lipschitz},
	$$
		\abs{
			\frac1nM_\rho(\bb x,\lambda)
			-
			\frac1nM_\rho(\bb x',\lambda')
		}
		\leq
		\frac1{\sqrt n}\vnorm{\bb x-\bb x'}_2
		+
		\frac12\abs{\lambda-\lambda'}\,.
	$$
	Furthermore,
	$$
		\bb x-\bb x'
		=
		\Delta_\theta\bb1
		+
		\bb Z\Delta_{\bu}
		+
		\Delta_\sigma\bg
		+
		(\lambda-\lambda')\bYtil\,,
	$$
	and hence
	$$
		\frac1{\sqrt n}\vnorm{\bb x-\bb x'}_2
		\leq
		\abs{\Delta_\theta}
		+
		\frac{\mnorm{\bb Z}_2}{\sqrt n}\vnorm{\Delta_{\bu}}_2
		+
		\frac{\vnorm{\bg}_2}{\sqrt n}\abs{\Delta_\sigma}
		+
		\frac{\vnorm{\bYtil}_2}{\sqrt n}\abs{\lambda-\lambda'}\,.
	$$
	Combining the last three displays gives the Moreau part of the Lipschitz bound.

	Next,
	$$
		\abs{
			-\sigma r\frac{\vnorm{\bh}_2}{\sqrt n}
			+
			\sigma'r'\frac{\vnorm{\bh}_2}{\sqrt n}
		}
		\leq
		\frac{\vnorm{\bh}_2}{\sqrt n}
		\left(
			\abs{\Delta_\sigma}
			+
			C_\sigma\abs{\Delta_r}
		\right)\,.
	$$
	Also,
	$$
		\abs{
			\frac{rt}{2}
			-
			\frac{r't'}{2}
		}
		\leq
		\frac{C_t}{2}\abs{\Delta_r}
		+
		\frac12\abs{\Delta_t}\,.
	$$
	For the linear term,
	$$
		\begin{aligned}
		&\abs{
			\frac1n
			\bYtil^\top
			\left[
				(\theta-\theta')\bb1
				+
				\bb Z(\bu-\bu')
				+
				(\sigma-\sigma')\bg
			\right]
		}\leq
		\frac{\abs{\bYtil^\top\bb1}}{n}\abs{\Delta_\theta}
		+
		\frac{\vnorm{\bb Z^\top\bYtil}_2}{n}\vnorm{\Delta_{\bu}}_2
		+
		\frac{\abs{\bYtil^\top\bg}}{n}\abs{\Delta_\sigma}\,.
		\end{aligned}
	$$
	Finally,
	$$
		\abs{
			\frac{t}{2rn}\vnorm{\bYtil}_2^2
			-
			\frac{t'}{2r'n}\vnorm{\bYtil}_2^2
		}
		\leq
		\frac{\vnorm{\bYtil}_2^2}{n}
		\left(
			\frac{\abs{\Delta_t}}{2c_r}
			+
			\frac{C_t}{2c_r^2}\abs{\Delta_r}
		\right)\,.
	$$
	Collecting the coefficients of $\abs{\Delta_\sigma}$, $\abs{\Delta_r}$, $\vnorm{\Delta_{\bu}}_2$, $\abs{\Delta_\theta}$, and $\abs{\Delta_t}$ gives the stated bound.
\end{proof}

\begin{lemma}[Bounded Lipschitz constants]
	\label{lemma:ao_lipschitz_constants}
	Let
	$$
		L_n
		=
		\max
		\left\{
			L_{\sigma,n},
			L_{r,n},
			L_{u,n},
			L_{\theta,n},
			L_{t,n}
		\right\}\,.
	$$
	There exist constants $\bar C<\infty$, $C,c>0$ and $N\in\mathbb N$ such that, for all $n>N$,
	$$
		\Pr\left(
			L_n>\bar C
			\mid
			\mathcal G_n
		\right)
		\leq
		C\exp\{-cn\}\,.
	$$
	In particular,
	$$
		\mathds{1}_{\mathcal S_n(t_0)}
		\Pr\left(
			L_n>\bar C
			\mid
			\mathcal G_n
		\right)
		\leq
		C\exp\{-cn\}\,.
	$$
	The constants depend only on $\kappa,C_\sigma,C_\beta,C_\theta,C_t,c_r$
	and on universal Gaussian concentration constants.
\end{lemma}
\begin{proof}
	Since $0<\Ytil_i<1$, we have $\abs{\bYtil^\top\bb1}/n\leq1$, $\vnorm{\bYtil}_2/\sqrt n\leq1$, and $\vnorm{\bYtil}_2^2/n\leq1\,.$
	Therefore, $L_{\theta,n}\leq2$
	and
	$$
		L_{t,n}
		\leq
		\frac12
		+
		\frac1{c_r}
		\left(
			1+\frac12
		\right)
		+
		\frac1{2c_r}\,.
	$$
	Also,
	$$
		L_{r,n}
		\leq
		C_\sigma\frac{\vnorm{\bh}_2}{\sqrt n}
		+
		\frac{C_t}{2}
		+
		\frac{3C_t}{2c_r^2}
		+
		\frac{C_t}{2c_r^2}\,.
	$$
	We next control the Gaussian terms. Conditional on $\mathcal G_n$, the vector $\bg$ and the matrix $\bb Z$
	have standard Gaussian entries, while
	$$
		\bh
		\sim
		\mathrm N\left(
			\bb0_{p-\rank(\bb B)},
			\bb I_{p-\rank(\bb B)}
		\right)\,.
	$$
	These Gaussian objects are conditionally independent. Standard chi-square concentration
	(see, for example, \citealt[][Theorem~3.1.1]{vershynin:2018}) gives
	constants $C,c>0$ such that, for all sufficiently large $n$,
	$$
		\Pr\left(
			\frac{\vnorm{\bg}_2}{\sqrt n}>2
			\mid
			\mathcal G_n
		\right)
		\leq
		Ce^{-cn}\,,
	$$
	and
	$$
		\Pr\left(
			\frac{\vnorm{\bh}_2}{\sqrt n}>2\sqrt{\kappa}+1
			\mid
			\mathcal G_n
		\right)
		\leq
		Ce^{-cn}\,.
	$$
	Moreover, since $\bb Z$ is $n\times2$ with i.i.d. standard Gaussian entries, the Gaussian singular-value bound gives
	$$
		\Pr\left(
			\frac{\mnorm{\bb Z}_2}{\sqrt n}>2
			\mid
			\mathcal G_n
		\right)
		\leq
		Ce^{-cn}\,,
	$$
	for all sufficiently large $n$. Finally,
	$$
		\frac{\vnorm{\bb Z^\top\bYtil}_2}{n}
		\leq
		\frac{\mnorm{\bb Z}_2\vnorm{\bYtil}_2}{n}
		\leq
		\frac{\mnorm{\bb Z}_2}{\sqrt n}\,,
	$$
	and
	$$
		\frac{\abs{\bYtil^\top\bg}}{n}
		\leq
		\frac{\vnorm{\bYtil}_2\vnorm{\bg}_2}{n}
		\leq
		\frac{\vnorm{\bg}_2}{\sqrt n}\,.
	$$
	On the intersection of the above Gaussian events, all quantities $L_{\sigma,n}$, $L_{r,n}$, $L_{u,n}$, $L_{\theta,n}$, and $L_{t,n}$ are bounded by a deterministic constant $\bar C<\infty$ depending only on $\kappa,C_\sigma,C_t,c_r\,.$
	A union bound gives
	$$
		\Pr\left(
			L_n>\bar C
			\mid
			\mathcal G_n
		\right)
		\leq
		C\exp\{-cn\}\,.
	$$
	This proves the claim.
\end{proof}

\begin{lemma}[Lipschitz bound for the limiting AO objective]
	\label{lemma:limit_ao_lipschitz}
	There exists a constant $\bar C_F<\infty$, depending only on
	$$
		\kappa,\alpha,\bb\Gamma,\theta_0^*,\theta_P^*,
		C_\sigma,C_\beta,C_\theta,C_t,c_r\,,
	$$
	such that, for all $\bb\omega,\bb\omega'\in\bb\Omega$,
	$$
		\abs{F(\bb\omega)-F(\bb\omega')}
		\leq
		\bar C_F
		\vnorm{\bb\omega-\bb\omega'}_2\,.
	$$
	\end{lemma}
	\begin{proof}
		Let $\bb\omega=(\sigma,r,\bu,\theta,t), \quad \bb\omega'=(\sigma',r',\bu',\theta',t')$
		be points in $\bb\Omega$, and write $\lambda=t/r, \quad \lambda'=t'/r'\,.$
	Since $r,r'\in[c_r,1]$ and $t,t'\in[0,C_t]$,
	$$
		\abs{\lambda-\lambda'}
		\leq
		\frac{\abs{t-t'}}{c_r}
		+
		\frac{C_t}{c_r^2}\abs{r-r'}\,.
	$$
	Also, if
	$$
		\xi=\theta+\bz^\top\bu+\sigma G,
		\quad
		\xi'=\theta'+\bz^\top\bu'+\sigma'G\,,
	$$
	then
	$$
		\abs{\xi-\xi'}
		\leq
		\abs{\theta-\theta'}
		+
		\vnorm{\bz}_2\vnorm{\bu-\bu'}_2
		+
		\abs{G}\abs{\sigma-\sigma'}\,.
	$$
	Because $0<\Ytil<1$, Lemma~\ref{lemma:moreau_lipschitz} gives
	$$
		\begin{aligned}
		&\abs{
			m_\rho(\xi+\lambda \Ytil,\lambda)
			-
			\Ytil\xi
			-
			\frac{\lambda}{2}\Ytil^2
			-
			m_\rho(\xi'+\lambda'\Ytil,\lambda')
			+
			\Ytil\xi'
			+
			\frac{\lambda'}{2}\Ytil^2
		}\leq
		2\abs{\xi-\xi'}
		+
		2\abs{\lambda-\lambda'}\,.
		\end{aligned}
	$$
	Furthermore,
	$$
		\abs{\sigma r-\sigma'r'}
		\leq
		\abs{\sigma-\sigma'}
		+
		C_\sigma\abs{r-r'}\,,
	$$
	and
	$$
		\abs{rt-r't'}
		\leq
		\abs{t-t'}
		+
		C_t\abs{r-r'}\,.
	$$
	Since
	$$
		\expect[\vnorm{\bz}_2]
		+
		\expect[\abs G]
		<
		\infty\,,
	$$
	taking expectations and collecting the coordinatewise coefficients
	gives a deterministic constant $\bar C_F<\infty$ with the stated
	dependence such that
	$$
		\abs{F(\bb\omega)-F(\bb\omega')}
		\leq
		\bar C_F
		\vnorm{\bb\omega-\bb\omega'}_2\,.
	$$
\end{proof}

\subsection{Lemmas for Section~\ref{sec:limiting-AO-properties}}

\subsubsection{Well-posedness and natural-domain presentation}

\begin{lemma}[Population residual gap for the limiting pseudo-response]
\label{lemma:limiting_AO_population_gap}
Let $\bb Q = [Q_1, Q_2]^\top$ as in Section~\ref{sec:limiting-ao-focs}, and set
$$
    q_0=\rho'(\theta_0^*+Q_1),
    \quad
    q_P=\rho'(\theta_P^*+Q_2),
    \quad
    y
    =
    \alpha\mathds{1}\{\varepsilon<q_0\}+(1-\alpha)q_P
    =
    \Ytil\,,
$$
so that $y$ is the limiting pseudo-response of
Section~\ref{subsec:mer}. Let
$$
    y_0=(1-\alpha)q_P,
    \quad
    y_1=\alpha+(1-\alpha)q_P\,,
$$
and define $\eta_y=\logit{y}\,.$
Then $0<y<1$ almost surely, $\eta_y\in L^2$, and, with
\begin{equation}
\label{eq:limiting_AO_d0}
    d_0^2
    =
    \underset{\substack{
    \theta\in\Re\\
    \bu\in\range(\bb\Gamma)
    }}{\inf}
    \expect\left[
        \left(
            \eta_y-\theta-\bz^\top\bu
        \right)^2
    \right]\,,
\end{equation}
one has
\begin{equation}
\label{eq:limiting_AO_d0_lower}
    d_0^2
    \geq
    16\alpha^2\expect[q_0(1-q_0)]
    >0\,.
\end{equation}
Moreover, for every $\sigma\geq0$,
\begin{equation}
\label{eq:limiting_AO_sigma_distance}
    \underset{\substack{
    \theta\in\Re\\
    \bu\in\range(\bb\Gamma)
    }}{\inf}
    \left\|
        \eta_y-\theta-\bz^\top\bu-\sigma G
    \right\|_{L^2}
    =
    \sqrt{d_0^2+\sigma^2}\,,
\end{equation}
and consequently
\begin{equation}
\label{eq:limiting_AO_uniform_radius_gap}
    \underset{\sigma\geq0}{\inf}
    \left
    \{
        \sqrt{d_0^2+\sigma^2}
        -\sigma\sqrt\kappa
    \right\}
    =
    d_0\sqrt{1-\kappa}
    >0\,.
\end{equation}
\end{lemma}

\begin{proof}
Conditional on $\bz$, the pseudo-response $y$ equals $y_1$ with
probability $q_0$ and $y_0$ with probability $1-q_0$. Since
$0<q_0,q_P<1$ almost surely and $\alpha\in(0,1)$, both conditional atoms
have positive probability and belong to $(0,1)$.

Put $A=\theta_P^*+Q_2$. Direct calculation gives $\logit{y_0} = A+\log(1-\alpha)-\log(1+\alpha e^A)$
and $\logit{y_1} = \log(\alpha+e^A)-\log(1-\alpha)\,.$
Hence there is a constant $C_\alpha<\infty$ such that
$$
    \max\left
    \{
        \abs{\logit{y_0}},
        \abs{\logit{y_1}}
    \right\}
    \leq
    \abs A+C_\alpha\,.
$$
Since $A$ is Gaussian, $\eta_y\in L^2$.

Conditional on $\bz$,
$$
    \var(\eta_y\mid\bz)
    =
    q_0(1-q_0)
    \left(
        \logit{y_1}-\logit{y_0}
    \right)^2\,.
$$
Because $y_1-y_0=\alpha$ and $d\logit{v}/dv = 1/\{v(1-v)\} \geq4, \quad v\in(0,1)\,,$
one has $\logit{y_1}-\logit{y_0} \geq4\alpha\,.$
For every $\theta\in\Re$ and
$\bu\in\range(\bb\Gamma)$, the variable
$\theta+\bz^\top\bu$ is $\sigma(\bz)$-measurable. Therefore
$$
	\expect\left[ \left( \eta_y-\theta-\bz^\top\bu \right)^2 \right] \geq \expect[\var(\eta_y\mid\bz)] \geq 16\alpha^2\expect[q_0(1-q_0)]\,.
$$
The final expectation is strictly positive because
$q_0(1-q_0)>0$ almost surely. Taking the infimum proves
\eqref{eq:limiting_AO_d0_lower}.

Finally, $G$ is independent of $(\eta_y,\bz)$, has mean zero and variance
one. Thus, for every admissible $(\theta,\bu)$,
$$
    \begin{aligned}
    &\expect\left[
        \left(
            \eta_y-\theta-\bz^\top\bu-\sigma G
        \right)^2
    \right]=
    \expect\left[
        \left(
            \eta_y-\theta-\bz^\top\bu
        \right)^2
    \right]
    +\sigma^2\,.
    \end{aligned}
$$
Taking the infimum gives \eqref{eq:limiting_AO_sigma_distance}. A direct
one-dimensional minimisation gives
$$
    \inf_{\sigma\geq0}
    \left
    \{
        \sqrt{d_0^2+\sigma^2}
        -\sigma\sqrt\kappa
    \right\}
    =
    d_0\sqrt{1-\kappa}\,,
$$
which proves \eqref{eq:limiting_AO_uniform_radius_gap}.
\end{proof}

\begin{lemma}[Geometry and smoothness of the limiting AO]
\label{lemma:limiting_AO_geometry}
For $y\in(0,1)$, write
$$
    \ell_y(\eta)
    =
    \rho(\eta)-y\eta
$$
and, for $\lambda\geq0$, recall from the paragraph preceding Lemma~\ref{lemma:nullspace_no_help_barFn} that
$$
    e_\lambda(x,y)
    =
    m_\rho(x+\lambda y,\lambda)
    -
    yx
    -
    \frac{\lambda}{2}y^2\,,
$$
with $m_\rho(x,0)=\rho(x)$. For $r>0$ and $t>0$,
\begin{equation}
\label{eq:limiting_AO_perspective}
    e_{t/r}(x,y)
    =
    \underset{\eta\in\Re}{\inf}
    \left\{
        \ell_y(\eta)
        +\frac{r}{2t}(\eta-x)^2
    \right\}\,.
\end{equation}
At $t=0$, give the quadratic perspective its closed extension
$$
    \left[
        \frac{r}{2t}(\eta-x)^2
    \right]_{t=0}
    =
    \begin{cases}
        0,&\eta=x,\\
        +\infty,&\eta\neq x\,.
    \end{cases}
$$
Then the right-hand side of \eqref{eq:limiting_AO_perspective} equals
$\ell_y(x)$ at $t=0$, consistently with the continuous extension of the
compact limiting objective. The resulting closed function $F$ is jointly
continuous on
$$
    [0,\infty)
    \times(0,\infty)
    \times\range(\bb\Gamma)
    \times\Re
    \times[0,\infty)\,,
$$
where the coordinates are ordered as $(\sigma,r,\bu,\theta,t)$.

For every fixed $r>0$, the closed map $(\sigma,\bu,\theta,t) \longmapsto F(\sigma,r,\bu,\theta,t)$
is strictly convex on
$$
    [0,\infty)
    \times\range(\bb\Gamma)
    \times\Re
    \times[0,\infty)\,.
$$
For every fixed $(\sigma,\bu,\theta,t)$ with $t>0$, the map $r\longmapsto F(\sigma,r,\bu,\theta,t)$
is strictly concave on $(0,\infty)$.

On the relative positive domain $r,t>0$, the function $F$ is $C^2$. Put
$$
    \lambda=\frac tr,
    \quad
    \xi=\theta+\bz^\top\bu+\sigma G\,,
$$
$$
    \eta
    =
    \prox{\lambda\rho}{\xi+\lambda y},
    \quad
    R=\rho'(\eta)-y,
    \quad
    W=
    \frac{\rho''(\eta)}{1+\lambda\rho''(\eta)}\,.
$$
Then
\begin{equation}
\label{eq:limiting_AO_gradient}
\begin{aligned}
    \partial_\theta F
    &=\expect[R],
    &
    \nabla_{\bu}F
    &=\expect[\bz R],
    \\
    \partial_\sigma F
    &=\expect[GR]-r\sqrt\kappa,
    &
    \partial_tF
    &=-\frac{1}{2r}\expect[R^2]+\frac r2,
    \\
    \partial_rF
    &=\frac{t}{2r^2}\expect[R^2]
      -\sigma\sqrt\kappa+\frac t2\,.
\end{aligned}
\end{equation}
For a direction $h=(h_\sigma,h_{\bu},h_\theta,h_t)$
in the minimising coordinates, set $A_h=h_\theta+\bz^\top h_{\bu}+Gh_\sigma\,.$
Then
\begin{equation}
\label{eq:limiting_AO_min_hessian}
    D_x^2F(x,r)[h,h]
    =
    \expect\left[
        W
        \left(
            A_h-\frac{R}{r}h_t
        \right)^2
    \right]\,,
\end{equation}
and
\begin{equation}
\label{eq:limiting_AO_r_hessian}
    \partial_{rr}F(x,r)
    =
    -\frac{t}{r^3}
    \expect\left[
        \frac{R^2}
        {1+(t/r)\rho''(\eta)}
    \right] <0\,.
\end{equation}
\end{lemma}

\begin{proof}
Expanding the square in the Moreau definition gives
\eqref{eq:limiting_AO_perspective}. Since $\ell_y'(\eta)=\rho'(\eta)-y$
has absolute value at most one, $\ell_y$ is one-Lipschitz. Testing
$\eta=x$ gives $e_\lambda(x,y)\leq\ell_y(x)\,,$
whereas $\ell_y(\eta) \geq \ell_y(x)-\abs{\eta-x}$
gives
$$
    e_\lambda(x,y)
    \geq
    \ell_y(x)-\frac\lambda2\,.
$$
Thus $e_\lambda(x,y)\to\ell_y(x)$ as $\lambda\to 0$, locally
uniformly when $r$ is bounded away from zero. Together with the
one-Lipschitz property of $\ell_y$, dominated convergence gives joint
continuity of the closed expectation and hence of $F$.

For fixed $r>0$, the map
$$
    (\eta,x,t)
    \longmapsto
    \ell_y(\eta)+\frac{r}{2t}(\eta-x)^2\,,
$$
with its closed perspective at $t=0$ is jointly convex. Partial
minimisation over $\eta$ preserves convexity in $(x,t)$. Since
$x=\theta+\bz^\top\bu+\sigma G$ is affine in
$(\sigma,\bu,\theta)$, expectation preserves convexity, and the remaining
terms
$$
    -\sigma r\sqrt\kappa+\frac{rt}{2}\,,
$$
are affine in the minimising coordinates. For fixed minimising
coordinates, the expression inside the infimum in
\eqref{eq:limiting_AO_perspective} is affine in $r$. Its infimum is
therefore concave in $r$, and the remaining $r$-dependence is affine.

Fix $\lambda>0$, $x\in\Re$, and $y\in(0,1)$, and set
$$
    \eta
    =
    \prox{\lambda\rho}{x+\lambda y},
    \quad
    R
    =
    \rho'(\eta)-y,
    \quad
    W
    =
    \frac{\rho''(\eta)}
    {1+\lambda\rho''(\eta)}\,.
$$
The Moreau-envelope derivatives of
\citet[Lemma~2, equation~(29)]{salehi+et+al:2019} and the proximal-map
derivatives of
\citet[Proposition~6.3]{donoho+montanari:2016}, together with the
proximal first-order condition and the chain rule, give
\begin{equation}
\label{eq:limiting_AO_shifted_calculus}
\begin{aligned}
    \eta
    &=
    x-\lambda R,
    &
    \partial_xe_\lambda(x,y)
    &=
    R,
    &
    \partial_\lambda e_\lambda(x,y)
    &=
    -\frac12R^2,
    \\
    \partial_xR
    &=
    W,
    &
    \partial_\lambda R
    &=
    -WR,
    &
    \partial_yR
    &=
    -\frac{1}{1+\lambda\rho''(\eta)}\,.
\end{aligned}
\end{equation}
Here the $\lambda$-derivatives are taken with $x$ and $y$ fixed. Since
$0<\rho'(\eta)<1$, $0<y<1$, and $0<\rho''(\eta)\leq1/4$,
\begin{equation}
\label{eq:limiting_AO_shifted_bounds}
    \abs{R}<1,
    \quad
    0<W\leq\frac14,
    \quad
    1-\lambda W
    =
    \frac{1}{1+\lambda\rho''(\eta)}
    >
    0\,.
\end{equation}
Consequently,
\begin{equation}
\label{eq:limiting_AO_shifted_hessian}
    D^2_{(x,\lambda)}e_\lambda(x,y)[(a,b),(a,b)]
    =
    W(a-Rb)^2
\end{equation}
for every $(a,b)\in\Re^2$.

Applying \eqref{eq:limiting_AO_shifted_calculus} pointwise at $x=\xi$
and using the chain rule gives \eqref{eq:limiting_AO_gradient}. On every
compact subset of the relative domain $r,t>0$, the first derivatives are
dominated by a constant times
$$
    1+\vnorm{\bz}_2+\abs{G}\,,
$$
and the second derivatives by a constant times
$$
    \left(1+\vnorm{\bz}_2+\abs{G}\right)^2\,.
$$
The bounds in \eqref{eq:limiting_AO_shifted_bounds} make these envelopes
integrable. Dominated convergence therefore justifies differentiation
under the expectation and gives continuous second derivatives.

For fixed $r$, an increment $h$ changes the two arguments
$(\xi,\lambda)$ by
$$
    d\xi=A_h,
    \quad
    d\lambda=\frac{h_t}{r}\,.
$$
Equation~\eqref{eq:limiting_AO_shifted_hessian} gives
\eqref{eq:limiting_AO_min_hessian}. Also,
$$
    \frac{\partial R}{\partial r}
    =
    \frac{t}{r^2}WR\,.
$$
Differentiating the last line of
\eqref{eq:limiting_AO_gradient} gives
$$
	\partial_{rr}F =-\frac{t}{r^3}\expect[R^2] +\frac{t^2}{r^4}\expect[WR^2] =-\frac{t}{r^3} \expect\left[R^2(1-\lambda W)\right]\,,
$$
which is \eqref{eq:limiting_AO_r_hessian} by
\eqref{eq:limiting_AO_shifted_bounds}.

It remains to prove strictness. Suppose first that $t>0$ and that the
quadratic form in \eqref{eq:limiting_AO_min_hessian} vanishes for a
direction $h$. Since $W>0$ almost surely,
$$
    A_h=\frac{R}{r}h_t
    \quad\text{almost surely}\,.
$$
If $h_t=0$, then the affine Gaussian variable $A_h$ vanishes almost
surely, which implies
$$
    h_\sigma=0,
    \quad
    h_{\bu}=\bb0,
    \quad
    h_\theta=0\,.
$$
If $h_t\neq0$, boundedness of $R$ implies that $A_h$ is bounded almost
surely. Hence again $h_\sigma=0$ and $h_{\bu}=\bb0$, so $R$ would have
to be almost surely constant. This is impossible: conditional on
$(\bz,G)$, the response takes the two distinct values $y_0$ and $y_1$
with positive probabilities, while
\eqref{eq:limiting_AO_shifted_calculus} gives
$$
    \partial_yR
    =
    -\frac{1}{1+\lambda\rho''(\eta)}
    <
    0\,.
$$
Thus the Hessian is positive definite in every nonzero minimising
direction when $t>0$.

To extend strict convexity to the closed domain, take two distinct points
in the minimising domain and restrict $F(\cdot,r)$ to the line segment
joining them. If at least one endpoint has positive $t$, every interior
point of the segment has positive $t$, and the preceding positive second
directional derivative gives strict convexity along the segment. If both
endpoints have $t=0$, the restriction is
$$
    \expect[\ell_y(\theta+\bz^\top\bu+\sigma G)]
    -\sigma r\sqrt\kappa\,,
$$
and its second directional derivative is
$$
    \expect[\rho''(\xi)A_h^2]>0\,,
$$
for every nonzero coefficient direction. This proves strict convexity on
the closed minimising domain.

Finally, the expectation in \eqref{eq:limiting_AO_r_hessian} is strictly
positive. Conditional on $(\bz,G)$, the two possible residual values are
distinct, so they cannot both be zero: both atoms have positive
conditional probability. Hence \eqref{eq:limiting_AO_r_hessian} is
strictly negative, proving strict concavity in $r$ when $t>0$.
\end{proof}

\begin{lemma}[Coercivity and existence of a closed limiting-AO saddle]
\label{lemma:limiting_AO_closed_saddle}
There exist deterministic constants $\bar r\in(0,1)$, $a>0$, $B_a<\infty$, and $B_t<\infty$
such that
\begin{equation}
\label{eq:limiting_AO_coercive_bound}
    F(\sigma,\bar r,\bu,\theta,t)
    \geq
    a
    \sqrt{\theta^2+\vnorm{\bu}_2^2+\sigma^2}\,,
\end{equation}
for every $\sigma\geq0$, $\bu\in\range(\bb\Gamma)$, $\theta\in\Re$, and $t\geq0\,.$
For every $c\in(0,\bar r]$, define
$$
    \mathcal X_{\mathrm{cl}}
    =
    [0,\infty)
    \times\range(\bb\Gamma)
    \times\Re
    \times[0,\infty)\,,
$$
and
$$
    J_c(x)
    =
    \underset{r\in[c,1]}{\max}F(x,r),
    \quad
    x=(\sigma,\bu,\theta,t)\in\mathcal X_{\mathrm{cl}}\,.
$$
Then $J_c$ is continuous and coercive, and 
\begin{equation}
\label{eq:limiting_AO_closed_minimax}
    \underset{x\in\mathcal X_{\mathrm{cl}}}{\min}
    \underset{r\in[c,1]}{\max}F(x,r)
    =
    \underset{r\in[c,1]}{\max}
    \underset{x\in\mathcal X_{\mathrm{cl}}}{\inf}F(x,r)\,.
\end{equation}
The common value is attained by a saddle point $(x_c,r_c)$. Every
minimiser $x=(\sigma,\bu,\theta,t)$ of $J_c$ satisfies
\begin{equation}
\label{eq:limiting_AO_Ba}
    \sqrt{\theta^2+\vnorm{\bu}_2^2+\sigma^2}
    \leq B_a\,,
\end{equation}
and
\begin{equation}
\label{eq:limiting_AO_Bt}
    t\leq B_t\,.
\end{equation}
The bounds are uniform over $c\in(0,\bar r]$.
\end{lemma}

\begin{proof}
Choose $\delta\in(0,1-\sqrt\kappa)\,.$
For a unit coefficient direction
$$
    d=(d_\theta,d_{\bu},d_\sigma),
    \quad
    d_\theta^2+\vnorm{d_{\bu}}_2^2+d_\sigma^2=1\,,
$$
write $D_d=d_\theta+\bz^\top d_{\bu}+d_\sigma G$
and set $m(y)=\min(y,1-y)>0 \quad\text{almost surely}\,.$
For every unit direction, $\expect[D_d^2]=1$, and the family
$\{D_d:\vnorm{d}_2=1\}$ has uniformly bounded fourth moments.
For $a>0$ and $r>0$,
$$
    \begin{aligned}
    &\expect\left[
        D_d^2\mathds{1}
        \left\{
            r\abs{D_d}>m(y)
        \right\}
    \right]\leq
    \left(\expect[D_d^4]\right)^{1/2}
    \Pr(m(y)<a)^{1/2}
    +
    \frac{r^2}{a^2}\expect[D_d^4]\,.
    \end{aligned}
$$
The bound is uniform in $d$. First choose $a$ small and then choose
$\bar r\in(0,1)$ small enough that
\begin{equation}
\label{eq:limiting_AO_dual_mass}
    \underset{\vnorm{d}_2=1}{\inf}
    \expect\left[
        D_d^2\mathds{1}
        \left\{
            \bar r\abs{D_d}\leq m(y)
        \right\}
    \right]
    \geq1-\delta\,.
\end{equation}

For $v\in[0,1]$, let $\rho^*(v) = v\log v+(1-v)\log(1-v)\,,$
with $0\log0=0$. Thus $\rho^*(v)\leq0$. Fenchel's inequality \citep[Chapter~4]{beck:2017} gives
\begin{equation}
\label{eq:limiting_AO_fenchel}
    \ell_y(\eta)
    \geq
    w\eta-\rho^*(y+w)\,,
\end{equation}
whenever $-y\leq w\leq1-y$.
For $X\in L^2$, define
$$
    K_r(X)
    =
    \underset{\eta\in L^2}{\inf}
    \left\{
        \expect[\ell_y(\eta)]
        +r\vnorm{\eta-X}_{L^2}
    \right\}\,.
$$
For $t>0$, apply \eqref{eq:limiting_AO_perspective} at the pointwise
proximal minimiser and use $a^2/(2t)+t/2\geq a$ for $a\geq0$ to obtain
\begin{equation}
\label{eq:limiting_AO_profile_lower}
    F(\sigma,r,\bu,\theta,t)
    \geq
    K_r(\theta+\bz^\top\bu+\sigma G)
    -r\sigma\sqrt\kappa\,.
\end{equation}
At $t=0$, the same bound follows by taking
$\eta=\theta+\bz^\top\bu+\sigma G$ in the definition of $K_r$.

Put $s=\sqrt{\theta^2+\vnorm{\bu}_2^2+\sigma^2}\,.$
If $s=0$, the desired lower bound follows from
$\ell_y\geq-\rho^*(y)\geq0$. If $s>0$, set
$$
    d=\frac{(\theta,\bu,\sigma)}s\,,
$$
and define
$$
    w
    =
    \bar rD_d\mathds{1}
    \left\{
        \bar r\abs{D_d}\leq m(y)
    \right\}\,.
$$
Then $-y\leq w\leq1-y$ almost surely and
$\vnorm{w}_{L^2}\leq\bar r$. For every $\eta\in L^2$,
$$
\begin{aligned}
    \expect[\ell_y(\eta)]
    +\bar r\vnorm{\eta-sD_d}_{L^2}\geq
    \expect[w\eta]
    -\expect[\rho^*(y+w)]
    +\bar r\vnorm{\eta-sD_d}_{L^2}
  \geq
    s\expect[wD_d]
    -\expect[\rho^*(y+w)]\geq
    s\bar r(1-\delta)\,,
\end{aligned}
$$
where the second inequality uses Cauchy--Schwarz \citep[Theorem~1.37(d)]{rudin+etal:1976} and
$\vnorm{w}_{L^2}\leq\bar r$, and the final inequality uses
\eqref{eq:limiting_AO_dual_mass} and $\rho^*\leq0$. Taking the infimum
over $\eta$ and using $\sigma\leq s$ in
\eqref{eq:limiting_AO_profile_lower} gives
$$
    F(\sigma,\bar r,\bu,\theta,t)
    \geq
    \bar r(1-\delta-\sqrt\kappa)s\,.
$$
Thus \eqref{eq:limiting_AO_coercive_bound} holds with
$$
    a
    =
    \bar r(1-\delta-\sqrt\kappa)>0\,.
$$
Since $J_c\geq F(\cdot,\bar r)$, the coefficient block is coercively
controlled. If that block remains bounded while $t\to\infty$, then
$$
    e_{t/\bar r}(\xi,y)
    \geq
    \underset{\eta\in\Re}{\inf}\ell_y(\eta)
    =
    -\rho^*(y)
    \geq0\,,
$$
so
$$
    F(\sigma,\bar r,\bu,\theta,t)
    \geq
    -\bar r\sigma\sqrt\kappa
    +\frac{\bar r t}{2}
    \to+\infty\,.
$$
Hence $J_c$ is coercive. Continuity follows from
Lemma~\ref{lemma:limiting_AO_geometry}. Therefore $J_c$ attains its
minimum.

The interval $[c,1]$ is compact and convex, the closed natural minimising
domain is convex, and Lemma~\ref{lemma:limiting_AO_geometry} gives
continuity, convexity in $x$, and concavity in $r$. Sion's min--max theorem
(see, for example, \citealt[][Theorem~3]{simons:1995}) gives
\eqref{eq:limiting_AO_closed_minimax}. Define
$$
    h_c(r)
    =
    \underset{x\in\mathcal X_{\mathrm{cl}}}{\inf}F(x,r)\,.
$$
As an infimum of continuous functions, $h_c$ is upper semicontinuous and
therefore attains its maximum on $[c,1]$. Let $x_c$ minimise $J_c$, let
$r_c$ maximise $h_c$, and let $v_c$ denote the common minimax value. Then
$$
    v_c=h_c(r_c)
    \leq F(x_c,r_c)
    \leq J_c(x_c)=v_c\,.
$$
Thus equality holds throughout and $(x_c,r_c)$ is a saddle point.

Finally, the feasible point $(\sigma,\bu,\theta,t)=(0,\bb0,0,0)$ has
value $\log2$ for every $r$, so every minimiser of $J_c$ satisfies
$J_c(x)\leq\log2$. The coercive bound gives
$$
    \sqrt{\theta^2+\vnorm{\bu}_2^2+\sigma^2}
    \leq
    \frac{\log2}{a}
    =B_a\,.
$$
Also,
$$
    F(\sigma,\bar r,\bu,\theta,t)
    \geq
    -\bar r\sigma\sqrt\kappa+\frac{\bar r t}{2}\,.
$$
Since $F(\sigma,\bar r,\bu,\theta,t)\leq J_c(x)\leq\log2$ and
$\sigma\leq B_a$,
$$
    t
    \leq
    \frac{2}{\bar r}
    \left\{
        \log2+\bar rB_a\sqrt\kappa
    \right\}
    =B_t\,.
$$
This proves the uniform bounds.
\end{proof}

\begin{lemma}[Exclusion of the lower-radius boundary]
\label{lemma:limiting_AO_lower_radius}
Let $\bar r$ be as in
Lemma~\ref{lemma:limiting_AO_closed_saddle}. There exists
$c_0\in(0,\bar r)$ such that, for every $c\in(0,c_0]$, every saddle point
$(x_c,r_c)$ of \eqref{eq:limiting_AO_closed_minimax} satisfies $r_c>c$.
\end{lemma}

\begin{proof}
Suppose otherwise. Then there are $c_j\downarrow0$ and saddle points
$(\sigma_j,c_j,\bu_j,\theta_j,t_j)$ of the closed games over $[c_j,1]$.
By Lemma~\ref{lemma:limiting_AO_closed_saddle}, their minimising
coordinates are uniformly bounded.

Set
$$
    \xi_j
    =
    \theta_j+\bz^\top\bu_j+\sigma_jG\,.
$$
Let $\eta_j$ be the minimiser in the perspective representation of
Lemma~\ref{lemma:limiting_AO_geometry} at
$(\sigma_j,c_j,\bu_j,\theta_j,t_j)$, and set
$$
    M_j
    =
    \vnorm{\eta_j-\xi_j}_{L^2}\,.
$$
Joint minimisation in $(\theta,\bu,t)$ gives
$$
    t_j
    =
    M_j
    =
    \underset{\substack{
        \theta\in\Re\\
        \bu\in\range(\bb\Gamma)
    }}{\inf}
    \vnorm{
        \eta_j-\theta-\bz^\top\bu-\sigma_jG
    }_{L^2}\,.
$$
By Lemma~\ref{lemma:limiting_AO_population_gap}, strict convexity and
compactness of the uniformly bounded predictor family give
$\Delta_B>0$ such that
$$
    \expect[\ell_y(\xi_j)]
    \geq
    \expect[\ell_y(\eta_y)]
    +
    \Delta_B\,.
$$
Testing the joint minimisation with $\eta=\eta_y$, its $L^2$ projection,
and $t=\sqrt{d_0^2+\sigma_j^2}$, using
\eqref{eq:limiting_AO_sigma_distance}, gives
$$
    \expect[\ell_y(\eta_j)]
    +
    c_jM_j
    \leq
    \expect[\ell_y(\eta_y)]
    +
    c_j\sqrt{d_0^2+\sigma_j^2}\,.
$$
Since $\ell_y$ is one-Lipschitz,
$$
    \expect[\ell_y(\eta_j)]
    \geq
    \expect[\ell_y(\eta_y)]
    +
    \Delta_B
    -
    M_j\,,
$$
and therefore
$$
    (1-c_j)M_j
    \geq
    \Delta_B
    -
    c_j\sqrt{d_0^2+\sigma_j^2}\,.
$$
Thus $M_j\geq\Delta_B/2>0$ for all sufficiently large $j$; we restrict attention to such $j$, so that $t_j=M_j$ is interior and the first-order condition in $\eta$ below applies.
The first-order condition in $\eta$ is
$$
    \rho'(\eta_j)-y
    +
    \frac{c_j}{M_j}
    (\eta_j-\xi_j)
    =
    0
    \quad\text{almost surely}\,.
$$
Since $\rho'$ is increasing and $\rho'(\eta_y)=y$,
$$
	\{\rho'(\eta_j)-y\}(\eta_j-\eta_y)
	\geq
	0
	\quad\text{almost surely}\,.
$$
Substituting the first-order condition gives
$$
	(\eta_j-\xi_j)(\eta_j-\eta_y)
	\leq
	0
	\quad\text{almost surely}\,,
$$
so $\eta_j$ lies between $\xi_j$ and $\eta_y$. In either case,
$$
	\abs{\eta_j-\eta_y}
	\vee
	\abs{\eta_j-\xi_j}
	\leq
	\abs{\xi_j-\eta_y}
	\quad\text{almost surely}\,.
$$
Consequently, the same first-order condition gives
$$
    \abs{\rho'(\eta_j)-y}
    \leq
    \frac{c_j}{M_j}
    \abs{\xi_j-\eta_y}\,.
$$
The coefficient bounds imply
$$
    \abs{\xi_j-\eta_y}
    \leq
    B_a
    \left(
        1+\vnorm{\bz}_2+\abs G
    \right)
    +
    \abs{\eta_y}
    \in
    L^2\,.
$$
Since $c_j/M_j\to0$, it follows that
$\eta_j\to\eta_y$ almost surely and, by dominated convergence, in $L^2$.
Now $M_j$ and $\sqrt{d_0^2+\sigma_j^2}$ are the distances of
$\eta_j$ and $\eta_y$, respectively, from the same closed affine set.
Hence
$$
    \abs{
        M_j-\sqrt{d_0^2+\sigma_j^2}
    }
    \leq
    \vnorm{\eta_j-\eta_y}_{L^2}
    =
    o(1)\,.
$$
By \eqref{eq:limiting_AO_uniform_radius_gap},
$$
    M_j-\sigma_j\sqrt\kappa
    \geq
    d_0\sqrt{1-\kappa}
    +
    o(1)
    >
    0
$$
for all sufficiently large $j$.
The first-order condition and $t_j=M_j$ also give
$$
    \expect
    \left[
        \left\{
            \rho'(\eta_j)-y
        \right\}^2
    \right]
    =
    c_j^2\,.
$$
Therefore, by \eqref{eq:limiting_AO_gradient},
$$
    \left.
    \partial_r
    F(\sigma_j,r,\bu_j,\theta_j,t_j)
    \right|_{r=c_j}
    =
    M_j-\sigma_j\sqrt\kappa
    >
    0\,.
$$
This contradicts maximality at the left endpoint $r=c_j$.
\end{proof}

\begin{lemma}[Compact and natural limiting-AO saddle]
\label{lemma:limiting_AO_natural_domain}
Let $\bar r,B_a,B_t$ be the deterministic constants from
Lemma~\ref{lemma:limiting_AO_closed_saddle}, and let $c_0$ be the
constant from Lemma~\ref{lemma:limiting_AO_lower_radius}. Suppose that
$$
    C_\sigma>B_a,
    \quad
    C_\beta>B_a,
    \quad
    C_\theta>B_a,
    \quad
    C_t>B_t,
    \quad
    0<c_r\leq\min\{\bar r,c_0\}\,.
$$
Then the compact limiting AO has a unique optimiser chain
$$
    (\sigma^*,r^*,\bu^*,\theta^*,t^*)\in\bb\Omega^*\,,
$$
in the sense that
$$
    \underset{\sigma\in[0,C_\sigma]}{\arg\min}\,
    V(\sigma)
    =
    \{\sigma^*\}\,,
$$
$$
    \underset{r\in[c_r,1]}{\arg\max}\,
    H(\sigma^*,r)
    =
    \{r^*\}\,,
$$
and
$$
    \underset{(\bu,\theta,t)\in\bb K^*}{\arg\min}\,
    F(\sigma^*,r^*,\bu,\theta,t)
    =
    \{(\bu^*,\theta^*,t^*)\}\,.
$$
Moreover,
$$
    0<\sigma^*<C_\sigma,
    \quad
    c_r<r^*<1,
    \quad
    0<t^*<C_t,
    \quad
    \vnorm{\bu^*}_2<C_\beta,
    \quad
    \abs{\theta^*}<C_\theta\,.
$$
In particular, $\bu^*$ lies in the relative interior of $U^*$ within
$\range(\bb\Gamma)$. The same tuple is the unique saddle point of the
natural-domain problem in \eqref{eq:limit_ao_natural}, whose value is
$\bar\phi$.
\end{lemma}

\begin{proof}
For $x=(\sigma,\bu,\theta,t)$, write
$$
    F(x,r)
    =
    F(\sigma,r,\bu,\theta,t)\,.
$$
Apply Lemma~\ref{lemma:limiting_AO_closed_saddle} with $c=c_r$.
It gives a saddle $(x^*,r^*)$, where
$$
    x^*
    =
    (\sigma^*,\bu^*,\theta^*,t^*)\,,
$$
such that
\begin{equation}
\label{eq:closed_AO_saddle_inequalities}
    F(x^*,r)
    \leq
    F(x^*,r^*)
    \leq
    F(x,r^*)
\end{equation}
for every $r\in[c_r,1]$ and every
$x\in\mathcal X_{\mathrm{cl}}$. Moreover,
Lemma~\ref{lemma:limiting_AO_lower_radius} gives
$r^*>c_r$.
Suppose that $\sigma^*=0$. If $t^*>0$, let $R^*$ denote the residual
appearing in \eqref{eq:limiting_AO_gradient}. Since $R^*$ is independent
of $G$ when $\sigma^*=0$,
$$
    \left.
    \partial_\sigma^+
    F(\sigma,r^*,\bu^*,\theta^*,t^*)
    \right|_{\sigma=0}
    =
    -r^*\sqrt\kappa
    <
    0\,.
$$
If $t^*=0$, direct differentiation of the closed extension gives the
same right derivative. Both conclusions contradict the minimality in
\eqref{eq:closed_AO_saddle_inequalities}. Hence $\sigma^*>0$
If $t^*=0$, then
$$
    F(x^*,r)
    =
    \expect\left[
        \ell_{\Ytil}
        \left(
            \theta^*+\bz^\top\bu^*+\sigma^*G
        \right)
    \right]
    -
    \sigma^*r\sqrt\kappa
$$
is strictly decreasing in $r$. Its maximiser on $[c_r,1]$ would therefore
be $c_r$, contradicting $r^*>c_r$. Thus $t^*>0$.
The $t$-coordinate is now interior in the closed minimising domain. Hence
$\partial_tF(x^*,r^*)=0$, and
\eqref{eq:limiting_AO_gradient} gives
$$
    (r^*)^2
    =
    \expect[(R^*)^2]\,.
$$
Since $\Ytil\in(0,1)$ and $\rho'(\eta)\in(0,1)$ almost surely,
$\abs{R^*}<1$ almost surely. Therefore $r^*<1$.
The bounds
\eqref{eq:limiting_AO_Ba}--\eqref{eq:limiting_AO_Bt} and the strict
choices of the compactification constants give
$$
    \sigma^*<C_\sigma,
    \quad
    \vnorm{\bu^*}_2<C_\beta,
    \quad
    \abs{\theta^*}<C_\theta,
    \quad
    t^*<C_t\,.
$$
Thus $(x^*,r^*)$ belongs to the stated relative interior.
Set
$$
    v^*
    =
    F(x^*,r^*)\,.
$$
The right inequality in
\eqref{eq:closed_AO_saddle_inequalities} gives
$$
    H(\sigma,r^*)
    \geq
    v^*
$$
for every $\sigma\in[0,C_\sigma]$, and hence
$V(\sigma)\geq v^*$. Conversely, the left inequality gives
$$
    H(\sigma^*,r)
    \leq
    F(x^*,r)
    \leq
    v^*
$$
for every $r\in[c_r,1]$, while
$H(\sigma^*,r^*)=v^*$. Therefore
$$
    V(\sigma^*)=v^*=\bar\phi\,.
$$

Suppose that $V(\sigma)=v^*$. Then
$H(\sigma,r^*)=v^*$. By compactness and continuity, the inner minimum is
attained at some $k=(\bu,\theta,t)\in\bb K^*$, so
$$
    F(\sigma,r^*,k)
    =
    v^*\,.
$$
Lemma~\ref{lemma:limiting_AO_geometry} makes the map
$x\mapsto F(x,r^*)$ strictly convex, while
\eqref{eq:closed_AO_saddle_inequalities} shows that $x^*$ is its
minimiser. Hence $(\sigma,k)=x^*$, and therefore $\sigma=\sigma^*$.
This proves uniqueness of the outer minimiser.
If $H(\sigma^*,r)=v^*$, then the preceding inequalities imply
$F(x^*,r)=v^*$. Since $t^*>0$,
Lemma~\ref{lemma:limiting_AO_geometry} makes
$r\mapsto F(x^*,r)$ strictly concave, so $r=r^*$. Thus $r^*$ is the
unique maximiser of $H(\sigma^*,\cdot)$. Strict convexity of
$x\mapsto F(x,r^*)$ also gives the stated uniqueness of
$(\bu^*,\theta^*,t^*)$.

It remains to pass to the natural domain. Since
$r^*\in(c_r,1)$,
$$
    \partial_rF(x^*,r^*)=0\,.
$$
Because $t^*>0$, the map $r\mapsto F(x^*,r)$ is strictly concave on
$(0,\infty)$. Hence
$$
    F(x^*,r)
    \leq
    F(x^*,r^*)
$$
for every $r>0$. The right inequality in
\eqref{eq:closed_AO_saddle_inequalities} already holds for every
$x\in\mathcal X_{\mathrm{cl}}$. Consequently, $(x^*,r^*)$ is a saddle
of the natural-domain problem in \eqref{eq:limit_ao_natural}, with value
$\bar\phi$.
Finally, let $(\widetilde x,\widetilde r)$ be any other natural-domain
saddle. The two saddle inequalities give
$$
\begin{aligned}
    F(x^*,r^*)
    &\leq
    F(\widetilde x,r^*)
    \leq
    F(\widetilde x,\widetilde r)
    \leq
    F(x^*,\widetilde r)
    \leq
    F(x^*,r^*)\,.
\end{aligned}
$$
Equality therefore holds throughout. Strict convexity in $x$ gives
$\widetilde x=x^*$, and strict concavity in $r$ gives
$\widetilde r=r^*$. This proves uniqueness on the natural domain.
\end{proof}

\subsubsection{Ambient-domain relaxation}

\begin{lemma}[Ambient-domain relaxation of the limiting AO]
\label{lemma:ambient_limiting_saddle_chain}
Let $(\sigma^*,r^*,\bu^*,\theta^*,t^*)$ be the unique compact saddle
chain from Lemma~\ref{lemma:limiting_AO_natural_domain}. For every
$(\sigma,r)\in[0,C_\sigma]\times[c_r,1]$,
$$
    H_0(\sigma,r)=H(\sigma,r)\,.
$$
For every $\sigma\in[0,C_\sigma]$,
$$
    V_0(\sigma)=V(\sigma)\,.
$$
Every minimiser defining $H_0(\sigma,r)$ belongs to
$\range(\bb\Gamma)$ and is therefore a minimiser defining
$H(\sigma,r)$. Consequently,
$$
    \operatorname{val}(F;\bb\Omega)
    =
    \operatorname{val}(F;\bb\Omega^*)
    =
    \bar\phi\,.
$$
Moreover, the ambient profiles have the same unique saddle chain: $V_0$
has the unique minimiser $\sigma^*$, $H_0(\sigma^*,\cdot)$ has the unique
maximiser $r^*$, and $(\bu^*,\theta^*,t^*)$ is the unique minimiser of
$$
    (\bu,\theta,t)
    \longmapsto
    F(\sigma^*,r^*,\bu,\theta,t)\,,
$$
over $\bb K$.
\end{lemma}

\begin{proof}
By Lemma~\ref{lemma:limit_ao_lipschitz}, the function $F$ is continuous
on the ambient compact domain $\bb\Omega$. Hence all inner minima and
outer maxima defining $H$, $V$, $H_0$, and $V_0$ are attained.

Fix $(\sigma,r)\in[0,C_\sigma]\times[c_r,1]$. Since
$\bb K^*\subseteq\bb K$, the ambient inner minimisation is over a larger
set, and therefore $H_0(\sigma,r)\leq H(\sigma,r)$.

For $(\bu,\theta,t)\in\bb K$, write
$$
    \bu_R=\bb P_{\bb \Gamma}\bu,
    \quad
    \bu_N=\bu-\bu_R\,,
$$
where $\bb P_{\bb \Gamma}$ is the orthogonal projector onto
$\range(\bb\Gamma)$. Since $\bb P_{\bb \Gamma}$ is an orthogonal projector,
$\vnorm{\bu_R}_2\leq\vnorm{\bu}_2\leq C_\beta$, so
$(\bu_R,\theta,t)\in\bb K^*$. By
Lemma~\ref{lemma:nullspace_no_help_F},
$$
    F(\sigma,r,\bu_R,\theta,t)
    \leq
    F(\sigma,r,\bu,\theta,t)\,.
$$
Taking the infimum over $\bb K$ gives
$H(\sigma,r)\leq H_0(\sigma,r)$, and hence
$H_0(\sigma,r)=H(\sigma,r)$. Taking the maximum over
$r\in[c_r,1]$ gives $V_0(\sigma)=V(\sigma)$, and taking the minimum over
$\sigma\in[0,C_\sigma]$ gives the displayed equality of values.
Now let $(\tilde\bu,\tilde\theta,\tilde t)$ be any minimiser defining
$H_0(\sigma,r)$ and write
$$
    \tilde\bu_R=\bb P_{\bb \Gamma}\tilde\bu,
    \quad
    \tilde\bu_N=\tilde\bu-\tilde\bu_R\,.
$$
If $\tilde\bu_N\neq\bb0$, the final assertion of
Lemma~\ref{lemma:nullspace_no_help_F} gives
$$
    F(\sigma,r,\tilde\bu_R,\tilde\theta,\tilde t)
    <
    F(\sigma,r,\tilde\bu,\tilde\theta,\tilde t)\,,
$$
contradicting ambient minimality. Thus
$\tilde\bu\in\range(\bb\Gamma)$, so the ambient minimiser belongs to
$\bb K^*$ and, because $H_0(\sigma,r)=H(\sigma,r)$, also minimises the
starred profile.
Since $V_0=V$ pointwise and $V$ has the unique minimiser $\sigma^*$ by
Lemma~\ref{lemma:limiting_AO_natural_domain}, $V_0$ has the same unique
minimiser. Likewise,
$H_0(\sigma^*,\cdot)=H(\sigma^*,\cdot)$, so
$H_0(\sigma^*,\cdot)$ has the unique maximiser $r^*$. Finally, every
ambient inner minimiser at $(\sigma^*,r^*)$ is a starred inner minimiser
by the preceding argument, and the latter is uniquely
$(\bu^*,\theta^*,t^*)$ by
Lemma~\ref{lemma:limiting_AO_natural_domain}. This proves the ambient
uniqueness conclusions.
\end{proof}

\subsubsection{First-order equations}

\begin{lemma}[First-order equations for the limiting AO]
	\label{lemma:limiting_AO_focs}
	Let $(\sigma^*,r^*,\bu^*,\theta^*,t^*)$
	be the unique interior saddle given by
	Lemma~\ref{lemma:limiting_AO_natural_domain}, and set
	$$
		\lambda^*=\frac{t^*}{r^*}\,.
	$$
	For $\lambda>0$ and
	$(\sigma,\bu,\theta)\in(0,\infty)\times\range(\bb\Gamma)\times\Re$,
	define $\xi=\theta+\bz^\top\bu+\sigma G\,.$ With $q_0$ and $q_P$ as in
	Lemma~\ref{lemma:limiting_AO_population_gap}, set
	$$
		\Ytil^{(0)}=(1-\alpha)q_P,
		\quad
		\Ytil^{(1)}=\alpha+(1-\alpha)q_P\,.
	$$
	For $b\in\{0,1\}$, let $\eta_b = \prox{\lambda\rho}{\xi+\lambda \Ytil^{(b)}}\,,$
	and define $R_b=\rho'(\eta_b)-\Ytil^{(b)}$
	and
	$$
		W_b
		=
		\frac{\rho''(\eta_b)}
		{1+\lambda\rho''(\eta_b)}\,.
	$$
	Define
	$
		\bar R
		=
		q_0R_1+(1-q_0)R_0$, $
		\bar R_2
		=
		q_0R_1^2+(1-q_0)R_0^2$, and $
		\bar W
		=
		q_0W_1+(1-q_0)W_0
	$. 
	Then $(\sigma^*,\lambda^*,\bu^*,\theta^*)$ satisfies
	\begin{equation}
		\label{eq:FOCs_derivation}
		\begin{aligned}
			\expect[\bar R] &= 0,\\
			\bb P_{\bb \Gamma}
			\expect[\bz\bar R] &= \bb0,\\
			\lambda\expect[\bar W] &= \kappa,\\
			\frac{\sigma^2\kappa}{\lambda^2}
			&=
			\expect[\bar R_2]\,,
		\end{aligned}
	\end{equation}
	where all expectations in \eqref{eq:FOCs_derivation} are with respect to $(\bz,G)$, and the equations are evaluated at $(\sigma,\lambda,\bu,\theta) = (\sigma^*,\lambda^*,\bu^*,\theta^*)\,.$
	Moreover,
	$$
		r^*=\frac{\sigma^*\sqrt\kappa}{\lambda^*}, \quad 
		t^*=r^*\lambda^*=\sigma^*\sqrt\kappa\,.
	$$
Conversely, let
$$
    (\sigma,\lambda,\bu,\theta)
    \in
    (0,\infty)^2
    \times
    \range(\bb\Gamma)
    \times
    \Re\,,
$$
and suppose that it satisfies the four equations in \eqref{eq:FOCs}, with
the starred quantities replaced by $(\sigma,\lambda,\bu,\theta)$. Then
$$
    0
    <
    \frac{\sigma\sqrt\kappa}{\lambda}
    <
    1 \,.
$$
Setting
$$
    r
    =
    \frac{\sigma\sqrt\kappa}{\lambda},
    \quad
    t
    =
    r\lambda
    =
    \sigma\sqrt\kappa\,,
$$
recovers the unique natural-domain saddle, and hence
$$
    (\sigma,r,\bu,\theta,t)
    =
    (\sigma^*,r^*,\bu^*,\theta^*,t^*) \,.
$$
In particular,
$(\sigma^*,\lambda^*,\bu^*,\theta^*)$ is the unique solution of the four
equations in \eqref{eq:FOCs} in
$$
    (0,\infty)^2
    \times
    \range(\bb\Gamma)
    \times
    \Re \,.
$$
\end{lemma}

\begin{proof}
	By Lemma~\ref{lemma:limiting_AO_natural_domain}, all first derivatives
	of $F$ vanish at the limiting saddle, with the $\bu$-derivative
	understood relative to $\range(\bb\Gamma)$.
	By \eqref{eq:F_limit_def}, after setting
	$$
		\lambda=\frac tr,
		\quad
		t=r\lambda\,,
	$$
	the limiting AO objective can be written in the variables
	$$
		a=(\sigma,\bu,\theta),
		\quad
		r,
		\quad
		\lambda\,,
	$$
	as
	\begin{equation}
		\label{eq:limiting_AO_transformed_objective}
		\mathcal F(a,r,\lambda)
		=
		\expect[
			e_\lambda(\xi,\Ytil)
		]
		-
		\sigma r\sqrt\kappa
		+
		\frac{r^2\lambda}{2}\,,
	\end{equation}
	where $\xi=\theta+\bz^\top\bu+\sigma G$
	and
	$$
		e_\lambda(x,v)
		=
		m_\rho(x+\lambda v,\lambda)-vx-\frac{\lambda}{2}v^2\,.
	$$
	Because $r^*>0$, the change of variables $(r,t)\mapsto(r,\lambda=t/r)$ is smooth and invertible in a neighbourhood of the saddle. 
    Hence, 
	$$
		\partial_\lambda\mathcal F
		=
		r\,\partial_tF,
		\quad
		\partial_r\mathcal F
		=
		\partial_rF+\lambda\,\partial_tF\,,
	$$
	and the first-order conditions in $(r,t)$ obtained above are equivalent
	to the first-order conditions in $(r,\lambda)$ at the limiting saddle.
	The conditions in $\sigma$, $\theta$, and the admissible
	$\bu$-directions are unchanged by this coordinate transformation.
	Let $\eta = \prox{\lambda\rho}{\xi+\lambda \Ytil}$
	and set $R=\rho'(\eta)-\Ytil\,.$
	By \eqref{eq:limiting_AO_shifted_calculus},
	$$
		\eta
		=
		\xi-\lambda R,
		\quad
		\partial_xe_\lambda(\xi,\Ytil)
		=
		R\,.
	$$
	Therefore $\partial_\theta\mathcal F = \expect[R]\,.$
	The first-order condition in $\theta$ gives
	\begin{equation}
		\label{eq:foc_theta_raw}
		\expect[R]=0\,.
	\end{equation}
	For a direction $\bb d\in\range(\bb\Gamma)$,
	$$
		\left.
		\frac{d}{ds}
		\mathcal F(\sigma,r,\bu+s\bb d,\theta,\lambda)
		\right|_{s=0}
		=
		\expect[(\bz^\top\bb d)R]
		=
		\bb d^\top\expect[\bz R]\,.
	$$
	Since this derivative vanishes for every $\bb d\in\range(\bb\Gamma)$,
	\begin{equation}
		\label{eq:foc_u_raw}
		\bb P_{\bb \Gamma}
		\expect[\bz R]=\bb0\,.
	\end{equation}
	The first-order condition in $\sigma$ is
	\begin{equation}
		\label{eq:foc_sigma_start}
		\partial_\sigma\mathcal F
		=
		\expect[GR]-r\sqrt\kappa
		=
		0\,.
	\end{equation}
	Conditional on $(\bz,\varepsilon)$, the random variable $\Ytil$ is fixed and independent of $G$, and the dependence of $R$ on $G$ is only through $\xi=\theta+\bz^\top\bu+\sigma G\,.$
	By \eqref{eq:limiting_AO_shifted_calculus},
	$$
		\frac{\partial R}{\partial \xi}
		=
		\frac{\rho''(\eta)}
		{1+\lambda\rho''(\eta)}\,.
	$$
	Hence
	$$
		\frac{\partial R}{\partial G}
		=
		\sigma
		\frac{\rho''(\eta)}
		{1+\lambda\rho''(\eta)}\,.
	$$
	Since $0\leq\Ytil\leq1$, $0<\rho'(\eta)<1$, and
	$0\leq\rho''(\eta)\leq1/4$, one has $\abs{R}\leq1$
	and
	$$
		0
		\leq
		\frac{\partial R}{\partial\xi}
		=
		\frac{\rho''(\eta)}{1+\lambda\rho''(\eta)}
		\leq
		\frac14\,.
	$$
	Thus, 
	Stein's lemma \citep[Lemma~1]{stein:1981} yields
	$$
		\expect[GR]
		=
		\sigma
		\expect\left[
			\frac{\rho''(\eta)}
			{1+\lambda\rho''(\eta)}
		\right]\,.
	$$
	Substituting into \eqref{eq:foc_sigma_start} gives
	\begin{equation}
		\label{eq:foc_sigma_raw}
		\sigma
		\expect\left[
			\frac{\rho''(\eta)}
			{1+\lambda\rho''(\eta)}
		\right]
		=
		r\sqrt\kappa\,.
	\end{equation}
	Next, \eqref{eq:limiting_AO_shifted_calculus} gives
	$$
		\partial_\lambda e_\lambda(\xi,\Ytil)
		=
		-\frac12R^2\,.
	$$
	Thus
	$$
		\partial_\lambda\mathcal F
		=
		-\frac12\expect[R^2]
		+
		\frac{r^2}{2}\,.
	$$
	The first-order condition in $\lambda$ gives
	\begin{equation}
		\label{eq:foc_lambda_raw}
		r^2=\expect[R^2]\,.
	\end{equation}
	The first-order condition in $r$ gives
	$$
		\partial_r\mathcal F
		=
		-\sigma\sqrt\kappa+r\lambda
		=
		0\,,
	$$
	and therefore
	\begin{equation}
		\label{eq:foc_r_raw}
		r\lambda=\sigma\sqrt\kappa\,.
	\end{equation}
	Since $\sigma>0$ by Lemma~\ref{lemma:limiting_AO_natural_domain}, combining \eqref{eq:foc_sigma_raw} and \eqref{eq:foc_r_raw} yields
	\begin{equation}
		\label{eq:foc_curvature_raw}
		\lambda
		\expect\left[
			\frac{\rho''(\eta)}
			{1+\lambda\rho''(\eta)}
		\right]
		=
		\kappa\,.
	\end{equation}
	Combining \eqref{eq:foc_lambda_raw} and \eqref{eq:foc_r_raw} gives
	\begin{equation}
		\label{eq:foc_residual_raw}
		\frac{\sigma^2\kappa}{\lambda^2}
		=
		\expect[R^2]\,.
	\end{equation}
	Also \eqref{eq:foc_r_raw} gives
	$$
		r=\frac{\sigma\sqrt\kappa}{\lambda}\,,
	$$
	and, since $t=r\lambda$, $t=\sigma\sqrt\kappa$.
	It remains to integrate out the Bernoulli randomness. Conditional on $(\bz,G)$, the pseudo-response is $\Ytil^{(1)}$ with probability $q_0$ and $\Ytil^{(0)}$ with probability $1-q_0$. Hence
	$$
		\expect[R\mid\bz,G]
		=
		q_0R_1+(1-q_0)R_0
		=
		\bar R\,,
	$$
	$$
		\expect[R^2\mid\bz,G]
		=
		q_0R_1^2+(1-q_0)R_0^2
		=
		\bar R_2\,,
	$$
	and
	$$
		\expect\left[
			\left.
			\frac{\rho''(\eta)}
			{1+\lambda\rho''(\eta)}
			\right|
			\bz,G
		\right]
		=
		q_0W_1+(1-q_0)W_0
		=
		\bar W\,.
	$$
	Substituting these conditional expectations into
	\eqref{eq:foc_theta_raw}, \eqref{eq:foc_u_raw},
	\eqref{eq:foc_curvature_raw}, and \eqref{eq:foc_residual_raw}
	gives exactly \eqref{eq:FOCs_derivation}. Evaluating at $(\sigma,\lambda,\bu,\theta) = (\sigma^*,\lambda^*,\bu^*,\theta^*)$
	proves the necessity statements.

For the converse, let
$$
    (\sigma,\lambda,\bu,\theta)
    \in
    (0,\infty)^2
    \times
    \range(\bb\Gamma)
    \times
    \Re\,,
$$
satisfy the four equations in \eqref{eq:FOCs}, with the starred quantities
replaced by $(\sigma,\lambda,\bu,\theta)$. Let
$$
    \xi
    =
    \theta+\bz^\top\bu+\sigma G,
    \quad
    \eta
    =
    \prox{\lambda\rho}{\xi+\lambda\Ytil}\,,
$$
and set
$$
    R
    =
    \rho'(\eta)-\Ytil,
    \quad
    W
    =
    \frac{\rho''(\eta)}
    {1+\lambda\rho''(\eta)} \,.
$$
Since $\alpha\in(0,1)$, the definition of the limiting pseudo-response in Section~\ref{sec:limiting-ao-focs}
gives $ 0<\Ytil<1$ and since $0<\rho'(\eta)<1$, we have $\abs{R}<1$.
The fourth equation in \eqref{eq:FOCs} therefore gives
$$
    0
    <
    \frac{\sigma^2\kappa}{\lambda^2}
    =
    \expect[R^2]
    <
    1 \,.
$$
Consequently,
$
    r
    =
    \sigma\sqrt\kappa / \lambda 
    \in(0,1)\,.
$
Setting
$
    t
    =
    r\lambda
    =
    \sigma\sqrt\kappa
$ we see that 
$(\sigma,r,\bu,\theta,t)$ belongs to the natural domain.
The first two equations in \eqref{eq:FOCs} and
\eqref{eq:limiting_AO_gradient} give
$$
    \partial_\theta F
    =
    0,
    \quad
    \bb P_{\bb \Gamma}
    \nabla_{\bu}F
    =
    \bb0 \,.
$$
Stein's identity used above gives
$$
    \expect[GR]
    =
    \sigma\expect[W] \,.
$$
Using the third equation in \eqref{eq:FOCs} and the definition of $r$,
$$
\begin{aligned}
    \expect[GR]-r\sqrt\kappa
    =
    \sigma\expect[W]
    -
    \frac{\sigma\sqrt\kappa}{\lambda}
    \sqrt\kappa
    =
    \frac{\sigma\kappa}{\lambda}
    -
    \frac{\sigma\kappa}{\lambda}
    =
    0 \,.
\end{aligned}
$$
Hence
$
    \partial_\sigma F=0
$.
The fourth equation in \eqref{eq:FOCs} gives
$$
    \expect[R^2]
    =
    \frac{\sigma^2\kappa}{\lambda^2}
    =
    r^2 \,.
$$
The formula for $\partial_tF$ in
\eqref{eq:limiting_AO_gradient} therefore yields
$$
    \partial_tF
    =
    -\frac{1}{2r}\expect[R^2]
    +
    \frac r2
    =
    0 \,.
$$
Finally,
\eqref{eq:limiting_AO_gradient} gives
$$
\begin{aligned}
    \partial_rF
    =
    \frac{t}{2r^2}\expect[R^2]
    -
    \sigma\sqrt\kappa
    +
    \frac t2
    =
    t-\sigma\sqrt\kappa
    =
    0 \,.
\end{aligned}
$$
Thus all relative first-order conditions in the minimising coordinates and
the first-order condition in the maximising coordinate are satisfied.

By Lemma~\ref{lemma:limiting_AO_geometry}, for this fixed $r$ the objective
is strictly convex in $(\sigma,\bu,\theta,t)$, while for the fixed
minimising coordinates it is strictly concave in $r$. Hence
$(\sigma,r,\bu,\theta,t)$ is the unique natural-domain saddle.
Lemma~\ref{lemma:limiting_AO_natural_domain} identifies that saddle with
$$
    (\sigma^*,r^*,\bu^*,\theta^*,t^*) \,.
$$
It follows that
$$
    (\sigma,\lambda,\bu,\theta)
    =
    (\sigma^*,\lambda^*,\bu^*,\theta^*) \,,
$$
which proves the converse and uniqueness on the stated positive domain.
\end{proof}

\subsection{Lemmas for Section~\ref{sec:cgmt-cost-comparisons}}

\subsubsection{Finite convex cover and limiting gaps}

\begin{lemma}[Finite convex cover of the signal-intercept bad set]
\label{lemma:finite_convex_bad_cover}
	For every $\epsilon>0$,
	$$
		\mathcal S_{\bv,\theta}^{\epsilon,c}
		\subseteq
		\bigcup_{j\in\mathcal J}C_{p,j}^{\epsilon}\,.
	$$
	Moreover, each nonempty $C_{p,j}^{\epsilon}$ is compact and convex.
\end{lemma}

\begin{proof}
For $i\in\mathcal I_{\mathcal A}$, the set
$C_{p,i}^{\epsilon}$ is the intersection of
$V_p\times[-C_\theta,C_\theta]$ with a closed affine half-space. For
$j\in\{+,-\}$,
$C_{p,j}^{\epsilon}=V_p\times\Theta_j^\epsilon$.
The set $V_p$ is compact and convex, and the
intervals are closed, possibly empty. Hence every nonempty
$C_{p,j}^{\epsilon}$ is compact and convex.

Let $(\bv,\theta)\in\mathcal S_{\bv,\theta}^{\epsilon,c}$. Then either
$$
    \vnorm{\bb\Gamma_p^{1/2}\bv-\bu^*}_2\geq\epsilon\,,
$$
or $\abs{\theta-\theta^*}\geq\epsilon$. If
$\theta-\theta^*\geq\epsilon$, then
$(\bv,\theta)\in C_{p,+}^{\epsilon}$. If
$\theta^*-\theta\geq\epsilon$, then
$(\bv,\theta)\in C_{p,-}^{\epsilon}$. Otherwise, set
$$
    \bb a_0
    =
    \frac{\bb\Gamma_p^{1/2}\bv-\bu^*}
    {\vnorm{\bb\Gamma_p^{1/2}\bv-\bu^*}_2}
    \in\mathbb S^1 \,.
$$
Choose $i\in\mathcal I_{\mathcal A}$ such that
$\vnorm{\bb a_i-\bb a_0}_2\leq1/2$. Then
$$
\begin{aligned}
    \bb a_i^\top
    (\bb\Gamma_p^{1/2}\bv-\bu^*)
    \geq
    \bb a_0^\top
    (\bb\Gamma_p^{1/2}\bv-\bu^*)
    -
    \vnorm{\bb a_i-\bb a_0}_2
    \vnorm{\bb\Gamma_p^{1/2}\bv-\bu^*}_2
    \geq
    \frac{\epsilon}{2} \,.
\end{aligned}
$$
Thus $(\bv,\theta)\in C_{p,i}^{\epsilon}$, proving the cover.
\end{proof}

\begin{lemma}[Finite-piece limiting AO value gaps]
\label{lemma:limiting_AO_gap}
Let $(\sigma^*,r^*,\bu^*,\theta^*,t^*)$ be the unique compact saddle
chain from Proposition~\ref{prop:limiting_AO_properties}. Then, for every
$\epsilon>0$, there exists a deterministic $\eta_\epsilon>0$, independent
of $n$, such that
$$
    \bar\phi_{\sigma}^{\epsilon,c}\ge\bar\phi+3\eta_\epsilon,
    \quad
    \bar\phi_j^\epsilon\ge\bar\phi+3\eta_\epsilon,
    \quad j\in\mathcal J\,.
$$
\end{lemma}

\begin{proof}
We first prove a strict gap for every nonempty restricted domain.
If $\bb\Omega_{\sigma}^{\epsilon,c}$ is empty, then
$\bar\phi_{\sigma}^{\epsilon,c}=+\infty$. Otherwise,
\eqref{eq:limiting_restricted_values} gives
$$
    \bar\phi_{\sigma}^{\epsilon,c}
    =
    \underset{\sigma\in\mathcal S_{\sigma}^{\epsilon,c}}{\min}
    V_0(\sigma) \,.
$$
The minimum is attained because $V_0=V$ by
Lemma~\ref{lemma:ambient_limiting_saddle_chain} and $V$ is continuous by
Lemma~\ref{lemma:limiting_AO_natural_domain}. Suppose that
$\bar\phi_{\sigma}^{\epsilon,c}=\bar\phi$. Then there is
$\tilde\sigma\in\mathcal S_{\sigma}^{\epsilon,c}$ such that
$V_0(\tilde\sigma)=\bar\phi$. Since
$$
    \bar\phi
    =
    \underset{\sigma\in[0,C_\sigma]}{\min}V_0(\sigma)\,,
$$
and $V_0$ has the unique minimiser $\sigma^*$ by
Lemma~\ref{lemma:ambient_limiting_saddle_chain}, we obtain
$\tilde\sigma=\sigma^*$, contradicting
$\tilde\sigma\in\mathcal S_{\sigma}^{\epsilon,c}$. Hence
$\bar\phi_{\sigma}^{\epsilon,c}>\bar\phi$.
Now fix $j\in\mathcal J$. If $\bb\Omega_j^\epsilon$ is empty, then
$\bar\phi_j^\epsilon=+\infty$. Otherwise, define
$$
    H_j^\epsilon(\sigma,r)
    =
    \underset{(\bu,\theta,t)\in\bb K_j^\epsilon(U)}{\min}
    F(\sigma,r,\bu,\theta,t),
    \quad
    V_j^\epsilon(\sigma)
    =
    \underset{r\in[c_r,1]}{\max}
    H_j^\epsilon(\sigma,r) \,.
$$
The set $\bb K_j^\epsilon(U)$ is nonempty and compact. Since $F$ is
continuous, Berge's maximum theorem gives continuity and attainment of
$H_j^\epsilon$ on
$[0,C_\sigma]\times[c_r,1]$. A second application gives continuity and
attainment of $V_j^\epsilon$ on $[0,C_\sigma]$. Because
$\bb K_j^\epsilon(U)\subseteq\bb K$,
$$
    H_j^\epsilon(\sigma,r)
    \geq
    H_0(\sigma,r),
    \quad
    V_j^\epsilon(\sigma)
    \geq
    V_0(\sigma) \,.
$$
Suppose that $\bar\phi_j^\epsilon=\bar\phi$. By attainment, there is
$\tilde\sigma\in[0,C_\sigma]$ such that
$V_j^\epsilon(\tilde\sigma)=\bar\phi$. Therefore
$$
    \bar\phi
    =
    V_j^\epsilon(\tilde\sigma)
    \geq
    V_0(\tilde\sigma)
    \geq
    \bar\phi \,.
$$
Thus $V_0(\tilde\sigma)=\bar\phi$, and
Lemma~\ref{lemma:ambient_limiting_saddle_chain} gives
$\tilde\sigma=\sigma^*$. Hence
$V_j^\epsilon(\sigma^*)=\bar\phi$. At $r=r^*$,
$$
    H_j^\epsilon(\sigma^*,r^*)
    \geq
    H_0(\sigma^*,r^*)
    =
    \bar\phi\,,
$$
by Lemma~\ref{lemma:ambient_limiting_saddle_chain}, while
$$
    H_j^\epsilon(\sigma^*,r^*)
    \leq
    V_j^\epsilon(\sigma^*)
    =
    \bar\phi \,.
$$
Consequently,
$H_j^\epsilon(\sigma^*,r^*)=\bar\phi$. The inner minimum is attained, so
there is
$(\tilde\bu,\tilde\theta,\tilde t)\in\bb K_j^\epsilon(U)$ with
$$
    F(\sigma^*,r^*,\tilde\bu,\tilde\theta,\tilde t)
    =
    \bar\phi \,.
$$
By Lemma~\ref{lemma:ambient_limiting_saddle_chain}, the unique minimiser
of $F(\sigma^*,r^*,\cdot)$ over $\bb K$ is
$(\bu^*,\theta^*,t^*)$. Hence
$(\tilde\bu,\tilde\theta,\tilde t)=(\bu^*,\theta^*,t^*)$. This point lies
in none of the bad pieces: if $j=i\in\mathcal I_{\mathcal A}$, then
$$
    \bb a_i^\top(\bu^*-\bu^*)
    =
    0
    <
    \frac{\epsilon}{2}\,,
$$
while $\theta^*\notin\Theta_+^\epsilon$ and
$\theta^*\notin\Theta_-^\epsilon$. This is a contradiction. Therefore $\bar\phi_j^\epsilon>\bar\phi$ for every nonempty
piece.
Define
$$
    \mathfrak D_\epsilon
    =
    \left\{
        \bar\phi_{\sigma}^{\epsilon,c}-\bar\phi:
        \bb\Omega_{\sigma}^{\epsilon,c}\neq\emptyset
    \right\}
    \cup
    \left\{
        \bar\phi_j^\epsilon-\bar\phi:
        j\in\mathcal J,\,
        \bb\Omega_j^\epsilon\neq\emptyset
    \right\} \,.
$$
If $\mathfrak D_\epsilon=\emptyset$, set $\eta_\epsilon=1$. Otherwise set
$$
    \eta_\epsilon
    =
    \frac13\min\mathfrak D_\epsilon \,.
$$
Every element of $\mathfrak D_\epsilon$ is strictly positive, and
$\mathcal J$ is finite. Hence $\eta_\epsilon>0$, and the stated bounds
follow. Empty restricted domains also satisfy them by convention.
\end{proof}

\subsubsection{Restricted scalarisation and AO barriers}

\begin{lemma}[Large-$t$ exclusion on restricted convex pieces]
\label{lemma:restricted_large_t_pieces}
	Fix $\epsilon>0$, $c_r>0$, and $C_t>0$. For $\bb\eta\in\mathcal B^n_{C_\eta}$, $(\bv,\theta)\in V_p\times[-C_\theta,C_\theta]$, $\sigma\in[0,C_\sigma]$, $r\in[c_r,1]$, and $t>0\,,$
	define
	$$
		\mathcal T_n(\sigma,r,\bv,\theta,\bb\eta,t)
		=
		\frac{r}{2t}
		\frac1n
		\vnorm{
			\bb\eta
			-
			\theta\bb1
			-
			\frac1{\sqrt p}\bH_1\bv
			-
			\sigma\bg
		}_2^2
		+
		\frac{rt}{2}
		+
		\frac1n
		\left\{
			\bb1^\top\bb\rho(\bb\eta)
			-
			\bYtil^\top\bb\eta
		\right\}\,.
	$$
	Let $\mathcal E_{n,\mathrm{res}}^{t,\epsilon}(c_r,C_t)$ be the event
	that, simultaneously for every $j\in\mathcal J_n^\epsilon$, every
	$\sigma\in[0,C_\sigma]$, and every $r\in[c_r,1]$,
	$$
	\begin{aligned}
		&\underset{\substack{(\bv,\theta)\in C_{p,j}^{\epsilon}\\
		\bb\eta\in\mathcal B^n_{C_\eta}\\
		t>0}}{\inf}
		\, \mathcal T_n(\sigma,r,\bv,\theta,\bb\eta,t)
		 =
		\underset{\substack{(\bv,\theta)\in C_{p,j}^{\epsilon}\\
		\bb\eta\in\mathcal B^n_{C_\eta}\\
		t\in(0,C_t]}}{\inf}
		\, \mathcal T_n(\sigma,r,\bv,\theta,\bb\eta,t)\,.
	\end{aligned}
	$$
	Recall $K_t$ from Lemma~\ref{lemma:t_ub}.
	Suppose that
	$$
		C_t>
		\max
		\left\{
			1,\,
			\frac{2(K_t+\log2)}{c_r}
		\right\}\,.
	$$
	Then there exist constants $C,c>0$ and $N\in\mathbb N$ such that, for
	all $n>N$\,,
	$$
		\Pr\left(
			\lnot\mathcal E_{n,\mathrm{res}}^{t,\epsilon}(c_r,C_t)
			\mid
			\mathcal G_n
		\right)
		\leq
		C\exp\{-cn\}\,.
	$$
	In particular, for every fixed $t_0>0$,
	$$
		\mathds{1}_{\mathcal S_n(t_0)}
		\Pr\left(
			\lnot\mathcal E_{n,\mathrm{res}}^{t,\epsilon}(c_r,C_t)
			\mid
			\mathcal G_n
		\right)
		\leq
		C\exp\{-cn\}\,.
	$$
	The constants depend only on
	$c_r$, $C_\beta$, $C_\theta$, $C_\sigma$
	and universal Gaussian concentration constants.
\end{lemma}

\begin{proof}
	Define the event
	$$
		\mathcal I_{n}^{t,\mathrm{res}}
		=
		\left\{
			\frac{\vnorm{\bg}_2}{\sqrt n}\leq 2
		\right\}
		\cap
		\left\{
			\underset{\bv\in V_p}{\sup}
			\frac1{\sqrt n}
			\vnorm{
				\frac1{\sqrt p}\bH_1\bv
			}_2
			\leq
			3C_\beta
		\right\}\,.
	$$
	We first show that
	$$
		\Pr\left(
			\lnot\mathcal I_{n}^{t,\mathrm{res}}
			\mid
			\mathcal G_n
		\right)
		\leq
		C\exp\{-cn\}\,.
	$$
	The first event follows from chi-square concentration
	\citep[Equation~(3.1)]{vershynin:2018} for $\bg\sim N(0,I_n)$.

	For the second event, use \eqref{eq:ao_H_repar}. On the common
	extension,
	$$
		\frac1{\sqrt p}\bH_1
		=
		\bb Z\bb\Gamma_p^{1/2}\,,
	$$
	where $\bb Z\in\mathbb R^{n\times2}$ has i.i.d. standard normal entries
	conditional on $\mathcal G_n$.
	For $\bv\in V_p$, set $\bu=\bb\Gamma_p^{1/2}\bv \,.$
	By the definition of $V_p$,
	$$
		\vnorm{\bu}_2^2
		=
		\bv^\top\bb\Gamma_p\bv
		=
		\frac1p\vnorm{\bB\bv}_2^2
		\leq
		C_\beta^2 \,.
	$$
	Therefore, on the common extension,
	$$
		\underset{\bv\in V_p}{\sup}
		\frac1{\sqrt n}
		\vnorm{
			\frac1{\sqrt p}\bH_1\bv
		}_2
		\leq
		C_\beta
		\frac{\mnorm{\bb Z}_2}{\sqrt n}\,.
	$$
	The standard Gaussian operator-norm bound for an $n\times2$ matrix gives
	$$
		\Pr\left(
			\frac{\mnorm{\bb Z}_2}{\sqrt n}>3
			\right)
			\leq
			C\exp\{-cn\}\,,
		$$
	for all large enough $n$. This proves the probability bound for
	$\mathcal I_{n}^{t,\mathrm{res}}$.
	Now work on $\mathcal I_{n}^{t,\mathrm{res}}$. Fix $j\in\mathcal J_n^\epsilon$, $\sigma\in[0,C_\sigma]$, and $r\in[c_r,1]\,.$
	Since $j\in\mathcal J_n^\epsilon$, the set $C_{p,j}^{\epsilon}$ is
	nonempty. Choose any $(\bv_j,\theta_j)\in C_{p,j}^{\epsilon}\,.$
	Then $\abs{\theta_j}\leq C_\theta, \quad \bv_j\in V_p\,.$
	Hence, on $\mathcal I_{n}^{t,\mathrm{res}}$,
	$$
		\frac1{\sqrt n}
		\vnorm{
			\theta_j\bb1
			+
			\frac1{\sqrt p}\bH_1\bv_j
			+
			\sigma\bg
		}_2
		\leq
		C_\theta+3C_\beta+2C_\sigma \,.
	$$
	Test the feasible point
	$$
		\bv=\bv_j,
		\quad
		\theta=\theta_j,
		\quad
		\bb\eta=0,
		\quad
		t=1\,.
	$$
	Since $C_t>1$, this point is feasible also for the truncated problem.
	At this point,
	$$
		\frac1n
		\left\{
			\bb1^\top\bb\rho(0)-\bYtil^\top0
		\right\}
		=
		\log2\,.
	$$
	Using $r\leq1$, we obtain
	$$
	\begin{aligned}
		\mathcal T_n(\sigma,r,\bv_j,\theta_j,0,1)
		&=
		\frac r2
		\frac1n
		\vnorm{
			\theta_j\bb1
			+
			\frac1{\sqrt p}\bH_1\bv_j
			+
			\sigma\bg
		}_2^2
		+
		\frac r2
		+
		\log2
		\\
		&\leq
		\frac12
		\left(
			C_\theta+3C_\beta+2C_\sigma
		\right)^2
		+
		\frac12
		+
		\log2
		\\
		&=
		K_t\,.
	\end{aligned}
	$$
	Thus the truncated infimum is at most $K_t$.
	On the other hand, for every $y\in(0,1)$ and every $x\in\mathbb R$,
	$$
		\rho(x)-yx
		\geq
		\underset{z\in\mathbb R}{\inf}\{\rho(z)-yz\}
		=
		-y\log y -(1-y)\log(1-y) \in [0, \log(2)]\,.
	$$
	Therefore, for every $\bb\eta\in\mathbb R^n$,
	$$
		\frac1n
		\left\{
			\bb1^\top\bb\rho(\bb\eta)
			-
			\bYtil^\top\bb\eta
		\right\}
		\geq
		0 \,. 
	$$
	Since the square term in $\mathcal T_n$ is nonnegative and
	$r\geq c_r$, every feasible point with $t>C_t$ satisfies
	$$
		\mathcal T_n(\sigma,r,\bv,\theta,\bb\eta,t)
		\geq
		\frac{c_r t}{2}
		\geq
		\frac{c_r C_t}{2}
		>
		K_t\,.
	$$
	Such a point cannot contribute to the infimum, because the feasible point
	$(\bv_j,\theta_j,0,1)$ has value at most $K_t$.
	Consequently
	$$
	\begin{aligned}
		&\underset{\substack{(\bv,\theta)\in C_{p,j}^{\epsilon}\\
		\bb\eta\in\mathcal B^n_{C_\eta}\\
		t>0}}{\inf}
		\, \mathcal T_n(\sigma,r,\bv,\theta,\bb\eta,t) =
		\underset{\substack{(\bv,\theta)\in C_{p,j}^{\epsilon}\\
		\bb\eta\in\mathcal B^n_{C_\eta}\\
		t\in(0,C_t]}}{\inf}
		\, \mathcal T_n(\sigma,r,\bv,\theta,\bb\eta,t)\,.
	\end{aligned}
	$$
	The argument was uniform in $j\in\mathcal J_n^\epsilon$,
	$\sigma\in[0,C_\sigma]$, and $r\in[c_r,1]$. Hence
	$$
		\mathcal I_{n}^{t,\mathrm{res}}
		\subseteq
		\mathcal E_{n,\mathrm{res}}^{t,\epsilon}(c_r,C_t)\,.
	$$
	The claimed conditional probability bound follows from the probability
	bound for $\mathcal I_{n}^{t,\mathrm{res}}$.
\end{proof}

\begin{lemma}[Simultaneous restricted scalarisation]
\label{lemma:restricted_piece_scalarisation}
Fix $\epsilon>0$ and $t_0\in(0,t_\star]$. Let
$
    d_r$, 
    $
    a$, 
    $
    K$, 
    $
    C_L$, 
    $
    c_r$, 
    $
    C_t$, 
    $
    C_M$, $
    C_\eta
$
be the constants fixed from Proposition~\ref{prop:scalarisation_events}
with structural tolerance $t_\star$, and write
$$
    \mathcal E_n^{\mathrm{sr}}
    =
    \mathcal E_n^{\mathrm{sr}}(d_r,a,K,C_L)\,.
$$
For every $j\in\mathcal J_n^\epsilon$, define the restricted predictor set
$$
    R_j^\epsilon
    =
    \left\{
        \theta\bb1+\frac1{\sqrt p}\bH_1\bv:
        (\bv,\theta)\in C_{p,j}^\epsilon
    \right\}\,.
$$
Then there exists an event $\mathcal R_{n,\mathrm{res}}^\epsilon$ and
constants $C_\epsilon,c_\epsilon>0$ and
$N_\epsilon\in\mathbb N$ such that, for all $n>N_\epsilon$,
$$
    \mathds{1}_{\mathcal S_n(t_0)}
    \Pr\left(
        \lnot\mathcal R_{n,\mathrm{res}}^\epsilon
        \mid
        \mathcal G_n
    \right)
    \leq
    C_\epsilon\exp\{-c_\epsilon n\}\,.
$$
On $\mathcal R_{n,\mathrm{res}}^\epsilon$,
$$
    \phi(\bg,\bh)
    =
    \operatorname{val}
    \left(
        F_n;
        \bb\Omega_n
    \right)\,,
$$
$$
    \phi_{\bw}^{\epsilon,c}(\bg,\bh)
    =
    \operatorname{val}
    \left(
        F_n;
        \bb\Omega_{n,\sigma}^{\epsilon,c}
    \right)\,,
$$
with both sides interpreted as $+\infty$ when
$\mathcal S_{\bw}^{\epsilon,c}$ is empty, and, simultaneously for every
$j\in\mathcal J_n^\epsilon$,
$$
    \phi_j^\epsilon(\bg,\bh)
    =
    \operatorname{val}
    \left(
        F_n;
        \bb\Omega_{n,j}^\epsilon
    \right)\,.
$$
\end{lemma}

\begin{proof}
Let
$$
\begin{aligned}
    \mathcal R_{n,\mathrm{res}}^\epsilon
    =
    \mathcal R_n
    \cap
    \mathcal E_n^{\mathrm{sr}}
    \cap
    \mathcal I_n^t
    \cap
    \mathcal E_{n,\mathrm{res}}^{t,\epsilon}(c_r,C_t)
    \cap
    \mathcal E_n^{\mathrm{ctr}}(C_M)\,.
\end{aligned}
$$
The bounds of Proposition~\ref{prop:scalarisation_events}, obtained with
structural tolerance $t_\star$, remain valid after multiplication by
$\mathds{1}_{\mathcal S_n(t_0)}$ because
$\mathcal S_n(t_0)\subseteq\mathcal S_n(t_\star)$. The probability bounds
for the remaining three events follow from Lemmas~\ref{lemma:t_ub},
\ref{lemma:restricted_large_t_pieces}, and
\ref{lemma:moreau_center_probability}. A union bound therefore gives the
claimed conditional estimate.

The unrestricted identity follows from $\mathcal R_n$. For the
$\bw$-restricted AO, the radial parametrisation
$\bw=\sqrt p\,\sigma\bb q$ changes only the outer $\sigma$-domain from
$[0,C_\sigma]$ to $\mathcal S_\sigma^{\epsilon,c}$. On
$\mathcal E_n^{\mathrm{sr}}$, Proposition~\ref{prop:small_r_exclusion}
provides the pointwise small-$r$ reduction for
$\mathcal R=\mathcal R_0$; on $\mathcal I_n^t$,
Lemma~\ref{lemma:t_ub} provides the pointwise $t$-reduction uniformly in
$\sigma$; and $\mathcal E_n^{\mathrm{ctr}}(C_M)$ gives pointwise Moreau
feasibility. Hence the same scalarisation argument gives
$$
    \phi_{\bw}^{\epsilon,c}(\bg,\bh)
    =
    \operatorname{val}
    \left(
        F_n;
        \bb\Omega_{n,\sigma}^{\epsilon,c}
    \right)\,.
$$

Fix $j\in\mathcal J_n^\epsilon$. By
Lemma~\ref{lemma:finite_convex_bad_cover}, the set
$Z_{p,j}^\epsilon=C_{p,j}^\epsilon\times\mathcal B^n_{C_\eta}$ is
nonempty, compact, and convex. Since the restricted predictor set
$R_j^\epsilon$ defined in the statement is the image of
$C_{p,j}^\epsilon$ under an affine map, it is also nonempty, compact, and
convex, and its definition gives $R_j^\epsilon\subseteq\mathcal R_0$.
Lemma~\ref{lemma:psi_conv_conc}, applied on this restricted set, therefore
gives the same ball-to-radius reduction as
in the unrestricted scalarisation. The resulting distance profile is
$M_{R_j^\epsilon}$.

On $\mathcal E_n^{\mathrm{sr}}$, Proposition~\ref{prop:small_r_exclusion}
applies simultaneously to every $R_j^\epsilon$. Since the final radius
$c_r$ fixed in Proposition~\ref{prop:scalarisation_events} satisfies
$c_r\leq r_0$, it permits the restriction
$r\in[c_r,1]$ for every piece. The event
$\mathcal E_{n,\mathrm{res}}^{t,\epsilon}(c_r,C_t)$ then permits the
restriction $t\in(0,C_t]$ simultaneously over all pieces.

Finally, on $\mathcal E_n^{\mathrm{ctr}}(C_M)$, the condition
$$
    C_\eta
    \geq
    C_M+\frac{C_t}{c_r}
$$
and Lemma~\ref{lemma:moreau_feasibility} show that the unconstrained
coordinatewise Moreau minimisers are feasible for
$\mathcal B^n_{C_\eta}$. This pointwise feasibility is inherited by every
restricted predictor set. Hence the completion-of-squares and
Moreau-envelope calculation of Section~\ref{sec:ao-scalar} applies without
change.

Under the Gaussian representation in \eqref{eq:ao_H_repar}, the map
$\bu=\bb\Gamma_p^{1/2}\bv$ sends
$C_{p,j}^\epsilon\times[0,C_t]$ exactly onto
$\bb K_j^\epsilon(U_n)$. The resulting objective is $F_n$ on
$\bb\Omega_{n,j}^\epsilon$, and therefore
$$
    \phi_j^\epsilon(\bg,\bh)
    =
    \operatorname{val}
    \left(
        F_n;
        \bb\Omega_{n,j}^\epsilon
    \right)\,.
$$
The construction is simultaneous over $j\in\mathcal J_n^\epsilon$.
\end{proof}

\begin{lemma}[Restricted AO value bounds on finite convex pieces]
\label{lemma:restricted_AO_value_bounds}
Recall the restricted scalar domains in
\eqref{eq:restricted_scalar_domains}, the restricted limiting values in
\eqref{eq:limiting_restricted_values}, and the restricted AO values defined
in Section~\ref{sec:cgmt-cost-comparisons}. For every $\epsilon>0$ and
$s>0$, define
$$
\begin{aligned}
    \mathcal A_n^{\epsilon}(s)
    =
    \left\{
        \abs{\phi(\bg,\bh)-\bar\phi}\leq s
    \right\}
    \cap
    \left\{
        \phi_{\bw}^{\epsilon,c}(\bg,\bh)
        \geq
        \bar\phi_{\sigma}^{\epsilon,c}-s
    \right\}
    \cap
    \bigcap_{j\in\mathcal J_n^\epsilon}
    \left\{
        \phi_j^\epsilon(\bg,\bh)
        \geq
        \bar\phi_j^\epsilon-s
    \right\} \,.
\end{aligned}
$$
Then there exist $t_0=t_0(\epsilon,s)>0$, constants
	$C_{\epsilon,s},c_{\epsilon,s}>0$, and
	$N_{\epsilon,s}\in\mathbb N$ such that, for all
	$n>N_{\epsilon,s}$,
	$$
		\mathds{1}_{\mathcal S_n(t_0)}
		\Pr\left(
			\lnot\mathcal A_n^{\epsilon}(s)
			\mid
			\mathcal G_n
		\right)
		\leq
		C_{\epsilon,s}\exp\{-c_{\epsilon,s}n\}\,.
	$$
\end{lemma}

\begin{proof}
	Set
	$$
		\delta=\frac{s}{3}\,.
	$$
	Let $\mathcal R_{n,\mathrm{res}}^\epsilon$ be the simultaneous
	scalarisation event from Lemma~\ref{lemma:restricted_piece_scalarisation}.
	On this event,
	$$
		\phi(\bg,\bh)
		=
		\operatorname{val}
		\left(
			F_n;
			\bb\Omega_n
		\right)\,,
	$$
	$$
		\phi_{\bw}^{\epsilon,c}(\bg,\bh)
		=
		\operatorname{val}
			\left(
				F_n;
				\bb\Omega_{n,\sigma}^{\epsilon,c}
			\right)\,,
		$$
	whenever $\mathcal S_{\bw}^{\epsilon,c}$ is nonempty, and
	$$
		\phi_j^\epsilon(\bg,\bh)
		=
		\operatorname{val}
		\left(
			F_n;
			\bb\Omega_{n,j}^{\epsilon}
		\right),
		\quad j\in\mathcal J_n^\epsilon\,.
	$$
	By Lemma~\ref{lemma:restricted_piece_scalarisation}, after possibly
	decreasing $t_0$ and adjusting constants,
	$$
		\mathds{1}_{\mathcal S_n(t_0)}
		\Pr\left(
			\lnot\mathcal R_{n,\mathrm{res}}^\epsilon
			\mid
			\mathcal G_n
		\right)
		\le
		C\exp\{-cn\}\,.
	$$

	By Proposition~\ref{prop:uniform_AO_control}, after replacing
	$t_0$ by the minimum of its current value and the structural
	tolerance furnished there, and after adjusting constants,
	$$
		\mathds{1}_{\mathcal S_n(t_0)}
		\Pr\left(
			\left\{\lnot\mathcal U_n(\delta)\right\}
			\cup
			\left\{\lnot\mathcal N_n(\delta)\right\}
			\mid
			\mathcal G_n
		\right)
		\le
		C_\delta\exp\{-c_\delta n\}\,.
	$$

	Work on
	$$
		\mathcal R_{n,\mathrm{res}}^\epsilon
		\cap
		\mathcal U_n(\delta)
		\cap
		\mathcal N_n(\delta)\,.
	$$
	For the unrestricted AO value, the standard implication from
	\eqref{eq:AO_as_val_Fn}, \eqref{eq:value_lipschitz}, and
	\eqref{event:N_n} gives $\abs{\phi(\bg,\bh)-\bar\phi} \le 3\delta=s\,.$
	In particular, $\phi(\bg,\bh)\le\bar\phi+s\,.$

	Next consider the $\bw$-restricted AO. If
	$\mathcal S_{\bw}^{\epsilon,c}$ is empty, there is nothing to prove.
	Otherwise,
	$$
		\phi_{\bw}^{\epsilon,c}(\bg,\bh)
		=
		\operatorname{val}
		\left(
			F_n;
			\bb\Omega_{n,\sigma}^{\epsilon,c}
		\right) \,.
	$$
	Since
	$$
		\bb\Omega_{n,\sigma}^{\epsilon,c}
		\subseteq
		\bb\Omega_{\sigma}^{\epsilon,c}\,,
		$$
	only through the innermost minimisation domain in $\bu$, enlarging the
	$\bu$-domain can only decrease the min--max--min value. Hence
	$$
		\operatorname{val}
		\left(
			F_n;
			\bb\Omega_{n,\sigma}^{\epsilon,c}
		\right)
		\ge
		\operatorname{val}
		\left(
			F_n;
			\bb\Omega_{\sigma}^{\epsilon,c}
		\right) \,.
	$$
	On $\mathcal U_n(\delta)$, \eqref{eq:value_lipschitz} gives
$$
    \operatorname{val}
    \left(
        F_n;
        \bb\Omega_{\sigma}^{\epsilon,c}
    \right)
    \geq
    \operatorname{val}
    \left(
        F;
        \bb\Omega_{\sigma}^{\epsilon,c}
    \right)
    -
    \delta
    =
    \bar\phi_{\sigma}^{\epsilon,c}-\delta \,.
$$
Therefore
$$
    \phi_{\bw}^{\epsilon,c}(\bg,\bh)
    \geq
    \bar\phi_{\sigma}^{\epsilon,c}-\delta
    \geq
    \bar\phi_{\sigma}^{\epsilon,c}-s \,.
$$

Now fix $j\in\mathcal J_n^\epsilon$. On the same event,
	$$
		\phi_j^\epsilon(\bg,\bh)
		=
		\operatorname{val}
		\left(
			F_n;
			\bb\Omega_{n,j}^{\epsilon}
		\right) \,.
	$$
	Since $\bb\Omega_{n,j}^{\epsilon} \subseteq \bb\Omega_j^\epsilon$
	only through the innermost minimisation domain in $(\bu,\theta,t)$,
	$$
		\operatorname{val}
		\left(
			F_n;
			\bb\Omega_{n,j}^{\epsilon}
		\right)
		\ge
		\operatorname{val}
		\left(
			F_n;
			\bb\Omega_j^\epsilon
		\right) \,.
	$$
	On $\mathcal U_n(\delta)$, \eqref{eq:value_lipschitz} gives
	$$
		\operatorname{val}
		\left(
			F_n;
			\bb\Omega_j^\epsilon
		\right)
		\ge
		\operatorname{val}
		\left(
			F;
			\bb\Omega_j^\epsilon
		\right)
		-
		\delta
		=
		\bar\phi_j^\epsilon-\delta\,.
	$$
	Hence
	$$
		\phi_j^\epsilon(\bg,\bh)
		\ge
		\bar\phi_j^\epsilon-s\,,
		$$
	for every $j\in\mathcal J_n^\epsilon$.

	By the definition of $\mathcal A_n^\epsilon(s)$ in the statement,
	we have shown that
	$$
		\mathcal R_{n,\mathrm{res}}^\epsilon
		\cap
		\mathcal U_n(\delta)
		\cap
		\mathcal N_n(\delta)
		\subseteq
		\mathcal A_n^\epsilon(s)\,.
	$$
	Taking complements and using the three conditional probability bounds
	above gives the claim, after adjusting constants.
\end{proof}

\subsubsection{Restricted PO barriers and value concentration}

\begin{lemma}[Restricted PO value barriers]
\label{lemma:restricted_PO_value_barriers}
For every $\epsilon>0$, with $\eta_\epsilon$ from
Lemma~\ref{lemma:limiting_AO_gap}, there are $t_0=t_0(\epsilon)>0$,
constants $C_\epsilon,c_\epsilon>0$, and $N_\epsilon\in\mathbb N$ such
that, for all $n>N_\epsilon$,
$$
    \mathds{1}_{\mathcal S_n(t_0)}
    \Pr\left(
        \Psi>\bar\phi+\eta_\epsilon
        \mid
        \mathcal G_n
    \right)
    \leq
    C_\epsilon\exp\{-c_\epsilon n\}\,,
$$
$$
    \mathds{1}_{\mathcal S_n(t_0)}
    \Pr\left(
        \Psi_{\bw}^{\epsilon,c}
        <
        \bar\phi_{\sigma}^{\epsilon,c}-\eta_\epsilon
        \mid
        \mathcal G_n
    \right)
    \leq
    C_\epsilon\exp\{-c_\epsilon n\}\,,
$$
where the second display is trivial when
$\mathcal S_{\bw}^{\epsilon,c}=\emptyset$, and
$$
    \mathds{1}_{\mathcal S_n(t_0)}
    \Pr\left(
        \exists j\in\mathcal J_n^\epsilon:
        \Psi_j^\epsilon
        <
        \bar\phi_j^\epsilon-\eta_\epsilon
        \mid
        \mathcal G_n
    \right)
    \leq
    C_\epsilon\exp\{-c_\epsilon n\} \,.
$$
Consequently, the event $\mathcal P_n^\epsilon$ in
\eqref{event:P_n_eps} satisfies
$$
    \mathds{1}_{\mathcal S_n(t_0)}
    \Pr\left(
        \lnot\mathcal P_n^\epsilon
        \mid
        \mathcal G_n
    \right)
    \leq
    C_\epsilon\exp\{-c_\epsilon n\} \,.
$$
\end{lemma}

\begin{proof}
Set $s_\epsilon=\eta_\epsilon/2$. By
Lemma~\ref{lemma:restricted_AO_value_bounds} with $s=s_\epsilon$, after
possibly decreasing $t_0$ and adjusting constants,
$$
    \mathds{1}_{\mathcal S_n(t_0)}
    \Pr\left(
        \abs{\phi(\bg,\bh)-\bar\phi}>s_\epsilon
        \mid
        \mathcal G_n
    \right)
    \leq
    C_\epsilon\exp\{-c_\epsilon n\}\,,
$$
$$
    \mathds{1}_{\mathcal S_n(t_0)}
    \Pr\left(
        \phi_{\bw}^{\epsilon,c}(\bg,\bh)
        <
        \bar\phi_{\sigma}^{\epsilon,c}-s_\epsilon
        \mid
        \mathcal G_n
    \right)
    \leq
    C_\epsilon\exp\{-c_\epsilon n\}\,,
$$
and
$$
    \mathds{1}_{\mathcal S_n(t_0)}
    \Pr\left(
        \exists j\in\mathcal J_n^\epsilon:
        \phi_j^\epsilon(\bg,\bh)
        <
        \bar\phi_j^\epsilon-s_\epsilon
        \mid
        \mathcal G_n
    \right)
    \leq
    C_\epsilon\exp\{-c_\epsilon n\} \,.
$$
For the unrestricted value, apply
Proposition~\ref{thm:pathwise_cgmt_localisation} with
$
    W_0=W_p$, $
    D_0=V_p\times[-C_\theta,C_\theta]
$. 
Since $W_p$ is convex, the upper-tail comparison gives, almost surely,
$$
    \Pr\left(
        \Psi>\bar\phi+\eta_\epsilon
        \mid
        \mathcal C_n
    \right)
    \leq
    2
    \Pr\left(
        \phi(\bg,\bh)\geq\bar\phi+\eta_\epsilon
        \mid
        \mathcal C_n
    \right) \,.
$$
Because $s_\epsilon<\eta_\epsilon$,
$$
    \left\{
        \phi(\bg,\bh)\geq\bar\phi+\eta_\epsilon
    \right\}
    \subseteq
    \left\{
        \abs{\phi(\bg,\bh)-\bar\phi}>s_\epsilon
    \right\} \,.
$$
Taking conditional expectations with respect to $\mathcal G_n$ gives the
first bound.
For the radial comparison, define the $\mathcal C_n$-measurable
nonempty compact-valued fallback set 
$$
    \widetilde W_\epsilon
    =
    \begin{cases}
        \mathcal S_{\bw}^{\epsilon,c},
        & \mathcal S_{\bw}^{\epsilon,c}\neq\emptyset,\\
        W_p,
        & \mathcal S_{\bw}^{\epsilon,c}=\emptyset
    \end{cases}\,.
$$
Apply Proposition~\ref{thm:pathwise_cgmt_localisation} with
$W_0=\widetilde W_\epsilon$ and
$D_0=V_p\times[-C_\theta,C_\theta]$. On the event
$\{\mathcal S_{\bw}^{\epsilon,c}\neq\emptyset\}$, the resulting PO and
AO are exactly $\Psi_{\bw}^{\epsilon,c}$ and
$\phi_{\bw}^{\epsilon,c}$. On the complementary event the desired lower-tail
bound is trivial because $\Psi_{\bw}^{\epsilon,c}=+\infty$. Since
$\mathcal S_{\bw}^{\epsilon,c}$ is compact when nonempty and no convexity is
required for the lower-tail comparison, almost surely,
$$
\begin{aligned}
    \mathds{1}_{\{\mathcal S_{\bw}^{\epsilon,c}\neq\emptyset\}}
    \Pr\left(
        \Psi_{\bw}^{\epsilon,c}
        <
        \bar\phi_{\sigma}^{\epsilon,c}-\eta_\epsilon
        \mid
        \mathcal C_n
    \right)
  \leq
    2\mathds{1}_{\{\mathcal S_{\bw}^{\epsilon,c}\neq\emptyset\}}
    \Pr\left(
        \phi_{\bw}^{\epsilon,c}(\bg,\bh)
        \leq
        \bar\phi_{\sigma}^{\epsilon,c}-\eta_\epsilon
        \mid
        \mathcal C_n
    \right) \,.
\end{aligned}
$$
The AO event on the right is contained in
$$
    \left\{
        \phi_{\bw}^{\epsilon,c}(\bg,\bh)
        <
        \bar\phi_{\sigma}^{\epsilon,c}-s_\epsilon
    \right\} \,.
$$
Taking conditional expectations with respect to $\mathcal G_n$ gives the
radial bound.

For each fixed $j\in\mathcal J$, define the
$\mathcal C_n$-measurable nonempty compact convex fallback set 
$$
    \widetilde C_{p,j}^{\epsilon}
    =
    \begin{cases}
        C_{p,j}^{\epsilon},
        & j\in\mathcal J_n^\epsilon,\\
        V_p\times[-C_\theta,C_\theta],
        & j\notin\mathcal J_n^\epsilon
    \end{cases}\,.
$$
Apply Proposition~\ref{thm:pathwise_cgmt_localisation} with
$W_0=W_p$ and $D_0=\widetilde C_{p,j}^{\epsilon}$. On the event
$\{j\in\mathcal J_n^\epsilon\}$, the resulting PO and AO are exactly
$\Psi_j^\epsilon$ and $\phi_j^\epsilon$. Hence, almost surely,
$$
\begin{aligned}
    \mathds{1}_{\{j\in\mathcal J_n^\epsilon\}}
    \Pr\left(
        \Psi_j^\epsilon
        <
        \bar\phi_j^\epsilon-\eta_\epsilon
        \mid
        \mathcal C_n
    \right)
   \leq
    2\mathds{1}_{\{j\in\mathcal J_n^\epsilon\}}
    \Pr\left(
        \phi_j^\epsilon(\bg,\bh)
        \leq
        \bar\phi_j^\epsilon-\eta_\epsilon
        \mid
        \mathcal C_n
    \right) \,.
\end{aligned}
$$
The AO event on the right is contained in
$$
    \left\{
        \phi_j^\epsilon(\bg,\bh)
        <
        \bar\phi_j^\epsilon-s_\epsilon
    \right\} \,.
$$
Summing over the deterministic finite set $\mathcal J$, taking conditional
expectations with respect to $\mathcal G_n$, and using
$\abs{\mathcal J_n^\epsilon}\leq\abs{\mathcal J}$ gives
$$
\begin{aligned}
    \Pr\left(
        \exists j\in\mathcal J_n^\epsilon:
        \Psi_j^\epsilon
        <
        \bar\phi_j^\epsilon-\eta_\epsilon
        \mid
        \mathcal G_n
    \right)
  \leq
    2\abs{\mathcal J}
    \Pr\left(
        \exists j\in\mathcal J_n^\epsilon:
        \phi_j^\epsilon(\bg,\bh)
        <
        \bar\phi_j^\epsilon-s_\epsilon
        \mid
        \mathcal G_n
    \right) \,.
\end{aligned}
$$
The simultaneous AO bound yields the piecewise PO estimate after absorbing
$2\abs{\mathcal J}$ into $C_\epsilon$. A final union bound proves the
estimate for $\lnot\mathcal P_n^\epsilon$.
\end{proof}

\subsection{Lemmas for Section~\ref{sec:test-function-convergence}}

\begin{lemma}[Spherical representation of the PO solution]

\label{lemma:compact_PO_spherical_representation}

Let $\bb d$ be as in \eqref{eq:compact_PO_direction}. Conditionally on
$\mathcal C_n$, the direction $\bb d$ is uniform on
$\mathbb S_{\bb B}^\perp$ and independent of $\bb{\mathfrak Z}$.

\end{lemma}

\begin{proof}
Condition on $\mathcal C_n$ and work on a fixed rank stratum. Write
$
    s=\rank(\bb B)$, $
    k=p-s\geq1$, $
    \mathcal E_\sigma=\{\hat\sigma>0\}$, $
    \bb d'=\bb E^\top\bb d
$. 
By Lemma~\ref{lemma:compact_PO_uniqueness}, outside a conditional null
event the compact problem in \eqref{eq:lagrangian4} has a unique primal
minimiser.
Let $\bb R$ be any $k\times k$ orthogonal matrix. Replacing $\bH_2$ by
$\bH_2\bb R^\top$ maps the primal minimiser to
$
    (\hat\bv,\bb R\hat\bw,\hat\theta,\hat{\bb\eta})
$. 
Indeed,
$
    \bH_2\bb R^\top\bb R\hat\bw
    =
    \bH_2\hat\bw$, $
    \vnorm{\bb R\hat\bw}_2
    =
    \vnorm{\hat\bw}_2
$, 
so both the objective and the feasible domain are preserved. Consequently,
$\bb{\mathfrak Z}$ and $\mathcal E_\sigma$ are unchanged, while on
$\mathcal E_\sigma$,
$
    \bb d'
    =
    \frac{\hat\bw}{\vnorm{\hat\bw}_2}
$
is mapped to $\bb R\bb d'$.
Conditional on $\mathcal C_n$, the matrix $\bH_2$ has independent standard
Gaussian entries, and hence
$$
    \bH_2\bb R^\top
    \overset{\mathrm d}{=}
    \bH_2\,,
$$
conditionally on $\mathcal C_n$. 
It follows that, for every bounded Borel function
$f:\mathbb S^{k-1}\to\Re$ and every bounded Borel function $g$ of
$\bb{\mathfrak Z}$,
\begin{equation}
\label{eq:compact_PO_fixed_rotation_identity}
\begin{aligned} 
    \expect\left[
        f(\bb d')g(\bb{\mathfrak Z})
        \mathds{1}\{\mathcal E_\sigma\}
       \mid
        \mathcal C_n
    \right]
   =
    \expect\left[
        f(\bb R\bb d')g(\bb{\mathfrak Z})
        \mathds{1}\{\mathcal E_\sigma\}
       \mid
        \mathcal C_n
    \right] \,. 
\end{aligned}
\end{equation}
Let $\mu$ be the conditional law of $\bb d'$ given
$\sigma(\mathcal C_n,\bb{\mathfrak Z})$. Since
$
    \mathcal E_\sigma\in\sigma(\bb{\mathfrak Z})
$,
\eqref{eq:compact_PO_fixed_rotation_identity} and the almost-sure uniqueness
of regular conditional distributions \citep[Chapter~2, Section~7]{shiryaev:2016} imply that, for each fixed
$k\times k$ orthogonal matrix $\bb R$,
\begin{equation}
\label{eq:compact_PO_fixed_rotation_kernel}
    \mu(\bb R^{-1}A)
    =
    \mu(A) \,, 
\end{equation}
for every Borel set $A\subseteq\mathbb S^{k-1}$, on
$\mathcal E_\sigma$, outside a null set that may depend on $\bb R$.

The collection of $k\times k$ orthogonal matrices is a compact metric
subset of $\Re^{k\times k}$ and therefore has a countable dense subset
$D$ \citep[Lemma~3.26]{aliprantis+et+al:2006}. Intersect the probability-one events in
\eqref{eq:compact_PO_fixed_rotation_kernel} over $\bb R\in D$. On the
resulting probability-one event, $\mu$ is invariant under every
$\bb R\in D$ simultaneously on $\mathcal E_\sigma$.
Now let $\bb R$ be any $k\times k$ orthogonal matrix. Since $D$ is dense,
there exists a sequence $(\bb R_m)_{m\in\mathbb N}\subseteq D$ such that
$$
    \mnorm{\bb R_m-\bb R}_2
    \longrightarrow
    0\,.
$$
Therefore, on
$\mathcal E_\sigma$, and continuous $f$, 
\begin{equation}
\begin{aligned}
    \int_{\mathbb S^{k-1}}
    f(\bb Rx)\,\mu(\mathrm dx)
    =
    \lim_{m\to\infty}
    \int_{\mathbb S^{k-1}}
    f(\bb R_mx)\,\mu( \mathrm dx)
=
    \int_{\mathbb S^{k-1}}
    f(x)\,\mu(\mathrm  dx) \,. 
\end{aligned} 
\end{equation}
Since continuous functions determine Borel probability measures on the
compact sphere \citep[Lemma~9.3.2]{dudley:2018}, $\mu$ is invariant under every $k\times k$
orthogonal matrix. The unique probability measure with this invariance is
the uniform probability measure on $\mathbb S^{k-1}$ (see
\citealt[Theorem~4.4]{eaton:1989}). Hence, on $\mathcal E_\sigma$,
$\bb d'$ is conditionally uniform on $\mathbb S^{k-1}$. 
Since $\bb E^\top\bb E=\bb I_k$ and
$\range(\bb E)=\range(\bb B)^\perp$, the maps
$x\mapsto\bb E x$ and $y\mapsto\bb E^\top y$ are inverse isometries between
$\mathbb S^{k-1}$ and $\mathbb S_{\bb B}^\perp$. Hence, for every Borel set
$A\subseteq\mathbb S_{\bb B}^\perp$,
$$
    \mathds{1}\{\mathcal E_\sigma\}
    \Pr\left(
        \bb d\in A
        \mid
        \sigma(\mathcal C_n,\bb{\mathfrak Z})
    \right)
    =
    \mathds{1}\{\mathcal E_\sigma\}
    \Pr\left(
        \bb d'\in\bb E^\top A
        \mid
        \sigma(\mathcal C_n,\bb{\mathfrak Z})
    \right) \,. 
$$
The right-hand side is the uniform spherical probability of
$\bb E^\top A$, which equals that of $A$ because $\bb E^\top$ is an
isometry. Thus $\bb d$ is conditionally uniform on
$\mathbb S_{\bb B}^\perp$ on $\mathcal E_\sigma$.
On $\lnot \mathcal E_\sigma$, \eqref{eq:compact_PO_direction} gives, outside a
null event,
$$
    \bb d
    =
    \frac{
        \bb P^\perp\bb\xi_n
    }{
        \vnorm{\bb P^\perp\bb\xi_n}_2
    }
    =
    \bb E
    \frac{
        \bb E^\top\bb\xi_n
    }{
        \vnorm{\bb E^\top\bb\xi_n}_2
    }\,.
$$
The vector $\bb\xi_n$ is independent of
$\sigma(\mathcal C_n,\bb{\mathfrak Z})$, and $\bb E$ is
$\mathcal C_n$-measurable. Hence, conditionally on
$\sigma(\mathcal C_n,\bb{\mathfrak Z})$,
$$
    \bb E^\top\bb\xi_n
    \sim
    \mathrm N(\bb0_k,\bb I_k)\,.
$$
The Gaussian polar decomposition
\citep[Exercise~3.3.7]{vershynin:2018} therefore shows that $\bb d$ is
uniform on $\mathbb S_{\bb B}^\perp$ and independent of
$\bb{\mathfrak Z}$ on $\lnot \mathcal E_\sigma$.
Since $\mathcal E_\sigma\in\sigma(\bb{\mathfrak Z})$, the two cases give
$$
        \bb d
       \mid
        \mathcal C_n,\bb{\mathfrak Z}
    \overset{\mathrm d}{=}
        \bb d
       \mid
        \mathcal C_n
    \overset{\mathrm d}{=}
    \operatorname{Unif}
    \left(
        \mathbb S_{\bb B}^\perp
    \right)\,,
$$
almost surely as required.
\end{proof}

\subsection{Minor auxiliary lemmas}


\subsubsection{Concentration results}

\begin{lemma}[Gaussian design spectral-norm bound]
	\label{lemma:sv_ub}
	Let
	$
		\widetilde{\bH}
		=
		[\bb1,\bb G]
		\in
		\Re^{n\times(p+1)}
	$,
	where $\bb G$ has i.i.d.\ $\mathrm{N}(0,1)$ entries. There exists a
	universal constant $c_{\mathrm{sv}}>0$ such that, for every $t>0$,
	$$
		\Pr\left(
			\mnorms{\widetilde{\bH}}_2
			>
			2\sqrt n+\sqrt p+t
		\right)
		\leq
		2\exp\{-c_{\mathrm{sv}}t^2\}\,.
	$$
\end{lemma}

\begin{proof}
	The triangle inequality gives
	$$
		\mnorms{\widetilde{\bH}}_2
		\leq
		\sqrt n+\mnorm{\bb G}_2\,,
	$$
	and the claim follows from
	\citealt[][Corollary~7.3.3]{vershynin:2018}.
\end{proof}

\begin{lemma}[Uniform row-projection bound]
	\label{lemma:set_bound}
	Let
	$
		\widetilde{\bH}
		=
		[\bb1,\bb G]
		\in
		\Re^{n\times(p+1)}
	$.
	For $\epsilon\in(0,1)$ and $\bv\in\Re^{p+1}$, define
	$$
		S_\epsilon(\bv)
		=
		\left\{
			i\in[n]:
			\abs{\bv^\top\widetilde{\bh}_i}
			\leq
			\frac{3}{\sqrt\epsilon}\vnorm{\bv}_2
		\right\}\,,
	$$
	where $\widetilde{\bh}_i^\top$ is the $i$th row of
	$\widetilde{\bH}$. If
	$\mnorms{\widetilde{\bH}}_2\leq3\sqrt n$, then
	$$
		\underset{\bv\in\Re^{p+1}}{\inf}
		\abs{S_\epsilon(\bv)}
		\geq
		(1-\epsilon)n\,.
	$$
	If, moreover, $\bb G$ has i.i.d.\ $\mathrm{N}(0,1)$ entries and
	$n>p$, then
	$$
		\Pr\left(
			\underset{\bv\in\Re^{p+1}}{\inf}
			\abs{S_\epsilon(\bv)}
			<
			(1-\epsilon)n
		\right)
		\leq
		2\exp\left\{
			-c_{\mathrm{sv}}n
			\left(
				1-\sqrt{\frac pn}
			\right)^2
		\right\}\,.
	$$
\end{lemma}

\begin{proof}
    The deterministic argument is similar to
\citet[][Lemma~7]{sur:2019}.
	On the event
	$\{\mnorms{\widetilde{\bH}}_2\leq3\sqrt n\}$, for every
	$\bv\neq\bb0_{p+1}$,
	$$
	\begin{aligned}
		9n\vnorm{\bv}_2^2
		&\geq
		\vnorm{\widetilde{\bH}\bv}_2^2
		\\
		&\geq
		\sum_{i\notin S_\epsilon(\bv)}
		(\bv^\top\widetilde{\bh}_i)^2
		\\
		&\geq
		\left\{
			n-\abs{S_\epsilon(\bv)}
		\right\}
		\frac{9\vnorm{\bv}_2^2}{\epsilon}\,.
	\end{aligned}
	$$
	Hence
	$
		\abs{S_\epsilon(\bv)}
		\geq
		(1-\epsilon)n
	$.
	For $\bv=\bb0_{p+1}$,
	$
		S_\epsilon(\bv)=[n]
	$,
	so the deterministic claim follows. The probability bound follows from
	Lemma~\ref{lemma:sv_ub} with
	$
		t=\sqrt n-\sqrt p
	$.
\end{proof}

\subsubsection{Moreau-envelope analytic facts}

\begin{lemma}[Lipschitz continuity of the Moreau envelope]
\label{lemma:moreau_lipschitz}
Let $\rho:\Re\to\Re$ be convex and $L_\rho$-Lipschitz. For $\lambda>0$,
define
$$
    m_\rho(x,\lambda)
    =
    \underset{\eta\in\Re}{\min}
    \left\{
        \rho(\eta)
        +
        \frac{(\eta-x)^2}{2\lambda}
    \right\},
    \quad
    m_\rho(x,0)=\rho(x)\,.
$$
Then:
\begin{enumerate}[label=(\roman{*})]
    \item for every $\lambda\geq0$, the map
    $x\mapsto m_\rho(x,\lambda)$ is $L_\rho$-Lipschitz on $\Re$;

    \item for every $x\in\Re$, the map
    $\lambda\mapsto m_\rho(x,\lambda)$ is
    $L_\rho^2/2$-Lipschitz on $[0,\infty)$.
\end{enumerate}
In particular, for $\rho(x)=\log(1+\exp\{x\})$, define, for $\lambda>0$,
$$
    M_\rho(\bx,\lambda)
    =
    \underset{\bb\eta\in\Re^n}{\min}
    \left\{
        \sum_{i=1}^n\rho(\eta_i)
        +
        \frac{\vnorm{\bb\eta-\bx}_2^2}{2\lambda}
    \right\},
    \quad
    M_\rho(\bx,0)
    =
    \bb1^\top\bb\rho(\bx)\,.
$$
Then $M_\rho$ is $\sqrt n$-Lipschitz in $\bx$,
$n/2$-Lipschitz in $\lambda$, and
$\sqrt{n+n^2/4}$-Lipschitz jointly in $(\bx,\lambda)$ on
$\Re^n\times[0,\infty)$.
\end{lemma}

\begin{proof}
We prove each claim in turn.
\begin{enumerate}[label=(\roman{*})]
    \item
    Fix $\lambda\geq0$. For $\lambda=0$, the claim follows from
    $m_\rho(\cdot,0)=\rho$. Suppose that $\lambda>0$, and set
    $\eta=\prox{\lambda\rho}{x}$. The standard
    Moreau-envelope derivative formula
    \citep[Lemma~2, equation~(29)]{salehi+et+al:2019}, applied with
    $d=1$, $\Phi=\rho$, and $t=\lambda$, gives
    $$
        \partial_xm_\rho(x,\lambda)
        =
        \frac{x-\eta}{\lambda}\,.
    $$
    The proximal optimality condition
    \citep[cf.][Theorem~6.39]{beck:2017} gives
    $$
        \frac{x-\eta}{\lambda}
        \in
        \partial\rho\{\eta\}\,.
    $$
    Since every subgradient of the $L_\rho$-Lipschitz function $\rho$
    has absolute value at most $L_\rho$
    \citep[cf.][Theorem~3.61]{beck:2017},
    $$
        \abs{\partial_xm_\rho(x,\lambda)}
        \leq
        L_\rho\,.
    $$
    The mean-value theorem proves that
    $x\mapsto m_\rho(x,\lambda)$ is $L_\rho$-Lipschitz.

    \item
    Fix $x\in\Re$. For $\lambda>0$, write again
    $\eta=\prox{\lambda\rho}{x}$. The second
    Moreau-envelope derivative formula in
    \citet[Lemma~2, equation~(29)]{salehi+et+al:2019} gives
    $$
        \partial_\lambda m_\rho(x,\lambda)
        =
        -
        \frac{
            \{x-\eta\}^2
        }{
            2\lambda^2
        }\,.
    $$
    By the same proximal-optimality and bounded-subgradient argument as
    above,
    $$
        \abs{x-\eta}
        \leq
        L_\rho\lambda\,,
    $$
    and therefore
    $$
        \abs{\partial_\lambda m_\rho(x,\lambda)}
        \leq
        \frac{L_\rho^2}{2}\,.
    $$
    The mean-value theorem proves the claimed bound when both smoothing
    parameters are positive.

    It remains to include $\lambda=0$. For $\lambda>0$, testing the
    defining minimum at $\eta=x$ gives
    $m_\rho(x,\lambda)\leq\rho(x)$, while
    $\rho(\eta)\geq\rho(x)-L_\rho\abs{\eta-x}$ gives
    $$
        m_\rho(x,\lambda)
        \geq
        \rho(x)
        +
        \underset{r\geq0}{\min}
        \left\{
            \frac{r^2}{2\lambda}
            -
            L_\rho r
        \right\}
        =
        \rho(x)
        -
        \frac{L_\rho^2}{2}\lambda\,.
    $$
    Since $m_\rho(x,0)=\rho(x)$,
    $$
        \abs{
            m_\rho(x,\lambda)-m_\rho(x,0)
        }
        \leq
        \frac{L_\rho^2}{2}\lambda\,.
    $$
    Hence $\lambda\mapsto m_\rho(x,\lambda)$ is
    $L_\rho^2/2$-Lipschitz on $[0,\infty)$.
\end{enumerate}

For logistic $\rho$, one has $L_\rho=1$, and separability
\citep[Lemma~6, equation~(35)]{salehi+et+al:2019} gives
$M_\rho(\bx,\lambda)=\sum_{i=1}^nm_\rho(x_i,\lambda)$. Consequently,
for $\bx,\bz\in\Re^n$ and $\lambda,\mu\geq0$,
$$
\begin{aligned}
    \abs{
        M_\rho(\bx,\lambda)
        -
        M_\rho(\bz,\mu)
    }
    &\leq
    \vnorm{\bx-\bz}_1
    +
    \frac n2\abs{\lambda-\mu}
    \\
    &\leq
    \sqrt n\,\vnorm{\bx-\bz}_2
    +
    \frac n2\abs{\lambda-\mu}
    \\
    &\leq
    \sqrt{
        n+\frac{n^2}{4}
    }\,
    \vnorm{
        \begin{bmatrix}
            \bx\\
            \lambda
        \end{bmatrix}
        -
        \begin{bmatrix}
            \bz\\
            \mu
        \end{bmatrix}
    }_2\,.
\end{aligned}
$$
The separate vector bounds follow by taking $\lambda=\mu$ or $\bx=\bz$.
\end{proof}

\clearpage
\putbib
\end{bibunit}

\end{document}